%% file: cn-proofs.tex
\documentclass{gtpart}
\agtart

\usepackage{tikz}
\usetikzlibrary{calc,arrows,decorations.pathmorphing,backgrounds,fit} %Added options except calc
\usetikzlibrary{patterns}
\usetikzlibrary{positioning}
\usepackage{ifthen}
\usepackage{amsmath}
\usepackage{amsfonts}
\usepackage{amssymb}
\usepackage{amsthm}
\usepackage{comment}
\usepackage{epsfig}
\usepackage{psfrag}
\usepackage{mathrsfs}
\usepackage{amscd}
\usepackage[all,rotate]{xy}
\usepackage{rotating}
\usepackage{lscape}
\usepackage{amsbsy}
\usepackage{verbatim}
\usepackage{moreverb}
\usepackage{xspace}
\usepackage{calc}
\usepackage{multicol} % Used for appendix
\usepackage{needspace}

\title{Conformal nets IV: The $3$-category} 

       \author[A. Bartels]{Arthur Bartels}
       \givenname{Arthur}
       \surname{Bartels}
      \address{Westf\"alische Wilhelms-Universit\"at M\"unster\\
               Mathematisches Institut\\
               Einsteinstr.~62,
               D-48149 M\"unster, Deutschland}
        \email{bartelsa@wwu.de}
        
       \author[C. L. Douglas]{Christopher L. Douglas} 
       \givenname{Christopher L.}
       \surname{Douglas}
      \address{Mathematical Institute\\ Radcliffe Observatory Quarter\\ Woodstock Road\\ Oxford\\ OX2 6GG\\ United Kingdom}
        \email{cdouglas@maths.ox.ac.uk}
        
       \author[A. Henriques]{Andr{\'e} Henriques}
       \givenname{Andr{\'e}}
       \surname{Henriques}
      \address{Mathematisch Instituut\\
               Universiteit Utrecht, Postbus 80.010\\
               3508 TA Utrecht, The Netherlands}
        \email{a.g.henriques@uu.nl}

\arxivreference{1605.00662}
\arxivpassword{6npwu}

\newtheorem*{theorem*}{Theorem}
\newtheorem*{proposition*}{Proposition}
\newtheorem*{lemma*}{Lemma}
\newtheorem*{corollary*}{Corollary}
\newtheorem*{construction*}{Construction}
\newtheorem{theorem}{Theorem}
\newtheorem{proposition}[theorem]{Proposition}
\newtheorem{lemma}[theorem]{Lemma}

\newtheoremstyle{lemmaletter}{}{}{\itshape}{}{\bfseries}{.}{.5em}{#1\thmnote{ #3}}
\theoremstyle{lemmaletter}
\newtheorem*{lemmaletter}{Lemma}

\theoremstyle{definition}
\newtheorem*{definition}{Definition}

\newtheorem*{axprop}{Proposition} % To avoid italics in the diagrams.
\newtheorem*{axlem}{Lemma} % To avoid italics in the diagrams.

\theoremstyle{remark}
\newtheorem{remark}[theorem]{Remark}
\newtheorem*{remark*}{Remark}

\renewcommand{\theenumi}{(\roman{enumi})}

\newcommand{\ra}{\rightarrow}

\newcommand{\cb}{\raisebox{.6ex-.5\height}}

\DeclareRobustCommand{\SkipTocEntry}[5]{}

       \newcommand{\cala}{\mathcal{A}}
       \newcommand{\calb}{\mathcal{B}}
  \newcommand{\IC}{\mathbb{C}}     \newcommand{\calc}{\mathcal{C}}

  \newcommand{\IR}{\mathbb{R}}

  \newcommand{\IZ}{\mathbb{Z}}

  \newcommand{\bfB}{{\mathbf B}}

  \DeclareMathOperator{\id}{id}
  
  \DeclareMathOperator{\Diff}{Diff}
  \DeclareMathOperator{\Hom}{Hom}  

  \newcommand{\x}{{\times}}
  \newcommand{\ox}{{\otimes}}
  
  \newcommand{\e}{{\varepsilon}}

  \newcommand{\forget}{\mathit{forget}}
  \newcommand{\op}{\mathrm{op}}

  \newcommand{\byDef}{D}
  \newcommand{\byNat}{N}
  \newcommand{\byComOp}{C}
  \newcommand{\byLemmaL}{J}
  \newcommand{\byLemmaM}{K}
  \newcommand{\byLemmaX}{L}
  \newcommand{\byLemmaZ}{M}
  
  \newcommand{\INT}{{\mathsf{INT}}}
  \newcommand{\VN}{{\mathsf{VN}}}
  \newcommand{\CN}{{\mathsf{CN}}}

  \def\tworarrow{\hspace{.1cm}{\setlength{\unitlength}{.50mm}\linethickness{.09mm}
	\begin{picture}(8,8)(0,0)\qbezier(0,4)(4,7)(8,4)\qbezier(0,1)(4,-2)(8,1)\qbezier(3.5,4)(3.5,3)(3.5,1.5)
	\qbezier(4.5,4)(4.5,3)(4.5,1.5)\qbezier(4,0.8)(4.5,1.7)(5.5,2)\qbezier(4,0.8)(3.5,1.7)(2.5,2)
	\qbezier(8,1)(7.4,.2)(7.7,-.7)\qbezier(8,1)(7,1)(6.5,1.5)\qbezier(8,4)(7.4,4.8)(7.7,5.7)
	\qbezier(8,4)(7,4)(6.5,3.5)\end{picture}\hspace{.1cm}}}

  \newcounter{commentcounter}

  \newcounter{newpicturecounter}

\graphicspath{{FinalPix/miscPix/}{../../../Archived-projects/Completed/IB-early/LQFTPix/FinalPix/ci3Pix/}{FinalPix/cdcPix/}{FinalPix/cbcPix/}{FinalPix/c2cPix/newnewNumberScheme/}{FinalPix/MoritaPix/}{FinalPix/dualPix/}{FinalPix/triPix/}{FinalPix/dbPix/}}
 
\def\cC{\mathcal C} 
\def\Rho{\mathrm{P}}

\newcommand{\ingeps}{\includegraphics[scale=1]}

\newcommand{\ing}{\includegraphics[scale=.3]}
\newcommand{\ingt}{\includegraphics[scale=.25]}

\newcommand{\xys}{\xymatrix@-15pt@M+5pt}
\newcommand{\arr}{\ar@{=>}}
\newcommand{\scs}{.35in}

\newcommand{\sca}{.5in}
\newcommand{\scaa}{.55in}
\newcommand{\scb}{.6in}
\newcommand{\scbb}{.65in}
\newcommand{\scc}{.7in}
\newcommand{\sccc}{.75in}

\newcommand{\sce}{.9in}
\newcommand{\scf}{1in}

\newcommand{\hsqza}{-25pt}
\newcommand{\hsqzb}{-65pt}

\newcommand{\nid}{\noindent}

\newcommand{\db}{\raisebox{3ex-\height}}

\newcommand{\bb}{\raisebox{2ex-.5\height}}

\newcommand\ignore[1]{}

\definecolor{AHcolor}{rgb}{0.5,0.0,0.5}   % Textcolor for AB
\definecolor{CDcolor}{rgb}{0.7,0.0,0.3}   % Textcolor for CD
\definecolor{ABcolor}{rgb}{0.2,0.8,0.2}   % Textcolor for AH
\include{ib-defs-forcniv}

\include{ib-code-forcniv}

\include{cnii-defs}

\include{cnii-code}

\begin{document}

\begin{abstract}
Conformal nets are a mathematical model for conformal field theory, and defects between conformal nets are a model for an interaction or phase transition between two conformal field theories.  We previously introduced a notion of composition, called \emph{fusion}, between defects.  We also described a notion of sectors between defects, modeling an interaction among or transformation between phase transitions, and defined fusion composition operations for sectors.  In this paper we prove that altogether the collection of conformal nets, defects, sectors, and intertwiners, equipped with the fusion of defects and fusion of sectors, forms a symmetric monoidal 3-category.  This 3-category encodes the algebraic structure of the possible interactions among conformal field theories.
\end{abstract}

\maketitle

%%% Added for appendix
\xyoption{poly}
\xyoption{2cell}
\objectmargin{5pt}
%%%

%\setcounter{tocdepth}{2}
%\tableofcontents

%\pagebreak

\setcounter{page}{1} 

\section*{Introduction}

%\addtocontents{toc}{\SkipTocEntry }
\addtocontents{toc}{\SkipTocEntry }
\subsection*{Background and results}
Conformal nets are a mathematical formalization of the notion of conformal field theory~\cite{Brunetti-Guido-Longo(1993modular+duality-in-CQFT),Buchholz-Mack-Todorov(1988current-alg),Buchholz-Schulz-Mirbach(Haag-duality-in-conformal-quantum-field-theory),Gabbiani-Froehlich(OperatorAlg-CFT),Kawahigashi-Longo(2004classification),Kawahigashi-Longo-Mueger(2001multi-interval),Wassermann(Operator-algebras-and-conformal-field-theory)}.  This paper is the fourth in a series investigating \emph{coordinate-free conformal nets} and their \emph{defects}~\cite{BDH(nets), BDH(modularity), BDH(1*1)}.  The notion of coordinate-free conformal nets is a modification of the usual notion of conformal nets, in which one does not demand the positive-energy condition and does not require the existence of a vacuum vector in the vacuum sector.\footnote{We refer to the first paper of the series~\cite[\S4]{BDH(nets)} for a detailed comparison between the usual notion and the coordinate-free notion of conformal nets.}  Here we will use the term `conformal net' to refer to coordinate-free conformal nets.

In the preceding paper~\cite{BDH(1*1)}, we introduced defects between conformal nets and gave examples thereof.  (Conjecturally, such defects can also be obtained by taking spatial slices of defects on two-dimensional Minkowski space-time~\cite{Bischoff-Kawahigashi-Longo-Rehren(Phase-boundaries-in-algebraic-conformal-QFT)}.)  In this paper, defects will provide a notion of 1-morphism between conformal nets.  Crucially, in~\cite[Def.\,3.1]{BDH(1*1)}, we introduced the operation of \emph{fusion of defects}, when the conformal nets have finite-index---this fusion will provide the composition operation on 1-morphisms between nets.  We defined the associator isomorphism for this fusion operation, and proved it satisfies the pentagon equation.  We investigated \emph{sectors between defects}---in what follows these will form the 2-morphisms between 1-morphisms between nets---and introduced the horizontal and vertical fusion of sectors, along with their respective associators.  We also constructed the fundamental \emph{interchange isomorphism} relating the vertical and horizontal fusion of sectors.  All these structures and more are summarized in the first appendix of that paper~\cite{BDH(1*1)}.

Given all the work done in~\cite{BDH(1*1)} constructing the notions of higher morphisms between nets and their composition operations, one might imagine that it would be straightforward to establish that finite-index conformal nets, defects, sectors, and intertwiners form a \emph{tricategory}.  But unfortunately, tricategories are complicated beasts.  This entire paper is, in effect, devoted to the proof of the following theorem:
\begin{theorem}\label{thm-tricat}
The collections of all finite-index conformal nets, defects, sectors, and intertwiners form the objects, $1$-morphisms, $2$-morphisms, and $3$-morphisms of a tricategory.
\end{theorem}

Let $\mathsf{CAT}$ denote the 2-category of categories; previously~\cite{Douglas-Henriques(Internal-bicategories)}, we introduced a notion of tricategory called a \emph{dicategory object in $\mathsf{CAT}$}.\footnote{A dicategory object differs from a bicategory object in that the associativity structures (but not the unital structures) are strict.}  A dicategory object in $\mathsf{CAT}$ consists of a category of objects, a category of 1-morphisms, and a category of 2-morphisms, together with various functors encoding identity and composition operations and various natural transformations encoding compatibility relationships between the operations, all satisfying various axioms encoding coherences between the compatibility relationships.  We proved in~\cite{Douglas-Henriques(Internal-bicategories)} that every dicategory object in $\mathsf{CAT}$ has an associated tricategory---Theorem~\ref{thm-tricat} is therefore a corollary of the following result:
\begin{theorem}\label{thm-dicat}
The groupoid of all finite-index conformal nets, the groupoid of all defects, and the category of all sectors form, respectively, the category of objects, the category of $1$-morphisms, and the category of $2$-morphisms for a dicategory object in $\mathsf{CAT}$.
\end{theorem}

The tensor product of two conformal nets is again a conformal net, the tensor product of two defects is again a defect, and similarly for sectors and for intertwiners.  One might therefore naturally conjecture that finite-index conformal nets form not just a tricategory but in fact a \emph{symmetric monoidal tricategory}.  A useful feature of the formalism of dicategory objects is that it makes it particularly simple to incorporate the symmetric monoidal structure.  Let $\mathsf{SMC}$ denote the $2$-category of symmetric monoidal categories, strong (as opposed to lax) symmetric monoidal functors, and symmetric monoidal natural transformations.  Our approach to Theorem~\ref{thm-dicat} immediately also proves the follow symmetric monoidal strengthening:
\begin{theorem}\label{thm-dicatsmc}
The symmetric monoidal groupoid of finite-index conformal nets, the symmetric monoidal groupoid of defects, and the symmetric monoidal category of sectors form, respectively, the object of objects, the object of $1$-morphisms, and the object of $2$-morphisms of a dicategory object in $\mathsf{SMC}$.
\end{theorem}
\noindent This theorem is one precise formulation (indeed the only one available at present) of the statement that finite-index conformal nets form a symmetric monoidal tricategory: finite-index conformal nets form a dicategory object in $\mathsf{SMC}$, and the notion of dicategory object in $\mathsf{SMC}$ may be interpreted as a notion of symmetric monoidal tricategory.

A dicategory in $\mathsf{SMC}$ consists of: [0-data] a symmetric monoidal category of objects, a symmetric monoidal category of 1-morphisms, and a symmetric monoidal category of 2-morphisms; [1-data] various symmetric monoidal functors between (fiber products of) these categories, encoding identity and composition operations; [2-data] various symmetric monoidal natural transformations between (products and composites of) these functors, encoding compatibility relationships between the operations; and [3-axioms] various axioms for these transformations, encoding coherences between the compatibility relationships.  These pieces of structure are tracked by labels of the form [D0-x], [D1-x], [D2-x], [D3-x], following the same numbering scheme as in~\cite[Def.\,3.3]{Douglas-Henriques(Internal-bicategories)}; in the body of this paper we will abbreviate these labels to [0-x], [1-x], [2-x], [3-x] respectively.  For ease of reference, the complete definition of a dicategory object in a 2-category is compiled here in Appendix~\ref{appendix}, and an abbreviated definition is depicted there in Table~\ref{table-internaldicategory}.

%\addtocontents{toc}{\SkipTocEntry }
\addtocontents{toc}{\SkipTocEntry }
\subsection*{Overview}
In our case, the category of objects [D0-0] is the (symmetric monoidal) category $\CN_0$ of finite conformal nets (finite direct sums of irreducible conformal nets with finite-index), together with isomorphisms between them; the category of 1-morphisms [D0-1] is the (symmetric monoidal) category $\CN_1$ of defects (between finite conformal nets), together with isomorphisms of defects; and the category of 2-morphisms [D0-2] is the (symmetric monoidal) category $\CN_2$ of sectors (between defects between finite conformal nets), together with homomorphism of sectors (also called intertwiners) that cover isomorphisms of defects and of conformal nets; the intertwiners play the role of 3-morphisms in the overall 3-category.  These categories of nets, defects, and sectors are discussed in Section~\ref{sec-0data} below.

The most important operation in the 3-category is the composition of 1-morphisms [D1-2], which here is a (symmetric monoidal) functor $$\CN_1 \times_{\CN_0} \CN_1 \to \CN_1$$ given by the fusion of defects: $$({}_\cala D_\calb, {}_\calb E_\calc) \mapsto D \circledast_\calb E.$$  The existence of this fusion operation is proven in~\cite[Thm~1.44]{BDH(1*1)}.  The vertical composition of 2-morphisms [D1-4] is a functor $$\CN_2 \times_{\CN_1} \CN_2 \to \CN_2,$$ given by the vertical fusion of sectors: $$({}_D H_E, {}_E K_F) \mapsto H \boxtimes_E K.$$  This operation is defined in~\cite[Sec 2.C]{BDH(1*1)}.  Horizontal composition of 2-morphisms is encoded indirectly using the vertical composition of 2-morphisms together with left and right whisker operations [D1-5, D1-6], which are functors
\begin{align*}
\CN_2 \times_{\CN_0} \CN_1 &\to \CN_2, \\
\CN_1 \times_{\CN_0} \CN_2 &\to \CN_2. 
\end{align*}
These whisker functors are given by horizontal fusion with a vacuum sector (identity 2-morphism): for instance the right whisker is 
$$({}_{({}_{\cala} (D_1)_{\calb})} H_{({}_{\cala} (D_2)_{\calb})}, {}_{\calb} E_{\calc}) \mapsto H \boxtimes_\calb H_0(E),$$ where $H_0(E)$ is the vacuum sector of the defect $E$.  These various composition operations, and others, are presented in Section~\ref{sec-1data}.

Compatibility transformations encode various relationships among the composition operations, for instance the associativity of vertical composition of 2-morphisms [D2-3], the interaction of the horizontal whiskering operation and vertical composition [D2-6, D2-7], and the associativity of the horizontal whiskering [D2-9, D2-10, D2-11] and of the fusion of defects [D2-12].  The most important compatibility transformation is the switch transformation [D2-8], which is a (symmetric monoidal) natural transformation
$$\CN_2 \times_{\CN_0} \CN_2 \tworarrow \CN_2,$$ between the two functors
\begin{align*}
\left({}_{({}_{\cala} (D_1)_{\calb})} H_{({}_{\cala} (D_2)_{\calb})}, {}_{({}_{\calb} (E_1)_{\calc})} K_{({}_{\calb} (E_2)_{\calc})}\right)
&\mapsto
\left(H \boxtimes_{\calb} H_0(E_1)\right)
\boxtimes_{(D_2 \circledast_\calb E_1)}
\left(H_0(D_2) \boxtimes_{\calb} K)\right), \\
&\mapsto \left(H_0(D_1) \boxtimes_{\calb} K\right)
\boxtimes_{(D_1 \circledast_\calb E_2)}
\left(H \boxtimes_{\calb} H_0(E_2)\right).
\end{align*}
These and other compatibility transformations are constructed in Section~\ref{sec-2data}.

The compatibility transformations are subject to various coherence axioms, for instance pentagon conditions for vertical composition of 2-morphisms [D3-4], and for horizontal whiskering [D3-15, D3-16] and horizontal fusion of defects [D3-17], along with conditions governing the interaction of vertical associativity with horizontal whiskering [D3-7].  Crucial coherence conditions are the one controlling the interaction of the switch transformation with vertical composition [D3-8], and the ones controlling the interaction of the switch operation with horizontal whiskering [D3-13, D3-14].  These and other conditions are proven in Section~\ref{sec-3axioms}.

\setcounter{section}{-1}

\section{Nets, defects, and sectors} \label{sec-0data} %{$0$-data}
%\addtocontents{toc}{\setcounter{tocdepth}{1} }

%\addtocontents{toc}{\SkipTocEntry }
\addtocontents{toc}{\SkipTocEntry }
\subsection*{[0-0] Conformal nets}

   By an \emph{interval}, we shall mean a smooth oriented $1$-manifold that is 
   diffeomorphic to the standard interval $[0, 1]$.
   We let $\INT$ denote the category whose objects are intervals and 
   whose morphisms are embeddings (not necessarily orientation-preserving
   and not necessarily boundary preserving).
   Let $\VN$ be the category whose objects are von Neumann algebras, 
   and whose morphisms are $\IC$-linear $*$-homomorphisms, and $\IC$-linear 
   $*$-antihomomorphisms. 
   A \emph{net} is a covariant functor 
   \[ 
      \cala \colon \INT \to \VN 
   \] 
   taking orientation-preserving embeddings to homomorphisms and 
   orientation-revers\-ing embeddings to antihomomorphisms. 
   It is said to be \emph{isotonic} if the induced maps for embeddings are 
   injective. 
   In this case, given a subinterval $I \subseteq K$, 
   we will often not distinguish between $\cala(I)$ and its image in $\cala(K)$.
   A \emph{conformal} net $\cala$ is an isotonic net subject to a number of 
   axioms~\cite[Def 1.1]{BDH(nets)}.
   Conformal nets form a symmetric monodial category, whose morphisms
   are natural transformations and whose tensor product is the tensor product of 
   $\VN$ applied objectwise.
   There is also the operation of direct sum of conformal nets;
   it is also defined objectwise. 
   A conformal net $\cala$ is said to be \emph{irreducible} 
   if every algebra $\cala(I)$ is a  factor. 
   A direct sum of finitely many irreducible conformal nets is 
   called \emph{semisimple}.
   There is a 
   notion of a \emph{finite} semisimple conformal 
   net (direct sum of conformal nets with finite $\mu$-index \cite{Kawahigashi-Longo-Mueger(2001multi-interval)}), defined utilizing the minimal index from subfactor theory~\cite[Sec 3]{BDH(nets)}. 
   The object category $\CN_0$ of our $3$-category $\CN$ is the 
   subcategory of the category of conformal nets whose objects are
   finite semisimple conformal nets and whose morphisms are
   natural isomorphisms. 
   This subcategory is a symmetric monoidal subcategory~\cite[Sec 3]{BDH(nets)}.  
   From now on all nets will be finite and semisimple 
   and will be simply referred to as conformal nets.
  
%\addtocontents{toc}{\SkipTocEntry }
\addtocontents{toc}{\SkipTocEntry }
\subsection*{[0-1] Defects}

  A \emph{bicolored interval} is an interval $I$ (always oriented), 
  equipped with a covering by two closed, connected, possibly empty subsets 
  $I_{\circ}, I_{\bullet}\subset I$ with disjoint interiors, 
  along with a local coordinate in the neighborhood of $I_{\circ} \cap I_{\bullet}$.
  We disallow the cases where $I_{\circ}$ or $I_{\bullet}$ consists of a single point.
  The local coordinate does {\it not} need to preserve the orientation, 
  but is required to send $(-\e,0]$ into $I_{\circ}$ and $[0,\e)$ into $I_{\bullet}$.
  If either $I_\circ$ or $I_\bullet$ is empty, then there is no local coordinate specified.
  An embedding $f \colon J\hookrightarrow I$ of bicolored intervals is 
  called color preserving if $f^{-1}(I_{\circ})= J_{\circ}$ and 
  $f^{-1}(I_{\bullet}) = J_{\bullet}$.
  The bicolored intervals form a category $\INT_{\circ\bullet}$, 
  whose morphisms are the color preserving embeddings that preserve the 
  local coordinate.
  Let $\cala$ and $\calb$ be conformal nets.
  A \emph{defect} from $\cala$ to $\calb$ is a functor
  \[
  D \colon \INT_{\circ\bullet} \to \VN
  \] that extends $\cala$ and $\calb$
  in the following sense: if $I = I_\circ$ then $D(I) = \cala(I)$;
  if $I = I_\bullet$ then $D(I) = \calb(I)$.
  Moreover, $D$ is subject to a number of axioms~\cite[Def 1.7]{BDH(1*1)}.  
  We often say $D$ is an $\cala$-$\calb$-defect and write $D = {}_\cala D_\calb$.
  Direct sum and tensor product for defects can be defined objectwise, 
  as for nets.
  As morphisms between defects we have again natural transformations. 
  Such a natural transformation ${}_\cala D_\calb \to {}_{\cala'}{D'}_{\calb'}$
  restricts to natural transformations $\cala \to \cala'$ and
  $\calb \to \calb'$.
%  Similarly to nets, there are \emph{semisimple} and \emph{finite} semisimple defects. % and \emph{finite} natural transformation between defects.
  The $1$-morphism category of our $3$-category $\CN$ is the 
  symmetric monodial category $\CN_1$ whose objects are 
  %finite semisimple 
  defects between finite semisimple nets,
  and whose morphisms are natural isomorphisms. %finite natural transformations. 
%  From now on all defects will be finite and semisimple. 
%  All morphism between defects will be finite natural transformations.  
  There are forgetful source and target functors 
  $s, t \colon \CN_1 \to \CN_0$ defined by $s({}_\cala D_\calb) = \cala$
  and $t({}_\cala D_\calb) = \calb$. 
 
  \begin{proposition*}
    The symmetric monoidal functor $s \x t \colon \CN_1 \to \CN_0 \x \CN_0$ is a 
    fibration in the sense of~\cite[Def 2.1, Def 2.2]{Douglas-Henriques(Internal-bicategories)}.
  \end{proposition*}

  \begin{proof}
  Observe as follows that the underlying (non-monoidal) functor $s \x t$ is a fibration of categories.  Given finite semisimple conformal nets $\cala$, $\calb$, $\cala'$, $\calb'$, natural isomorphisms $\phi: \cala \to \cala'$ and $\psi:\calb \to \calb'$, and a finite semisimple defect ${}_{\cala'} D'_{\calb'}$, we must construct a defect ${}_{\cala} D_{\calb}$ and a natural isomorphism ${}_\cala D_\calb \to {}_{\cala'}{D'}_{\calb'}$.  We may take $D(I) = \cala(I)$ when $I$ is white, $D(I) = \calb(I)$ when $I$ is black, and $D(I) = D'(I)$ when $I$ is genuinely bicolored, together with $D(I \hookrightarrow J) = D'(I \hookrightarrow J) \circ \phi(I)$ when $I$ is white and $J$ is genuinely bicolored, and $D(I \hookrightarrow J) = D'(I \hookrightarrow J) \circ \psi(I)$ when $I$ is black and $J$ is genuinely bicolored.  The isomorphism $D \ra D'$ is the identity on genuinely bicolored intervals and is $\phi$, respectively $\psi$ on white and black intervals.  That $s \x t$ is in fact a fibration of symmetric monoidal categories is similarly straightforward.
  \end{proof}

  Throughout this paper we will depend heavily on graphical notation.
  Defects will often be represented by a picture, thought of as a bicolored interval, as follows:
  \[
  \phantom{.}
  \cb{
  \displaytikz{
  \draw (0,0) -- (3,0); \bull(0,0);\bull(1,0);\bull(2,0);
                         \bull(3,0);\dfct(1.5,0);
  } 
  } .
  \] 
  The four bullets on this interval indicate that this
  interval is of length three.
  The marked point $\mathrm{x}$ denotes the point where the color of the
  interval changes.
  We often call this marked point the defect point.
  Strictly speaking we should include an orientation of
  our interval, for example from left to right.
  (Later our intervals will often sit on the boundary
  of $2$-manifolds embedded in the plane. 
  Such a $2$-manifold inherits its orientation from the plane
  and the interval from the boundary of the $2$-manifold.)  
  For a defect ${}_\cala D_\calb$ we think of the above interval
  as representing a collection of von Neumann algebras indexed
  by subintervals of our interval.
  If $I$ is a subinterval to 
  the left of $\mathrm{x}$, then it represents the algebra $\cala(I) = D(I)$; 
  if $I$ is a subinterval containing $\mathrm{x}$, then it represents
  $D(I)$; if $I$ is a subinterval to
  the right of $\mathrm{x}$, then it represents  $\calb(I) = D(I)$.
  Sometimes we will simplify our graphical notation and drop the 
  marked point from the interval.
  If we need coordinates on the above interval, then we will
  identify it with $[0,3]$, where $0$ corresponds to the left
  boundary point and $3$ to the right boundary point.
  The defect point then has the coordinate $1.5$. 

%\addtocontents{toc}{\SkipTocEntry }
\addtocontents{toc}{\SkipTocEntry }
\subsection*{[0-2] Sectors}

  Consider the standard regular hexagon with side length 1:
  \[\cb{
  \displaytikz{
  \sh(0,0);
  \hexb(1,1)1; \dfct(0.5,1); \dfct(0.5,0);
  } 
  }  
  \]
  In this paper $S^1$ is \emph{defined} to be the 
  boundary of this hexagon.
  Sometimes we emphasize this and write $S_6^1$ for the 
  boundary of  this hexagon.
  Later we will also need the regular octagon with side length 1;
  we denote its boundary as $S^1_8$.
  The hexagon inherits an orientation from the plane; this also orients
  its boundary.
  We will pick the clockwise orientation of the plane; thus the circle is
  also clockwise oriented. %, but this will not be of any importance.
  We think of the above circle as bicolored: 
  the left hand side is white $\circ$, while
  the right hand side is black $\bullet$.
  The two marked points are the points where the color changes.
  In particular, every subinterval of $S^1_6$ that contains at most
  one of the marked points (and none on the boundary) inherits a bicoloring.  
  
  Let ${}_\cala D_\calb$ and ${}_\cala E_\calb$ be defects.
  A $D$-$E$-sector is a Hilbert space $H$ equipped with actions of
  algebras represented by the bicolored subintervals of $S^1_6$ as follows:
  for every white subinterval $I \subseteq S^1_6$, the algebra $\cala(I)$
  acts; for every black subinterval $I \subseteq S^1_6$, the algebra
  $\calb(I)$ acts; for every bicolored subinterval containing
  the upper defect point, the algebra $D(I)$ acts;
  for every bicolored subinterval containing
  the lower defect point, the algebra $E(I)$ acts.
  These actions are subject to compatibility axioms~\cite[Def 2.2]{BDH(1*1)}.
  We often write $H = {}_D H_E$ to emphasize that $H$ is a $D$-$E$
  sector.
  For fixed $D$ and $E$, the $D$-$E$-sectors form a category whose morphisms
  are the bounded linear maps that commute with the actions associated
  to the bicolored subintervals of $S^1_6$.  
  There is also a natural notion of morphism 
  $f \colon {}_D H_E \to {}_{D'} {H'}_{E'}$.
  In this case $f$ includes morphisms $\cala \to \cala'$, $\calb \to \calb'$, 
  ${}_\cala D_\calb \to {}_{\cala'} {D'}_{\calb'}$, ${}_\cala E_{\calb'} \to {}_{\cala'} {E'}_{\calb'}$ and an operator $T \colon H \to H'$ that commutes with the induced maps $D(I) \to D'(I)$ and
  $E(I) \to E'(I)$.
  The tensor product of Hilbert space yields a symmetric monodial
  structure on sectors.
  Thus we obtain the symmetric monodial
  category $\CN_2$ of sectors~\cite[Def 2.7]{BDH(1*1)}. 
  This is the $2$-morphism category of our $3$-category $\CN$.
 
  In the graphical notation we think of a sector ${}_D H_E$ as represented by
  the above hexagon.
  We then think of the upper defect point as the $D$ point and the
  lower defect point as the $E$ point.
  By definition every bicolored subinterval of $S^1_6$ corresponds then 
  to a von Neumann algebra that acts on the sector $H$.
  (Later, other $2$-manifolds will also be thought of as representing
  certain sectors.)   
  Often we will drop the marked points from the picture.
  Moreover, we will often draw the hexagon in a rectilinear
  fashion, for example as one of the following:
  \[
  \phantom{.} \quad
  \displaytikz{
  \sh(0,0);
  \bxb(1,1)1;
  \dfct(.5,0);\dfct(.5,1);
  \sh(2,0);
  \bxbur(1,1)1;
  \dfct(2.5,0);\dfct(2.75,1);
  \sh(4,0);
  \bxbul(1,1)1;
  \dfct(4.5,0);\dfct(4.25,1);
  \sh(6,0);
  \bxbdr(1,1)1;
  \dfct(6.75,0);\dfct(6.5,1);
  \sh(8,0);
  \bxbdl(1,1)1;
  \dfct(8.25,0);\dfct(8.5,1);
  \sh(10,0);
  \bxbt(1,1)1;
  \dfct(10.5,0);\dfct(10.5,1);
  \sh(12,0);
  \bxbb(1,1)1;
  \dfct(12.5,0);\dfct(12.5,1);
  } \quad .
  \] 
  Despite their appearance, all these pictures refer to 
  \emph{the} standard regular hexagon  with its 
  two marked points  as drawn.

  Later it will sometimes be convenient to have coordinates 
  on $S^1_6$ and $S^1_8$. 
  Then we will identify $S^1 \cong \IR / 6 \IZ$ such that
  the coordinates of the corner points are $0, 1, \dots, 5$,
  where we start on the left and proceed clockwise from there.
  The coordinate of the upper defect point is then $1.5$, and
  the coordinate of the lower defect point is $-1.5 \equiv 4.5$.
  In a similar fashion we will identify $S^1_8$ with $\IR / 8\IZ$. 

  \begin{proposition*}
    The symmetric monoidal functor $s \x t \colon \CN_2 \to \CN_1 \x_{\CN_0 \x \CN_0} \CN_1$ is a 
    fibration in the sense of~\cite[Def 2.1, Def 2.2]{Douglas-Henriques(Internal-bicategories)}.
  \end{proposition*}

  \begin{proof}
  Observe as follows that the underlying (non-monoidal) functor $s \x t$ is a fibration of categories.  Given finite semisimple conformal nets $\cala$, $\calb$, $\cala'$, $\calb'$, finite semisimple defects ${}_\cala D_\calb$, ${}_\cala E_\calb$, ${}_{\cala'} D'_{\calb'}$, ${}_{\cala'} E'_{\calb'}$, natural isomorphisms $\phi : D \to D'$ and $\psi: E \to E'$, and a sector ${}_{D'} H'_{E'}$, we must construct a sector ${}_D H_E$ and an isomorphism ${}_D H_E \to {}_{D'} H'_{E'}$.  We may take the Hilbert space $H$ to be $H'$, together with $d \in D(I)$ acting by $\phi(d) \in D'(I)$ when $I$ does not contain the lower defect point and with $e \in E(I)$ acting by $\psi(e) \in E'(I)$ when $I$ does not contain the upper defect point.  The isomorphism ${}_D H_E \to {}_{D'} H'_{E'}$ is the identity on the Hilbert space, and is $\phi$ and $\psi$ on $D$ and $E$ respectively.  That $s \x t$ is in fact a fibration of symmetric monoidal categories is similarly straightforward.
  \end{proof}

%\addtocontents{toc}{\SkipTocEntry }
\addtocontents{toc}{\SkipTocEntry }
\subsection*{Implementation of diffeomorphisms}

  Let $\cala$ be a net, $I$ an interval, and $\varphi \colon I \to I$
  a diffeomorphism that is the identity in a neighborhood of
  the boundary of $I$.  To this diffeomorphism, there is an associated automorphisms $\cala(\varphi)$ of $\cala(I)$,
  and it is one of the requirements for conformal nets that this
  automorphism is inner: there is a unitary $U_\varphi \in \cala(I)$
  such that $\cala(\varphi)(a) = U_\varphi a {U_{\varphi}}^*$.
  We then say that $U_\varphi$ implements $\varphi$ on $\cala(I)$.
  Of course, $U_\varphi$ is not unique.
  If $\alpha \colon \cala \to \calb$ is a morphism of $\CN_0$,
  then $\alpha(I)(U_\varphi)$ is an implementation
  of $\varphi$ on $\calb(I)$.

  Let ${}_\cala D_\calb$ and ${}_\cala E_\calb$ be defects and let
  ${}_D H_E$ be a sector.
  Let $\varphi \colon S^1_6 \to S^1_6$ be a diffeomorphism that
  fixes a neighborhood of both defect points.
  We can then pick subintervals $I_L$ of the left half of $S^1_6$
  and $I_R$ of the right half of the circle $S^1_6$ such that $\varphi$
  is the identity on a neighborhood of the complement of
  $I_L \cup I_R$.
  In particular, $\varphi$ restricts to diffeomorphisms
  $\varphi_L$ and $\varphi_R$ of $I_L$ and $I_R$. 
  We obtain automorphisms $\cala(\varphi_L)$ of $\cala(I_L)$
  and $\calb(\varphi_R)$ of $\calb(I_R)$.
  A unitary $U \colon H \to H$ is said to \emph{implement} $\varphi$
  if  
  \begin{equation*}
    \cala(\varphi_L)(a) \circ U = U \circ a  \quad \text{and} \quad
    \calb(\varphi_R)(b) \circ U = U \circ b 
  \end{equation*}
  as operators on $H$ for all $a \in \cala(I_L)$, $b \in \calb(I_R)$.
  Such an implementation always exits; for example we can set
  $U := U_L \circ U_R$ where $U_L$ implements $\varphi_L$ on
  $\cala(I_L)$ and $U_R$ implements $\varphi_R$ on $\calb(I_R)$.
  (It is part of the axioms for sectors that the actions of
  $\cala(I_L)$ and $\calb(I_R)$ on $H$ commute; in particular
  $U_L$ is $\calb(I_R)$-linear and $U_R$ is $\cala(I_L)$-linear.)
  
%g%\pagebreak  % This is a hack because it was breaking after the section headings

\section{Composition and identity operations} \label{sec-1data} %{$1$-data}
%\addtocontents{toc}{\setcounter{tocdepth}{2}}

\subsection{Horizontal identity and composition} %{Identity and composition of defects}
     
\subsubsection*{[1-1] Horizontal identity.}
  Let $\cala \colon \INT \to \VN$ be a conformal net.
  Then the identity defect $\id_\cala$ for $\cala$ is defined by
  \[
  \id_\cala = \cala \circ \forget
  \] 
  where $\forget \colon \INT_{\circ\bullet} \to \INT$ is the
  functor that forgets the bicoloring. 
  The $1$-cell identity $\CN_0 \to \CN_1$ is defined to be the functor 
  $\cala \mapsto \id_\cala$.
  We will draw the identity defect as
  \[
  \phantom{\text{[1-1]}} \qquad
  \cb{
  \displaytikz{
  \draw (0,0) -- (3,0); \bull(0,0);\bull(1,0);\bull(2,0);\bull(3,0);\idfct(1.5,0);
  } 
  } \qquad \text{[1-1]}
  \]
  where we use an equal sign (rotated) in the place of the usual
  $x$ at the defect point. 
  Sometimes we simplify this by dropping the defect marker altogether:
  \[
  \phantom{.} \quad
  \cb{
  \displaytikz{
  \draw (0,0) -- (3,0); \bull(0,0);\bull(1,0);\bull(2,0);\bull(3,0);
  } 
  } \quad .
  \]

\subsubsection*{[1-2] Horizontal composition.}
  The horizontal composition is defined as the
  horizontal composition of defects as introduced   
  in~\cite[Sec 1.E]{BDH(1*1)}.
  We write this horizontal composition functor
  $\CN_1 \x_{\CN_0} \CN_1 \to \CN_1$ as 
  $({}_\cala D_\calb, {}_\calb E_\calc) \mapsto D \circledast_\calb E$
  and draw the composite of two defects as
  \[
  \phantom{\text{[1-2]}} \qquad
  \cb{
  \displaytikz{
  \draw (0,0) -- (4,0); 
  \bull(0,0);\bull(1,0);\bull(2,0);\bull(3,0);\bull(4,0);
  \dfct(1.5,0);\dfct(2.5,0);
  } 
  } \qquad \text{[1-2]} 
  .
  \]
  Here the left defect point is associated to to $D$ and the right defect point 
  is associated to $E$. We will review horizontal composition of defects and explain the picture in more detail below in Section~\ref{sec-horizwhiskers}.
  
\subsection{Vertical identity and composition} %{Identity and composition of sectors}

\subsubsection*{[1-3] Vertical identity.}
  Let $D = {}_\cala D_\calb$ be a defect, and let $S$ be a circle with a bicoloring-preserving automorphism that exchanges the two color change points---we refer to such an automorphism as a reflection.  The vacuum sector
  $H_0(D,S) = {}_D H_0(D,S)_D$ of $D$ on $S$ was introduced 
  in~\cite[Sec 1.B]{BDH(1*1)}. 
  If $S$ is the standard circle $S^1_6$, then the reflection along the horizontal axis is a canonical
  choice of a bicoloring-preserving reflection. 
  (In coordinates the reflection is $t \mapsto 6-t$).
  We call $H_0(D) := H_0(D,S^1_6)$ the vacuum sector of $D$. 
  The functor $D \mapsto H_0(D)$ defines
  the $2$-cell identity $\CN_1 \to \CN_2$.
  The underlying Hilbert space of the vacuum sector 
  is the standard form $L^2(D(I))$ of the von Neumann algebra $D(I)$, where $I$ is the upper half of
  the circle $S^1_6$.
  (This is also the interval of length $3$ used earlier.)
  Pictorially we denote the vacuum sector as
  \[
  \phantom{\text{[1-2]} .} \qquad 
  \cb{
  \displaytikz{
  \sh(0,0)
  \hexb(1,1)2;
    \draw(\x,\y+.5) -- +(.15,0) (\x+1,\y+.5) -- +(-.15,0);
	\dfct(0.5,0);\dfct(0.5,1);
  }
  }  \qquad \text{[1-2]} .
  \] 
  In this picture, the gray shading indicates that the sector is a vacuum sector; an arbitrary sector would have no interior shading.
  The upper and the lower half of $S^1_6$ are
  both copies of our bicolored interval $I$ and correspond
  to the two actions of $D(I)$ on $L^2(D(I))$.
  Sometimes we will drop the defect points from
  our pictures.
  Moreover, we might draw the picture in a rectilinear 
  fashion such as
  \[
  \phantom{.} \quad 
  \displaytikz{
  \sh(0,0)
  \bxb(1,1)2;
  } \quad .
  \]
  We point out that whenever the boundary of the circle is split into
  two intervals each of which contains a defect point in the interior,
  then the corresponding algebras are commutants of each other~\cite[Prop 1.16]{BDH(1*1)}. 

\subsubsection*{[1-4] Vertical composition.}
  The vertical composition $\CN_2 \x_{\CN_1} \CN_2 \to \CN_2$ 
  is defined as the vertical fusion from~\cite[Sec 2.C]{BDH(1*1)}.
  Our picture for the vertical fusion is
  \[\phantom{\text{[1-4]}} \qquad
  \cb{
  \displaytikz{
  \sh(0,0)
  \hexb(1,1)1;\draw (0,.5) -- (1,.5);\bull(.33,.5);\bull(.67,.5);
  \dfct(.5,0);\dfct(.5,.5);\dfct(.5,1);
  }
  } \qquad \text{[1-4]}
  .
  \]
  Note that the boundary of this picture is canonically $S^1_6$.
  In particular, no boundary reparametriation is needed in the definition
  of vertical fusion.
  Often the picture is simplified by omitting defect points 
  and is drawn as a rectilinear equivalent
  \[
  \phantom{.} \quad
  \displaytikz{
  \sh(0,0)
  \bxbt(1,1)1;
  \sh(0,1)
  \bxbb(1,1)1;
  } \quad .
  \]
  The underlying Hilbert space of the vertical fusion of sectors 
  is the Connes fusion of the Hilbert spaces for the two sectors
  over the algebra associated to the horizontal interval of
  length 3 in the middle of the pictures. 
  Sometimes will will draw this pictures in the following different, 
  but equivalent, forms
  \[
  \phantom{.} \quad
  \cb{
  \displaytikz{
  \path
  (0,0) coordinate (1)
  (2,0) coordinate (2);
  \draw (1) node[littlenode] {
  \innertikz{
		\sh(0,0);
		\hex(1.25,1)1;
		\sh(0,1);
		\dhkd(1.25,1*\fudge)1,0,1;
		\bull(0,.5);\bull(1.25,.5);\bull(0,1+1*\fudge);\bull(1.25,1+1*\fudge);
		\bull(.42,0); \bull(1.25-.42,0);\bull(.42,1);\bull(1.25-.42,1);%!% These bullets will be wrong is hexangleA changes.
  } 
  };
  \draw (2) node[littlenode] {
  \innertikz{
		\sh(0,0);
		\dhku(1.25,1*\fudge)1,0,1;
		\sh(0,1*\fudge);
		\hex(1.25,1)1;
		\bull(0,0);\bull(1.25,0);\bull(0,.5+\fudge);\bull(1.25,.5+\fudge);
				\bull(.42,\fudge); \bull(1.25-.42,\fudge);\bull(.42,1+\fudge);\bull(1.25-.42,1+\fudge);%!% These bullets will be wrong is hexangleA changes.
  } 
  }; 
  }
  } \quad.
  \]
  These versions will be helpful when we discuss the vertical
  fusion of three sectors.
  
\subsection{Horizontal whiskers} \label{sec-horizwhiskers}%{Whiskering sectors by defects} 
  
\subsubsection*{Horizontal fusion.}
  The definition of a $3$-category that we are using in this 
  paper does not (for reasons of efficiency) directly include a notion of horizontal composition of $2$-morphisms.
  Nevertheless, there is such a composition for our sectors called
  \emph{horizontal fusion} of sectors and this operation will be the basis for many pieces of structure in our $3$-category. Horizontal fusion is a functor
  $\CN_2 \x_{\CN_0} \CN_2 \to \CN_2$ and is defined 
  in~\cite[Sec 2.B]{BDH(1*1)}.
  In symbols, given defects ${}_\cala (D_1)_\calb$, ${}_\cala (D_2)_\calb$, ${}_\calb (E_1)_\calc$, and ${}_\calb (E_2)_\calc$, we will write the horizontal fusion functor as
  \[
    ({}_{D_1} H_{D_2}, {}_{E_1} K_{E_2}) \mapsto 
     ({}_{D_1 \circledast_\calb E_1} 
                   (H \boxtimes_{\calb} K)_{D_2 \circledast_\calb E_2})
        .
  \] 
  We draw $H \boxtimes_\calb K$ as 
  \[
  \phantom{.}
  \cb{
  \displaytikz{
  \sh(0,0);
  \halfhexl(1,1)1;
  \sh(1,0);
  \halfhexr(1,1)1;
  \bull(.42,0);\bull(1.58,0);\bull(.42,1);\bull(1.58,1);
  \bull(0,.5);\bull(2,.5);\bull(1,0);\bull(1,.5);\bull(1,1);
  \dfct(0.71,0);\dfct(1.29,0);\dfct(0.71,1);\dfct(1.29,1); %!% These will be in the wrong place if hexangleC changes.
  }
  }. 
  \]
  The underlying Hilbert space is the Connes fusion of $H$ and 
  $K$ along the algebra associated by $\calb$ to the vertical 
  interval $I$ of length $2$ in the middle of the picture. 
  Note that $I$ inherits two different orientations
  from the two (deformed) hexagons.
  If we orient $I$ using the right hexagon (corresponding to $K$), 
  then $\calb(I)$ acts on $K$, while $\calb(-I) = \calb(I)^\op$ acts
  on $H$.
  It is exactly this situation that allows the use of
  Connes fusion (just as one may take the tensor product of a right module and a left module).
  Sometimes we drop  defect points from the picture.
  Moreover, we often draw a rectilinear version of the picture,
  \[
  \phantom{.} \quad
  \displaytikz{
  \sh(0,0);
  \bxb(1,1)1;
  \sh(1,0);
  \bxb(1,1)1;
  \dfct(0.5,0);\dfct(0.5,1);\dfct(1.5,0);\dfct(1.5,1);
  } \quad .
  \]
  We can now give a brief summary of the composition of defects. 
  Let ${}_\cala D_\calb$ and ${}_\calb E_\calc$ be defects.
  Consider $H_0(D) \boxtimes_\calb H_0(E)$, 
  the horizontal composition of the vacuum sectors
  for $D$ and $E$,
  \[
  \phantom{.} \quad
  \cb{
  \displaytikz{
  \sh(0,0);
  \halfhexl(1,1)2;
  \sh(1,0);
  \halfhexr(1,1)2;
  \bull(.42,0);\bull(1.58,0);\bull(.42,1);\bull(1.58,1);
  \bull(0,.5);\bull(2,.5);\bull(1,0);\bull(1,.5);\bull(1,1);
    \dfct(0.71,0);\dfct(1.29,0);\dfct(0.71,1);\dfct(1.29,1); %!% These will be in the wrong place if hexangleC changes.
  }
  } \quad. 
  \]
  The boundary of this picture is drawn as an irregular hexagon,
  but its boundary has length 8.
  Thus it can be identified with the octagon $S^1_8$, and the upper four segments of the boundary 
  can be identified with the interval $I_4$ of length $4$
  \[\phantom{.} \quad
  \cb{
  \displaytikz{
  \draw (0,0) -- (4,0); 
  \bull(0,0);\bull(1,0);\bull(2,0);\bull(3,0);\bull(4,0);
  \dfct(1.5,0);\dfct(2.5,0);
  } 
  } \quad . 
  \]
  The evaluation of the composed defect 
  $D \circledast_\calb E$ on this interval is generated in the
  algebra of bounded operators on $H_0(D) \boxtimes_\calb H_0(E)$ by
  the evaluation of $D$ on the first two segments and by the evaluation
  of $E$ on the last two segments. 
  Similarly, we obtain an algebra acting on 
  $H_0(D) \boxtimes_\calb H_0(E)$ for any subinterval of $I_4$.
  The  interval $I_4$ is not bicolored, but there is a
  map onto a bicolored interval $I_3$ of length 3 that collapses 
  the two half segments between
  the two defect points to a single defect point, and this collapse map is
  used to view $D \circledast_\calb E$ as a functor on 
  $\INT_{\circ\bullet}$. 
  The evaluation of $D \circledast_\calb E$ on a subinterval of $I_3$ 
  is defined via its preimage in $I_4$ under the collapse map. 
  In our pictures we never indicate this collapse map in any way. 
  Thus the pictures remember
  more than just the structure of $D \circledast_\calb E$ as a defect:
  we see more subintervals to which we can associate algebras, for example 
  we could consider a little neighborhood of the left (say) defect point.
  In a similar fashion our picture for the horizontal fusion
  remembers more than just the structure of a sector; it also encodes
  the actions of some additional algebras.
  If we compose more than two defects, then we obtain intervals of yet longer 
  length with yet more defect points.
  
  To formally define the horizontal fusion of sectors 
  a similar collapse
  map $\pi \colon S^1_8 \to S^1_6$ is used.
  It collapses four half segments to two points.  On the upper half this is just the collapse map used before, on 
  the lower half this is the reflection of that collapse map.
 
\subsubsection*{[1-5] Right whisker.}
  The right composition of a $1$-cell with a $2$-cell 
  $\CN_2 \x_{\CN_0} \CN_1 \to \CN_2$ is defined using horizontal
  fusion and the vacuum sector.
  Let $\cala$, $\calb$, and $\calc$ be nets, let ${}_\cala (D_1)_\calb$, ${}_\cala (D_2)_\calb$, and 
  ${}_\calb E_{\calc}$ be defects, and let ${}_{D_1} H_{D_2}$ be a sector. 
  The right composition  of $H$ with $E$ is defined as the
  horizontal fusion $H \boxtimes_\calb H_0(E)$.
  We draw this as
  \[
  \phantom{\text{[1-5]}} \qquad
  \cb{
  \displaytikz{
  \sh(0,0);
  \halfhexl(1,1)1;
  \sh(1,0);
  \halfhexr(1,1)2;
  \bull(.42,0);\bull(1.58,0);\bull(.42,1);\bull(1.58,1);
  \bull(0,.5);\bull(2,.5);\bull(1,0);\bull(1,.5);\bull(1,1);
  }
  } \qquad \text{[1-5]}  .
  \]
  Again,  often this is drawn as
  \[
  \phantom{.} \quad
  \displaytikz{
  \sh(0,0);
  \bxb(1,1)1;
  \sh(1,0);
  \bxb(1,1)2;
  } \quad .
  \]
  Here we omitted defect points from the pictures, but sometimes we will
  include them for clarity.    

\subsubsection*{[1-6] Left whisker.}
  The left composition of a $1$-cell with a $2$-cell is
  defined similarly to the right composition, 
  and is drawn as
  \[
  \phantom{\text{[1-6]}} \qquad
  \cb{
  \displaytikz{
  \sh(0,0);
  \halfhexl(1,1)2;
  \sh(1,0);
  \halfhexr(1,1)1;
  \bull(.42,0);\bull(1.58,0);\bull(.42,1);\bull(1.58,1);\bull(0,.5);
  \bull(2,.5);\bull(1,0);\bull(1,.5);\bull(1,1);
  } \qquad  \qquad
  \displaytikz{
  \sh(0,0);
  \bxb(1,1)2;
  \sh(1,0);
  \bxb(1,1)1;
  }
  } \qquad \text{[1-6]}
  .
  \] 
  
  \medskip \medskip \medskip
  
  The data discussed so far is part of the definition of both
  a $2$-category object and a  dicategory object
  (in the $2$-category of symmetric monoidal
  categories). 
  We remind the reader that our $3$-category of conformal nets
  is a dicategory object; the next two pieces of data labeled here [1-7]
  and [1-8] are  labeled [D1-7] and [D1-8] in the appendix and correspondingly in~\cite{Douglas-Henriques(Internal-bicategories)}.
  
%g%\pagebreak % hack bad  
\subsection{Directional identity cells} %{Half-identity sectors}
\subsubsection*{[1-7] Left identity.}
  The (upper) left $2$-cell identity is a functor $\CN_1 \to \CN_2$.  Its role is to show that the composition of a defect with an identity defect is, at least in a weak sense, equivalent to the original defect.

  Let  ${}_\cala D_\calb$ be a defect.
  There is no canonical isomorphism between 
  $\id_\cala \circledast_\cala D$ and $D$ in $\CN_1$.
  (For this reason our $3$-category of conformal nets
  is  not a $2$-category object in symmetric monoidal categories.)
  There is however a canonical 
  $(\id_\cala \circledast_\cala D)$-$D$-sector, the left identity for
  $D$.
  Our picture of this left identity is
  %%% Here is the length two version
  %\[\cb{
  %\displaytikz{
  %\sh(0,0);
  %\boxt(.25,1,.75,1)2;
  %\sh(0,1);
  %\bubble(.25,1,.75)2;
  %\begin{pgfonlayer}{foreground}
  %	\idfct(.5,1+\cw);\bull(1,1+\cw);\dfct(1.5,1);\dfct(1,0);\bull(0,.5);\bull(2,.5);
  %\end{pgfonlayer}
  %}
  %} \qquad \text{[1-7]}
  %\]
  %%% Here is the length 1.5 version
  \[
  \phantom{ \text{[1-7]}} \qquad 
  \cb{
  \displaytikz{
  \sh(0,0);
  \boxt(.25*.75,.75,.75*.75,1)2;
  \sh(0,1);
  \bubble(.25*.75,.75,.75*.75)2;
  \begin{pgfonlayer}{foreground}
	\idfct(.375,1+\cw);\bull(.75,1+\cw);\dfct(1.125,1);\dfct(.75,0);\bull(0,.5);\bull(1.5,.5);
  \end{pgfonlayer}
  }
  } \qquad \text{[1-7]}.
  \]
  This sector is the vacuum sector $H_0(D)$ for $D$ (this is the box 
  part of the picture), twisted by  a diffeomorphism 
  (indicated by the balloon in the picture).
  Details of this construction follow. 
 
  We begin by reviewing the defect $\id_\cala \circledast_\cala D$.
  Consider the collapse map $\pi \colon S^1_8 \to S^1_6$ used earlier.
  This map is symmetric with respect to the reflection along the horizontal
  axis (in coordinates the reflection is given by $x \mapsto -x$).
  The restriction of $\pi$ to the upper half of $S^1_8$ is a map $I_4 \to I_3$, which
  collapses $[1.5,2.5] \subset I_4$ to  $1.5 \in I_3$,
  sends $x \in [0,1.5] \subset I_4$ to $x \in [0,1.5]$,
  and sends $x \in [2.5,4] \subset I_4$ to  
  $x-1 \in [1.5,3] \subset I_3$.  
  The evaluation of $\id_\cala \circledast_\cala D$
  on a bicolored subinterval $I$ of $I_3$ 
  is the algebra $D(\pi^{-1}(I))$.
  (This algebra is isomorphic to $D(I)$, but there is no canonical
  isomorphism if $I$ contains the upper defect point $1.5$.)
  
  To construct the left identity we pick a 
  diffeomorphism $\Phi_L \colon S^1_8 \to S^1_6$
  such that $\Phi_L(x) = x$ on a neighborhood of $[0,1]$
  and $\Phi_L(x) = x - 1$ on a neighborhood of $[2.5,4]$
  and is symmetric with respect to the vertical reflection of the circles; 
  that is, we require that $\Phi_L( -x ) = - \Phi_L(x)$. 
  In particular, $\Phi_L$ coincides with $\pi$ on $[0,1] \cup [2.5,4]$.
  Now we start with the vacuum sector $H_0(D)$ for $D$:
\[
\phantom{.} \quad 
\cb{
\displaytikz{
\sh(0,0);
\bxb(1.5,1)2;
\begin{pgfonlayer}{foreground}
	\dfct(.75,1);\dfct(.75,0);
\end{pgfonlayer}
}
} \quad.
\]
  It is a $D$-$D$-sector.
  We can twist the upper $D$-action by the restriction of $\Phi_L$ to
  the upper half $I_3$ of $S^1_6$,
  turning $H_0(D)$ into an $(\id_\cala \circledast_\cala D)$-$D$-sector.
  More precisely, if $I \subseteq S^1_6$ is a bicolored interval
  containing the upper defect point,
  then the action of
  $(\id_\cala \circledast_\cala D)(I)$ on $H_0(D)$ is defined via
  the isomorphism
  \[
  (\id_\cala \circledast_\cala D)(I) = D(\pi^{-1}(I)) 
       \xrightarrow{D(\Phi_L|_{\pi^{-1}(I)})} D(\Phi_L(\pi^{-1}(I))).
  \]
  Because 
  $\Phi_L(x) = x$ on a neighborhood of $[0,1]$
  and $\Phi_L(x) = x - 1$ on a neighborhood of $[2.5,4]$ 
  it follows that this construction indeed defines a sector. 
  We define the left identity for $D$ to be this sector.
  (The left identity functor $\CN_1 \to \CN_2$ constructed in this section, and similarly the right identity functor in the next section, depends on the choice of the diffeomorphism $\Phi_L$.  However, as discussed in the later Remark~\ref{remark-equivalence}, distinct choices here will result in isomorphic 3-categories.)
  
  The small balloon in the above picture  
  represents the restriction of $\Phi_L$ to $[1,2.5] \to [1,1.5]$. 
  Occasionally we use the abbreviated notation
  \[
  \displaytikz{
  \sh(0,0);
  \bxbul(1.5,1)2;\bull(1.5,.5);
  \halfidl(1.5,1){.2};
  }
  \]
  or
  \[
  \displaytikz{
  \sh(0,0);
  \sbx(1.5,1)2;
  \halfidl(1.5,1){.2};
  \tick(.75,1);
  }
  \]
  in which a small vertical tick indicates that the top is a composition 
  of two defects, or, when there could be no confusion, simply by
  \[
  \phantom{.} \quad
  \displaytikz{
  \sh(0,0);
  \sbx(1.5,1)2;
  \halfidl(1.5,1){.2};
  } \quad .
  \]   

\subsubsection*{[1-8] Right identity.}
  The right identity is defined similarly to the left identity; specifically
  the right identity is a horizontal reflection of
  the left identity.
  Thus we replace the diffeomorphism $\Phi_L$ by the diffeomorphism
  $\Phi_R \colon S^1_8 \to S^1_6$, $\Phi_R(x) := 3 - \Phi_L(4-x)$.
  The pictures for the right identity are: 
\[
\phantom{\text{[1-8]}} \qquad
\cb{
\displaytikz{
\sh(0,0);
\boxt(.75*.75,.75,.25*.75,1)2;
\sh(0,1);
\bubble(.75*.75,.75,.25*.75)2;
\begin{pgfonlayer}{foreground}
	\dfct(.375,1);\bull(.75,1+\cw);\idfct(1.125,1+\cw);\dfct(.75,0);\bull(0,.5);\bull(1.5,.5);
\end{pgfonlayer}
}

  \qquad
  \displaytikz{
  \sh(0,0);
  \bxbur(1.5,1)2;\bull(0,.5);
  \halfidr(1.5,1){.8};
  }
  \qquad
  \displaytikz{
  \sh(0,0);
  \sbx(1.5,1)2;
  \halfidr(1.5,1){.8};
  \tick(.75,1);
  }
  \qquad
  \displaytikz{
  \sh(0,0);
  \sbx(1.5,1)2;
  \halfidr(1.5,1){.8};
  }
  } \qquad \text{[1-8]}.
  \]

  \begin{lemma*}
    The left and the right identity sectors are invertible with respect to vertical fusion of sectors, as required
    in the definition of a dicategory object. %~\cite{Douglas-Henriques(Internal-bicategories)}. 
  \end{lemma*}
  
  \begin{proof}
    An inverse for the left identity is given by a vertical reflection
    of the left identity.
    Similarly, an inverse for the right identity is given by
    a vertical reflection of the right identity. 
  \end{proof}

  The procedure of twisting with a diffeomorphism 
  as in the construction of the left identity [1-7]  can be applied to other
  defects than the vacuum sector.
  We can twist any ${}_\cala D_\calb$-${}_\cala E_\calb$-sector by a diffeomorphism 
  to obtain an $(\id_\cala \circledast_\cala D)$-$E$-sector.
  Varying the position of the diffeomorphism we can also produce a 
  $(D \circledast_\calb \id_\calb)$-$E$-sector or a
  $D$-$(\id_\cala \circledast_\cala E)$-sector or a
  $D$-$(E \circledast_\calb \id_\calb)$-sector.
  Moreover this process can be reversed.
  For example given a $(\id_\cala \circledast_\cala D)$-$E$-sector
  we can twist by the inverse of $\Phi_L$ to obtain
  a $D$-$E$-sector.
  These constructions are inverse to each other.
  Also note that in vertical compositions we will often move the diffeomorphism from 
  one sector to another, when this does not affect the resulting composite sector; for example the following pictures are interchangeable:
\[\phantom{.} \quad
 \cb{
  \displaytikz{
  \path
  (0,0) coordinate (1)
  (3,0) coordinate (2);
  \draw (1.2,.04) -- (1.8,.04); \draw (1.2,-.04) -- (1.8,-.04);
  \draw (1) node[littlenode] {
  \displaytikz{
		\sh(0,0);
		\boxt(.625,.5,.875,1)2;
		\sh(0,1);
		\bubble(.625,.5,.875)2;
		\sh(0,1+\gap);
		\boxslopedr(.625,.375,1)2;
		\sh(1,1+\gap);
		\boxslopedl(.875,.125,1)1;
		\begin{pgfonlayer}{foreground}
			\bull(.5,1);\bull(1.5,1);\bull(.5,1+\gap);
                        \bull(1,1+\cw);
			\bull(1.5,1+\gap);
                        \bull(1,1.5+\gap);
                        \idfct(.5,2+\gap);\idfct(.75,1+\gap+\cw);\idfct(.75,1+\cw);
                        \dfct(1,0);\dfct(1.25,1);\dfct(1.25,1+\gap);\dfct(1.5,2+\gap);
		\end{pgfonlayer}
  }
  };
  \draw (2) node[littlenode] {
  \displaytikz{
		\sh(0,0);
		\boxt(.625,.5,.875,1)2;
		\sh(0,1+\gap);
		\bubble(.625,.5,.875)2;
		\sh(0,1+\gap);
		\boxslopedr(.625,.375,1)2;
		\sh(1,1+\gap);
		\boxslopedl(.875,.125,1)1;
		\begin{pgfonlayer}{foreground}
			\bull(.5,1);\bull(1.5,1);\bull(.5,1+\gap);
			\bull(1.5,1+\gap);
                        \bull(1,1.5+\gap);
                        \idfct(.5,2+\gap);\idfct(.75,1+\gap+\cw);
                        \dfct(1,0);\dfct(1.25,1);\dfct(1.25,1+\gap);\dfct(1.5,2+\gap);
		\end{pgfonlayer}
  }
  };
}
}\quad .
\]

  %==================================================================

\section{Compatibility transformations for composition and identity operations} \label{sec-2data}%{2-Data}

\subsection{Transformations for vertical identity and composition}

\subsubsection*{[2-1] Top identity}
  There is a canonical natural isomorphism 
  \[
  \phantom{\text{[2-1]}} \quad
  \cb{
  \displaytikz{
  \path
  (0,0) coordinate (1)
  (2,0) coordinate (2);
  \draw (1) -- (2);
  \draw (1) node[littlenode] {
  \innertikzb{
  \sh(0,0);
  \bxbt(1,1)1;
  \sh(0,1);
  \bxbb(1,1)2;
  }
  };
  \draw (2) node[littlenode] {
  \innertikzb{
  \sh(0,0);
  \bxb(1,1)1;
  }
  };
  }
  }  \qquad \text{[2-1]}
  \]
  because the underlying Hilbert space of the identity defect on
  the left hand side is the standard form
  of the algebra associated to the interval of length 3
  in the middle of the picture on the left hand side.  (Here we have denoted the isomorphism between the left and right Hilbert spaces simply as a horizontal dash, without an arrowhead, to indicate that the morphism may be read in either direction.)  This is the top identity.

\subsubsection*{[2-2] Bottom identity}
  The bottom identity is similarly depicted
  \[
  \phantom{\text{[2-2]}.} \qquad
  \cb{
  \displaytikz{
  \path
  (0,0) coordinate (1)
  (2,0) coordinate (2);
  \draw (1) -- (2);
  \draw (1) node[littlenode] {
  \innertikzb{
  \sh(0,0);
  \bxbt(1,1)2;
  \sh(0,1);
  \bxbb(1,1)1;
  }
  };
  \draw (2) node[littlenode] {
  \innertikzb{
  \sh(0,0);
  \bxb(1,1)1;
  }
  };
  }
  } \qquad \text{[2-2]}
  .
  \]

\subsubsection*{[2-3] Vertical associator}
  Connes fusion of bimodules over von Neumann algebras is not strictly
  associative, but there is a coherent associator for this operation
  (similar to the associator for the algebraic tensor product of
  bimodules over rings).
  Because vertical fusion is defined using fusion along 
  the algebra corresponding to the upper, respectitively lower,
  half of our standard circle, the associator for Connes 
  fusion over von Neumann
  algebras is also an associator for vertical fusion of sectors.
  We will draw this associator as
  \[
  \phantom{\text{[2-3]}.} \qquad
  \cb{
  \displaytikz{
  \path
  (0,0) coordinate (1)
  (2,0) coordinate (2);
  \draw (1) -- (2);
  \draw (1) node[littlenode] {
  \innertikzb{
		\sh(0,0);
		\dhku(1.25,1*\fudge)1,1,1;
		\sh(0,1*\fudge+\gap);
		\hex(1.25,1)1;
		\sh(0,1+1*\fudge+\gap);
		\dhkd(1.25,1*\fudge)1,0,1;
  } 
  };
  \draw (2) node[littlenode] {
  \innertikzb{
		\sh(0,0);
		\dhku(1.25,1*\fudge)1,0,1;
		\sh(0,1*\fudge);
		\hex(1.25,1)1;
		\sh(0,1+1*\fudge+\gap);
		\dhkd(1.25,1*\fudge)1,1,1;
  } 
  }; 
  }
  } \qquad \text{[2-3]}.
  \]
  The little gap on the left hand side illustrates that here 
  we first do the Connes fusion along the upper algebra;
  on the right hand side the gap illustrates that we first
  do the Connes fusion along the lower algebra.
  Because this associator just comes from the fact that
  Connes fusion over von Neumann algebras is not strictly associative
  we will henceforth very often suppress this isomorphism and
  treat the right hand and left hand side of the above picture as equal;
  we therefore simply draw this vertical fusion as
  \[\cb{
  \displaytikz{
    \path (0,0) coordinate (1);
  \draw (1) node[littlenode] {
  \innertikzb{
		\sh(0,0);
		\dhku(1.25,1*\fudge)1,0,1;
		\sh(0,1*\fudge);
		\hex(1.25,1)1;
		\sh(0,1+1*\fudge);
		\dhkd(1.25,1*\fudge)1,0,1;
		}
		};
  } 
}
  \]

\subsection{Transformations for horizontal composition and whiskers}  
  
\subsubsection*{$1 \boxtimes 1$-isomorphism}
  Crucial for the construction of our $3$-category is the 
  $1 \boxtimes 1$-isomorphism from~\cite[Thm 6.2]{BDH(1*1)}.
  The $1 \boxtimes 1$-isomorphism provides a natural isomorphism
  between two functors $\CN_1 \x_{\CN_0} \CN_1 \to \CN_2$
  defined as follows.
  Let ${}_\cala D_\calb$ and ${}_\calb E_\calc$ be defects. 
  The first functor sends  $(D,E)$ to 
  $H_0(D) \boxtimes_\calb H_0(E)$ and the second functor
  sends $(D,E)$ to $H_0(D \circledast_\calb E)$.
  Thus the $1 \boxtimes 1$ isomorphism shows in particular that
  the horizontal fusion of two vacuum sectors is again a vacuum
  sector.
  In pictures the $1 \boxtimes 1$-isomorphism is denoted
  \[
  \phantom{.} \quad
  \cb{
  \displaytikz{
  \path
  (0,0) coordinate (1)
  (3,0) coordinate (2);
  \draw (1) -- (2);
  \draw (1) node[littlenode] {  
  \displaytikz{
		\sh(0,0);
		\bxb(1,1)2;
		\sh(1,0);
		\bxb(1,1)2;
		\dfct(.5,0);\dfct(.5,1);\dfct(1.5,0);\dfct(1.5,1);
  }
  };
  \draw (2) node[littlenode] {
  \displaytikz{
		\sh(0,0);
		\bxb(2,1)2;\bull(1,0);\bull(1,1);
				\dfct(.5,0);\dfct(.5,1);\dfct(1.5,0);\dfct(1.5,1);
  }
  };
  }
  }   \quad .
  \]  
  
\subsubsection*{[2-4] and [2-5] Right and left vertical 
       identity expansion.}
  The right and left vertical identity expansions coincide and are both
  given  by the $1 \boxtimes 1$-isomorphism. 
  In pictures we have
  \[
  \phantom{\text{[2-4], [2-5]},} \qquad 
  \cb{
  \displaytikz{
  \path
  (0,0) coordinate (1)
  (3,0) coordinate (2);
  \draw (1) -- (2);
  \draw (1) node[littlenode] {
  \displaytikz{
		\sh(0,0);
		\bxb(1,1)2;
		\sh(1,0);
		\bxb(1,1)2;
		\dfct(.5,0);\dfct(.5,1);\dfct(1.5,0);\dfct(1.5,1);
  }
  };
  \draw (2) node[littlenode] {
  \displaytikz{
		\sh(0,0);
		\bxb(2,1)2;\bull(1,0);\bull(1,1);
		\dfct(.5,0);\dfct(.5,1);\dfct(1.5,0);\dfct(1.5,1);
  }
  };
  }
  } \qquad \text{[2-4], [2-5]},
  \] 
  but often we will drop the defect points from the notation.

\subsubsection*{The categories ${\CN_\sim}_2$, 
           ${\CN^\sim}_2$ and ${\CN^\sim_\sim}_2$.}

  We will later need variants of $\CN_2$ that have more
  morphisms. 
 
  The following  notation will be helpful.
  If $\varphi \colon A \to B$ is a map of von Neumann algebras, 
  $H$ is an $A$-module,
  and $K$ is a $B$-module, then we denote by $\Hom_\varphi(H,K)$ the space of
  all bounded linear maps $T: H \to K$ that are \emph{$\varphi$-linear}, 
  that is such that for all $a \in A$ and
  $\xi \in H$ we have $T(a \xi) = \varphi(a) T(\xi)$.

  We start by recalling the precise definition of the morphisms in $\CN_2$.
  For defects ${}_\cala D_\calb$, ${}_\cala E_\calb$,
  ${}_{\cala'} {D'}_{\calb'}$, and ${}_{\cala'} {E'}_{\calb'}$, 
  and sectors ${}_D H_E$ and ${}_{D'} {H'}_{E'}$ a morphism 
  $f \colon {}_D H_E \to {}_{D'} {H'}_{E'}$ is a triple $f = (T,\delta,\e)$ where 
  $T \colon H \to H'$ is a bounded linear map and
  $\delta \colon D \to D'$ and $\e \colon E \to E'$ are 
  morphisms from $\CN_1$ such that $s(\delta) = s(\e) \colon \cala \to \cala'$ and
  $t(\delta) = t(\e) \colon \calb \to \calb'$.
  Moreover, $T$ is required to be  $\delta(I)$-linear for all
  bicolored subintervals $I$ of $S^1_6$ not containing the lower defect point,
  and to be $\e(I)$-linear for all bicolored subintervals
  not containing the upper defect points.
  (On subintervals not containing a defect point these two requirements
  coincide.)
  
  Informally, the categories 
  ${\CN_\sim}_2$, ${\CN^\sim}_2$, ${\CN^\sim_\sim}_2$
  are obtained, respectively, by relaxing the linearity of morphisms
  around the lower defect point (for ${\CN_\sim}_2$),
  around the upper defect point (for ${\CN^\sim}_2$),
  or around both defect points (for ${\CN^\sim_\sim}_2$).  
  In all three cases the objects coincide with the objects
  of $\CN_2$.
  Morphisms are defined more formally as follows.
  We use the following subintervals of our standard circle $S^1_6$:
  \begin{equation*}
    I_\sim := [4,5], \; {I_\sim}^c := [-1,4], \; I^\sim := [1,2], \; 
    {I^\sim}^c := [2,7], \; I_l := [-1,1], \; I_r := [2,4].
  \end{equation*} 
  Let ${}_\cala D_\calb$, ${}_\cala E_\calb$,
  ${}_{\cala'} {D'}_{\calb'}$, and ${}_{\cala'} {E'}_{\calb'}$
  be defects and let ${}_D H_E$ and ${}_{D'} {H'}_{E'}$ be sectors.
  A morphism $f \colon H \to H'$ in ${\CN^\sim}_2$ is a pair
  $f = (T,\e)$ where $\e \colon E \to E'$ is a morphism of
  $\CN_1$ and $T \in \Hom_{\e({I^\sim}^c)}(H,H')$. 
  A morphism $f \colon H \to H'$ in ${\CN_\sim}_2$ is a pair $f = (T,\delta)$
  where $\delta \colon D \to D'$ is a morphism of $\CN_1$
  and $T \in \Hom_{\delta({I_\sim}^c)}(H,H')$.
  Finally, a morphism $f \colon H \to H'$ in ${\CN_\sim^\sim}_2$ is a triple 
  $f = (T,\alpha,\beta)$ where $\alpha \colon \cala \to \cala'$,
  $\beta \colon \calb \to \calb'$ are morphisms of $\CN_0$ and
  $T \colon H \to H'$ belongs to both
  $\Hom_{\alpha(I_l)}(H,H')$ and  $\Hom_{\beta(I_r)}(H,H')$.
  Note that there are forgetful functors
  ${\CN^\sim}_2 \to {\CN^\sim_\sim}_2$ and 
  ${\CN_\sim}_2 \to {\CN^\sim_\sim}_2$ that are the identity on objects
  and for morphisms are defined by
  $(T,\e) \mapsto (T,s(\e),t(\e))$ and
  $(T,\delta) \mapsto (T,s(\delta),t(\delta))$. 
  
  We remark that many of our previous constructions extend to
  these variants of $\CN_2$.
  For example the vertical fusion [1-4] is natural for these more
  general morphisms and thus extends canonically to a functor
  ${\CN^\sim}_2 \x_{\CN_1} {\CN_\sim}_2 \to {\CN_\sim^\sim}_2$.
  As a rule of thumb: whenever we have a neighborhood of 
  a defect point on the boundary of the picture describing 
  one of our functors, we can canonically extend that functor, 
  adding an appropriate $\sim$ to source and target of the functor.  

  Homomorphisms in ${\CN^\sim}_2$ (or ${\CN_\sim}_2$) between vacuum sectors 
  can be more concretely described, as follows:

  \begin{lemma} \label{lem:sim-morphism-for-vaccum}
    Let ${}_\cala D_{\calb}$, ${}_{\cala'} D'_{\calb'}$ be defects, and
    let $\delta \colon D \to D'$ be an (iso)morphism of defects.  
    Then there is an isomorphism of vector spaces from $D'(I^\sim)$ to 
    $\Hom_{\delta({I^\sim}^c)}(H_0(D), \allowbreak H_0(D'))$ given by 
    $b \mapsto  b \circ H_0(\delta)$.  
    Here we use ``\;\!$b$'' to refer both to an algebra element and to the 
    linear map given by multiplying by that element.
    
    In particular, if we view $H_0(D)$ and $H_0(D')$ as objects 
    of ${\CN^\sim}_2$, then given a bounded linear map 
    $F \colon H_0(D) \to H_0(D')$,
    the pair $f := (F,\delta)$ defines
    a morphism $H_0(D) \to H_0(D')$ in ${\CN^\sim}_2$ if and only
    if $F = b \circ H_0(\delta)$ for some $b \in D'(I^\sim)$.
  \end{lemma} 

  \begin{proof}
    Haag duality for defects~\cite[Prop 1.16]{BDH(1*1)} 
    implies that 
    $\Hom_{D'({I^\sim}^c)} ( H_0( D'),H_0(D')) \allowbreak = D'(I^\sim)$,
    i.e., every $D'({I^\sim}^c)$-linear operator on $H_0(D')$ 
    is given by the action of  a unique element in $D'(I^\sim)$.
  
    Now $H_0(\delta) \colon H_0(D) \to H_0(D')$ is a $\delta(K)$-linear
    isomorphism for all intervals $K \subseteq S^1_6$.
    In particular it is $\delta({I^\sim}^c)$-linear and 
    induces an isomorphism
    $$D'(I^\sim) = \Hom_{D'({I^\sim}^c)} ( H_0( D'), H_0(D'))
         \to \Hom_{\delta({I^\sim}^c)}(H_0(D), H_0(D')).$$    

    The second statement follows from the first and 
    the definition of ${\CN^\sim}_2$. 
  \end{proof} 

\subsubsection*{The categories ${\CN^\simeq}_2$ and ${\CN_\simeq}_2$.}

  We define ${\CN^{\simeq}}_2$ as the full subcategory of ${\CN^\sim}_2$
  on objects of the form $H_0(D) \ox \ell$, where 
  ${}_\cala D_\calb$ is a defect and $\ell$ is a separable Hilbert space.
  It is a monodial subcategory.
  Similarly we obtain a monodial subcateory ${\CN_{\simeq}}_2$ of ${\CN_\sim}_2$.

  \begin{proposition} \label{prop:CNsim-as-subobjects-of-CNsimeq}
    Every object from ${\CN_2}$ is isomorphic in ${\CN^\sim}_2$ 
    (resp.\ in ${\CN_\sim}_2$)
    to a direct summand of an object from ${\CN^\simeq}_2$ 
    (resp.\ from ${\CN_\simeq}_2$).
  \end{proposition}

  \begin{proof}
    Let ${}_\cala D_\calb$, ${}_\cala E_\calb$ be defects and 
    let ${}_D K_E$ be a sector (in other words, an object of $\CN_2$).
    In particular, $K$ is an $E(I)$-module, where $I = [2,7]$ is the complement
    of the interval $(1,2)$ of our standard circle $S^1_6$.
    Observe that the vacuum sector $H_0(E)$ is faithful as an $E(I)$-module 
    (in fact $H_0(E)$ is isomorphic to the standard form $L^2(E(I))$).
    Recall that whenever $A$ is a separable von Neumann algebra acting on a separable Hilbert space $H$ and acting faithfully on a separable Hilbert space $H'$, then there is an $A$-linear isometric embedding $H \ox \ell \to H' \ox \ell$, and in particular an isometric embedding $H \to H' \ox \ell$, where $\ell$ is an infinite-dimensional Hilbert space.  We can therefore find an $E(I)$-linear isometric 
    embedding of $K$ into $H_0(E) \ox \ell$ for 
    some separable Hilbert space $\ell$.
    By the definition of ${\CN^\sim}_2$ this embedding defines a
    morphism
    $K \to H_0(E) \ox \ell$ in ${\CN^\sim}_2$, as desired.
  \end{proof}

  \begin{proposition} \label{prop:morphisms-CN-simeq}
    Let ${}_\cala D_\calb$ and ${}_{\cala'} {D'}_{\calb'}$
    be  defects and let $\ell$ and $\ell'$ be separable Hilbert spaces.
    Consider $X := H_0(D) \ox \ell$ and $X' := H_0(D') \ox \ell'$ 
    as objects  of ${\CN^\simeq}_2$. 
    Let $f = (\tilde F, \delta) \colon X \to X'$ 
    be a morphism in ${\CN^\sim}_2$.
    Let $Y := H_0(D') \ox \ell$.
    Then $f$ can be factored through $Y$ as $f = f_2 \circ f_1$ in ${\CN^\simeq}_2$ 
    where $f_1$ is induced by $\delta$ and $f_2$ is the identity
    on $D'$.
    More precisely we have
    \begin{enumerate}
    \item 
      $f_1 = (H_0(\delta) \ox \id_{\ell_1}, \delta) \colon X \to Y$;
    \item
      $f_2 = (T , \id_{D'}) \colon Y \to X'$
      where $T \in D'(I) \ox \bfB(\ell,\ell')$, with $I=[1,2]$.
    \end{enumerate}
  \end{proposition}

  Here, the space $\bfB(\ell,\ell')$ of bounded
  linear maps $\ell \to \ell'$ is a corner in the von Neumann algebra $\bfB(\ell \oplus \ell')$, and 
  the tensor product $D'(I) \ox \bfB(\ell,\ell')$ is defined to be
  the closure of the corresponding algebraic tensor product in the 
  von Neumann algebra $D'(I) \ox \bfB(\ell \oplus \ell')$.

  \begin{proof}[Proof of Proposition~\ref{prop:morphisms-CN-simeq}]
    We have to find $T \in D'(I^\sim) \ox \bfB(\ell,\ell')$
    such that $\tilde F = T \circ (H_0(\delta) \ox \id_{\ell})$.

    By the definition of ${\CN^\sim}_2$ we have
    \[
     \tilde F \in \Hom_{\delta({I^\sim}^c)}(H_0(D) \ox \ell, H_0(D') \ox \ell') 
        \cong  \Hom_{\delta({I^\sim}^c)}(H_0(D), H_0(D')) \ox \bfB(\ell,\ell').
    \] 
    By Lemma~\ref{lem:sim-morphism-for-vaccum}, the map 
    $T_0 \mapsto T_0 \circ H_0(\delta)$ gives an isomorphism 
    \[
    D'(I^\sim) \cong \Hom_{\delta({I^\sim}^c)}(H_0(D), H_0(D')).
    \]
    Therefore $T' \mapsto \tilde T \circ (H_0(\delta) \ox \id_{\ell})$
    yields an isomorphism
    \[
    D'(I^\sim) \ox \bfB(\ell,\ell') \xrightarrow{\cong} 
     \Hom_{\delta({I^\sim}^c)}(H_0(D) \ox \ell, H_0(D') \ox \ell').
    \]    
    The inverse image of $\tilde F$ under this isomorphism 
    provides the desired factorization.
  \end{proof}

\subsubsection*{[2-6] Right dewhisker.}
  The right dewhisker is an isomorphism
  \[
  \phantom{\text{[2-6]}.} \qquad  
  \cb{
  \displaytikz{
  \path
  (0,0) coordinate (1)
  (3,0) coordinate (2);
  \draw (1) -- (2);
  \draw (1) node[littlenode] {
  \innertikzb{
		\sh(0,0);
		\bxbt(1,1)1;
		\sh(0,1);
		\bxbb(1,1)1;
		\sh(1+\gap,0);
		\bxb(1,2)2;
  }
  };
  \draw (2) node[littlenode] {
  \innertikzb{
		\sh(0,0);
		\bxbur(1,1)1;
		\sh(1,0);
		\bxbul(1,1)2;
		\sh(0,1+\gap);
		\bxbdr(1,1)1;
		\sh(1,1+\gap);
		\bxbdl(1,1)2;
  }
  };
  }
  } \qquad \text{[2-6]}.
  \] 
  The left and right sides of the above picture describe functors
  \[
  L,R \colon (\CN_2 \x_{\CN_1} \CN_2) \x_{\CN_0} \x \CN_1 \to \CN_2
  \]
  and the right dewhisker is a natural isomorphism $\tau \colon L \to R$.
  Its construction will be a bit involved.
  We will show that in order to construct $\tau$
  it suffices to define $\tau$ on the image
  of the functor 
  $I \colon \CN_1 \x_{\CN_0} \CN_1 \to 
    (\CN_2 \x_{\CN_1} \CN_2) \x_{\CN_0} \x \CN_1$ defined
  by $I({}_\cala D_\calb,{}_\calb E_\calc) = (H_0(D),H_0(D),E)$.
  Here the natural isomorphism $\tau_0 \colon L \circ I \to R \circ I$ can be 
  constructed as the following composition:
  \[
  \phantom{.} \quad  
  \displaytikz{
  \path (0,0) coordinate (1) (3,0) coordinate (2)
        (6,0) coordinate (3) (9,0) coordinate (4)
        (12,0) coordinate (5);
  \draw (1) -- (2) -- (3) -- (4) -- (5);
  \draw (1) node[littlenode] {
  \innertikzb{
		\sh(0,0);
		\bxbt(1,1)3;
		\sh(0,1);
		\bxbb(1,1)3;
		\sh(1+\gap,0);
		\bxb(1,2)2;
  }
  };
  \draw (2) node[littlenode] {
  \innertikzb{
		\sh(0,0);
		\bxb(1,1)2;
		\sh(1,0);
		\bxb(1,1)2;
  }
  };
  \draw (3) node[littlenode] {
  \innertikzb{
		\sh(0,0);
		\bxb(2,1)2;\bull(1,0);\bull(1,1);
  }
  };
  \draw (4) node[littlenode] {
  \innertikzb{
		\sh(0,0);
		\sbx(2,1)2;
		\sh(0,1);
		\sbx(2,1)2;
		\bull(0,0);\bull(1,0);\bull(2,0);\bull(0,1);
                \bull(.5,1);\bull(1,1);\bull(1.5,1);\bull(2,1);
                \bull(0,2);\bull(1,2);\bull(2,2);
  }
  };
  \draw (5) node[littlenode] {
  \innertikzb{
		\sh(0,0);
		\bxbur(1,1)3;
		\sh(1,0);
		\bxbul(1,1)2;
		\sh(0,1+\gap);
		\bxbdr(1,1)3;
		\sh(1,1+\gap);
		\bxbdl(1,1)2;
  }
  };
  } \quad .
  \] 
  The darker shading indicates that those sectors are assumed to be 
  vacuum sectors in the definition of $\tau_0$, 
  but will need to be replaced by arbitrary sectors in order to define $\tau$. 
  The first and third isomorphisms are given by the 
  isomorphisms [2-1] or [2-2] (which are equivalent by axiom [3-1]).  
  The second and fourth isomorphisms are given by the 
  $1 \boxtimes 1$-isomorphism.
  
  In order to promote $\tau_0$ to $\tau$ we use the following diagram
  of functors.
  \begin{equation*}
    \xymatrix{
    \CN_1 \x_{\CN_0} \CN_1 \ar[r]^-{I} \ar[dd]^{i_0} &
    (\CN_2 \x_{\CN_1} \CN_2) \x_{\CN_0} \CN_1 
      \ar@/^1.5pc/[r]^{L} \ar@/_1.5pc/[r]_{R} \ar[dd]^{i_1} &
    \CN_2 \ar[dd]^{i_2}
    \\ \\
    ({\CN^{\simeq}}_2 \x_{\CN_1} {\CN_\simeq}_2) \x_{\CN_0} \CN_1
      \ar[r]^{\tilde I} &
    ({\CN^\sim}_2 \x_{\CN_1} {\CN_\sim}_2) \x_{\CN_0} \CN_1 
      \ar@/^1.5pc/[r]^{\tilde L} \ar@/_1.5pc/[r]_{\tilde R} &
    {\CN^\sim_\sim}_2 
    }
  \end{equation*}  %!% Around here is where the upper tilde started to mean nonlinear on top.
  Here we use the variations of $\CN_2$ introduced earlier.
  The functors $\tilde L$ and $\tilde R$ are the canonical extensions
  of $L$ and $R$.
  The functor $I$ applies the identity sector twice in the first entry
  and has already been defined.
  The vertical functors $i_1$ and $i_2$ are induced from the three
  inclusions of $\CN_2$  into ${\CN^\sim}_2$, ${\CN_\sim}_2$, 
  and ${\CN^\sim_\sim}_2$.
  The functor $\tilde I$ is induced from the inclusions
  ${\CN^\simeq}_2 \to {\CN^\sim}_2$ and ${\CN_\simeq}_2 \to {\CN_\sim}_2$.
  The composition $i_1 \circ I$ canonically factors as $\tilde I \circ i_0$.  
  
  In the next step we use $\tau_0$ to construct a natural isomorphism
  $\tilde \tau_0 \colon \tilde L \circ \tilde I \to \tilde R \circ \tilde I$.  Let $X_0 := (D,E)$ be an object of $\CN_1 \x_{\CN_0} \CN_1$ and let $X = (H_0(D) \ox \ell, H_0(D) \ox \ell' ,E)$ be an object
  from $({\CN^{\simeq}}_2 \x_{\CN_1} {\CN_\simeq}_2) \x_{\CN_0} \CN_1$.
  We have natural identifications
  $\tilde L(\tilde I (X)) = L(I(X_0)) \ox \ell \ox \ell'$ and 
  $\tilde R(\tilde I (X)) = R(I(X_0)) \ox \ell \ox \ell'$.
  We set $(\tilde \tau_0)_{X} := (\tau_0) \ox \id_{\ell} \ox \id_{\ell'}$.
  However, there are more morphisms in 
  $({\CN^{\simeq}}_2 \x_{\CN_1} {\CN_\simeq}_2) \x_{\CN_0} \CN_1$
  than there are in $\CN_1 \x_{\CN_0} \CN_1$, and we need to check that
  $\tilde \tau_0$ is natural with respect to these extra morphisms.
  Note that $\tilde \tau_0$ is natural for morphisms from $\CN_1$,
  because $\tau_0$ is.
  By Proposition~\ref{prop:morphisms-CN-simeq}, to check naturality with respect to morphisms in ${\CN^\simeq}_2$, it 
  suffices to consider morphisms   of the form
  \begin{enumerate}
  \item $f_1 = (H_0(\delta) \ox \id_{\ell},\delta)$, where $\delta: D \mapsto D'$, or
  \item $f_2 = (T,\id_{D})$
      where $T \in D(I) \ox \bfB(\ell,\ell')$, for $I=[1,2]$.
  \end{enumerate}
  Now $\tilde \tau_0$ is natural with respect to morphisms of the first kind
  because $\tau_0$ is natural for the morphisms 
  from $\CN_1$.
  As $\tau_0$ is equivariant for the action of $D(I)$, it follows that
  $\tilde \tau_0$ is also natural for morphisms of the second kind.
  Similarly, $\tilde \tau_0$ is natural for morphisms from ${\CN_\simeq}_2$.
  Thus $\tilde \tau_0$ is a natural transformation.

  By Proposition~\ref{prop:CNsim-as-subobjects-of-CNsimeq}
  every object of 
  $({\CN^\sim}_2 \x_{\CN_1} {\CN_\sim}_2) \x_{\CN_0} \CN_1$ 
  can be embedded as a direct summand in an object of
  $({\CN^\simeq}_2 \x_{\CN_1} {\CN_\simeq}_2) \x_{\CN_0} \CN_1$. %!%\CDcomm{Was that what was meant?}
  Thus we can extend $\tilde \tau_0$ canonically
  to a natural transformation $\tilde \tau \colon \tilde L \to \tilde R$.  %\CDcomm{Maybe this last sentence needs more explanation.}

  Let $X$ be an object from $(\CN_2 \x_{\CN_1} \CN_2) \x_{\CN_0} \CN_1$.
  Then 
  \[
  \tilde \tau_{i_1(X)} \colon \tilde L(i_1(X)) \to \tilde R(i_1(X))
  \]
  is a morphism in ${\CN^\sim_\sim}_2$.
  In fact we will see that $\tilde \tau_{i_1(X)}$ is a morphism in $\CN_2$, and so we can define $\tau_{X} := \tilde \tau_{i_1(X)}$ and obtain the desired
  natural isomorphism.
  To check that $\tilde \tau_{i_1(X)}$ is in $\CN_2$, we write $X = ({}_D H_{D'},{}_{D'} {H'}_{D''},E)$, where    
  $D$, $D'$, and $D''$ are $\cala$-$\calb$-defects and
  $E$ is a $\calb$-$\calc$-defect. For convenience we ignore the collapse map $S^1_8 \mapsto S^1_6$ and think of the $(D \circledast_\calb E)$-$(D'' \circledast_\calb E)$-sectors 
  $L(X)$ and $R(X)$ as defined on $S^1_8$ instead of on $S^1_6$.
  We have to show that $\tilde \tau_{i_1(X)}$ is equivariant for the
  actions of $D(I^\sim)$ and $D''(I_\sim)$, where we now view $I^\sim = [1,2]$
  and $I_\sim = [6,7]$ as subintervals of $S^1_8$.
  Elements of $D(I^\sim)$ and $D''(I_\sim)$ can be viewed as morphisms
  in ${\CN^\sim}_2$ and ${\CN_\sim}_2$, and therefore as morphisms in
  $({\CN^\sim}_2 \x_{\CN_1} {\CN_\sim}_2) \x_{\CN_0} \CN_1$.
  Therefore the required equivariance of $\tilde \tau_{i_1(X)}$ follows
  from the naturality of $\tilde \tau$.
  This finishes the construction of the right dewhisker.  
       
\subsubsection*{[2-7] Left dewhisker.}
  The left dewhisker is defined analogously to the right dewhisker 
  and is drawn as
  \[
  \phantom{\text{[2-7]}.} \qquad
  \cb{
  \displaytikz{
  \path
  (0,0) coordinate (1)
  (3,0) coordinate (2);
  \draw (1) -- (2);
  \draw (1) node[littlenode] {
  \innertikzb{
		\sh(1+\gap,0);
		\bxbt(1,1)1;
		\sh(1+\gap,1);
		\bxbb(1,1)1;
		\sh(0,0);
		\bxb(1,2)2;
  } 
  };
  \draw (2) node[littlenode] {
  \innertikz{
		\sh(0,0);
		\bxbur(1,1)2;
		\sh(1,0);
		\bxbul(1,1)1;
		\sh(0,1+\gap);
		\bxbdr(1,1)2;
		\sh(1,1+\gap);
		\bxbdl(1,1)1;
  }
  };
  }
  } \qquad \text{[2-7]}.
  \]

\subsubsection*{[2-8] Switch.}   
  The switch isomorphism is a composite of two isomorphisms
  \[
  \phantom{\text{[2-8]}.} \qquad
  \cb{
  \displaytikz{
  \path
  (0,0) coordinate (1)
  (3,0) coordinate (2)
  (6,0) coordinate (3);
  \draw (1) -- (2) -- (3);
  \draw (1) node[littlenode] {
  \innertikzb{
	\sh(0,0);
	\bxbur(1,1)2;
	\sh(1,0);
	\bxbul(1,1)1;
	\sh(0,1+\gap);
	\bxbdr(1,1)1;
	\sh(1,1+\gap);
	\bxbdl(1,1)2;
  }
  };
  \draw (2) node[littlenode] {
  \innertikzb{
		\sh(0,0);
		\bxb(1,1)1;
		\sh(1,0);
		\bxb(1,1)1;
  }
  };
  \draw (3) node[littlenode] {
  \innertikzb{
	\sh(0,0);
	\bxbur(1,1)1;
	\sh(1,0);
	\bxbul(1,1)2;
	\sh(0,1+\gap);
	\bxbdr(1,1)2;
	\sh(1,1+\gap);
	\bxbdl(1,1)1;
  }
  };
  }
  } \qquad \text{[2-8]}.
  \] 
  Each of those two isomorphisms is referred to as a half-switch. 
  Arguing as in the construction of the right dewhisker
  it suffices to construct these isomorphisms in the cases where 
  all sectors are vacuum sectors.
  In this case the first half switch is defined as 
  \[
  \phantom{.} \quad
  \displaytikz{
  \path
  (0,0) coordinate (1)
  (3,0) coordinate (2)
  (6,0) coordinate (3)
  (9,0) coordinate (4);
  \draw (1) -- (2) -- (3) -- (4);
  \draw (1) node[littlenode] {
  \innertikzb{
		\sh(0,0);
		\bxbur(1,1)2;
		\sh(1,0);
		\bxbul(1,1)3;
		\sh(0,1+\gap);
		\bxbdr(1,1)3;
		\sh(1,1+\gap);
		\bxbdl(1,1)2;
  }
  };
  \draw (2) node[littlenode] {
  \innertikzb{
		\sh(0,0);
		\sbx(2,1)2;
		\sh(0,1);
		\sbx(2,1)2;
		\bull(0,0);\bull(1,0);\bull(2,0);\bull(0,1);\bull(.5,1);
                \bull(1,1);\bull(1.5,1);\bull(2,1);\bull(0,2);
                \bull(1,2);\bull(2,2);
  }
  };
  \draw (3) node[littlenode] {
  \innertikzb{
		\sh(0,0);
		\bxb(2,1)2;\bull(1,0);\bull(1,1);
  }
  };
  \draw (4) node[littlenode] {
  \innertikzb{
		\sh(0,0);
		\bxb(1,1)3;
		\sh(1,0);
		\bxb(1,1)3;
  }
  };
  } \quad .
  \] 
  Here the first and third isomorphisms are both 
  the $1 \boxtimes 1$-isomorphism 
  and the second is [2-1] (which agrees in this case with [2-2]).  
  The second half-switch is defined analogously.
  
%g%\pagebreak  % This is a hack because it was breaking after the section headings

\subsection{Transformations for horizontal associators}  
\subsubsection*{[2-9], [2-10], and [2-11] Whisker associators.}
  The associators for twice whiskered sectors are given, using the $1 \boxtimes 1$-isomorphism and the 
  associativity of Connes fusion, as
  \[
  \phantom{\text{[2-9]},} \qquad
  \cb{
  \displaytikz{
  \path
  (0,0) coordinate (1)
  (3,0) coordinate (2)
  (6,0) coordinate (3);
  \draw (1) -- (2) -- (3);
  \draw (1) node[littlenode] {
  \innertikz{
	\sh(0,0);
	\bxb(1,1)1;
	\sh(1,0);
	\bxb(1,1)2;
	\sh(2+\gap,0);
	\bxbdr(1,1)2;
  }
  };
  \draw (2) node[littlenode] {
  \innertikz{
		\sh(0,0);
		\bxb(1,1)1;
		\sh(1+\gap,0);
		\bxb(1,1)2;
		\sh(2+\gap,0);
		\bxb(1,1)2;
  }
  };
  \draw (3) node[littlenode] {
  \innertikz{
	\sh(0,0);
	\bxb(1,1)1;
	\sh(1,0);
	\bxb(2,1)2;
	\bull(2,0);\bull(2,1);
  }
  };
  }
  } \qquad \text{[2-9]},
  \]
  \[\phantom{\text{[2-10]},} \qquad
  \cb{
  \displaytikz{
  \path
  (0,0) coordinate (1)
  (3,0) coordinate (2)
  (6,0) coordinate (3);
  \draw (1) -- (2) -- (3);
  \draw (1) node[littlenode] {
  \innertikz{
	\sh(0,0);
	\bxb(2,1)2;
	\sh(2,0);
	\bxb(1,1)1;
	\bull(1,0);\bull(1,1);
  }
  };
  \draw (2) node[littlenode] {
  \innertikz{
		\sh(0,0);
		\bxb(1,1)2;
		\sh(1,0);
		\bxb(1,1)2;
		\sh(2+\gap,0);
		\bxb(1,1)1;
  }
  };
  \draw (3) node[littlenode] {
  \innertikz{
	\sh(0,0);
	\bxb(1,1)2;
	\sh(1+\gap,0);
	\bxb(1,1)2;
	\sh(2+\gap,0);
	\bxb(1,1)1;
  }
  };
  }
  } \qquad \text{[2-10]},
  \] 
  \[\phantom{\text{[2-11]}.} \qquad
  \cb{
  \displaytikz{
  \path
  (0,0) coordinate (1)
  (3,0) coordinate (2);
  \draw (1) -- (2);
  \draw (1) node[littlenode] {
  \innertikz{
	\sh(0,0);
	\bxb(1,1)2;
	\sh(1,0);
	\bxb(1,1)1;
	\sh(2+\gap,0);
	\bxb(1,1)2;
  }
  };
  \draw (2) node[littlenode] {
  \innertikz{
		\sh(0,0);
		\bxb(1,1)2;
		\sh(1+\gap,0);
		\bxb(1,1)1;
		\sh(2+\gap,0);
		\bxb(1,1)2;
  }
  };
  }
  } \qquad \text{[2-11]}.
  \]   
  As in the case of the vertical associator [2-3] we will also here
  often suppress the associator for Connes fusion.
  In particular, we will suppress [2-11]. 

\subsubsection*{[2-12] Horizontal associator.}
  The associator for the composition of defects is induced from the
  associator for fusion (or fiber product) of von Neumann algebras
  and is discussed in~\cite[Eq 1.55]{BDH(1*1)}. 
  Here we will suppress this isomorphism.  

  In~\cite[Prop 4.32]{BDH(1*1)}
  we proved that the
  $1 \boxtimes 1$-isomorphism is associative in the following sense. 

  \begin{lemmaletter}[J] %formerly L
  The $1 \boxtimes 1$-isomorphism is associative for the
  composition of defects, that is the following diagram commutes:
  \[
  \displaytikz{
  \path
  (0,2) coordinate (1)
  (3,2) coordinate (2)
  (3,0) coordinate (3)
  (0,0) coordinate (4);
  \draw (1) -- (2) -- (3) -- (4) -- (1);
  \draw (1) node[littlenode] {
  \innertikz{
	\sh(0,0);
	\sbx(1,1)2;
	\sh(1,0);
	\sbx(1,1)2;
	\sh(2,0);
	\sbx(1,1)2;
  }
  };
  \draw (2) node[littlenode] {
  \innertikz{
	\sh(0,0);
	\sbx(2,1)2;
	\sh(2,0);
	\sbx(1,1)2;
  }
  };
  \draw (3) node[littlenode] {
  \innertikz{
	\sh(0,0);
	\sbx(3,1)2;
  }
  };
  \draw (4) node[littlenode] {
  \innertikz{
	\sh(0,0);
	\sbx(1,1)2;
	\sh(1,0);
	\sbx(2,1)2;
  }
  };
  }  \quad .
  \] 
  \end{lemmaletter}

\noindent  We point out that in the upper left corner of the diagram in this lemma,
  we have suppressed the morphism of sectors associated to the associator of horizontal fusion.
Similarly, in the lower right corner of the diagram 
  we have suppressed the morphism of sectors associated to the associator for fusion of defects.

%  \medskip \medskip

%%\pagebreak
\subsection{Transformations for horizontal identities}

  The $2$-data discussed so far is part of the definition of both
  a $2$-category object and a  dicategory object
  (in the $2$-category of symmetric monoidal
  categories). 
  We remind the reader that our $3$-category of conformal nets
  is a dicategory object; the remaining  pieces of data labeled here [2-13]
  to [2-18] are  labeled [D2-13] to [D2-18] in the appendix and correspondingly in~\cite{Douglas-Henriques(Internal-bicategories)}.

\subsubsection*{Left and right quasi-identity}
  Let ${}_\cala D_\calb$ and ${}_\cala E_\calb$ be defects
  and let ${}_D H_E$ be a sector.
  There is no canonical isomorphism 
  $H_0(\id_\cala) \boxtimes_\cala H \cong H$ in $\CN_2$.
  In fact, $H_0(\id_\cala) \boxtimes_\cala H$ is an
  $(\id_\cala \circledast_\cala D)$-$(\id_\cala \circledast_\cala E)$-sector  
  and not a $D$-$E$-sector.
  There is however such an isomorphism if 
  we are willing to twist $H$ on top and bottom by the diffeomorphism $\Phi_L$ introduced in the construction of the left identity [1-7].
  The left quasi-identity has been constructed in~\cite[Def 6.20]{BDH(1*1)} (where it was called the ``left unitor" $\hat{\Upsilon}^l$),
  and will be draw as
  \[\cb{
  \displaytikz{
  \path
  (0,0) coordinate (1)
  (3.5,0) coordinate (2);
  \draw (1) -- (2);
  \draw (1) node[littlenode] {
  \displaytikz{
		\sh(0,0);
		\bxb(1,1)2;
		\sh(1,0);
		\bxb(1,1)1;
		\idfct(.5,0);\idfct(.5,1);\dfct(1.5,0);\dfct(1.5,1);
  }
  };
  \draw (2) node[littlenode] {
  \displaytikz{
		\sh(0,0);
		\bxb(2,1)1;
		\bubble(.25,1,.75)1;
		\sh(0,1);
		\bubble(.25,1,.75)1;
		\bull(1,-\cw);\bull(1,1+\cw);
				\idfct(.5,-\cw);\idfct(.5,1+\cw);\dfct(1.5,0);\dfct(1.5,1);
  }
  };
  }
 }  \]   
  Here the balloons on the right box
  signal that the $D$-$E$-sector structure on $H$ has
  been twisted to an 
  $(\id_\cala \circledast_\cala D)$-$(\id_\cala \circledast_\cala E)$-sector structure
  using $\cala(\Phi_L)$ at the indicated portions of the picture.
  Of course we can move one or two of the balloons to the left box by composing with the appropriate inverse diffeomorphism.
  For example we obtain from the left quasi-identity  
  an isomorphism of 
  $(\id_\cala \circledast_\cala D)$-$E$-sectors 
  drawn as
\[\cb{
  \displaytikz{
  \path
  (0,0) coordinate (1)
  (3.5,0) coordinate (2);
  \draw (1) -- (2);
  \draw (1) node[littlenode] {
  \displaytikz{
		\sh(0,0);
		\boxslopedr(.25,.75,1)2;
		\bubble(.25,1,.75)2;
		\sh(1,0);
		\boxslopedl(.75,.25,1)1;
		\begin{pgfonlayer}{foreground}
			\bull(0,.5);\bull(1,.5);\bull(2,.5);
			\idfct(.5,\cw);\idfct(.5,1);\dfct(1.5,0);\dfct(1.5,1);
		\end{pgfonlayer}
  }
  };
  \draw (2) node[littlenode] {
  \displaytikz{
	\sh(0,0);
	\boxt(.25*.75,.75,.75*.75,1)1;
	\sh(0,1);
	\bubble(.25*.75,.75,.75*.75)1;
	\begin{pgfonlayer}{foreground}
		\bull(0,.5);\bull(1.5,.5);\bull(.75,1+\cw);
		\idfct(.375,1+\cw);\dfct(.75,0);\dfct(1.125,1);
	\end{pgfonlayer}
  }  
  };
}
}
\]
  Similarly, there is a right quasi-identity using $\Phi_R$ drawn as
  \[\cb{
  \displaytikz{
  \path
  (0,0) coordinate (1)
  (3.5,0) coordinate (2);
  \draw (1) -- (2);
  \draw (1) node[littlenode] {
  \displaytikz{
		\sh(0,0);
		\bxb(1,1)1;
		\sh(1,0);
		\bxb(1,1)2;
		\dfct(.5,0);\dfct(.5,1);\idfct(1.5,0);\idfct(1.5,1);
  }
  };
  \draw (2) node[littlenode] {
  \displaytikz{
		\sh(0,0);
		\bxb(2,1)1;
		\bubble(.75,1,.25)1;
		\sh(0,1);
		\bubble(.75,1,.25)1;
		\bull(1,-\cw);\bull(1,1+\cw);
				\idfct(1.5,-\cw);\idfct(1.5,1+\cw);\dfct(.5,0);\dfct(.5,1);
  }
  };
  }
 }  \]

  The next two lemmas are proved respectively in~\cite[Eq 6.23]{BDH(1*1)} and in~\cite[Lem 6.21]{BDH(1*1)}. 
 
  \begin{lemmaletter}[K] %formerly M  
     The following diagram commutes: 
     \[
     \begin{tikzpicture}
     \path 
     (0,0) coordinate (1)
     (3.5,0) coordinate (2)
     (0,-2) coordinate (3)
     (3.5,-2) coordinate (4);
	\draw (1) -- (2) -- (4) -- (3) -- (1);
	\draw (1) node[littlenode] {
	\displaytikz{
	\sh(0,0);
	\boxslopeur(.75,.25,1)2;
	\sh(1,0);
	\boxslopeul(.25,.75,1)2;
	\sh(0,1);
	\bubble(.75,1,.25)2;
	\dfct(.5,0);\dfct(.5,1);\idfct(1.5,0);\idfct(1.5,1-\cw);\bull(0,.5);\bull(1,.5);\bull(2,.5);
	}
	};
	\draw (2) node[littlenode] {
	\displaytikz{
	\sh(0,0);
	\bxb(2,1)2;
	\bubble(.75,1,.25)2;
	\dfct(.5,0);\dfct(1,1);\idfct(1.5,-\cw);\bull(1,-\cw);
	}
	};
	\draw (3) node[littlenode] {
	\displaytikz{
	\sh(0,0);
	\boxt(.75,1,.25,1)2;
	\sh(0,1);
	\bubble(.75,1,.25)2;
		\dfct(.5,0);\dfct(.5,1);\idfct(1.5,0);\idfct(1.5,1-\cw); \bull(1,0);\bull(1,1-\cw);\bull(0,.5);\bull(2,.5);
	}
	};
		\draw (4) node[littlenode] {
	\displaytikz{
			\sh(0,0);
		\boxt(.875,.5,.625,1)2;
		\sh(0,1+\gap);
		\bubble(.875,.5,.625)2;
		\sh(0,1+\gap);
		\boxb(.875,.5,.625,1)2;
		\begin{pgfonlayer}{foreground}
			\bull(.5,1);\bull(1.5,1);\bull(.5,1+\gap);
			\bull(1.5,1+\gap); \bull(1,1-\cw);\bull(1,1+\gap-\cw);
                        \dfct(1,2+\gap);\dfct(.75,1+\gap); \dfct(.75,1);
                        \dfct(.5,0);\idfct(1.5,0);\bull(1,0);\idfct(1.25,1-\cw);\idfct(1.25,1+\gap-\cw);
		\end{pgfonlayer}
	}
	};	     
	\path (1) -- (2) node[la] {\emph{\quasiidentity}}; %Note emph to counteract Lemma style.
	\path (3) -- (4) node[lb] {\emph{\identity}};
	\path (3) -- (1) node[la] {\onetimesone};
	\path (4) -- (2) node[lb] {\emph{\identity}};
     \end{tikzpicture}
     \]
  \end{lemmaletter}

  %\CDcomm{To change the labels here, adjust the definitions in cnii-defs.tex.}

  \begin{lemmaletter}[L] %formerly X
     The $1 \boxtimes 1$-isomorphism is, in the following sense,
     natural with respect to the left quasi-identity---this diagram commutes: 
     \[
     \begin{tikzpicture}
     \path 
     (0,0) coordinate (1)
     (3.5,0) coordinate (2)
     (0,-2) coordinate (3)
     (3.5,-2) coordinate (4);
	\draw (1) -- (2) -- (4) -- (3) -- (1);
	\draw (1) node[littlenode] {
	\displaytikz{
	\sh(0,0);
	\bxb(1,1)2;
	\sh(1,0);
	\bxb(1,1)2;
	\sh(2,0);
	\bxb(1,1)2;
	\idfct(.5,0);\idfct(.5,1);\dfct(1.5,0);\dfct(1.5,1);\dfct(2.5,0);\dfct(2.5,1);
	}
	};
	\draw (2) node[littlenode] {
	\displaytikz{
	\sh(0,0);
	\bxb(2,1)2;
	\bubble(.25,1,.75)2;
	\sh(0,1);
	\bubble(.25,1,.75)2;
	\sh(2,0);
	\bxb(1,1)2;
	\idfct(.5,-\cw);\idfct(.5,1+\cw);\dfct(1.5,0);\dfct(1.5,1);\dfct(2.5,0);\dfct(2.5,1);\bull(1,-\cw);\bull(1,1+\cw);
	}
	};
	\draw (3) node[littlenode] {
	\displaytikz{
	\sh(0,0);
	\bxb(1,1)2;
	\sh(1,0);
	\bxb(2,1)2;
	\bull(2,0);\bull(2,1);
		\idfct(.5,0);\idfct(.5,1);\dfct(1.5,0);\dfct(1.5,1);\dfct(2.5,0);\dfct(2.5,1);
	}
	};
		\draw (4) node[littlenode] {
	\displaytikz{
	\sh(0,0);
	\bxb(3,1)2;
	\bubble(.25,1,.75)2;
	\sh(0,1);
	\bubble(.25,1,.75)2;
	\idfct(.5,-\cw);\idfct(.5,1+\cw);\dfct(1.5,0);\dfct(1.5,1);\dfct(2.5,0);\dfct(2.5,1);\bull(1,-\cw);\bull(1,1+\cw); \bull(2,0);\bull(2,1);
	}
	};	     
	\path (1) -- (2) node[la] {\emph{\quasiidentity}}; %Note emph to counteract Lemma style.
	\path (3) -- (4) node[lb] {\emph{\quasiidentity}};
	\path (3) -- (1) node[la] {\onetimesone};
	\path (4) -- (2) node[lb] {\onetimesone};
     \end{tikzpicture}
     \]
  \end{lemmaletter}
  
  Similarly, there are ``right versions'' of Lemma~K and Lemma~L
  for the right quasi-identity.

\subsubsection*{[2-13] Left identity pass.}
  The pass through a left identity is defined as the following 
  composite
  \[
  \phantom{\text{[2-13]}.} \qquad
  \cb{
  \displaytikz{
  \path
  (0,0) coordinate (1)
  (3.5,0) coordinate (2)
  (7,0) coordinate (3)
  (10.5,0) coordinate (4);
  \draw (1) -- (2) -- (3) -- (4);
  \draw (1) node[littlenode] {
  \displaytikz{
		\sh(0,0);
		\boxt(.625,.5,.875,1)2;
		\sh(0,1);
		\bubble(.625,.5,.875)2;
		\sh(0,1+\gap);
		\boxslopedr(.625,.375,1)2;
		\sh(1,1+\gap);
		\boxslopedl(.875,.125,1)1;
		\begin{pgfonlayer}{foreground}
			\bull(.5,1);\bull(1.5,1);\bull(.5,1+\gap);
                        \bull(1,1+\cw);
			\bull(1.5,1+\gap);
                        \bull(1,1.5+\gap);
                        \idfct(.5,2+\gap);\idfct(.75,1+\gap+\cw);\idfct(.75,1+\cw);
                        \dfct(1,0);\dfct(1.25,1);\dfct(1.25,1+\gap);\dfct(1.5,2+\gap);
		\end{pgfonlayer}
  }
  };
  \draw (2) node[littlenode] {
  \displaytikz{
		\sh(0,0);
		\boxslopedr(.25,.75,1)2;
		\bubble(.25,1,.75)2;
		\sh(1,0);
		\boxslopedl(.75,.25,1)1;
		\begin{pgfonlayer}{foreground}
			\bull(0,.5);\bull(1,.5);\bull(2,.5);
			\idfct(.5,\cw);\idfct(.5,1);\dfct(1.5,0);\dfct(1.5,1);
		\end{pgfonlayer}
  }
  };
  \draw (3) node[littlenode] {
  \displaytikz{
	\sh(0,0);
	\boxt(.25*.75,.75,.75*.75,1)1;
	\sh(0,1);
	\bubble(.25*.75,.75,.75*.75)1;
	\begin{pgfonlayer}{foreground}
		\bull(0,.5);\bull(1.5,.5);\bull(.75,1+\cw);
		\idfct(.375,1+\cw);\dfct(.75,0);\dfct(1.125,1);
	\end{pgfonlayer}
  }
  };
  \draw (4) node[littlenode] {
  \displaytikz{
  	\sh(0,0);
	\bx(1.5,1)1;
	\sh(0,1);
	\boxt(.25*.75,.75,.75*.75,1)2;
	\sh(0,2);
	\bubble(.25*.75,.75,.75*.75)2;
		\begin{pgfonlayer}{foreground}
			\bull(0,0);\bull(1.5,0);\bull(.5,1);\bull(1,1);\bull(.75,2+\cw);
			\idfct(.375,2+\cw);\dfct(.75,0);\dfct(.75,1);\dfct(1.125,2);
		\end{pgfonlayer}
  }
  };
  }
  } \qquad \!\!\! \text{[2-13]}.
  \] 
  For the first map observe that we can move the diffeomorphism
  from the lower defect to the horizontal composition of
  the two upper defects.
  The first map is then obtained by the bottom identity [2-2].
  The second map is obtained from the left quasi-identity
  by applying the inverse of $\Phi_L$ on the
  lower half---this moves the lower bubble from the right hand picture
  to the left hand picture. 
  The third map is the top identity [2-1].

\subsubsection*{[2-14] Right identity pass.}
  The pass through a right identity is a reflection along a vertical axis of
  the pass through a left identity, and is defined similarly:
  \[
  \phantom{\text{[2-14]}.} \qquad
  \cb{
  \displaytikz{
  \path
  (0,0) coordinate (1)
  (3.5,0) coordinate (2)
  (7,0) coordinate (3)
  (10.5,0) coordinate (4);
  \draw (1) -- (2) -- (3) -- (4);
  \draw (1) node[littlenode] {
  \displaytikz{
		\sh(0,0);
		\boxt(.875,.5,.625,1)2;
		\sh(0,1);
		\bubble(.875,.5,.625)2;
		\sh(0,1+\gap);
		\boxslopedr(.875,.125,1)1;
		\sh(1,1+\gap);
		\boxslopedl(.625,.375,1)2;
		\begin{pgfonlayer}{foreground}
			\bull(.5,1);\bull(1.5,1);\bull(.5,1+\gap);
                        \bull(1,1+\cw);
			\bull(1.5,1+\gap);
                        \bull(1,1.5+\gap);
                        \dfct(.5,2+\gap);\dfct(.75,1+\gap);\dfct(.75,1);
                        \dfct(1,0);\idfct(1.25,1+\cw);\idfct(1.25,1+\gap+\cw);\idfct(1.5,2+\gap);
		\end{pgfonlayer}
  }
  };
  \draw (2) node[littlenode] {
  \displaytikz{
		\sh(0,0);
		\boxslopedr(.75,.25,1)1;
		\bubble(.75,1,.25)2;
		\sh(1,0);
		\boxslopedl(.25,.75,1)2;
		\begin{pgfonlayer}{foreground}
			\bull(0,.5);\bull(1,.5);\bull(2,.5);
			\dfct(.5,0);\dfct(.5,1);\idfct(1.5,\cw);\idfct(1.5,1);
		\end{pgfonlayer}
  }
  };
  \draw (3) node[littlenode] {
  \displaytikz{
	\sh(0,0);
	\boxt(.75*.75,.75,.25*.75,1)1;
	\sh(0,1);
	\bubble(.75*.75,.75,.25*.75)1;
	\begin{pgfonlayer}{foreground}
		\bull(0,.5);\bull(1.5,.5);\bull(.75,1+\cw);
		\dfct(.375,1);\dfct(.75,0);\idfct(1.125,1+\cw);
	\end{pgfonlayer}
  }
  };
  \draw (4) node[littlenode] {
  \displaytikz{
  	\sh(0,0);
	\bx(1.5,1)1;
	\sh(0,1);
	\boxt(.75*.75,.75,.25*.75,1)2;
	\sh(0,2);
	\bubble(.75*.75,.75,.25*.75)2;
		\begin{pgfonlayer}{foreground}
			\bull(0,0);\bull(1.5,0);\bull(.5,1);\bull(1,1);\bull(.75,2+\cw);
			\dfct(.375,2);\dfct(.75,0);\dfct(.75,1);\idfct(1.125,2+\cw);
		\end{pgfonlayer}
  }
  };
  }
  } \qquad \!\! \text{[2-14]}.
  \] 

\subsubsection*{[2-15] Swap.}
  The swap is an isomorphism
  \[
  \phantom{\text{[2-15]}.} \qquad
  \cb{
  \displaytikz{
  \path
  (0,0) coordinate (1)
  (3,0) coordinate (2);
  \draw (1) -- (2);
  \draw (1) node[littlenode] {
  \displaytikz{
	\sh(0,0);
	\boxt(.25*.75,.75,.75*.75,1)2;
	\sh(0,1);
	\bubble(.25*.75,.75,.75*.75)2;
	\begin{pgfonlayer}{foreground}
		\bull(0,.5);\bull(1.5,.5);\bull(.75,1+\cw);
		\idfct(.375,1+\cw);\idfct(.75,0);\idfct(1.125,1);
	\end{pgfonlayer}
  }
  };
  \draw (2) node[littlenode] {
  \displaytikz{
	\sh(0,0);
	\boxt(.75*.75,.75,.25*.75,1)2;
	\sh(0,1);
	\bubble(.75*.75,.75,.25*.75)2;
	\begin{pgfonlayer}{foreground}
		\bull(0,.5);\bull(1.5,.5);\bull(.75,1+\cw);
		\idfct(.375,1);\idfct(.75,0);\idfct(1.125,1+\cw);
	\end{pgfonlayer}
  }
  };
  }
  } \qquad \text{[2-15]}.
  \]
  The construction of the swap depends on the flip
  and will be given after the construction of the flip below.

\subsubsection*{[2-16] Left identity expansion.}
  The left identity expansion is obtained from the
  $1 \boxtimes 1$-isomorphism and is drawn as
  \[
  \phantom{\text{[2-16]}.} \qquad
  \cb{
  \displaytikz{
  \path
  (0,0) coordinate (1)
  (4,0) coordinate (2);
  \draw (1) -- (2);
  \draw (1) node[littlenode] {
  \displaytikz{
	\sh(0,0);
	\boxt(.25*.75,.75,.75*.75,1)2;
	\sh(0,1);
	\bubble(.25*.75,.75,.75*.75)2;
	\sh(1.5,0);
	\bxb(1,1)2;
	\begin{pgfonlayer}{foreground}
		\bull(0,.5);\bull(1.5,.5);\bull(.75,1+\cw);
		\idfct(.375,1+\cw);\dfct(.75,0);\dfct(1.125,1);
		\dfct(2,0);\dfct(2,1);
	\end{pgfonlayer}
  }
  };
  \draw (2) node[littlenode] {
  \displaytikz{
	\sh(0,0);
	\boxt(.25*.75,.75,.75*.75+1,1)2;
	\sh(0,1);
	\bubble(.25*.75,.75,.75*.75+1)2;
	\begin{pgfonlayer}{foreground}
		\bull(0,.5);\bull(.75,1+\cw);\bull(1.5,0);\bull(1.5,1);\bull(2.5,.5);
		\idfct(.375,1+\cw);\dfct(.75,0);\dfct(1.125,1);
		\dfct(2,0);\dfct(2,1);
	\end{pgfonlayer}
  }
  };
  }
  } \qquad \text{[2-16]}.
  \]
 
\subsubsection*{[2-17] Right identity expansion.}
  The right identity expansion is a reflection along a vertical axis of
  the left identity expansion, and is defined similarly:
  \[\phantom{\text{[2-17]}.} \qquad
  \cb{
  \displaytikz{
  \path
  (0,0) coordinate (1)
  (4,0) coordinate (2);
  \draw (1) -- (2);
  \draw (1) node[littlenode] {
  \displaytikz{
	\sh(1,0);
	\boxt(.75*.75,.75,.25*.75,1)2;
	\sh(1,1);
	\bubble(.75*.75,.75,.25*.75)2;
	\sh(0,0);
	\bxb(1,1)2;
	\begin{pgfonlayer}{foreground}
		\bull(0,.5);\bull(1.75,1+\cw); \bull(2.5,.5);
		\idfct(2.125,1+\cw);\dfct(1.75,0);\dfct(1.375,1);
		\dfct(.5,0);\dfct(.5,1);
	\end{pgfonlayer}
  }
  };
  \draw (2) node[littlenode] {
  \displaytikz{
	\sh(0,0);
	\boxt(.75*.75+1,.75,.25*.75,1)2;
	\sh(0,1);
	\bubble(.75*.75+1,.75,.25*.75)2;
	\begin{pgfonlayer}{foreground}
		\bull(0,.5);\bull(1,1);\bull(1,0);\bull(1.75,1+\cw);\bull(2.5,.5);
		\dfct(1.375,1);\dfct(1.75,0);\idfct(2.125,1+\cw);
		\dfct(.5,0);\dfct(.5,1);
	\end{pgfonlayer}
  }
  };
  }
  } \qquad \text{[2-17]}.
  \]

\subsubsection*{[2-18] Flip.}
  The flip is an isomorphism
  \[\phantom{\text{[2-18]}.} \qquad
  \cb{
  \displaytikz{
  \path
  (0,0) coordinate (1)
  (4,0) coordinate (2);
  \draw (1) -- (2);
  \draw (1) node[littlenode] {
  \displaytikz{
	\sh(0,0);
	\boxt(.75*.75,.75,.25*.75,1)2;
	\sh(0,1);
	\bubble(.75*.75,.75,.25*.75)2;
	\sh(1.5,0);
	\bxb(1,1)2;
	\begin{pgfonlayer}{foreground}
		\bull(0,.5);\bull(1.5,.5);\bull(.75,1+\cw);
		\dfct(.375,1);\dfct(.75,0);\idfct(1.125,1+\cw);
		\dfct(2,0);\dfct(2,1);
	\end{pgfonlayer}
  }
  };
  \draw (2) node[littlenode] {
  \displaytikz{
	\sh(1,0);
	\boxt(.25*.75,.75,.75*.75,1)2;
	\sh(1,1);
	\bubble(.25*.75,.75,.75*.75)2;
	\sh(0,0);
	\bxb(1,1)2;
	\begin{pgfonlayer}{foreground}
		\bull(0,.5);\bull(1.75,1+\cw); \bull(2.5,.5);
		\dfct(2.125,1);\dfct(1.75,0);\idfct(1.375,1+\cw);
		\dfct(.5,0);\dfct(.5,1);
	\end{pgfonlayer}
  }
  };
  }
  } \qquad \text{[2-18]}
  \] 
  and is defined as follows.
  Note that both the left and right hand side of the flip are obtained
  from the horizontal composition of two identity sectors
  \[
  \displaytikz{
		\sh(0,0);
		\bxb(1,1)2;
		\sh(1,0);
		\bxb(1,1)2;
		\dfct(.5,0);\dfct(.5,1);\dfct(1.5,0);\dfct(1.5,1);
  }
  \]
  by two different diffeomorphism $\phi_L$, $\phi_R$ from 
  an interval of length $2$ to an interval of length $1$---the intervals under consideration are those between the non-identity defect 
  points.  (The diffeomorphisms $\phi_L$ and $\phi_R$ are extensions by an identity of the restrictions of the diffeomorphisms $\Phi_L$ and $\Phi_R$ used in the construction of the left identity [1-7] and right identity [1-8].)
  Note that $\phi_L$ and $\phi_R$ coincide on a neighborhood
  of the boundary of this interval.  
  In order to define the flip we need to implement the 
  diffeomorphism $\phi := \phi_R \circ (\phi_L)^{-1}$;
  this is possible, because $\phi$ acts as the identity in a neighborhood of the
  boundary of the interval.
  However, there is no canonical implementation. 
  In order to choose these implementations consistently we 
  proceed as follows.
  The group $\Diff_0([0,1])$, of diffeomorphisms that are the identity on a 
  neighborhood of the boundary,
  is perfect by~\cite{Mather(1974commutators)}.
  Thus it admits a universal central extension 
  $\pi \colon \widetilde{\Diff_0}([0,1]) \to \Diff_0([0,1])$.
  For any net $\cala$ we can implement any $\varphi \in \Diff_0([0,1])$ 
  by a unitary  $U_\varphi \in U(\cala(I))$;
thus $U_\varphi a {U_\varphi}^* = \cala(\varphi)(a)$ for
  all $a \in \cala([0,1])$. 
  Moreover, $U_\varphi$ is unique modulo the center $Z(\cala([0,1]))$.
  This induces a unique homomorphism 
  $U \colon \widetilde{\Diff_0}([0,1]) \to \cala([0,1])$ such that
  $U_{\tilde{\varphi}} a {U_{\tilde{\varphi}}}^* = \cala(\pi(\tilde{\varphi}))(a)$
  for all $a \in \cala([0,1])$, and all
  $\tilde{\varphi} \in \widetilde{\Diff_0}([0,1])$.
  The uniqueness of this map implies that it is compatible with 
  tensor products of nets.  %\CDcomm{Maybe this sentence needs explanation.}
  Now we choose $\tilde{\phi}$ such that $\pi(\tilde{\phi}) = \phi$.
  The flip is the map induced by the action of $U_{\tilde{\phi}}$.  (Though the flip map depends on the choice of lift $\tilde{\phi}$, the overall resulting 3-category does not, up to isomorphism---see Remark~\ref{remark-equivalence} below.)
      %\CDcomm{I think we'll be able to polish the exposition in this paragraph at some point.}

  We note that $U_{\tilde \phi}$ also provides an isomorphism
  \[\cb{
  \displaytikz{
  \path
  (0,0) coordinate (1)
  (4,0) coordinate (2);
  \draw (1) -- (2);
  \draw (1) node[littlenode] {
  \displaytikz{
	\sh(0,0);
	\boxt(.75*.75,.75,.25*.75,1)1;
	\sh(0,1);
	\bubble(.75*.75,.75,.25*.75)1;
	\sh(1.5,0);
	\bxb(1,1)1;
	\begin{pgfonlayer}{foreground}
		\bull(0,.5);\bull(1.5,.5);\bull(.75,1+\cw);
		\dfct(.375,1);\dfct(.75,0);\idfct(1.125,1+\cw);
		\dfct(2,0);\dfct(2,1);
	\end{pgfonlayer}
  }
  };
  \draw (2) node[littlenode] {
  \displaytikz{
	\sh(1,0);
	\boxt(.25*.75,.75,.75*.75,1)1;
	\sh(1,1);
	\bubble(.25*.75,.75,.75*.75)1;
	\sh(0,0);
	\bxb(1,1)1;
	\begin{pgfonlayer}{foreground}
		\bull(0,.5);\bull(1.75,1+\cw); \bull(2.5,.5);
		\dfct(2.125,1);\dfct(1.75,0);\idfct(1.375,1+\cw);
		\dfct(.5,0);\dfct(.5,1);
	\end{pgfonlayer}
  }
  };
  }
  }
  \] 
  for general sectors.
  This generalization of the flip will be helpful in Lemma~M below.
   
\subsubsection*{Construction of the swap [2-15]}
  The domain and target of the swap isomorphism are obtained by twisting the
  vacuum sector of the same identity defect with different
  diffeomorphisms.
  Using implementation of diffeomorphisms for the net
  in question we can implement the difference between these
  diffeomorphism  (as in the construction of the flip [2-18]) 
  and see that domain and target of the swap are
  indeed isomorphic.
  However, there is a priori no preferred 
  implementation and therefore no canonical choice for
  the swap.
  Because every net can be canonically written as a direct sum
  of irreducible nets, it 
  suffices to determine the swap for irreducible nets.
  In this case there is up to phase a unique implementation.
  Therefore it remains to determine the phase of the swap
  in this case.
  Consider the diagram
  \[
     \begin{tikzpicture}
     \path 
     (0,0) coordinate (1)
     (3.5,0) coordinate (2)
     (0,-2) coordinate (3)
     (3.5,-2) coordinate (4);
	\draw (1) -- (2) -- (4) -- (3) -- (1);
	\draw (1) node[littlenode] {
	\displaytikz{
	\sh(1,0);
	\boxt(.25*.75,.75,.75*.75,1)2;
	\sh(1,1);
	\bubble(.25*.75,.75,.75*.75)2;
	\sh(0,0);
	\bxb(1,1)2;
	\begin{pgfonlayer}{foreground}
		\bull(0,.5);\bull(1.75,1+\cw); \bull(2.5,.5);
		\idfct(2.125,1);\idfct(1.75,0);\idfct(1.375,1+\cw);
		\idfct(.5,0);\idfct(.5,1);
	\end{pgfonlayer}
	}
	};
	\draw (2) node[littlenode] {
	\displaytikz{
	\sh(0,0);
	\boxt(.75*.75,.75,.25*.75,1)2;
	\sh(0,1);
	\bubble(.75*.75,.75,.25*.75)2;
	\sh(1.5,0);
	\bxb(1,1)2;
	\begin{pgfonlayer}{foreground}
		\bull(0,.5);\bull(1.5,.5);\bull(.75,1+\cw);
		\idfct(.375,1);\idfct(.75,0);\idfct(1.125,1+\cw);
		\idfct(2,0);\idfct(2,1);
	\end{pgfonlayer}
	}
	};
	\def\cuspwidthtmp{\cuspwidth};
	\def\cuspwidth{.28};	%c%
	\draw (3) node[littlenode] {
	\displaytikz{
\sh(0,1);
\bubble(.25,.5,1.75)2;
\sh(0,0);
\bxb(2.5,1)2;
\sh(0,1);
\draw \bend(0,1,.25,1.5,.75);
\def\cuspwidth{.18}; %c% coded by hand!!!
\sh(0,0);
\bubble(.25,1,1.25)2;
	\begin{pgfonlayer}{foreground}
		\bull(1,-\cw);\idfct(2,0);\idfct(.5,-\cw);
		\idfct(2,1); \bull(.615,1+.14);\idfct(.385,1+.14); %c%
		\idfct(1,1);\bull(1.5,1);
	\end{pgfonlayer}
	}
	};
\def\cuspwidth{.18}; %c%
		\draw (4) node[littlenode] {
	\displaytikz{
	\sh(0,0);
	\boxt(.25*.75,.75,.75*.75,1)2;
	\sh(0,1);
	\bubble(.25*.75,.75,.75*.75)2;
	\sh(1.5,0);
	\bxb(1,1)2;
	\begin{pgfonlayer}{foreground}
		\bull(0,.5);\bull(1.5,.5);\bull(.75,1+\cw);
		\idfct(.375,1+\cw);\idfct(.75,0);\idfct(1.125,1);
		\idfct(2,0);\idfct(2,1);
	\end{pgfonlayer}
	}
	};	     
	\path (1) -- (2) node[la] {\flip}; %Note emph to counteract Lemma style.
	\path (3) -- (4) node[lb] {\quasiidentity};
	\path (3) -- (1) node[la] {\quasiidentity};
	\path (4) -- (2) node[lb] {\swap};
     \end{tikzpicture}
  \]
  Here the lower horizontal map is obtained by twisting the quasi-identity 
  \[\cb{
  \displaytikz{
  \path
  (0,0) coordinate (1)
  (3.5,0) coordinate (2);
  \draw (1) -- (2);
  \draw (1) node[littlenode] {
  \displaytikz{
	\sh(0,0);
	\bxb(2,1)2;
	\bubble(.25,1,.75)2;
	\sh(0,1);
	\bubble(.25,1,.75)2;
	\begin{pgfonlayer}{foreground}
		\bull(1,-\cw); \bull(1,1+\cw);
		\idfct(.5,-\cw);\idfct(.5,1+\cw);\idfct(1.5,0);\idfct(1.5,1);
	\end{pgfonlayer}
  }
  };
  \draw (2) node[littlenode] {
  \displaytikz{
	\sh(0,0);
	\bxb(1,1)2;
	\sh(1,0);
	\bxb(1,1)2;
	\begin{pgfonlayer}{foreground}
		\idfct(.5,0);\idfct(.5,1);\idfct(1.5,0);\idfct(1.5,1);
	\end{pgfonlayer}
  }
  };
  }
  }
  \] 
  while the left vertical map is obtained by applying the quasi-identity
  to the left identity
\[\cb{
\displaytikz{
	\sh(0,0);
	\bxb(2,1)2;
	\sh(0,1);
	\bubble(.25,1,.75)2;
	\begin{pgfonlayer}{foreground}
		\bull(1,1+\cw);
		\idfct(.5,1+\cw);\idfct(1,0);\idfct(1.5,1);
	\end{pgfonlayer}
}
}
\]
  (In particular, for the lower horizontal map in the above square diagram, the inner diffeomorphism
  of the two upper diffeomorphisms in the lower left hand item 
  is added by  the quasi-identity; by contrast, for the left vertical map 
  the outer of those two diffeomorphisms is added by the quasi-identity.) 
  The phase of the swap is now fixed by requiring
  the above diagram to commute.
  
  %This construction of the swap produces a specific implementation of a 
  %diffeomorphism. 
  %Similar to the case of the flip, this also yields an isomorphism 
  %\ABcomm{Another picture}
  %\begin{equation*}
  %  \text{\tiny Picture of swap, no shading, general defects, sectors.}
  %\end{equation*}
  %for general sectors.

  The defining diagram for the swap generalizes as follows.  

  \begin{lemmaletter}[M] %formerly Z
%  \begin{lemma}[Z]
    The following diagram commutes
     \[
     \begin{tikzpicture}
     \path 
     (0,0) coordinate (1)
     (3,0) coordinate (2)
     (0,-2) coordinate (3)
     (3,-2) coordinate (4);
	\draw (1) -- (2) -- (4) -- (3) -- (1);
	\draw (1) node[littlenode] {
	\displaytikz{
	\sh(1,0);
	\boxt(.25*.75,.75,.75*.75,1)1;
	\sh(1,1);
	\bubble(.25*.75,.75,.75*.75)1;
	\sh(0,0);
	\bxb(1,1)2;
	\begin{pgfonlayer}{foreground}
		\bull(0,.5);\bull(1.75,1+\cw); \bull(2.5,.5);
		\dfct(2.125,1);\dfct(1.75,0);\idfct(1.375,1+\cw);
		\idfct(.5,0);\idfct(.5,1);
	\end{pgfonlayer}
	}
	};
	\draw (2) node[littlenode] {
	\displaytikz{
	\sh(0,0);
	\boxt(.75*.75,.75,.25*.75,1)2;
	\sh(0,1);
	\bubble(.75*.75,.75,.25*.75)2;
	\sh(1.5,0);
	\bxb(1,1)1;
	\begin{pgfonlayer}{foreground}
		\bull(0,.5);\bull(1.5,.5);\bull(.75,1+\cw);
		\idfct(.375,1);\idfct(.75,0);\idfct(1.125,1+\cw);
		\dfct(2,0);\dfct(2,1);
	\end{pgfonlayer}	
	}
	};
	\def\cuspwidthtmp{\cuspwidth};
	\def\cuspwidth{.28};	%c%
	\draw (3) node[littlenode] {
	\displaytikz{	
\sh(0,1);
\bubble(.25,.5,1.75)1;
\sh(0,0);
\bxb(2.5,1)1;
\sh(0,1);
\draw \bend(0,1,.25,1.5,.75);
\def\cuspwidth{.18}; %c% coded by hand!!!
\sh(0,0);
\bubble(.25,1,1.25)1;
	\begin{pgfonlayer}{foreground}
		\bull(1,-\cw);\dfct(2,0);\idfct(.5,-\cw);
		\dfct(2,1); \bull(.615,1+.14);\idfct(.385,1+.14); %c%
		\idfct(1,1);\bull(1.5,1);
	\end{pgfonlayer}
	}
	};
\def\cuspwidth{.18}; %c%
		\draw (4) node[littlenode] {
	\displaytikz{
	\sh(0,0);
	\boxt(.25*.75,.75,.75*.75,1)2;
	\sh(0,1);
	\bubble(.25*.75,.75,.75*.75)2;
	\sh(1.5,0);
	\bxb(1,1)1;
	\begin{pgfonlayer}{foreground}
		\bull(0,.5);\bull(1.5,.5);\bull(.75,1+\cw);
		\idfct(.375,1+\cw);\idfct(.75,0);\idfct(1.125,1);
		\dfct(2,0);\dfct(2,1);
	\end{pgfonlayer}	
	}
	};	     
%	\path (1) -- (2) node[la] {\emph{\flip}}; %Note emph to counteract Lemma style.
%	\path (3) -- (4) node[lb] {\emph{\quasiidentity}};
%	\path (3) -- (1) node[la] {\emph{\quasiidentity}};
%	\path (4) -- (2) node[lb] {\emph{\swap}};
     \end{tikzpicture}
     \]
  \end{lemmaletter}
%  \end{lemma}

  Later we will only need Lemma~M for vacuum sectors, but the proof
  of the more general statement is a bit cleaner. 

  \begin{proof}
    We denote the non-vacuum sector in the diagram by ${}_D H_E$.
    Here ${}_\cala D_\calb$ and ${}_\cala E_\calb$ are defects.
    If $D = \id_\cala$ and $H = H_0(\id_\cala)$, 
    then the diagram commutes by the construction
    of the swap.
    This also implies that the diagram commutes if
    $H = H_0(\id_\cala) \ox \ell$ for any Hilbert space $\ell$.
    For $\e > 0$ we now use the subinterval $I_\e = [- 1.5 + \e, 1.5 - \e]$
    of the circle $S^1_6$ bounding the non-vacuum sector.
    There exits an $\cala(I_\e)$-linear isometry
    $U_\e \colon H_0(\id_\cala) \ox \ell \to H \ox \ell$.
    If $\e$ is sufficiently small, then $U_\e$ will commute
    with all four sides of the above diagram.
    Thus the square also commutes for $H \ox \ell$ and thus for $H$ itself. 
  \end{proof}

\begin{remark}\label{remark-equivalence}
  The flip isomorphism [2-18] and the swap isomorphism [2-15] depended on the choice of lift $\tilde{\phi}$ of the diffeomorphism $\phi := \phi_R \circ (\phi_L)^{-1}$; changing the lift changes the flip by a scalar and the swap by the inverse scalar.  However, the symmetric monoidal 3-categories resulting from distinct choices are canonically isomorphic, as follows.  (Here the isomorphism will be given by a functor $A \ra B$ of dicategory objects in symmetric monoidal categories, which is a triple of functors $A_0 \ra B_0$, $A_1 \ra B_1$, $A_2 \ra B_2$, together with eight symmetric monoidal natural transformations, comparing each piece of 1-data for $A$ with the corresponding piece of 1-data for $B$, such that eighteen squares of natural transformations commute, one for each piece of 2-data.)  
  
  Recall that the left identity [1-7] and right identity [1-8] depended on a choice of diffeomorphism $\Phi_L$.  Assume for a moment the diffeomorphism $\Phi_L$ is fixed; let $\CN[\Phi_L,\tilde{\phi}]$ denote the dicategory object resulting from the choice $\tilde{\phi}$ of lift of $\phi$, and let $\CN[\Phi_L,\bar{\phi}]$ denote the dicategory object resulting from the alternative choice of lift $\bar{\phi}$ of $\phi$.  The ratio $\bar{\phi} / \tilde{\phi}$ is an invertible scalar.  There is therefore a canonical isomorphism of dicategory objects $\CN[\Phi_L,\tilde{\phi}] \rightarrow \CN[\Phi_L,\bar{\phi}]$ which is the identity functor on all 0-data, and the identity natural transformation on the 1-data [1-1] through [1-7], but is the natural transformation given by scalar multiplication by $(\bar{\phi} / \tilde{\phi})^{-1}$ on the 1-data [1-8].  These natural transformations commute with all the 2-data.
  
  We can moreover remove the dependence of the 3-category $\CN$ on the choice of the diffeomorphism $\Phi_L$.  Suppose $\Phi_L'$ is an alternative choice of diffeomorphism suitable for the construction of the left identity.  Choose a lift $\tilde{\Lambda} \in \widetilde{\Diff_0}[.5,1.5]$ of the diffeomorphism $\Lambda := \Phi_L' |_{[.5,2.5]} \circ (\Phi_L)^{-1} |_{[.5,1.5]}$.  Recall from the section on the right identity [1-8] that $\Phi_R$ is the reflection of the diffeomorphism $\Phi_L$; similarly we let $\Phi_R'$ be the corresponding reflection of $\Phi_L'$.  The lift $\tilde{\Lambda}$ determines by reflection a lift $\tilde{\Rho}$ of $\Rho := \Phi_R' |_{[1.5,3.5]} \circ (\Phi_R)^{-1} |_{[1.5,2.5]}$.  As before, let $\tilde{\phi}$ be a lift of $\phi := \phi_R \circ (\phi_L)^{-1}$.  Note that the diffeomorphisms $\phi_L : [1.5,3.5] \rightarrow [1.5,2.5]$ and $\phi_R : [1.5,3.5] \rightarrow [1.5,2.5]$ were, up to a shift, determined by restrictions of $\Phi_L$ and $\Phi_R$; more specifically $\phi_L(x) = \Phi_L(x - 1) + 1$ and $\phi_R(x) = \Phi_R(x)$.  Therefore, the product $\widetilde{\phi'} := \tilde{\Rho} \cdot \tilde{\phi} \cdot [\tilde{\Lambda}]^{-1}$ is a lift of $\phi' := \phi_R' \circ (\phi_L')^{-1}$, where $[-]$ denotes the canonical shift isomorphism from $\widetilde{\Diff_0}[.5,1.5]$ to $\widetilde{\Diff_0}[1.5,2.5]$.
  
  Now let $U_{\tilde{\Lambda}} \in \cala([.5,1.5])$ be the unitary, associated to $\tilde{\Lambda}$, implementing $\cala(\Lambda)$; here $U: \widetilde{\Diff_0}([.5, 1.5]) \rightarrow \cala([.5,1.5])$ is the homomorphism described in the section on the flip.  The action of $U_{\tilde{\Lambda}}$ provides a natural isomorphism from the left identity of $\CN[\Phi_L,\tilde{\phi}]$ (a functor $\CN_1 \rightarrow \CN_2$) to the left identity of $\CN[\Phi_L',\widetilde{\phi'}]$ (also a functor $\CN_1 \rightarrow \CN_2$).  Similarly, the action of the unitary $U_{\tilde{\Rho}} \in \calb([1.5,2.5])$ implementing $\calb(\Rho)$, provides a natural isomorphism from the right identity of $\CN[\Phi_L,\tilde{\phi}]$ to the right identity of $\CN[\Phi_L',\widetilde{\phi'}]$.  Altogether, this constructs an isomorphism of dicategory objects $\CN[\Phi_L,\tilde{\phi}] \rightarrow \CN[\Phi_L',\widetilde{\phi'}]$ which is the identity functor on all 0-data, is the identity natural transformation on the 1-data [1-1] through [1-6], and is the natural transformation $U_{\tilde{\Lambda}}$ on [1-7] and the natural transformation $U_{\tilde{\Rho}}$ on [1-8]---these natural transformations commute with all the 2-data.  %Note that the isomorphism we constructed $\CN[\Phi_L,\tilde{\phi}] \rightarrow \CN[\Phi_L',\widetilde{\phi'}]$ does depend on the choices of lift $\tilde{\Lambda}$, but it depends on it only up to canonical 2-isomorphism (given by scalar multiplication by the ratio of the lift $\tilde{\Lambda}$ to an alternative lift $\bar{\Lambda}$).
  \end{remark}

  %===================================================================

\pagebreak
\section{Coherence axioms for compatibility transformations} \label{sec-3axioms}%{Axioms}

\subsection{Axioms for vertical identity and composition}

  \begin{axprop} Axiom [3-1] is satisfied. \end{axprop}
  
  \begin{proof}
    Axiom [3-1] asserts that top and bottom identity agree
    in the case where both sectors are vacuum sectors.
    This holds because the corresponding statement is already true 
    for Connes fusion.
  \end{proof}

  \begin{axprop} Axioms [3-2] and [3-3] are satisfied. \end{axprop}

  \begin{proof}
    Axioms [3-2] and [3-3] assert that top and bottom identity
    are compatible with the vertical associator.
    This holds because the corresponding statement is already true 
    for Connes fusion. 
    %In fact, if we suppress the associator for Connes fusion, then axioms [3-2] and [3-3] are  tautologies.  
  \end{proof}

  \begin{axprop} Axiom [3-4] is satisfied. \end{axprop}
    
  \begin{proof}
    Axiom [3-4] asserts that the vertical associator satisfies
    the pentagon identity.
    This holds because the associator for Connes fusion
    satisfies the pentagon identity.
  \end{proof}

  \begin{remark*}
    It is possible to base the definition of sectors 
    on a square rather than a hexagon and to define
    vertical fusion using just a side of this square
    (rather than half of the hexagon).
    Then our pictures would become a little simpler,
    but vertical composition would require a diffeomorphism
    and the associator would then involve this diffeomorphism.
    Axioms [3-1] through [3-4] would be more cumbersome to prove
    in such a set-up.
  \end{remark*} 
  
\pagebreak  
\subsection{Axioms for horizontal composition and whiskers}

  \begin{axprop} Axiom [3-5] is satisfied. \end{axprop} 

  \begin{proof}
  The argument is summarized by the following diagram.
  \[
  \def\axnum{5}
  \begin{tikzpicture}[scale=.4ex]
  \axgrid{\axnum};
  \foreach \x in {1,...,6}
  \draw (\x) node[littlenode] {\innertikz{\ax{\axnum}-{\x}}};
  \comm{\axnum};
  \labelout{\axnum};
  \end{tikzpicture} \quad 
  \]
  Each of the four corners of the diagram denotes a functor 
  $\CN_2 \x_{\CN_0} \CN_1 \to \CN_2$.
  Each of the four lines on the boundary of the diagram
  denotes a natural isomorphism determined by its label.
  These four natural isomorphisms 
  (simply referred to as \emph{maps} for brevity) are 
  explained in more detail as follows.
  \begin{itemize}
  \item The map labeled ``whisker'' is the ``right dewhisker'' [2-6].
  \item The horizontal map labeled ``identity'' is obtained by applying
    the top identity [2-1] %(or equivalently the bottom identity [2-2])
    to the left half of the item in the top left corner of the diagram.
  \item The right vertical map labeled ``identity'' is obtained by
    applying the top identity [2-1].
  \item The horizontal map labeled ``expand'' is obtained by
    applying the left 
    (or equivalently right) vertical identity expansion [2-4]
    to the top half of the item in the lower left corner of the diagram.
  \end{itemize}
  Axiom [3-5] asserts that the boundary of this diagram commutes,
  that is, if we start at some corner of the diagram and compose the
  four maps along the boundary of the diagram then we should obtain
  the identity natural transformation on the functor corresponding
  to the corner where we started.
  (We remark that in~\cite{Douglas-Henriques(Internal-bicategories)}, the axiom is rotated by $-\pi/2$ from the version depicted above.)
  Now observe that each corner of the diagram can also be viewed as
  determining a functor ${\CN_{\sim}}_2 \x_{\CN_0} \CN_1 \to {\CN_{\sim}}_2$.
  And similarly each map on the boundary describes a natural isomorphism
  between these functors.
  Moreover, the question whether the diagram commutes or not is
  invariant under this change from $\CN_2$ to ${\CN_{\sim}}_2$.
  But the ${\CN_\sim}_2$-version has the advantage that because of
  Proposition~\ref{prop:CNsim-as-subobjects-of-CNsimeq}, it suffices 
  to check the commutativity of the diagram in the
  case when the darker shaded sector is an identity sector,
  not an arbitrary sector.  
  Therefore we  can and will assume
  that this sector is also an identity sector.
  Under this assumption, the internal maps and nodes of
  the diagram make sense.
  Axiom [3-5] will follow, once we have shown that the three cells
  in the interior commute.

  \emph{The cell labeled  ``D''.}
    The composition of the maps not labeled ``whisker'' around this cell 
    is the definition of the map labeled ``whisker''. 
    Thus this cell commutes by definition. 
    This is the reason for the label ``D''.

  \emph{The cell labeled ``N''.}
    The map that is counterclockwise after the map labeled ``identity''
    is the left 
    (or equivalently right) vertical identity expansion [2-4].
    The map clockwise after the map labeled ``identity''
    is the vertical fusion of the identity (on the top) and
    the left vertical identity expansion (on the bottom).
    The remaining map is a top identity (as is the map labeled identity). 
    Thus this cell commutes by the naturality of the top identity map.
    This is the reason for the label ``N''.

  \emph{The cell labeled ``C''.}
    Consider the item in the lower left corner of the diagram.
    Here we can apply  the left (or right) vertical identity 
    expansion [2-4] 
    to both the bottom and top half of this item.
    These applications do not interact with each other and can be done 
    in any order or simultaneously.
    All three maps on the boundary of this cells are obtained from these
    commuting operations. 
    The cell therefore commutes.
    We recorded this in the diagram by the label ``C'', 
    for commuting operations. 
  \end{proof}

  Formally the proofs of the remaining axioms will be very 
  similar to the proof of axiom [3-5].
  We will however not repeat the arguments in every case in detail.
  In particular, we will trust the reader to determine the correct 
  maps from our pictures.
  Moreover, the trick that allows us to assume that some sector is
  not an arbitrary sector but an identity sector 
  (by replacing $\CN_2$ temporarily with ${\CN_\sim}_2$ or ${\CN^\sim}_2$)
  will be used very often in the remainder of this paper.
  We will always refer to this as the \emph{corner trick} and
  indicate the sector to which it is applied by a darker shading. 

%  There will however be a number of new arguments that 
%  will be discussed in detail whenever they are used for the first time.

%\Needspace{3\baselineskip}

\pagebreak
  \begin{axprop} Axiom [3-6] is satisfied. \end{axprop}   
  
  \begin{proof}
    This axiom asserts that the following diagram commutes
    \[
    \def\axnum{6}
    \begin{tikzpicture}[scale=.45ex]
    \axgrid{\axnum};
    \foreach \x in {1,...,5}
    \draw (\x) node[littlenode] {\innertikz{\ax{\axnum}-{\x}}};
    \comm{\axnum};
    \labelout{\axnum};
    \end{tikzpicture} \quad .
    \]
    The diagonal maps are the two half-switches from the definition
    of the switch isomorphism [2-8]. 
    Thus the cell labeled ``D'' commutes by definition.
    There is a mirror symmetry between the two remaining cells.
    Thus it suffices to prove that [3-6a] commutes: this is the content 
    of the next lemma. 
  \end{proof}

%\Needspace{3\baselineskip}
%g%\pagebreak
  \begin{axlem} The diagram [3-6a] commutes. \end{axlem} 

  \begin{proof}
    The argument is similar to the proof of
    Axiom [3-5] and is summarized as
    \[
    \def\axnum{6a}
    \begin{tikzpicture}[scale=.25ex]
    \axgrid{\axnum};
    \foreach \x in {1,...,5}
    \draw (\x) node[littlenode] {\innertikz{\ax{\axnum}-{\x}}};
    \comm{\axnum};
    \labelout{\axnum};
    \end{tikzpicture}  \quad \qedhere
     \] 
  \end{proof}

\pagebreak
%\Needspace{3\baselineskip}
  \begin{axprop} Axiom [3-7] is satisfied. \end{axprop} 

  \begin{proof}
    The argument is summarized as
    \[
    \def\axnum{7}
    \begin{tikzpicture}[scale=.23ex]
    \axgrid{\axnum};
    \foreach \x in {1,...,16}
    \draw (\x) node[littlenode] {\innertikz{\ax{\axnum}-{\x}}};
    \comm{\axnum};
    \labelout{\axnum};
    \end{tikzpicture} \quad 
    \]
    The boundary of this diagram is a square, not a hexagon as in~\cite{Douglas-Henriques(Internal-bicategories)},
    because we suppress the vertical associativity 
    isomorphisms.
    To help the reader to decode the precise meaning of the items of this
    diagram we give a more detailed picture of the left top corner
    where the bullets are added (even though these can be reconstructed 
    from the form of the picture).
\[
%\begin{tikzpicture}
\displaytikz{
		\sh(0,0);
		\bx(1,1)3;
		\sh(0,1);
		\qhexdl(1,1)1;
		\sh(0,2+\gap);
		\hkdl(1+.5*\gap,1*\fudge)2,1,3;
		\sh(1+\gap,0);
		\bx(1,2)2;
		\sh(1+.5*\gap,2+\gap);
		\bx(1+.5*\gap,1*\fudge)2;
		\bull(0,0);\bull(1,0);\bull(1+\gap,0);\bull(2+\gap,0);
		\bull(0,1);\bull(.33,1);\bull(.67,1);\bull(1,1);\bull(1+\gap,1);
		\bull(1,2);\bull(1+\gap,2);\bull(1.5+\gap,2);\bull(2+\gap,2);
		\bull(0,2+\gap+\fudge);\bull(1+.5*\gap,2+\gap);\bull(1+.5*\gap,2.5+\gap);\bull(1+.5*\gap,2+\gap+\fudge);\bull(1.5+\gap,2+\gap);\bull(2+\gap,2+\gap);\bull(2+\gap,2+\gap+\fudge);
		\bull(.7,2);  % 1/tan(hexangleB) %!% This will break if hexangleB changes.
		\bull(0,1+\gapfactorB*\gap); 
		\bull(.597,2+\gap); % tan(90-hexangleB)*(1+(\gap)(1-\gapfactor)) %!% This will break if gapfactorB changes.
		\dfct(.5,0);\dfct(1.5+\gap,0);\dfct(.5,1);\dfct(.85,2);\dfct(1.25+\gap,2);\dfct(.85,2+\gap);\dfct(1.25+\gap,2+\gap);\dfct(.5+.25*\gap,2+\gap+\fudge);\dfct(1.5+.75*\gap,2+\gap+\fudge);
}
%\end{tikzpicture}
\]
    Each cell of the diagram commutes for the reason indicated in the 
    diagram. 
  \end{proof}
  
%\pagebreak
%\Needspace{3\baselineskip}

  \begin{axprop} Axiom [3-8] is satisfied. \end{axprop} 

  \begin{proof}
    Consider the diagram 
    \[
    \def\axnum{8}
    \begin{tikzpicture}[scale=.38ex]
    \axgrid{\axnum};
    \foreach \x in {1,...,8}
    \draw (\x) node[littlenode] 
        {\innertikz{\ax{\axnum}-{\x} \draw (1,0) node [above] {};}};
    \comm{\axnum};
    \labelout{\axnum};
    \end{tikzpicture}
    \]
    Here, the left isomorphism labeled $*$ is defined using the corner trick (as in the construction of 
    the right dewhisker [2-6]) to be the following composite:
    \[
    \def\axnum{8s}
    \begin{tikzpicture}[scale=.7ex]
    \axgrid{\axnum};
    \foreach \x in {1,...,4}
    \draw (\x) node[littlenode] {\innertikz{\ax{\axnum}-{\x}}};
    \comm{\axnum};
    \labelout{\axnum};
    \end{tikzpicture}
    \]
    The second isomorphism labeled $*$ is defined similarly.
    There is a horizontal symmetry between the cells labeled
    [3-8a] and [3-8c].
    Thus it suffices to prove that [3-8a] and [3-8b] commute.
    This is the content of the next two lemmas.
  \end{proof}

%\pagebreak
%\Needspace{3\baselineskip}

%g%\pagebreak
  \begin{axlem} The diagram [3-8a] commutes. \end{axlem} 

  \begin{proof}
    Applying the corner trick twice,
    the cell [3-8a] can be filled as follows:
    \[
    \def\axnum{8a}
    \begin{tikzpicture}[scale=.34ex]
    \axgrid{\axnum};
    \foreach \x in {1,...,14}
    \draw (\x) node[littlenode] {\innertikz{\ax{\axnum}-{\x} 
         \draw (1,0) node [above] {};}};
    \comm{\axnum};
    \labelout{\axnum};
    \end{tikzpicture} \quad \qedhere
    \]
  \end{proof}

%%\pagebreak
%\Needspace{3\baselineskip}

\pagebreak
  \begin{axlem} The diagram [3-8b] commutes. \end{axlem} 

  \begin{proof}
    Using the corner trick twice, we fill the cell [3-8b] as follows:
    \[
    \def\axnum{8b}
    \begin{tikzpicture}[scale=.4ex]
    %\draw[step=.5] (0,0) grid (6,6);
    \axgrid{\axnum};
    \foreach \x in {1,...,10}
    \draw (\x) node[littlenode] {\innertikz{\ax{\axnum}-{\x} 
     \draw (1,0) node [above] {};}};
    \comm{\axnum};
    \labelout{\axnum};
    \end{tikzpicture} \quad \qedhere
    \]
  \end{proof}

%%\pagebreak % This is good
\subsection{Axioms for horizontal associators}
  
%\Needspace{3\baselineskip}

  \begin{axprop} Axioms [3-9] and [3-10] are satisfied. \end{axprop}

  \begin{proof}
    Both are a consequence of the associativity of 
    the $1 \boxtimes 1$-isomorphism, in the form of Lemma~J.
  \end{proof}

%\Needspace{3\baselineskip}

%g%\pagebreak
  \begin{axprop} Axiom [3-11] is satisfied. \end{axprop}
  \begin{proof}
    Using the corner trick this axiom is proved by this diagram:
    \[
    \def\axnum{11}
    \begin{tikzpicture}[scale=.36ex]
    \axgrid{\axnum};
    \foreach \x in {1,...,12}
    \draw (\x) node[littlenode] {\innertikz{\ax{\axnum}-{\x} 
      \draw (1.5,0) node [above] {};}};
    \comm{\axnum};
    \labelout{\axnum};
    \end{tikzpicture} \quad 
    \]
    The cells labeled ``J'' commute by Lemma J.
  \end{proof}

%\Needspace{3\baselineskip}

  \begin{axprop} Axiom [3-12] is satisfied. \end{axprop} 

  \begin{proof}
    The formulation of axiom [3-12] simplifies from a hexagon to a square
    because we suppress the whisker associator [2-11].
    Using the corner trick we can fill in this square as follows:
    \[
    \def\axnum{12}
    \begin{tikzpicture}[scale=.57ex]
    \axgrid{\axnum};
    \foreach \x in {1,...,13}
    \draw (\x) node[littlenode] {\innertikz{\ax{\axnum}-{\x} 
     \draw (1.5,0) node [above] {};}};
    \comm{\axnum};
    \labelout{\axnum};
    \end{tikzpicture} \quad \qedhere
    \]
  \end{proof}

%\Needspace{3\baselineskip}
\pagebreak
  \begin{axprop} Axiom [3-13] is satisfied. \end{axprop} 

  \begin{proof}
    The axiom follows from the commutativity of the following diagram:
    \[
    \def\axnum{13}
    \begin{tikzpicture}[scale=.3ex]
    \axgrid{\axnum};
    \foreach \x in {1,...,7}
    \draw (\x) node[littlenode] {\innertikz{\ax{\axnum}-{\x} 
      \draw (1.5,-.3) node [above] {};}};
    \comm{\axnum};
    \labelout{\axnum};
    \end{tikzpicture} \quad 
    \]
    There is a horizontal symmetry between the cells labeled 
    [3-13a] and [3-13b].
    Thus it remains to prove that [3-13a] commutes. 
    This is the content of the next lemma.
  \end{proof}

%\Needspace{3\baselineskip}
\pagebreak
  \begin{axlem} The diagram [3-13a] commutes. \end{axlem}

  \begin{proof}
    Using the corner trick we can fill in the diamond as follows:
    \[
    \def\axnum{13a}
    \begin{tikzpicture}[scale=.5ex]
    \axgrid{\axnum};
    \foreach \x in {1,...,7}
    \draw (\x) node[littlenode] {\innertikz{\ax{\axnum}-{\x} 
       \draw (1.5,-.3) node [above] {};}};
    \comm{\axnum};
    \labelout{\axnum};
    \end{tikzpicture} \quad \qedhere
    \]
  \end{proof}

%\Needspace{3\baselineskip}

  \begin{axprop} Axiom [3-14] is satisfied. \end{axprop} 

  \begin{proof}
    The axiom follows from the commutativity of the diagram 
    \[
    \def\axnum{14}
    \begin{tikzpicture}[scale=.8ex]
    \axgrid{\axnum};
    \foreach \x in {1,...,5}
    \draw (\x) node[littlenode] {\innertikz{\ax{\axnum}-{\x} 
       \draw (1.5,-.3) node [above] {};}};
    \comm{\axnum};
    \labelout{\axnum};
    \end{tikzpicture} \quad .
    \]
    There is a horizontal symmetry between the cells [3-14a] and [3-14b].
    Thus it suffices to show that [3-14a] commutes.
    This is the content of the next lemma.
  \end{proof}

%%\pagebreak
%\Needspace{3\baselineskip}

%g%\pagebreak
  \begin{axlem} The diagram [3-14a] commutes. \end{axlem} 

  \begin{proof}
    Using the corner trick twice we can fill in [3-14a] as follows:
    \[    
    \def\axnum{14a}
    \begin{tikzpicture}[scale=.38ex]
    %\draw[step=.5] (0,0) grid (6,7);
    \axgrid{\axnum};
    \foreach \x in {1,...,12}
    \draw (\x) node[littlenode] {\innertikz{\ax{\axnum}-{\x} 
        \draw (1.5,-.3) node [above] {};}};
    \comm{\axnum};
    \labelout{\axnum};
    \end{tikzpicture}  \quad \qedhere
    \]
  \end{proof}

%\Needspace{3\baselineskip}

  \begin{axprop} Axiom [3-15] is satisfied. \end{axprop} 

  \begin{proof}
    This axiom follows from the associativitiy of the 
    $1 \boxtimes 1$-isomorphism (Lemma J):
    \[
    \def\axnum{15}
    \begin{tikzpicture}[scale=.5ex]
    \axgrid{\axnum};
    \foreach \x in {1,...,4}
    \draw (\x) node[littlenode] {\innertikz{\ax{\axnum}-{\x} 
        \draw (2,-.3) node [above] {};}};
    \comm{\axnum};
    \labelout{\axnum};
    \end{tikzpicture}  \quad \qedhere
    \]
  \end{proof}

%\Needspace{3\baselineskip}
%g%\pagebreak
  \begin{axprop} Axiom [3-16] is satisfied. \end{axprop}

  \begin{proof}
    Upon suppression of horizontal associators, 
    three of the nodes of the axiom 
    reduce to the left hand picture below. 
    The remaining two nodes reduce to the right hand picture; 
    both edges between these 
    sets of nodes are the indicated expansion: 
    \[
    \def\axnum{16}
    \begin{tikzpicture}[scale=.5ex]
    \axgrid{\axnum};
    \foreach \x in {1,...,2}
    \draw (\x) node[littlenode] {\innertikz{\ax{\axnum}-{\x} 
       \draw (2,-.3) node [above] {};}};
    \comm{\axnum};
    \labelout{\axnum};
    \end{tikzpicture} \quad \qedhere
    \]
  \end{proof}

%\Needspace{3\baselineskip}

  \begin{axprop} Axiom [3-17] is satisfied. \end{axprop}

  \begin{proof}
    This axiom asserts that the associator for defects [2-12]
    satisfies the pentagon identity.
    This holds because the corresponding statement is
    already true for fusion (or fiber product) of von Neumann
    algebras.   
%!%    \commenta{Do we have a reference?} 
  \end{proof}

  \medskip
  
\subsection{Axioms for horizontal identities}

  The axioms [3-1] to [3-17] do not involve identity defects.
  These axioms are part of the definition of both
  a $2$-category object and a  dicategory object
  (in the $2$-category of symmetric monoidal
  categories). 
  %We remind the reader that our $3$-category of conformal nets is a dicategory object, so 
  The remaining axioms labeled here [3-18]
  to [3-26] are the axioms labeled [D3-18] to [D3-26] in the appendix and correspondingly in~\cite{Douglas-Henriques(Internal-bicategories)}.

%\Needspace{3\baselineskip}
\pagebreak
  \begin{axprop} Axiom [3-18] is satisfied. \end{axprop} 

  \begin{proof}
    This axiom follows from the commutativity of this diagram:
%    \commenta{There is a missing bullet from the top of the item in        left lower corner of the inner triangle.}
%    \[
%    \def\axnum{18}
%    \begin{tikzpicture}[scale=.5ex]
%    \axgrid{\axnum};
%    \foreach \x in {1,...,7}
%    \draw (\x) node[littlenode] {\innertikzb{\ax{\axnum}-{\x} 
%           \draw (1,-.3) node [above] {};}};
%    \comm{\axnum};
%    \labelout{\axnum};
%    \end{tikzpicture} \quad .
%    \]
    \[
    \def\axnum{18new}
    \begin{tikzpicture}[scale=.5ex]
    %\draw[step=.5] (0,0) grid (3,3);
    \axgrid{\axnum};
    \foreach \x in {1,...,6}
    \draw (\x) node[littlenode] 
    	{\displaytikz{\ax{\axnum}-{\x} 
           \draw (1,-.3) node [above] {};}};
    \comm{\axnum};
    \labelout{\axnum};
    \end{tikzpicture} \quad 
    \] 
%    \ABcomm{In the file cnii-def and cnii-code everything I added has      the old name together with a ``new''. So the number used here is ``18new'' and I added the commands ``axgridRnew'', ``commRnew'',  ``labelRnew''. I did it the same way for the later grids.}
    The cells labeled ``N'' and ``D'' commute by naturality and by 
    definition.
    The  remaining cell commutes for a reason that we have not yet encountered, namely by Lemma~K.
  \end{proof}

%\Needspace{3\baselineskip}
%g%\pagebreak
  \begin{axprop} Axiom [3-19] is satisfied. \end{axprop} 

  \begin{proof}
    Axiom [3-19] reduces to a square, 
    because we are suppressing the vertical associator [2-3].
    We can partially fill the square as follows.  For readability and ease of comparison with a subsequent diagram, we use the abbreviated notation for the four corner configurations, and the full bullet and bubble notation for the interior configurations.
    \[
    \def\axnum{19new}
    \begin{tikzpicture}[scale=.6ex]
    %\draw[step=.5] (0,0) grid (3,3);
    \axgrid{\axnum};
    \foreach \x in {1,...,7}
    \draw (\x) node[littlenode] % {\innertikz{\ax{\axnum}-{\x}}};
        	{\displaytikz{\ax{\axnum}-{\x} 
           \draw (1,-.3) node [above] {};}};
    \comm{\axnum};
    \labelout{\axnum};
    \end{tikzpicture}  \quad 
    \]
    For the remaining cell we can use the corner trick 
    and assume that the top left sector is an identity sector.
    This reduces axiom [3-19] to the case where only one of the
    sectors is not an identity sector.
    Using this additional assumption we can fill in Axiom [3-19]
    as follows (using the simpler notation that  suppresses the
    bullets and bubbles for diffeomorphisms):
    \[
    \def\axnum{19}
    \begin{tikzpicture}[scale=.6ex]
    %\draw[step=.5] (0,0) grid (3,3);
    \axgrid{\axnum};
    \foreach \x in {1,...,8}
    \draw (\x) node[littlenode] {\innertikzb{\ax{\axnum}-{\x}}};
    \comm{\axnum};
    \labelout{\axnum};
    \end{tikzpicture}  \quad \qedhere
    \]
  \end{proof}

%%\pagebreak
%\Needspace{3\baselineskip}
\pagebreak
  \begin{axprop} Axiom [3-20] is satisfied. \end{axprop}

  \begin{proof}
    We can fill in axiom [3-20] partially as follows:
%    \[
%    \def\axnum{20}
%    \begin{tikzpicture}[scale=.8ex]
%      \pgftransformyscale{1.1} %!%
%      \axgrid{\axnum};
%      \foreach \x in {1,...,17}
%      \draw (\x) node[littlenode] 
%       {\innertikzb{\ax{\axnum}-{\x} \draw (1,-.3) node [above] {};}};
%      \comm{\axnum};
%      \labelout{\axnum};
%    \end{tikzpicture} \quad .
%    \]
    \[
    \def\axnum{20new}
    \begin{tikzpicture}[scale=.6ex]
      %\draw[step=.5] (0,0) grid (4,4);
      %\pgftransformyscale{1.1} %!%
      \axgrid{\axnum};
      \foreach \x in {1,...,11}
      \draw (\x) node[littlenode] 
       {\displaytikz{\ax{\axnum}-{\x} \draw (1,-.3) node [above] {};}};
      \comm{\axnum};
      \labelout{\axnum};
    \end{tikzpicture} \quad 
    \]
    Here we used the more precise notation using bullets and bubbles.
    The above diagram shows that axiom [3-20]  is equivalent to the 
    commutativity of the remaining hexagon.  
    By the corner trick this hexagon commutes if and only 
    if it commutes for the identity sector, 
    and the hexagon with identity sector commutes if and only if 
    axiom [3-20] commutes for the identity sector.  
    Thus it suffices to establishes axiom [3-20] for the identity sector.
    This follows from the following diagram (where we drop bullets and 
    bubbles from the notation).
    \[
    \def\axnum{20a}
    \begin{tikzpicture}[scale=.64ex]
      %\draw[step=.5] (0,0) grid (4,5);
      \axgrid{\axnum};
      \foreach \x in {1,...,16}
      \draw (\x) node[littlenode] 
        {\innertikz{\ax{\axnum}-{\x} \draw (1,-.3) node [above] {};}};
      \comm{\axnum};
      \labelout{\axnum};
    \end{tikzpicture} \quad \qedhere
    \]
  \end{proof}
  
%%\pagebreak
%\Needspace{3\baselineskip}
\pagebreak
  \begin{axprop} Axiom [3-21] is satisfied. \end{axprop} 

  \begin{proof}
    We can partially fill Axiom [3-21] as follows.
    \[
    \def\axnum{21new}
    \begin{tikzpicture}[scale=.65ex]
      %\draw[step=.5] (0,0) grid (4,6);
      %\pgftransformxscale{1.2} %!%
      \axgrid{\axnum};
      \foreach \x in {1,...,6}
      \draw (\x) node[littlenode] 
        {\displaytikz{\ax{\axnum}-{\x} \draw (1,-.3) node [above] {};}};
      \comm{\axnum};
      \labelout{\axnum};
    \end{tikzpicture}
    \]
    Thus it remains to prove the commutativity of the cell [3-21a].
    This is the content of the next Lemma.
  \end{proof}

\pagebreak
  \begin{axlem} Diagram [3-21a] commutes. \end{axlem} 

  \begin{proof}
    Using the corner trick, we can fill in [3-21a] as follows
    \[
    \def\axnum{21anew}
    \begin{tikzpicture}[scale=.6ex]
    \axgrid{\axnum};
    \foreach \x in {1,...,9}
    \draw (\x) node[littlenode] 
      {\displaytikz{\ax{\axnum}-{\x} \draw (1.25,-.3) node [above] {};}};
    \comm{\axnum};
    \labelout{\axnum};
    \end{tikzpicture} \quad 
    \]
    Here we use Lemma~L for the first time. 
    It ensures that the cell labeled ``L'' commutes. 
    The hexagon labeled N at the bottom of the diagram
    commutes by naturality of the bottom identity with respect
    to two applications of $1 \boxtimes 1$-isomorphisms.
    On one of the sides of this pentagon these two applications of 
    the $1 \boxtimes 1$-isomorphism are denoted by just one map.
  \end{proof}

%  \begin{proof}[Proof of 3-21a for the old version of the proof of 3-21.]
%    \newpicture{A bullet is missing on the middle of the top line in the upper 
%      two items.}
%    Consider the diagram
%    \[
%    \def\axnum{21a}
%    \begin{tikzpicture}[scale=.7ex]
%    \axgrid{\axnum};
%    \foreach \x in {1,...,4}
%    \draw (\x) node[littlenode] 
%      {\innertikzb{\ax{\axnum}-{\x} \draw (1.25,-.3) node [above] {};}};
%    \comm{\axnum};
%    \labelout{\axnum};
%    \end{tikzpicture} \quad .
%    \]
%    To prove commutativity of this diagram, we can use the corner trick
%    and assume that the white sector is an identity sector. 
%    Then commutativity follows from Lemma~K.  %% Added much later: old lemma K
%    The commutativity of [3-21a] can now be deduced.
%    Replace the white sector in the above diagram with the top line
%    from each item in [3-21a]. 
%    Then the above diagram is [3-21a], up to a reflection along a horizontal
%    axis. 
%    \newpicture{It would be better to draw the diagram here with the white
%      sector on top.}
%  \end{proof}

%\Needspace{3\baselineskip}
\pagebreak
  \begin{axprop} Axiom [3-22] is satisfied. \end{axprop} 

  \begin{proof}
    This follows from this diagram:
%    \[
%    \def\axnum{22}
%    \begin{tikzpicture}[scale=.7ex]
%    %\draw[step=.5] (0,0) grid (2,2);
%    \axgrid{\axnum};
%    \foreach \x in {1,...,9}
%    \draw (\x) node[littlenode] 
%      {\innertikzb{\ax{\axnum}-{\x} \draw (1.25,-.3) node [above] {};}};
%    \comm{\axnum};
%    \labelout{\axnum};
%    \end{tikzpicture} .
%    \]
    \[
    \def\axnum{22new}
    \begin{tikzpicture}[scale=.72ex]
    %\draw[step=.5] (0,0) grid (2,2);
    \axgrid{\axnum};
    \foreach \x in {1,...,8}
    \draw (\x) node[littlenode]
      {\displaytikz{\ax{\axnum}-{\x} \draw (1.25,-.3) node [above] {};}};
    \comm{\axnum};
    \labelout{\axnum};
    \end{tikzpicture} \quad \qedhere
    \]
  \end{proof}

  Almost all of the axioms of a dicategory object
  assert that a diagram and a number of variants of the
  diagram commute.
  So far we have ignored the variants---their commutativity
  can always be established by a straightforward variation
  of the argument for the original diagram.
  The only exception to this is axiom [3-23].
  Here our definition of the swap [2-15] was designed
  to ensure that [3-23L], the left hand version of [3-23],
  holds.
  For the right hand version [3-23R] we will have to 
  use a different argument.

%\Needspace{3\baselineskip}
\pagebreak
  \begin{axprop} Axiom [3-23L] is satisfied. \end{axprop}

  \begin{proof}
    Axiom [3-23L] can be filled as follows: %partially filled as follows:
%    \newpicture{In the old picture: the triangle labeled C should be 
%         labeled K.}
%    \newpicture{In the old picture: Consider the top left item of 
%     the inner hexagon. I think the
%     lower diffeomorphism bubble is misplaced. It should be as in the top
%     right item of this hexagon. This will also get rid of one bullet;
%     there should be 8 bullets on the boundary of each item. 
%     On the middle item to the left of this hexagon the middle bullet
%     on the bottom line should be removed.
%     On the bottom item to the left of this hexagon the following bullet 
%     should be removed: the bullet to the left of the top right corner.}
%    \[
%    \def\axnum{23}
%    \begin{tikzpicture}[scale=.6ex]
%    %\draw[step=.5] (0,0) grid (3,4);
%    \axgrid{\axnum};
%    \foreach \x in {1,...,10}
%    \draw (\x) node[littlenode] 
%        {\innertikzb{\ax{\axnum}-{\x} \draw (1.25,-.3) node [above] {};}};
%    \comm{\axnum};
%    \labelout{\axnum};
%    \end{tikzpicture}
%    \]
    \[
    \def\axnum{23new}
    \begin{tikzpicture}[scale=.6ex]
    %\draw[step=.5] (0,0) grid (3,4);
    \axgrid{\axnum};
    \foreach \x in {1,...,9}
    \draw (\x) node[littlenode]
        {\displaytikz{\ax{\axnum}-{\x} \draw (1.25,-.3) node [above] {};}};
    \comm{\axnum};
    \labelout{\axnum};
    \end{tikzpicture}
    \]
	The inner square commutes by Lemma~M.
%    It follows that [3-23L] is equivalent to the commutativity of the
%    inner square which commutes by Lemma~M.
  \end{proof}
  
%\Needspace{3\baselineskip}

  \begin{axprop} Axiom [3-23R] is  satisfied.  \end{axprop}

  \begin{proof}
    Consider again the diagram from the proof of [3-23L].
    This diagram reduced [3-23L] to Lemma~M.
    The proof of Lemma~M in turn reduced to the case
    where all defects are identity defects.
    The same argument can be applied to [3-23R] to reduce to the case of identity defects: 
    we therefore only need to prove [3-23R] in the case where all the
    defects are identity defects.
    In this case the following diagram reduces [3-23R]
    to [3-23L] (which is already proved) 
    and [3-24] (which we prove next): 
    \[
    \def\axnum{23a} 
    \begin{tikzpicture}[scale=.46ex]
    %\draw[step=.5] (0,0) grid (5,5);
    \axgrid{\axnum};
    \foreach \x in {1,...,12}
    \draw (\x) node[littlenode] 
        {\innertikz{\ax{\axnum}-{\x} \draw (1.25,-.3) node [above] {};}};
    \comm{\axnum};
    \labelout{\axnum};
    \end{tikzpicture} \quad \qedhere
    \]
  \end{proof}

%\Needspace{3\baselineskip}
\pagebreak
  \begin{axprop} Axiom [3-24] is satisfied. \end{axprop} 

  \begin{proof}
    Axiom [3-24] can be filled as follows:
%    \newpicture{For the old picture: There is abullet missing from 3-24-4.
%      (Consider the top right corner of the triangle in
%        upper half. There is a bullet missing in the
%        lower right corner of this item.)}
%    \newpicture{For the old picture: 
%       The two triangles labeled D should be labeled M.}
%    \[
%    \def\axnum{24}
%    \begin{tikzpicture}[scale=.6ex]
%    %\draw[step=.5] (0,0) grid (4,4);
%    \axgrid{\axnum};
%    \foreach \x in {1,...,11}
%    \draw (\x) node[littlenode] 
%      {\innertikzb{\ax{\axnum}-{\x} \draw (1.25,-.3) node [above] {};}};
%    \comm{\axnum};
%    \labelout{\axnum};
%    \end{tikzpicture}
%    \]
    \[
    \def\axnum{24new}
    \begin{tikzpicture}[scale=.57ex]
    %\draw[step=.5] (0,0) grid (4,4);
    \axgrid{\axnum};
    \foreach \x in {1,...,9}
    \draw (\x) node[littlenode]
      {\displaytikz{\ax{\axnum}-{\x} \draw (1.25,-.3) node [above] {};}};
    \comm{\axnum};
    \labelout{\axnum};
    \end{tikzpicture} \quad \qedhere
    \]
  \end{proof}

%\Needspace{3\baselineskip}
\pagebreak
  \begin{axprop} Axiom [3-25] is satisfied. \end{axprop} 

  \begin{proof}
    This follows from the associativity of $1 \boxtimes 1$:
    \[
    \def\axnum{25}
    \begin{tikzpicture}[scale=.6ex]
    \axgrid{\axnum};
    \foreach \x in {1,...,4}
    \draw (\x) node[littlenode] 
      {\displaytikz{\ax{\axnum}-{\x} \draw (1.25,-.3) node [above] {};}};
    \comm{\axnum};
    \labelout{\axnum};
    \end{tikzpicture} \quad \qedhere
    \]
  \end{proof}

%\Needspace{3\baselineskip}
%g%\pagebreak
  \begin{axprop} Axiom [3-26] is satisfied. \end{axprop} 

  \begin{proof}
    This follows from the naturality of the flip, applied to the
    expand isomorphism:
    \[
    \def\axnum{26}
    \begin{tikzpicture}[scale=.6ex]
    \axgrid{\axnum};
    \foreach \x in {1,...,4}
    \draw (\x) node[littlenode] 
      {\displaytikz{\ax{\axnum}-{\x} \draw (1.25,-.3) node [above] {};}};
    \comm{\axnum};
    \labelout{\axnum};
    \end{tikzpicture} \quad \qedhere
    \]
  \end{proof}

\newpage
%!%\newgeometry{margin=1.1in}
\appendix
\section{Internal dicategories} \label{appendix}

For ease of reference while reading the above proof, we compile here the definition of an internal dicategory, as introduced in~\cite{Douglas-Henriques(Internal-bicategories)}.  We also list, in Table~\ref{table-internaldicategory}, a single pictorial abbreviation for each piece of data and each axiom in the definition; the meaning of these pictograms is given in the full definition following. 
\def\tick{
	\coordinate (A) at ($ (\x,\y) + (0,.5*\tickv) $);
	\coordinate (B) at ($ (\x,\y) + (0,-.5*\tickv) $);
	\draw (A) to (B);
}

\begin{figure}[ht]
\centering
\begin{tikzpicture}[xscale=\tikztablescalex,yscale=\tikztablescaley]
\draw (0,0) node {0-data};
\draw (-.3,-1) node[anchor=north west] {
\def\tablenum{1}
\begin{tikzpicture}[scale=\tikztablescaleinner]
\tablegrid{\tablenum};
\foreach \x in {1,...,3}
\draw (\x) node[tablenode,anchor=west]
	{\innertikztikz{\table{\tablenum}-{\x};}};
\end{tikzpicture}
};
\draw (.95,0) node {1-data};
\draw (.65,-1) node[anchor=north west] {
\def\tablenum{2}
\begin{tikzpicture}[scale=\tikztablescaleinner]
\tablegrid{\tablenum};
\foreach \x in {1,...,8}
\draw (\x) node[tablenode,anchor=west]
	{\innertikztikz{\table{\tablenum}-{\x};}};
\end{tikzpicture}
};
\draw (2.05,0) node {2-data};
\draw (1.75,-1) node[anchor=north west] {
\def\tablenum{3}
\begin{tikzpicture}[xscale=\tikztablescaleinnerx,yscale=\tikztablescaleinnery]
\tablegrid{\tablenum};
\foreach \x in {1,...,2}
\draw (\x) node[tablenode-eps,anchor=west]
	{\innertikzeps{\table{\tablenum}-{\x};}};
\foreach \x in {3,...,11}
\draw (\x) node[tablenode,anchor=west]
	{\innertikztikz{\table{\tablenum}-{\x};}};
\foreach \x in {9}
\draw (12) node[tablenode,anchor=west]
	{\innertikztikz{\table{2}-{\x};}};
\foreach \x in {13,...,18}
\draw (\x) node[tablenode,anchor=west]
	{\innertikztikz{\table{\tablenum}-{\x};}};
\end{tikzpicture}
};
\draw (.1,-26) node {3-axioms};
\draw (-.2,-27) node[anchor=north west] {
\def\tablenum{4}
\begin{tikzpicture}[xscale=\tikztablescaleinnerx,yscale=\tikztablescaleinnery]
\tablegrid{\tablenum};
\foreach \x in {1,2,3,6,7,8,9,12,13}
\draw (\x) node[tablenode-eps,anchor=west]
	{\innertikzeps{\table{\tablenum}-{\x};}};
\foreach \x in {4,5,10,11}
\draw (\x) node[tablenode,anchor=west]
	{\innertikztikz{\table{\tablenum}-{\x};}};
\foreach \x in {14,...,28}
\draw (\x) node[tablenode,anchor=west]
	{\innertikztikz{\table{\tablenum}-{\x};}};
\foreach \x in {12}
\draw (29) node[tablenode,anchor=west]
	{\innertikztikz{\table{3}-{\x};}};
\foreach \x in {30,...,46}
\draw (\x) node[tablenode,anchor=west]
	{\innertikztikz{\table{\tablenum}-{\x};}};
\end{tikzpicture}
};
\end{tikzpicture}
%\vspace*{-20pt}
\renewcommand{\figurename}{Table}
\caption{Abbreviated definition of an internal bicategory.} \label{table-internaldicategory}
\end{figure}

\begin{definition}[{\cite[Def.\,3.3]{Douglas-Henriques(Internal-bicategories)}}] \label{def-dicat}
A dicategory object $C$ in the 2-category $\cC$ consists of the following three collections of data, subject to the listed axioms.

\begin{description}

\item[0-data]
There are three objects of $\cC$ as follows:
\begin{itemize}
\item[D0-0]: $C_0$, a groupoid object, denoted \cb{\ing{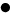}} and called the object of 0-cells.
\item[D0-1]: $C_1$, a groupoid object, denoted \cb{\ing{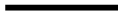}} and called the object of 1-cells.
\item[D0-2]: $C_2$ (typically not a groupoid object), denoted \cb{\ing{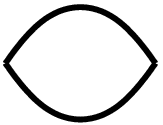}} and called the object of 2-cells.
\end{itemize}
In addition, there are morphisms $s,t: C_1 \ra C_0$ and $s,t: C_2 \ra C_1$, the source and target, such that $st = ss$ and $tt = ts$, and such that $s \times t: C_1 \ra C_0 \times C_0$ and $s \times t: C_2 \ra C_1 \times_{C_0 \times C_0} C_1$ are fibrations.

\item[1-data]
There are eight 1-morphisms of $\cC$ as follows:
\begin{itemize}
\item[D1-1]: $i: C_0 \ra C_1$, denoted \cb{\ing{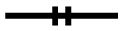}} and called the 1-cell identity.
\item[D1-2]: $m: C_1 \times_{C_0} C_1 \ra C_1$, denoted \cb{\ing{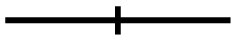}} and called the horizontal composition.
\item[D1-3]: $i_v: C_1 \ra C_2$, denoted \cb{\ing{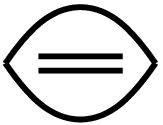}} and called the 2-cell identity.
\item[D1-4]: $m_v: C_2 \times_{C_1} C_2 \ra C_2$, denoted \cb{\ing{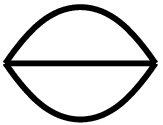}} and called the vertical composition.
\item[D1-5]: $w_r: C_2 \times_{C_0} C_1 \ra C_2$, denoted \cb{\ing{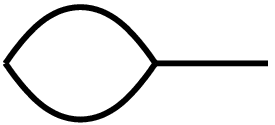}} and called the right composition or whisker of a 2-cell with a 1-cell.
\item[D1-6]: $w_l: C_1 \times_{C_0} C_2 \ra C_2$, denoted \cb{\ing{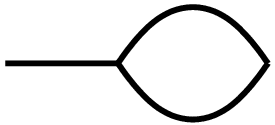}} and called the left composition or whisker of a 1-cell with a 2-cell.
\item[D1-7]: $i_l: C_1 \ra C_2$, denoted \cb{\ing{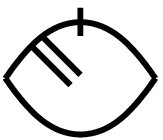}} and called the (upper) left 2-cell identity.
\item[D1-8]: $i_r: C_1 \ra C_2$, denoted \cb{\ing{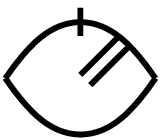}} and called the (upper) right 2-cell identity.
\end{itemize}
These morphisms are compatible with source and target maps. \\
\indent The morphisms [D1-7] and [D1-8] are required to be invertible, in the following sense.  There exists a morphism of $\cC$, denoted \cb{\ing{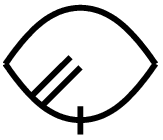}} (the lower left 2-cell identity), such that there are invertible 2-morphisms from \cb{\ing{d1stuff-halfid1}} to \cb{\ing{d1stuff-vertid1}} and from \cb{\ing{d1stuff-halfid2}} to \cb{\ing{d1stuff-vertid2}}, such that the two resulting 2-morphisms from \cb{\ing{d1stuff-triangle1}} to \cb{\ing{d1stuff-item7.eps}} are equal, and similarly the two 2-morphisms from \cb{\ing{d1stuff-triangle2}} to \cb{\ing{d1stuff-item9.eps}} are equal.  Similarly there exists a morphism \cb{\ing{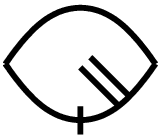}} (the lower right 2-cell identity) satisfying the corresponding conditions.

\item[2-data]
There are eighteen 2-isomorphisms of $\cC$ as follows: \vspace{-2pt} 
\setlength{\columnsep}{35pt}
\begin{itemize}
\begin{multicols}{2}
\item[D2-1]: \xymatrix@C-12pt{\cb{\ing{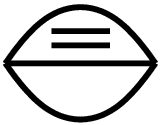}}   \ar@{=>}[r] &   \cb{\ing{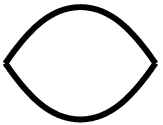}}} 
\item[D2-2]: \xymatrix@C-12pt{\cb{\ing{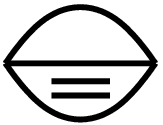}}   \ar@{=>}[r] &   \cb{\ing{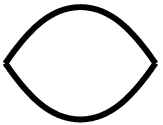}}} 
\item[D2-3]: \xymatrix@C-12pt{\cb{\ing{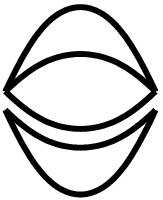}}   \ar@{=>}[r] &   \cb{\ing{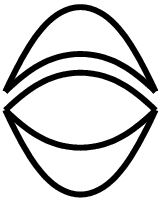}}} 
\item[D2-4]: \xymatrix@C-12pt{\cb{\ing{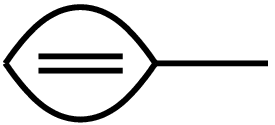}}   \ar@{=>}[r] &   \cb{\ing{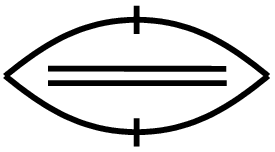}}} 
\item[D2-5]: \xymatrix@C-12pt{\cb{\ing{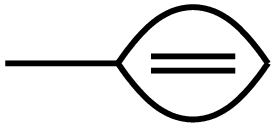}}   \ar@{=>}[r] &   \cb{\ing{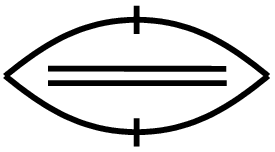}}} 
\item[D2-6]: \xymatrix@C-12pt{\cb{\ing{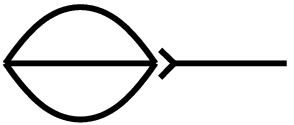}}   \ar@{=>}[r] &   \cb{\ing{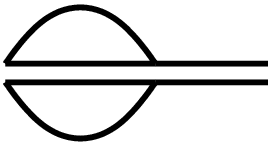}}}  
\item[D2-7]: \xymatrix@C-12pt{\cb{\ing{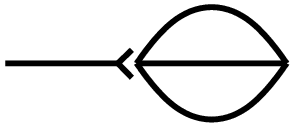}}   \ar@{=>}[r] &   \cb{\ing{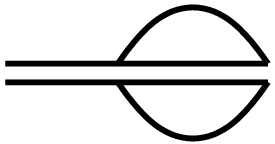}}} 
\item[D2-8]: \xymatrix@C-12pt{\cb{\ing{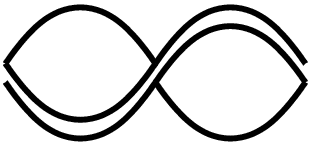}}   \ar@{=>}[r] &   \cb{\ing{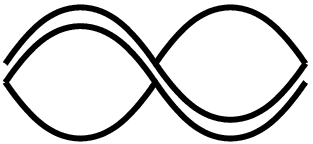}}} 
\item[D2-9]:\hspace{-5pt} \xymatrix@C-12pt{\cb{\ing{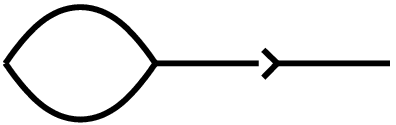}}   \ar@{=>}[r] &   \cb{\ing{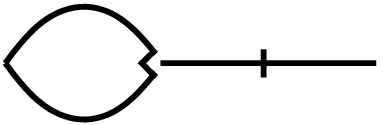}}}  
\item[D2-10]:\hspace{-5pt} \xymatrix@C-12pt{\cb{\ing{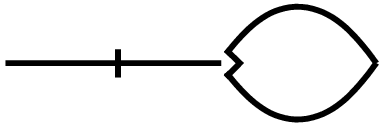}}   \ar@{=>}[r] &   \cb{\ing{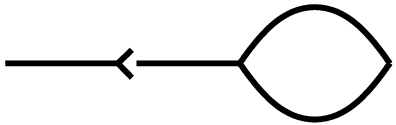}}} 
\item[D2-11]:\hspace{-5pt} \xymatrix@C-12pt{\cb{\ing{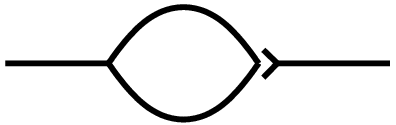}}   \ar@{=>}[r] &   \cb{\ing{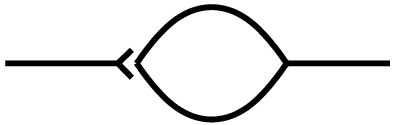}}}  
\item[D2-12]:\hspace{-5pt} \xymatrix@C-12pt{\cb{\ing{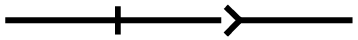}}   \ar@{=>}[r] &   \cb{\ing{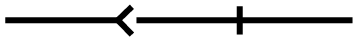}}} 
\item[D2-13]: \xymatrix@C-12pt{\cb{\ing{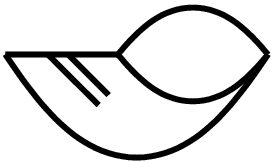}}   \ar@{=>}[r] &   \cb{\ing{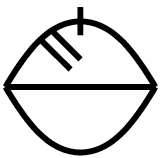}}} 
\item[D2-14]: \xymatrix@C-12pt{\cb{\ing{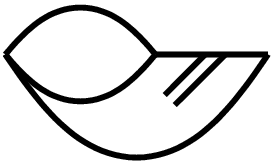}}   \ar@{=>}[r] &   \cb{\ing{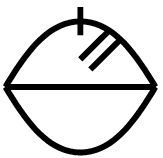}}} 
\item[D2-15]: \xymatrix@C-12pt{\cb{\ing{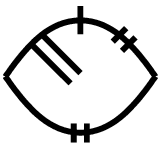}}   \ar@{=>}[r] &   \cb{\ing{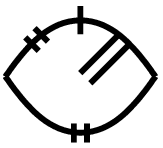}}} 
\item[D2-16]: \xymatrix@C-12pt{\cb{\ing{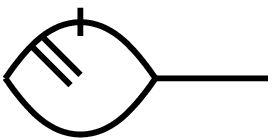}}   \ar@{=>}[r] &   \cb{\ing{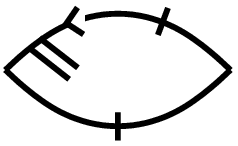}}} 
\item[D2-17]: \xymatrix@C-12pt{\cb{\ing{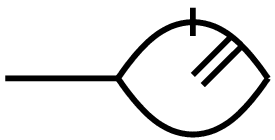}}   \ar@{=>}[r] &   \cb{\ing{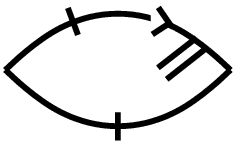}}} 
\item[D2-18]: \xymatrix@C-12pt{\cb{\ing{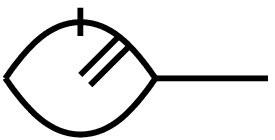}}   \ar@{=>}[r] &   \cb{\ing{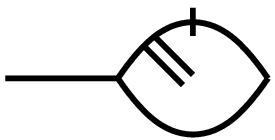}}} 
\end{multicols}
\end{itemize}
\nid These 2-isomorphisms are compatible with source and target maps in the sense that the sources and targets of [D2-1] through [D2-8] and of [D2-13] through [D2-15] are identity 2-morphisms, the sources and targets of [D2-9] through [D2-11] are the 2-isomorphism [D2-12], the sources of [D2-16] and [D2-18] are [D2-12], the source of [D2-17] is the inverse of [D2-12], and the targets of [D2-16] through [D2-18] are identity 2-isomorphisms.

\item[3-axioms]
The above data are such that the following twenty-six diagrams, as well as the variant diagrams abbreviated in parentheses, commute:

\end{description} % In order for axioms to have full page width.
\pagebreak

\setlength{\columnsep}{20pt}

\begin{itemize}
\begin{multicols}{2}
\item[D3-1]:
\db{
\xymatrix@C-12pt{\cb{\ingt{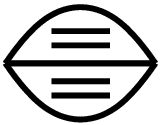}} \ar@/^10pt/ @{-}[r] \ar@/_10pt/@{-}[r] & \cb{\ingt{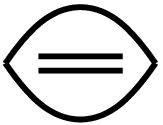}}}
}
\vspace{10pt}

\item[D3-2]: 
\def\alphanum{\ifcase\xypolynode\or \ingt{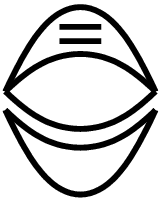} \or \ingt{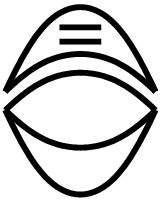} \or 
\ingt{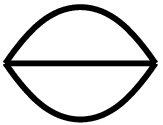} \fi}
\db{
\xy/r\scs/: \xypolygon3{~={0}~*{\alphanum}}
\endxy
}\nopagebreak
\\[-20pt]
\hspace*{\fill} $\big[\hspace{-1.5ex}\bb{\xys{\ingt{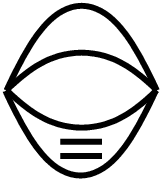}}}\hspace{-1.5ex}\big]$ \hspace*{30pt}

\item[D3-3]:  
\def\alphanum{\ifcase\xypolynode\or \ingt{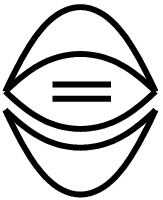} \or \ingt{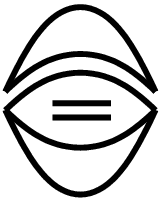} \or 
\ingt{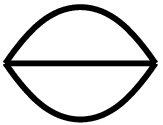} \fi}
\db{
\xy/r\scs/: \xypolygon3{~={0}~*{\alphanum}}
\endxy
}
\vspace{10pt}

\item[D3-4]: 
\def\alphanum{\ifcase\xypolynode\or \ingt{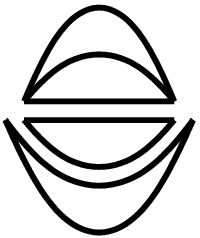} \or \ingt{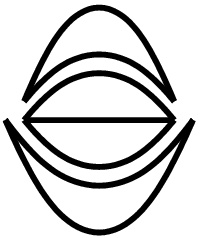} \or 
\ingt{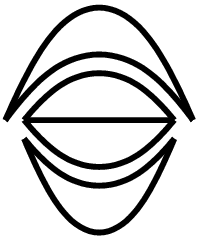} \or \ingt{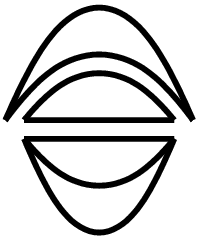} \or \ingt{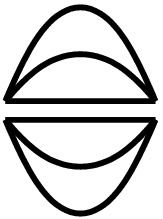} \fi}
\hspace{\hsqza}
\db{
\xy/r\sca/: \xypolygon5{~*{\alphanum}}
\endxy
}
\vspace{10pt}

\item[D3-5]: 
\def\alphanum{\ifcase\xypolynode\or \ingt{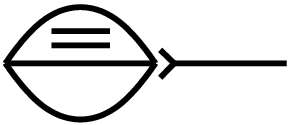} \or \ingt{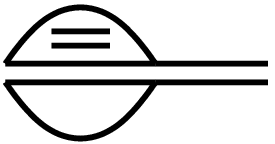} \or 
\ingt{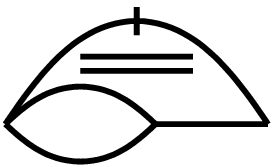} \or \ingt{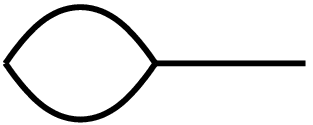} \fi}
\db{
\xy/r\sca/: \xypolygon4{~*{\alphanum}}
\endxy
}
\\
\hspace*{-10pt} $\big[\hspace{-1.5ex}\bb{\xys{\ingt{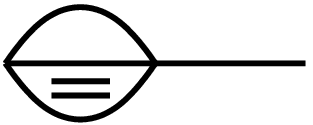}}} 
\hspace{-5pt} \bb{\xys{\ingt{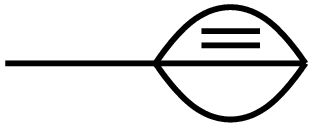}}} 
\hspace{-5pt} \bb{\xys{\ingt{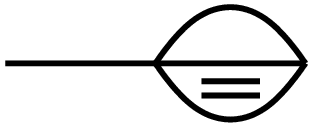}}}\hspace{-1.5ex}\big]$
\vspace{5pt}

\item[D3-6]: 
\def\alphanum{\ifcase\xypolynode\or \ingt{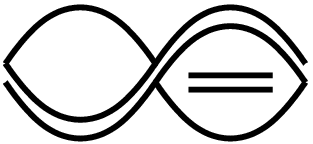} \or \ingt{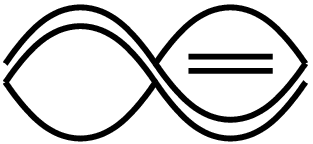} \or 
\ingt{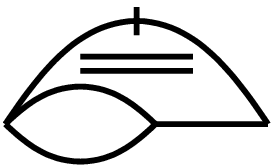} \or \ingt{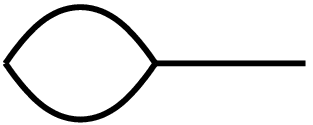} \or \ingt{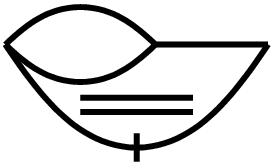} \fi}
\hspace{\hsqza}
\db{
\xy/r\scb/: \xypolygon5{~={0}~*{\alphanum}}
\endxy
}
\\[-10pt]
\hspace*{\fill} $\big[\hspace{-1.5ex}\bb{\xys{\ingt{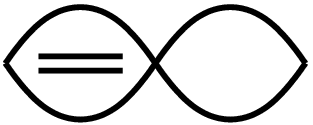}}}\hspace{-1.5ex}\big]$ \hspace*{10pt}
\vspace{5pt}

\item[D3-7]: 
\def\alphanum{\ifcase\xypolynode\or \ingt{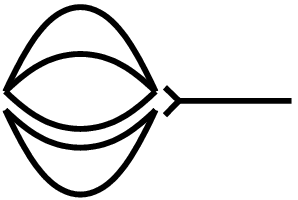} \or \ingt{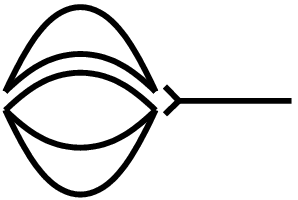} \or 
\ingt{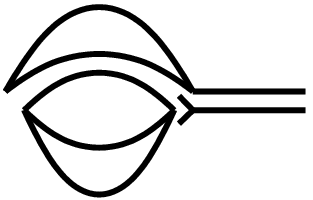} \or \ingt{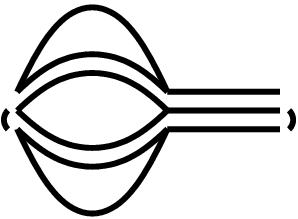} \or \ingt{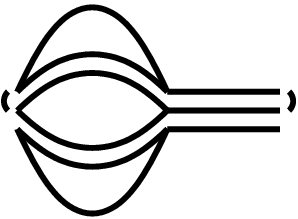} \or \ingt{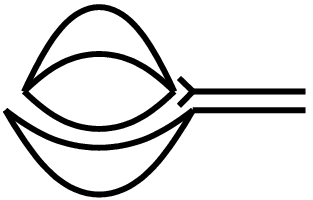} \fi}
\hspace{\hsqza}
\db{
\xy/r\scc/: \xypolygon6{~={30}~*{\alphanum}}
\endxy
}
\\[-15pt]
\hspace*{\fill} $\big[\hspace{-1.5ex}\bb{\xys{\ingt{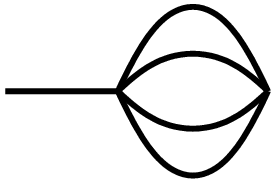}}}\hspace{-1.5ex}\big]$ \hspace*{10pt} 
\vspace{5pt}

%\end{multicols}

\item[D3-8]: 
\def\alphanum{\ifcase\xypolynode\or \ingt{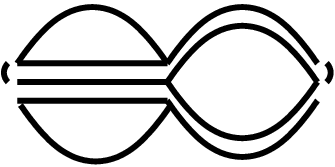} \or \ingt{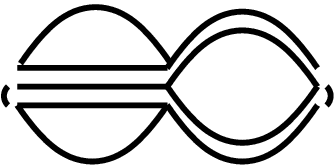} \or \ingt{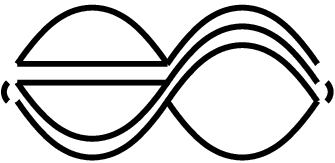} \or \ingt{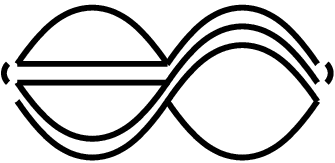} \or \ingt{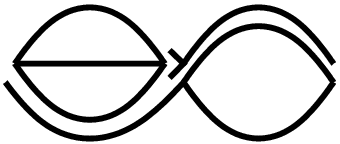} \or \ingt{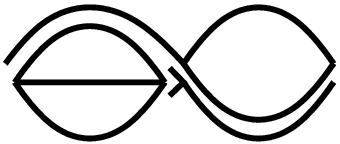} \or \ingt{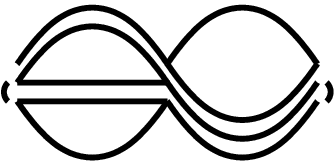} \or \ingt{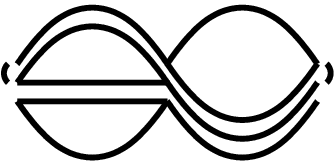} \fi}
%\hspace{\hsqza}
\hspace{-45pt}
\db{
\xy/r\scf/:  \xypolygon8{~*{\alphanum}}
\endxy
}
\\[-0pt]
\hspace*{\fill} $\big[\hspace{-1.5ex}\bb{\xys{\ingt{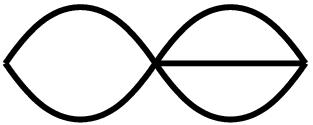}}}\hspace{-1.5ex}\big]$ \hspace*{10pt}
\vspace{5pt}

%\begin{multicols}{2}

\item[D3-9]: 
\hspace*{-25pt}
\def\alphanum{\ifcase\xypolynode\or \ingt{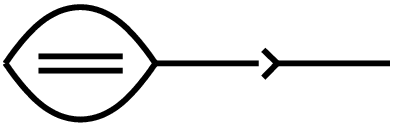} \or \ingt{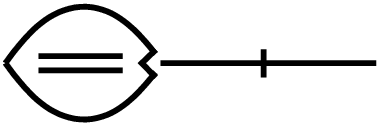} \or 
\ingt{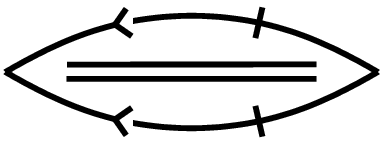} \or \ingt{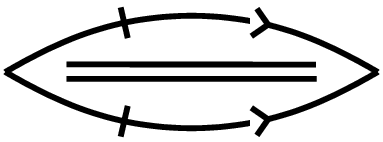} \or \ingt{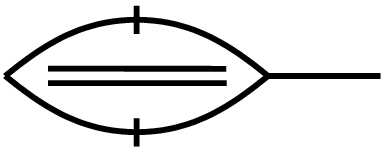} \fi}
\hspace{\hsqza}
\db{
\xy/r\scc/: \xypolygon5{~={0}~*{\alphanum}}
\endxy
}
\\[-10pt]
\hspace*{\fill} $\big[\hspace{-1.5ex}\bb{\xys{\ingt{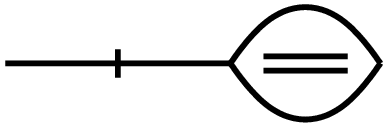}}}\hspace{-1.5ex}\big]$ \hspace*{10pt} 
\vspace{5pt}

\item[D3-10]: 
\hspace*{-25pt}
\def\alphanum{\ifcase\xypolynode\or \ingt{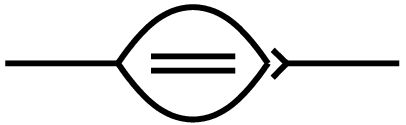} \or \ingt{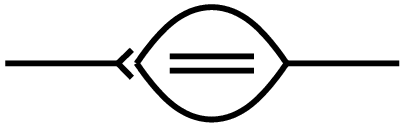} \or 
\ingt{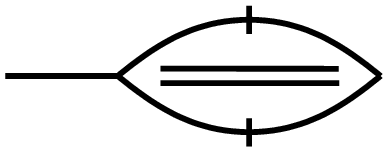} \or \ingt{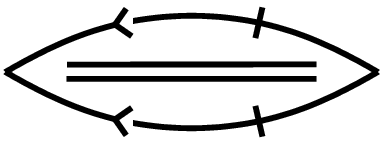} \or \ingt{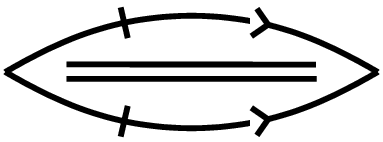} \or \ingt{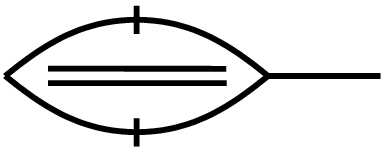} \fi}
\hspace{\hsqza}
\db{
\xy/r\scc/: \xypolygon6{~={30}~*{\alphanum}}
\endxy
}
\vspace{10pt}

\item[D3-11]: 
\hspace*{-25pt}
\def\alphanum{\ifcase\xypolynode\or \ingt{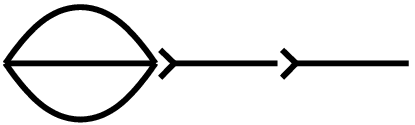} \or \ingt{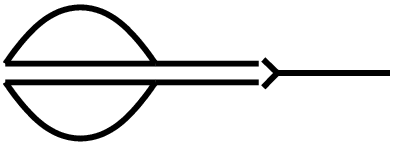} \or 
\ingt{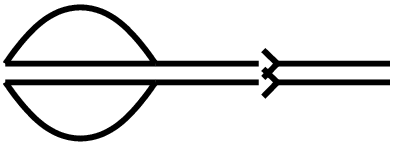} \or \ingt{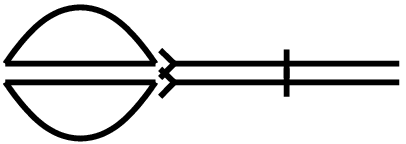} \or \ingt{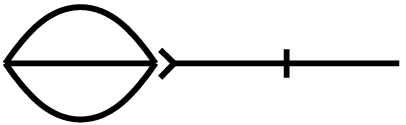} \fi}
\hspace{\hsqza}
\db{
\xy/r\scbb/: \xypolygon5{~={36}~*{\alphanum}}
\endxy
} \nopagebreak
\\[-10pt] 
\hspace*{\fill} $\big[\hspace{-1.5ex}\bb{\xys{\ingt{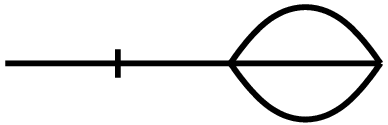}}}\hspace{-1.5ex}\big]$ \hspace*{20pt} 
\vspace{5pt}

\item[D3-12]: 
\hspace*{-25pt}
\def\alphanum{\ifcase\xypolynode\or \ingt{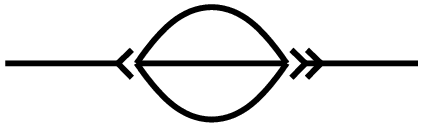} \or \ingt{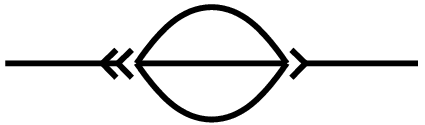} \or 
\ingt{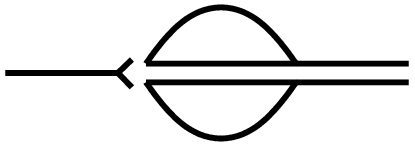} \or \ingt{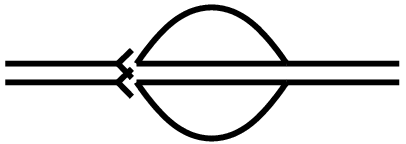} \or \ingt{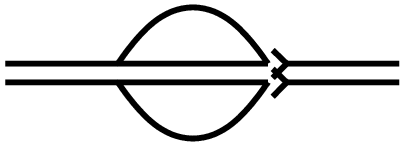} \or \ingt{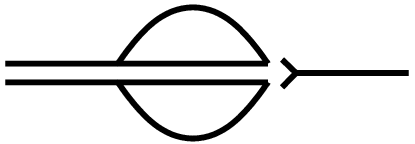} \fi}
\hspace{\hsqza}
\db{
\xy/r\sccc/: \xypolygon6{~={30}~*{\alphanum}}
\endxy
}
\vspace{10pt}

\item[D3-13]: 
\hspace*{-25pt}
\def\alphanum{\ifcase\xypolynode\or \ingt{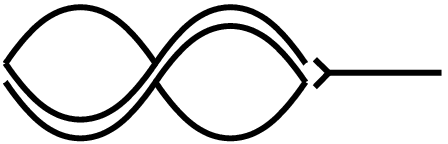} \or \ingt{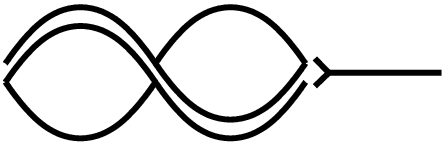} \or 
\ingt{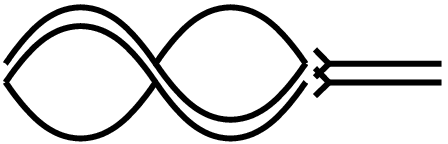} \or \ingt{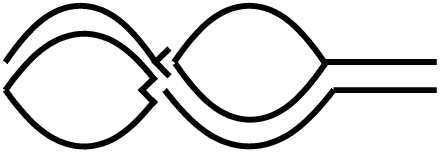} \or \ingt{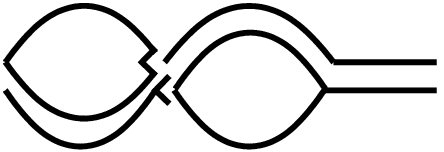} \or \ingt{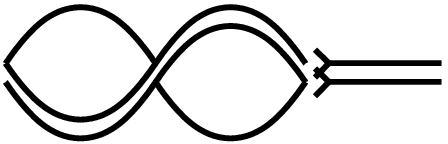} \fi}
\hspace{\hsqza}
\db{
\xy/r\sccc/: \xypolygon6{~={30}~*{\alphanum}}
\endxy
}
\\[-10pt]
\hspace*{\fill} $\big[\hspace{-1.5ex}\bb{\xys{\ingt{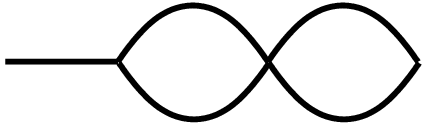}}}\hspace{-1.5ex}\big]$ \hspace*{10pt} 
\vspace{5pt}

\item[D3-14]: 
\def\alphanum{\ifcase\xypolynode\or \ingt{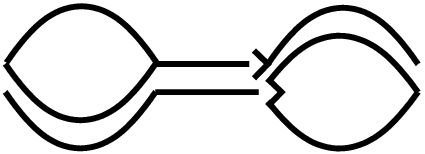} \or \ingt{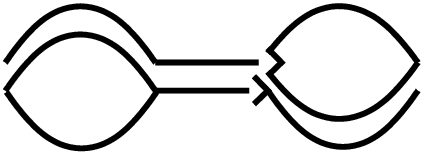} \or 
\ingt{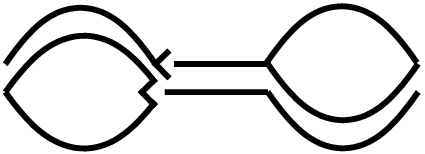} \or \ingt{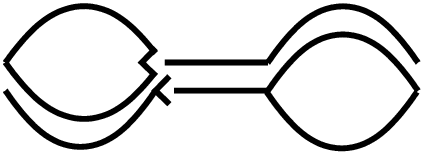} \fi}
\db{
\xy/r\scc/: \xypolygon4{~*{\alphanum}}
\endxy
}
\vspace{10pt}

\item[D3-15]: 
\hspace*{15pt}
\def\alphanum{\ifcase\xypolynode\or \ingt{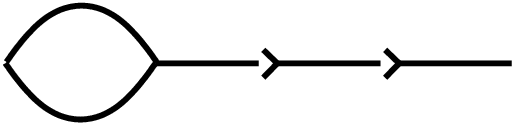} \or \ingt{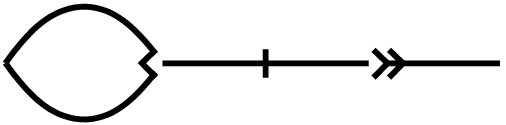} \or 
\ingt{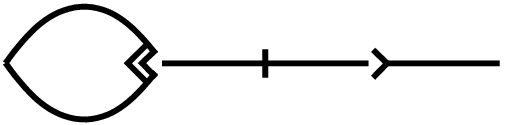} \or \ingt{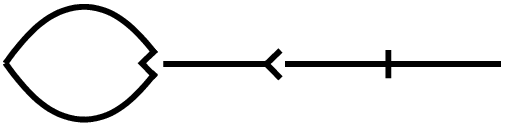} \or \ingt{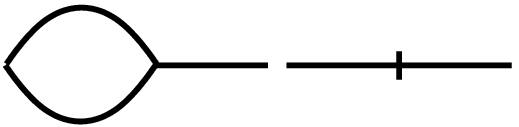} \fi}
\hspace{\hsqzb}
\db{
\xy/r\sccc/: \xypolygon5{~={36}~*{\alphanum}}
\endxy
}
\\[-5pt]
\hspace*{\fill} $\big[\hspace{-1.5ex}\bb{\xys{\ingt{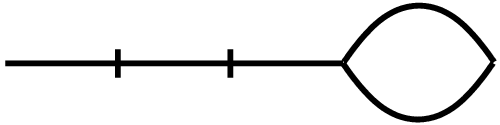}}}\hspace{-1.5ex}\big]$ \hspace*{35pt} 
\vspace{5pt}

\item[D3-16]: 
\hspace*{15pt}
\def\alphanum{\ifcase\xypolynode\or \ingt{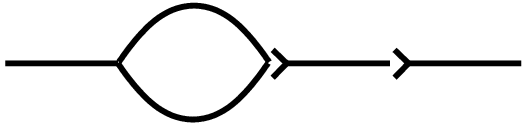} \or \ingt{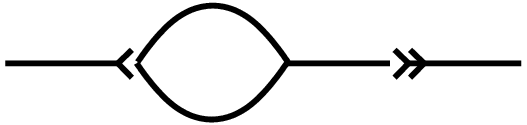} \or 
\ingt{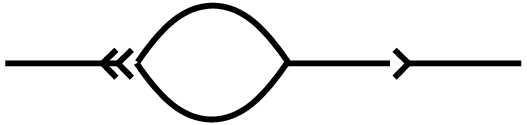} \or \ingt{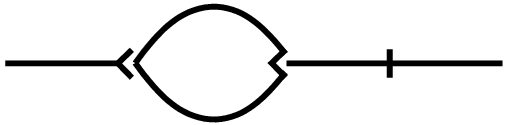} \or \ingt{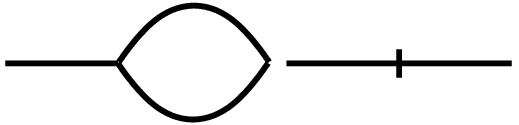} \fi}
\hspace{\hsqzb}
\db{
\xy/r\sccc/: \xypolygon5{~={36}~*{\alphanum}}
\endxy
}
\\[-5pt]
\hspace*{\fill} $\big[\hspace{-1.5ex}\bb{\xys{\ingt{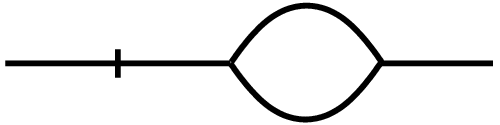}}}\hspace{-1.5ex}\big]$ \hspace*{35pt} 
\vspace{5pt}

\item[D3-17]: 
\hspace*{15pt}
\def\alphanum{\ifcase\xypolynode\or \ingt{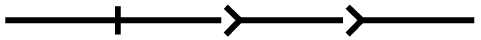} \or \ingt{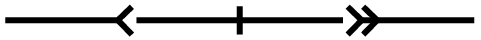} \or 
\ingt{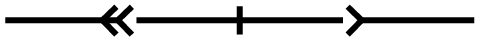} \or \ingt{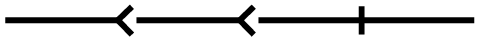} \or \ingt{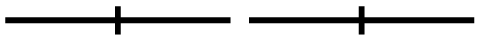} \fi}
\hspace{\hsqzb}
\db{
\xy/r\scc/: \xypolygon5{~={36}~*{\alphanum}}
\endxy
}
\vspace{50pt}

\item[D3-18]:
\def\alphanum{\ifcase\xypolynode\or \ingt{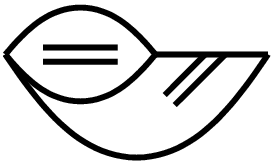} \or \ingt{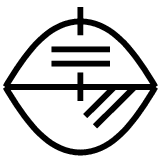} \or 
\ingt{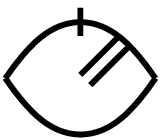} \or \ingt{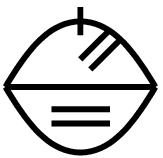} \fi}
\db{
\xy/r\sca/: \xypolygon4{~*{\alphanum}}
\endxy
}
\\[-5pt]
\hspace*{\fill} $\big[\hspace{-1.5ex}\bb{\xys{\ingt{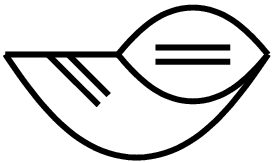}}}\hspace{-1.5ex}\big]$ \hspace*{10pt}
\vspace{5pt}

\item[D3-19]:
\def\alphanum{\ifcase\xypolynode\or \ingt{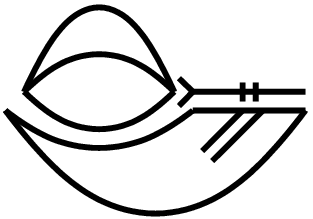} \or \ingt{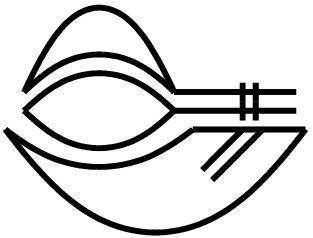} \or \ingt{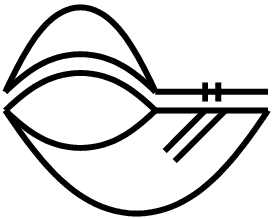} \or \ingt{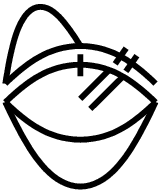} \or \ingt{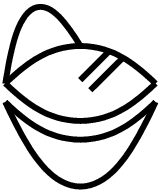} \or \ingt{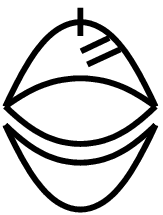} \or \ingt{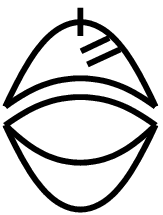} \fi}
\hspace{\hsqzb}
\db{
\xy/r\sce/:  \xypolygon7{~*{\alphanum}}
\endxy
}
\\[-20pt]
\hspace*{\fill} $\big[\hspace{-1.5ex}\bb{\xys{\ingt{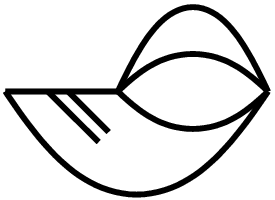}}}\hspace{-1.5ex}\big]$ \hspace*{10pt}
\vspace{5pt}

\item[D3-20]:
\hspace*{-15pt}
\def\alphanum{\ifcase\xypolynode\or \ingt{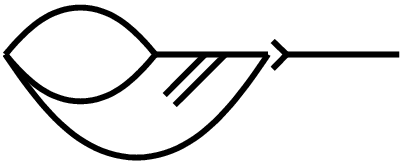} \or \ingt{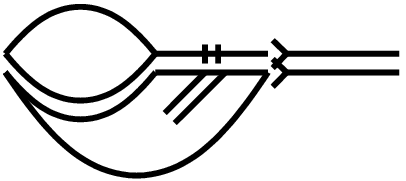} \or 
\ingt{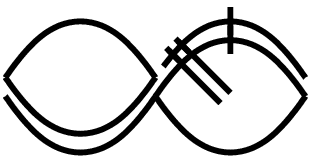} \or \ingt{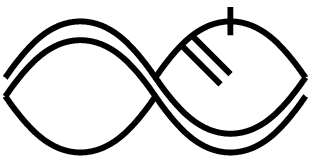} \or \ingt{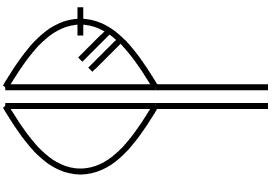} \or \ingt{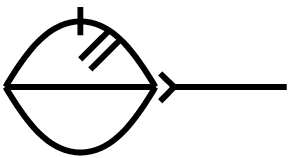} \fi}
\hspace{\hsqza}
\db{
\xy/r\sccc/: \xypolygon6{~={30}~*{\alphanum}}
\endxy
}
\\[-20pt]
\hspace*{\fill} $\big[\hspace{-1.5ex}\bb{\xys{\ingt{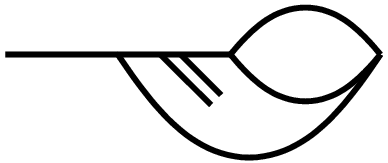}}}\hspace{-1.5ex}\big]$ \hspace*{10pt}
\vspace{5pt}

\item[D3-21]:
\def\alphanum{\ifcase\xypolynode\or \ingt{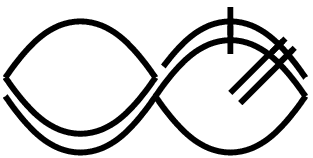} \or \ingt{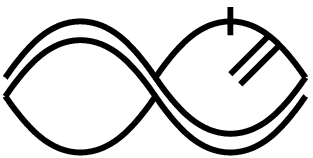} \or 
\ingt{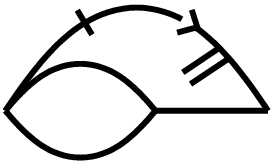} \or \ingt{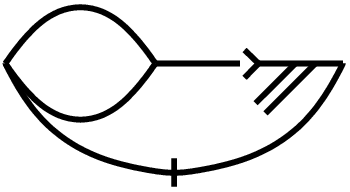} \fi}
\db{
\xy/r\scaa/: \xypolygon4{~*{\alphanum}}
\endxy
}
\\[-5pt]
\hspace*{\fill} $\big[\hspace{-1.5ex}\bb{\xys{\ingt{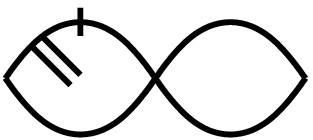}}}\hspace{-1.5ex}\big]$ \hspace*{25pt}
\vspace{5pt}

\item[D3-22]:
\hspace*{-15pt}
\def\alphanum{\ifcase\xypolynode\or \ingt{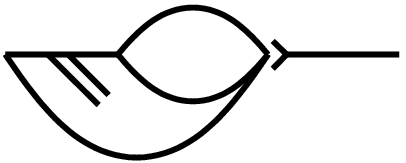} \or \ingt{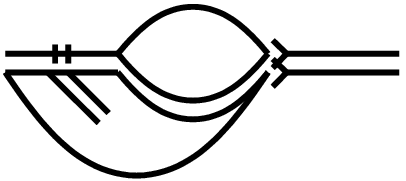} \or 
\ingt{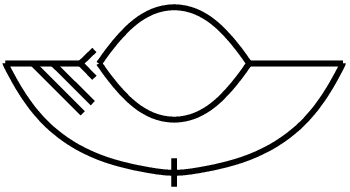} \or \ingt{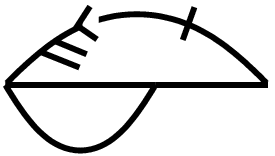} \or \ingt{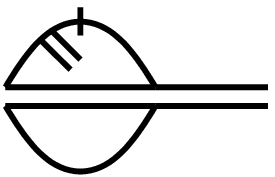} \or \ingt{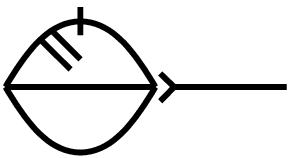} \fi}
\hspace{\hsqza}
\db{
\xy/r\scc/: \xypolygon6{~={30}~*{\alphanum}}
\endxy
} \nopagebreak
\\[-20pt]
\hspace*{\fill} $\big[\hspace{-1.5ex}\bb{\xys{\ingt{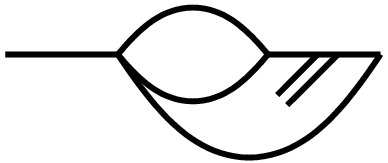}}}\hspace{-1.5ex}\big]$ \hspace*{25pt}
\vspace{5pt}

\item[D3-23]:
\def\alphanum{\ifcase\xypolynode\or \ingt{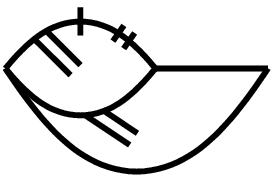} \or \ingt{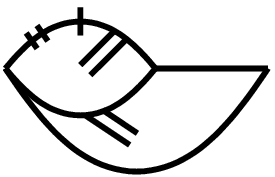} \or 
\ingt{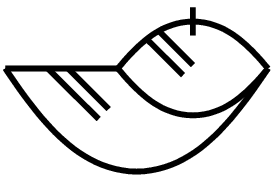} \or \ingt{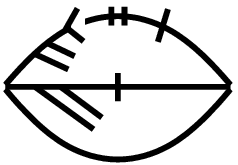} \fi}
\db{
\xy/r\sca/: \xypolygon4{~*{\alphanum}}
\endxy
}
\\[-5pt]
\hspace*{\fill} $\big[\hspace{-1.5ex}\bb{\xys{\ingt{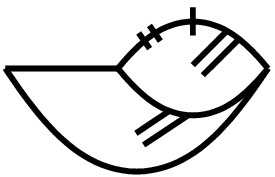}}}\hspace{-1.5ex}\big]$ \hspace*{25pt}
\vspace{5pt}

\item[D3-24]:
\def\alphanum{\ifcase\xypolynode\or \ingt{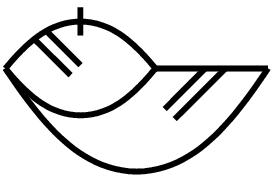} \or \ingt{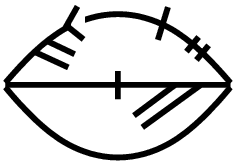} \or 
\ingt{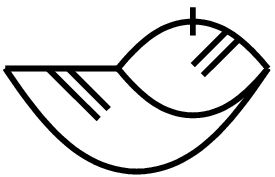} \or \ingt{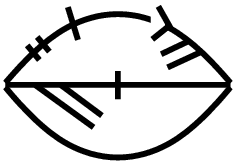} \fi}
\db{
\xy/r\sca/: \xypolygon4{~*{\alphanum}}
\endxy
}
\vspace{10pt}

\item[D3-25]:
\hspace*{-15pt}
\def\alphanum{\ifcase\xypolynode\or \ingt{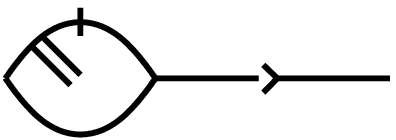} \or \ingt{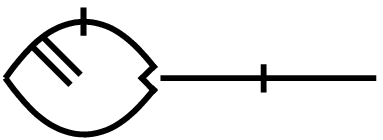} \or 
\ingt{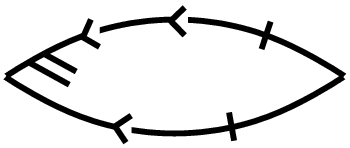} \or \ingt{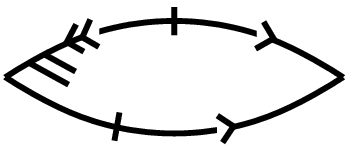} \or \ingt{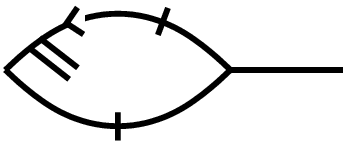} \fi}
\hspace{\hsqza}
\db{
\xy/r\scc/: \xypolygon5{~={36}~*{\alphanum}}
\endxy
}
\\[-10pt]
\hspace*{\fill} $\big[\hspace{-1.5ex}\bb{\xys{\ingt{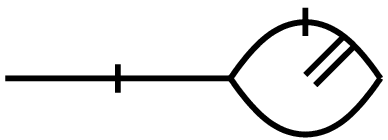}}}\hspace{-1.5ex}\big]$ \hspace*{10pt}
\vspace{5pt}

\item[D3-26]:
\hspace*{-15pt}
\def\alphanum{\ifcase\xypolynode\or \ingt{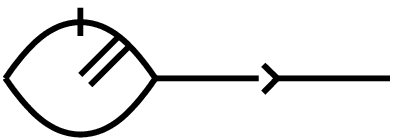} \or \ingt{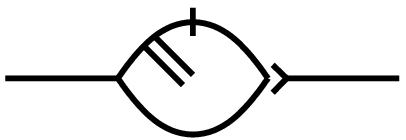} \or 
\ingt{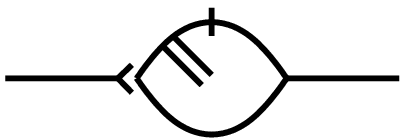} \or \ingt{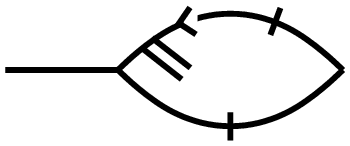} \or \ingt{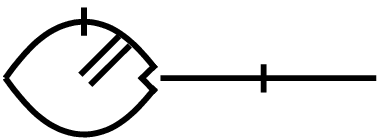} \fi}
\hspace{\hsqza}
\db{
\xy/r\scc/: \xypolygon5{~={36}~*{\alphanum}}
\endxy
}
\\[-10pt]
\hspace*{\fill} $\big[\hspace{-1.5ex}\bb{\xys{\ingt{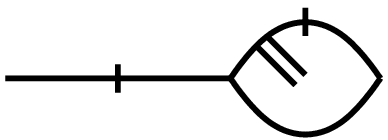}}}\hspace{-1.5ex}\big]$ \hspace*{10pt} 
\end{multicols}
\end{itemize}
\end{definition}

\bibliographystyle{gtart}
\bibliography{db-cn3}

\end{document}

%% file: ib-defs-forcniv.tex
\def\sectorh{4} % sector horizontal length
\def\dsectorh{7} % double sector h length
\def\tsectorh{10} % triple sector h length
\def\sectorv{3} % sector vertical length
\def\sectorvlow{1.65} % sector vertical length low
\def\sectorvmid{2.25} % sector vertical length mid
\def\sectorvhigh{4} % sector vertical length high
\def\dsectorvhigh{5.5} % vertical height of identity sectors of doulbe length in split situation.  Note, original was about 6.5, but because of different curve type, less is needed.
\def\tsectorvhigh{5.5} % vertical height of half identity sectors of triple length
\def\sectorhlong{5} % long sector length
\def\defecthsuplong{5} % defect super long

\def\sectorangle{60} % the angle at which the sector cookie begin.
\def\sectoranglelow{40}
\def\sectoranglemid{50}
\def\sectoranglehigh{65}
\def\dsectorangle{40}
\def\tsectorangle{28}
\def\dsectorhighangle{65}
\def\tsectorhighangle{65}
 
\def\vgap{.5} % distance between two sectors composed vertically
\def\hgap{.25} % distance associator ticks are away from picture
\def\exgap{.1} % extra gap when needed
\def\arcgap{.25} % gap at end of arc when an arrow arc is beginning there
\def\arcgapfudge{.1} % extension to cover the little end bits in the arc gap.
\def\arcgaparrowshiftv{.39} % vertical shift of caret placement in arrowarcgap
\def\arcgaparrowshifth{.24} % vertical shift of caret placement in arrowarcgap

\def\defecth{3} % defect (1-cell) horizontal length
\def\defectcgap{.5} % gap caused by caret in associator
\def\tickv{.75} % tick height
\def\tickvshort{.6} % short tick height
\def\tickvmini{.32} % mini tick height
\def\tickanglesharp{60} % angle of tick in a curving associator
\def\itickv{.5} % identity tick height
\def\itickvmini{.32} % mini id tick height
\def\tickvadd{.125} % extension of tick when near a curve
\def\iticksepv{.375} % identity tick separation for 2-cells
\def\itickseph{.36} % identity tick separation for 1-cells

\def\halfidsep{.375} % separation between lines in half identity sector
\def\halfidangle{45} % angle of half identity ticks.
\def\halfidangleshallow{35} % angle of half identity ticks when the main angle is too steep to fit.
\def\halfidanglesupshallow{18} % angle for the shallowest half id ticks
\def\hidticksangle{45} % angle for id ticks on an arc
\def\hidticksangleshallow{23} % angle for id ticks in triple sector
\def\hidshifth{1.27} % shift for id ticks on an arc % Note these are rotated shifts.
\def\hidshiftv{.18} % shift for id ticks on an arc
\def\hidshifthfact{1.21} % shift for id ticks on an arc when arc is \arcfactor stretched
\def\hidshiftvfact{.12} % shift for id ticks on an arc when arc is \arcfactor stretched
\def\hidshifttriph{1.6} % shift for id ticks on a triple arc % Note these are rotated shifts.
\def\hidshifttripv{.12} % shift for id ticks on an triple arc

\def\coordhhil{1.85}
\def\coordvhil{-.05}
\def\coordhhiltop{1.75}
\def\coordvhiltop{.5}
\def\coordhhilbotcut{2.6}
\def\coordvhilbotcut{-1.2}
\def\coordhhiltopthreecut{2}
\def\coordvhiltopthreecut{.75}
\def\coordhhilbotthreecut{2.6}
\def\coordvhilbotthreecut{-1.2}
\def\coordhhiltopthree{3}
\def\coordvhiltopthree{.5}
\def\coordhhilbot{3.5}
\def\coordvhilbot{-.75}
\def\coordhhilbotcutleft{3.35}%{5}
\def\coordvhilbotcutleft{-2.1}%{-3}

\def\ilength{2} % length of vertical identity lines
\def\ilengthd{4.5} % length of vertical identity lines in a two defect sector
\def\ilengtht{6.5} % length of vertical identity lines in a three defect sector
\def\ilengths{1.5} % length of vertical identity lines short
\def\assocv{1} % vertical height of associator

\def\caretangle{45} % caret angle
\def\caretlengthsect{.3} % caretlength in sectors
\def\caretlength{.375} % caretlength in arrows
\def\caretlengthmini{.2} % mini caretlength for tight spots

\def\arcfactor{.9} % factor by which arcs running along sectors are shrunk in vertical height.
\def\arcstretchfactor{1.1} % factor by which long arcs are stretched in vertical height.

\def\tickangle{13} % angle of tick in triple sector
\def\ttickh{1.31} % vertical position of tick and caret in triple sector, done by eye.
\def\dtickh{2.5} % vertical position of tick and caret in one exceptional half identity double width sector, namely hiltopthreecut.
\def\caretangle{11} %angle of caret in triple sector

\def\myinnersep{.15cm}
\def\myinnersepeps{.075cm}
\def\tableinnersep{.075cm}
\def\tableinnersepeps{.0375cm}
\def\linewidthvalue{.5985pt} % based on littlepictikzscale of .0315ex
\def\littlepicscale{.05ex}
\def\littlepicepsscale{.05ex}
\def\littlepictikzscale{.0315ex} % NB this is calibrated with respect to littlepicepsscale to produce the same size output.
\tikzset{littlepic/.style={scale=\littlepicscale, draw=black, line width=\linewidthvalue}}  %originally .08ex 
\tikzset{littlepiceps/.style={scale=\littlepicepsscale, draw=black, line width=\linewidthvalue}}  %originally .08ex 
\tikzset{littlepictikz/.style={scale=\littlepictikzscale, draw=black, line width=\linewidthvalue}}  %originally .08ex 
\tikzset{littlepicbull/.style={scale=.12ex, draw=black}}  %
\tikzset{littlenode/.style={fill=white,inner sep=\myinnersep}}
\tikzset{littlenode-eps/.style={fill=white,inner sep=\myinnersepeps}}
\tikzset{tablenode/.style={fill=white,inner sep=\tableinnersep}}
\tikzset{tablenode-eps/.style={fill=white,inner sep=\tableinnersepeps}}
\tikzset{displaypic/.style={scale=.15ex, draw=black}} 
\def\innertikz{\tikz[littlepic]}
\def\innertikzb{\tikz[littlepicbull]}
\def\innertikzeps{\tikz[littlepiceps]} % tricat eps
\def\innertikztikz{\tikz[littlepictikz]} % tricat tikz
\def\displaytikz{\tikz[displaypic]}

\def\tikztablescalex{.4ex}
\def\tikztablescaley{.05ex}
\def\tikztablescaleinner{.125ex}
\def\tikztablescaleinnerx{.35ex}
\def\tikztablescaleinnery{.125ex}

\def\sh(#1,#2){\def\x{#1}\def\y{#2}}

\def\sector{
	\coordinate (A) at (\x,\y);
	\coordinate (B) at ($ (A) + (.5*\sectorh,.5*\sectorv) $);
	\coordinate (C) at ($ (A) + (\sectorh,0) $);
	\coordinate (D) at ($ (A) + (.5*\sectorh,-.5*\sectorv) $);
	\draw (A) to [out=\sectorangle,in=180] (B) to [out=0,in=180-\sectorangle] (C);
	\draw (A) to [out=-\sectorangle,in=180] (D) to [out=0,in=180+\sectorangle] (C);
}

\def\arrowrarcu{
	\coordinate (A) at (\x,\y);
	\coordinate (B) at ($ (A) + (.5*\sectorh,\arcfactor*.5*\sectorv) $);
	\coordinate (C) at ($ (A) + (\sectorh,0) $);
	\draw (A) to [out=\sectorangle,in=180] (B) to [out=0,in=180-\sectorangle] (C);
	\caretr;
}

\def\arrowlarcu{
	\begin{scope}[cm={-1,0,0,1,(\x+\x+\sectorh,0)}]
		\arrowrarcu;
	\end{scope}
}

\def\arrowlarcd{
	\begin{scope}[cm={1,0,0,-1,(0,\y+\y)}]
		\arrowlarcu;
	\end{scope}
}

\def\arcu{
	\coordinate (A) at (\x,\y);
	\coordinate (B) at ($ (A) + (.5*\sectorh,\arcfactor*.5*\sectorv) $);
	\coordinate (C) at ($ (A) + (\sectorh,0) $);
	\draw (A) to [out=\sectorangle,in=180] (B) to [out=0,in=180-\sectorangle] (C);
}

\def\arcd{	
	\begin{scope}[cm={1,0,0,-1,(0,\y+\y)}]
		\arcu;
	\end{scope}
}

\def\arculong{
	\coordinate (A) at (\x,\y);
	\coordinate (B) at ($ (A) + (.5*\sectorh+.5,\arcstretchfactor*.5*\sectorv) $);
	\coordinate (C) at ($ (A) + (\sectorh+1,0) $);
	\draw (A) to [out=\sectorangle,in=180] (B) to [out=0,in=180-\sectorangle] (C);
}

\def\arcdlong{	
	\begin{scope}[cm={1,0,0,-1,(0,\y+\y)}]
		\arculong;
	\end{scope}
}

\def\arcugapl{
	\coordinate (D) at (\x-\arcgapfudge,\y-\arcgapfudge);
	\coordinate (E) at (\x+\arcgap,\y-\arcgapfudge);
	\coordinate (F) at (\x+\arcgap,\y+\arcgap+\arcgap+\arcgapfudge);
	\coordinate (G) at (\x-\arcgapfudge,\y+\arcgap+\arcgap+\arcgapfudge);
	\arcu;
	\draw[color=white,fill=white,line width=0pt] (D) -- (E) -- (F) -- (G) -- (D);
}

\def\arcugapr{
	\begin{scope}[cm={-1,0,0,1,(\x+\x+\sectorh,0)}]
		\arcugapl;
	\end{scope}
}

\def\arcdgapl{
	\begin{scope}[cm={1,0,0,-1,(0,\y+\y)}]
		\arcugapl;
	\end{scope}
}

\def\arcdgapr{
	\begin{scope}[cm={1,0,0,-1,(0,\y+\y)}]
		\arcugapr;
	\end{scope}
}

\def\arrowrarcugap{
	\arcugapl;
	\begin{scope}[shift={(\arcgaparrowshifth,\arcgaparrowshiftv)}]
		\caretr;
	\end{scope}
}

\def\arrowlarcugap{
	\begin{scope}[cm={-1,0,0,1,(\x+\x+\sectorh,0)}]
		\arrowrarcugap;
	\end{scope}
}

\def\arrowlarcdgap{
	\begin{scope}[cm={1,0,0,-1,(0,\y+\y)}]
		\arrowlarcugap;
	\end{scope}
}

\def\arrowrarcuturn{
	\arcugapl;
	\begin{scope}[shift={(.60,.57)}] %!% This is completely wrong, was tuned for a particular nonzero value of \y
	\begin{scope}[rotate=45]
		\caretrmini;
	\end{scope}
	\end{scope}
}

\def\arrowlarcuturn{
	\begin{scope}[cm={-1,0,0,1,(\x+\x+\sectorh,0)}]
		\arrowrarcuturn;
	\end{scope}
}

\def\arrowlarcdturn{
	\begin{scope}[cm={1,0,0,-1,(0,\y+\y)}]
		\arrowlarcuturn;
	\end{scope}
}

\def\splittopr{
	\coordinate (A) at (\x,\y);
	\coordinate (B) at ($ (A) + (.5*\sectorh,\arcfactor*.5*\sectorv) $);
	\coordinate (C) at ($ (A) + (\sectorh,0) $);
	\draw (A) to [out=\sectorangle,in=180] (B) to [out=0,in=180-\sectorangle] (C);
	\begin{scope}[shift={(.5*\sectorh,\arcfactor*.5*\sectorv)}]
		\caretrwhite;
	\end{scope}
}

\def\vid{
	\coordinate (A) at ($ (\x,\y) + (.5*\sectorh-.5*\ilength,-.5*\iticksepv) $);
	\coordinate (B) at ($ (\x,\y) + (.5*\sectorh+.5*\ilength,-.5*\iticksepv) $);
	\coordinate (C) at ($ (\x,\y) + (.5*\sectorh-.5*\ilength,.5*\iticksepv) $);
	\coordinate (D) at ($ (\x,\y) + (.5*\sectorh+.5*\ilength,.5*\iticksepv) $);
	\draw (A) to (B);
	\draw (C) to (D);
}

\def\vidd{
	\coordinate (A) at ($ (\x,\y) + (.5*\dsectorh-.5*\ilengthd,-.5*\iticksepv) $);
	\coordinate (B) at ($ (\x,\y) + (.5*\dsectorh+.5*\ilengthd,-.5*\iticksepv) $);
	\coordinate (C) at ($ (\x,\y) + (.5*\dsectorh-.5*\ilengthd,.5*\iticksepv) $);
	\coordinate (D) at ($ (\x,\y) + (.5*\dsectorh+.5*\ilengthd,.5*\iticksepv) $);
	\draw (A) to (B);
	\draw (C) to (D);
}

\def\vidt{
	\coordinate (A) at ($ (\x,\y) + (.5*\tsectorh-.5*\ilengtht,-.5*\iticksepv) $);
	\coordinate (B) at ($ (\x,\y) + (.5*\tsectorh+.5*\ilengtht,-.5*\iticksepv) $);
	\coordinate (C) at ($ (\x,\y) + (.5*\tsectorh-.5*\ilengtht,.5*\iticksepv) $);
	\coordinate (D) at ($ (\x,\y) + (.5*\tsectorh+.5*\ilengtht,.5*\iticksepv) $);
	\draw (A) to (B);
	\draw (C) to (D);
}

\def\vids{
	\coordinate (A) at ($ (\x,\y) + (.5*\sectorh-.5*\ilengths,-.5*\iticksepv) $);
	\coordinate (B) at ($ (\x,\y) + (.5*\sectorh+.5*\ilengths,-.5*\iticksepv) $);
	\coordinate (C) at ($ (\x,\y) + (.5*\sectorh-.5*\ilengths,.5*\iticksepv) $);
	\coordinate (D) at ($ (\x,\y) + (.5*\sectorh+.5*\ilengths,.5*\iticksepv) $);
	\draw (A) to (B);
	\draw (C) to (D);
}

\def\tripletopr{
	\coordinate (A) at (\x,\y);
	\coordinate (B) at ($ (A) + (.5*\tsectorh,.5*\sectorv) $);
	\coordinate (C) at ($ (A) + (\tsectorh,0) $);
	\draw (A) to [out=\tsectorangle,in=180] (B) to [out=0,in=180-\tsectorangle] (C);
	\begin{scope}[shift={(.31*\tsectorh,\ttickh)}]
		\tickcc;
	\end{scope}
	\begin{scope}[shift={(.69*\tsectorh,\ttickh)}]
		\caretrcwhite;
	\end{scope}
}

\def\tripletopl{
	\begin{scope}[cm={-1,0,0,1,(\x+\x+\tsectorh,0)}]
		\tripletopr;
	\end{scope}
}

\def\tripletoparc{
	\coordinate (A) at (\x,\y);
	\coordinate (B) at ($ (A) + (.5*\tsectorh,.5*\sectorv) $);
	\coordinate (C) at ($ (A) + (\tsectorh,0) $);
	\draw (A) to [out=\tsectorangle,in=180] (B) to [out=0,in=180-\tsectorangle] (C);
}

\def\triplebotr{
	\begin{scope}[cm={1,0,0,-1,(0,\y+\y)}]
		\tripletopr;
	\end{scope}
}

\def\triplebotl{
	\begin{scope}[cm={1,0,0,-1,(0,\y+\y)}]
		\tripletopl;
	\end{scope}
}

\def\triplebotarc{
	\begin{scope}[cm={1,0,0,-1,(0,\y+\y)}]
		\tripletoparc;
	\end{scope}
}

\def\tripletophigh{
	\coordinate (A) at (\x,\y);
	\coordinate (B) at ($ (A) + (.5*\tsectorh,.5*\tsectorvhigh) $);
	\coordinate (C) at ($ (A) + (\tsectorh,0) $);
	\draw (A) to [out=\tsectorhighangle,in=180] (B) to [out=0,in=180-\tsectorhighangle] (C);
}

\def\triplebothigh{
	\begin{scope}[cm={1,0,0,-1,(0,\y+\y)}]
		\tripletophigh;
	\end{scope}
}

\def\doubletop{
	\coordinate (A) at (\x,\y);
	\coordinate (B) at ($ (A) + (.5*\dsectorh,.5*\sectorv) $);
	\coordinate (C) at ($ (A) + (\dsectorh,0) $);
	\draw (A) to [out=\dsectorangle,in=180] (B) to [out=0,in=180-\dsectorangle] (C);
}

\def\doubletoptick{
	\doubletop;
	\begin{scope}[shift={(.5*\dsectorh,.5*\sectorv)}]
		\tick;
	\end{scope}
}

\def\doublebottick{
	\begin{scope}[cm={1,0,0,-1,(0,\y+\y)}]
		\doubletoptick;
	\end{scope}
}

\def\doubletopr{
	\doubletop;
	\begin{scope}[shift={(.5*\dsectorh,.5*\sectorv)}]
		\caretrwhite;
	\end{scope}
}

\def\doubletopl{
	\begin{scope}[cm={-1,0,0,1,(\x+\x+\dsectorh,0)}]
		\doubletopr;
	\end{scope}
}

\def\doublebotl{
	\begin{scope}[cm={1,0,0,-1,(0,\y+\y)}]
		\doubletopl;
	\end{scope}
}

\def\doubletophigh{
	\coordinate (A) at (\x,\y);
	\coordinate (B) at ($ (A) + (.5*\dsectorh,.5*\dsectorvhigh) $);
	\coordinate (C) at ($ (A) + (\dsectorh,0) $);
	\draw (A) to [out=\dsectorhighangle,in=180] (B) to [out=0,in=180-\dsectorhighangle] (C);
}

\def\doublebothigh{
	\begin{scope}[cm={1,0,0,-1,(0,\y+\y)}]
		\doubletophigh;
	\end{scope}
}

\def\caretrcwhite{
	\begin{scope}[rotate around={-\caretangle:(\x,\y)}]
		\caretrwhite;
	\end{scope}
}

\def\caretrwhite{
	\coordinate (A) at (\x,\y);
	\coordinate (B) at ($ (A) + (-\caretlength,\caretlength) $);
	\coordinate (C) at ($ (A) + (-\caretlength,-\caretlength) $);
	\coordinate (D) at ($ (A) + (-1.2*\caretlength,1.2*\caretlength) $);
	\coordinate (E) at ($ (A) + (-1.2*\caretlength,-1.2*\caretlength) $);
	\draw[color=white,fill=white,line width=0pt] (D) -- (A) -- (E) -- (D);
	\draw (B) to (A) to (C);
}

\def\caretr{
	\coordinate (A) at (\x,\y);
	\coordinate (B) at ($ (A) + (-\caretlength,\caretlength) $);
	\coordinate (C) at ($ (A) + (-\caretlength,-\caretlength) $);
	\draw (B) to (A) to (C);
}

\def\caretlwhite{
	\begin{scope}[cm={-1,0,0,1,(\x+\x+\caretlength,0)}]
		\caretrwhite;
	\end{scope}
}

\def\caretrmini{
	\coordinate (A) at (\x,\y);
	\coordinate (B) at ($ (A) + (-\caretlengthmini,\caretlengthmini) $);
	\coordinate (C) at ($ (A) + (-\caretlengthmini,-\caretlengthmini) $);
	\draw (B) to (A) to (C);
}

\def\tick{
	\coordinate (A) at ($ (\x,\y) + (0,.5*\tickv) $);
	\coordinate (B) at ($ (\x,\y) + (0,-.5*\tickv) $);
	\draw (A) to (B);
}

\def\tickaddup{
	\coordinate (A) at ($ (\x,\y) + (0,.5*\tickv+\tickvadd) $);
	\coordinate (B) at ($ (\x,\y) + (0,-.5*\tickv) $);
	\draw (A) to (B);
}

\def\tickadddown{
	\coordinate (A) at ($ (\x,\y) + (0,.5*\tickv) $);
	\coordinate (B) at ($ (\x,\y) + (0,-.5*\tickv-\tickvadd) $);
	\draw (A) to (B);
}

\def\tickshortupadddown{
	\coordinate (A) at ($ (\x,\y) + (0,.5*\tickv-\tickvadd) $);
	\coordinate (B) at ($ (\x,\y) + (0,-.5*\tickv-\tickvadd) $);
	\draw (A) to (B);
}

\def\tickshortdownaddup{
	\coordinate (A) at ($ (\x,\y) + (0,.5*\tickv+\tickvadd) $);
	\coordinate (B) at ($ (\x,\y) + (0,-.5*\tickv+\tickvadd) $);
	\draw (A) to (B);
}

\def\tickshort{
	\coordinate (A) at ($ (\x,\y) + (0,.5*\tickvshort) $);
	\coordinate (B) at ($ (\x,\y) + (0,-.5*\tickvshort) $);
	\draw (A) to (B);
}

\def\tickmini{
	\coordinate (A) at ($ (\x,\y) + (0,.5*\tickvmini) $);
	\coordinate (B) at ($ (\x,\y) + (0,-.5*\tickvmini) $);
	\draw (A) to (B);
}

\def\itickshere{
	\coordinate (A) at ($ (\x,\y) + (-.5*\itickseph,.5*\itickv) $);
	\coordinate (B) at ($ (\x,\y) + (-.5*\itickseph,-.5*\itickv) $);
	\coordinate (C) at ($ (\x,\y) + (.5*\itickseph,.5*\itickv) $);
	\coordinate (D) at ($ (\x,\y) + (.5*\itickseph,-.5*\itickv) $);	
	\draw (A) to (B);
	\draw (C) to (D);
}

\def\iticks{
	\coordinate (A) at ($ (\x,\y) + (.5*\defecth-.5*\itickseph,.5*\itickv) $);
	\coordinate (B) at ($ (\x,\y) + (.5*\defecth-.5*\itickseph,-.5*\itickv) $);
	\coordinate (C) at ($ (\x,\y) + (.5*\defecth+.5*\itickseph,.5*\itickv) $);
	\coordinate (D) at ($ (\x,\y) + (.5*\defecth+.5*\itickseph,-.5*\itickv) $);	
	\draw (A) to (B);
	\draw (C) to (D);
}

\def\iticksminihere{
	\coordinate (A) at ($ (\x,\y) + (-.5*\itickseph,.5*\itickvmini) $);
	\coordinate (B) at ($ (\x,\y) + (-.5*\itickseph,-.5*\itickvmini) $);
	\coordinate (C) at ($ (\x,\y) + (.5*\itickseph,.5*\itickvmini) $);
	\coordinate (D) at ($ (\x,\y) + (.5*\itickseph,-.5*\itickvmini) $);	
	\draw (A) to (B);
	\draw (C) to (D);
}

\def\iticksrotl{
	\begin{scope}[rotate around={\hidticksangle:(\x,\y)}]
	\begin{scope}[shift={(\hidshifth,\hidshiftv)}]
		\itickshere;
	\end{scope}
	\end{scope} %% Remark: the shifts end up being rotated; but if you try to do it the other way it doesn't work.
}

\def\iticksrotr{
	\begin{scope}[rotate around={-\hidticksangle:(\x,\y)}]
	\begin{scope}[shift={(.707*\sectorh-\hidshifth,.707*\sectorh+\hidshiftv)}]
		\itickshere;
	\end{scope}
	\end{scope} 
}

\def\iticksrotlminifact{
	\begin{scope}[rotate around={\hidticksangle:(\x,\y)}]
	\begin{scope}[shift={(\hidshifthfact,\hidshiftvfact)}]
		\iticksminihere;
	\end{scope}
	\end{scope} 
}

\def\iticksrotrminifact{
	\begin{scope}[rotate around={-\hidticksangle:(\x,\y)}]
	\begin{scope}[shift={(.707*\sectorh-\hidshifthfact,.707*\sectorh+\hidshiftvfact)}]
		\iticksminihere;
	\end{scope}
	\end{scope} 
}

\def\iticksrotlfact{
	\begin{scope}[rotate around={\hidticksangle:(\x,\y)}]
	\begin{scope}[shift={(\hidshifthfact,\hidshiftvfact)}]
		\itickshere;
	\end{scope}
	\end{scope} 
}

\def\iticksrotrfact{
	\begin{scope}[rotate around={-\hidticksangle:(\x,\y)}]
	\begin{scope}[shift={(.707*\sectorh-\hidshifthfact,.707*\sectorh+\hidshiftvfact)}]
		\itickshere;
	\end{scope}
	\end{scope} 
}

\def\itickstriprotlshallow{ %!%
	\begin{scope}[rotate around={\hidticksangleshallow:(\x,\y)}]
	\begin{scope}[shift={(\hidshifttriph,\hidshifttripv)}]
		\itickshere;
	\end{scope}
	\end{scope} 
}

\def\itickstriprotrshallow{ %!%
	\begin{scope}[rotate around={-\hidticksangleshallow:(\x,\y)}]
	\begin{scope}[shift={(-\hidshifttriph,\hidshifttripv)}]
		\itickshere;
	\end{scope}
	\end{scope} 
}

\def\tickcc{
	\begin{scope}[rotate around={\tickangle:(\x,\y)}]
		\tick;
	\end{scope}
}

\def\tickminicc{
	\begin{scope}[rotate around={\tickanglesharp:(\x,\y)}]
		\tickmini;
	\end{scope}
}

\def\defect{
	\coordinate (A) at (\x,\y);
	\coordinate (B) at ($ (A) + (\defecth,0) $);
	\draw (A) to (B);
}

\def\defectsemilong{
	\coordinate (A) at (\x,\y);
	\coordinate (B) at ($ (A) + (\defecth+.5,0) $); %!%
	\draw (A) to (B);
}

\def\defectlong{
	\coordinate (A) at (\x,\y);
	\coordinate (B) at ($ (A) + (\sectorh,0) $);
	\draw (A) to (B);
}

\def\defectsuplong{
	\coordinate (A) at (\x,\y);
	\coordinate (B) at ($ (A) + (\defecthsuplong,0) $);
	\draw (A) to (B);
}

\def\defectshortr{
	\coordinate (A) at (\x,\y);
	\coordinate (B) at ($ (A) + (\defecth-\defectcgap,0) $);
	\draw (A) to (B);
}

\def\defectshortl{
	\coordinate (A) at (\x+\defectcgap,\y);
	\coordinate (B) at ($ (A) + (\defecth-\defectcgap,0) $);
	\draw (A) to (B);
}

\def\defectshorterr{
	\coordinate (A) at (\x,\y);
	\coordinate (B) at ($ (A) + (\defecth-\defectcgap-\exgap,0) $);
	\draw (A) to (B);
}

\def\defectshorterl{
	\coordinate (A) at (\x+\defectcgap+\exgap,\y);
	\coordinate (B) at ($ (A) + (\defecth-\defectcgap-\exgap,0) $);
	\draw (A) to (B);
}

\def\strutv{
	\coordinate (A) at (\x,\y);
	\coordinate (B) at ($ (A) + (0,\assocv) $);
	\draw (A) to (B);
}

\def\strutlong{
	\coordinate (A) at (\x,\y);
	\coordinate (B) at ($ (A) + (0,\assocv+.5*\vgap) $);
	\draw (A) to (B);
}

\def\strutshort{
	\coordinate (A) at (\x,\y);
	\coordinate (B) at ($ (A) + (0,\assocv-.5*\vgap) $);
	\draw (A) to (B);
}

\def\strutshortdown{
	\coordinate (A) at (\x,\y);
	\coordinate (B) at ($ (A) - (0,\assocv-.5*\vgap) $);
	\draw (A) to (B);
}

\def\assoctickr{
	\coordinate (A) at (\x + \hgap,\y);
	\coordinate (B) at (\x + \hgap,\y + \vgap);
	\draw (A) to [out=45,in=-45] (B);
}

\def\assoctickl{
	\coordinate (A) at (\x - \hgap,\y);
	\coordinate (B) at (\x - \hgap,\y + \vgap);
	\draw (A) to [out=135,in=-135] (B);
}

\def\arrowr{
	\caretr;
	\defect;
}

\def\arrowrlong{
	\caretr;
	\defectlong;
}

\def\arrowl{
	\begin{scope}[cm={-1,0,0,1,(\x+\x+\defecth,0)}]
		\arrowr;
	\end{scope}
}

\def\arrowllong{
	\begin{scope}[cm={-1,0,0,1,(\x+\x+\sectorh,0)}]
		\arrowrlong;
	\end{scope}
}

\def\arrowrgap{
	\begin{scope}[cm={1,0,0,1,(\defectcgap,0)}]
		\caretr;
	\end{scope}
	\defectshortl;
}	

\def\arrowrgapbig{
	\begin{scope}[cm={1,0,0,1,(\defectcgap+\exgap,0)}]
		\caretr;
	\end{scope}
	\defectshorterl;
}

\def\arrowlgap{
	\begin{scope}[cm={-1,0,0,1,(\x+\x+\defecth-\defectcgap,0)}]
		\caretr;
	\end{scope}
	\defectshortr;
}

\def\arrowlgapbig{
	\begin{scope}[cm={-1,0,0,1,(\x+\x+\defecth-\defectcgap-\exgap,0)}]
		\caretr;
	\end{scope}
	\defectshorterr;
}

\def\sectorcaretr{
	\coordinate (A) at (\x,\y);
	\coordinate (B) at ($ (A) + (.5*\sectorh,.5*\sectorv) $);
	\coordinate (C) at ($ (A) + (\sectorh-\caretlengthsect,\caretlengthsect) $);
	\coordinate (CCC) at ($ (A) + (\sectorh-\caretlengthsect,-\caretlengthsect) $);
	\coordinate (CC) at ($ (A) + (\sectorh-2*\caretlengthsect,0) $);
	\coordinate (D) at ($ (A) + (.5*\sectorh,-.5*\sectorv) $);
	\draw (A) to [out=\sectorangle,in=180] (B) to [out=0,in=180-\sectorangle] (C) to (CC) to (CCC) to [out=180+\sectorangle,in=0] (D) to [out=180,in=-\sectorangle] (A);
}	

\def\sectorcaretl{
	\begin{scope}[cm={-1,0,0,1,(\x+\x+\sectorh,0)}]
		\sectorcaretr;
	\end{scope}
}

\def\sectortoparc{
	\coordinate (A) at (\x,\y);
	\coordinate (B) at ($ (A) + (.5*\sectorh,.5*\sectorv) $);
	\coordinate (C) at ($ (A) + (\sectorh,0) $);
	\draw (A) to [out=\sectorangle,in=180] (B) to [out=0,in=180-\sectorangle] (C);
}

\def\sectorbotarc{
	\coordinate (A) at (\x,\y);
	\coordinate (C) at ($ (A) + (\sectorh,0) $);
	\coordinate (D) at ($ (A) + (.5*\sectorh,-.5*\sectorv) $);
	\draw (A) to [out=-\sectorangle,in=180] (D) to [out=0,in=180+\sectorangle] (C);

}

\def\sectorbotarclong{
	\coordinate (A) at (\x,\y);
	\coordinate (C) at ($ (A) + (\sectorhlong,0) $);
	\coordinate (D) at ($ (A) + (.5*\sectorhlong,-.5*\sectorv) $);
	\draw (A) to [out=-\sectorangle,in=180] (D) to [out=0,in=180+\sectorangle] (C);

}

\def\sectortop{
	\sectortoparc;
	\draw (A) to (C);
}

\def\sectorbot{
	\sectorbotarc;
	\draw (A) to (C);
}

\def\sectortoparclow{
	\coordinate (A) at (\x,\y);
	\coordinate (B) at ($ (A) + (.5*\sectorh,.5*\sectorvlow) $);
	\coordinate (C) at ($ (A) + (\sectorh,0) $);
	\draw (A) to [out=\sectoranglelow,in=180] (B) to [out=0,in=180-\sectoranglelow] (C);
}

\def\sectortoparcmid{
	\coordinate (A) at (\x,\y);
	\coordinate (B) at ($ (A) + (.5*\sectorh,.5*\sectorvmid) $);
	\coordinate (C) at ($ (A) + (\sectorh,0) $);
	\draw (A) to [out=\sectoranglemid,in=180] (B) to [out=0,in=180-\sectoranglemid] (C);
}

\def\sectortoparchigh{
	\coordinate (A) at (\x,\y);
	\coordinate (B) at ($ (A) + (.5*\sectorh,.5*\sectorvhigh) $);
	\coordinate (C) at ($ (A) + (\sectorh,0) $);
	\draw (A) to [out=\sectoranglehigh,in=180] (B) to [out=0,in=180-\sectoranglehigh] (C);
}

\def\sectorbotarclow{
	\coordinate (A) at (\x,\y);
	\coordinate (B) at ($ (A) + (.5*\sectorh,-.5*\sectorvlow) $);
	\coordinate (C) at ($ (A) + (\sectorh,0) $);
	\draw (A) to [out=-\sectoranglelow,in=180] (B) to [out=0,in=180+\sectoranglelow] (C);
}

\def\sectorbotarcmid{
	\coordinate (A) at (\x,\y);
	\coordinate (B) at ($ (A) + (.5*\sectorh,-.5*\sectorvmid) $);
	\coordinate (C) at ($ (A) + (\sectorh,0) $);
	\draw (A) to [out=-\sectoranglemid,in=180] (B) to [out=0,in=180+\sectoranglemid] (C);
}

\def\sectorbotarchigh{
	\coordinate (A) at (\x,\y);
	\coordinate (B) at ($ (A) + (.5*\sectorh,-.5*\sectorvhigh) $);
	\coordinate (C) at ($ (A) + (\sectorh,0) $);
	\draw (A) to [out=-\sectoranglehigh,in=180] (B) to [out=0,in=180+\sectoranglehigh] (C);
}

\def\hil{ % half id left
	\coordinate (A) at (\x,\y);
	\coordinate (B) at ($ (A) + (.5*\sectorh,.5*\sectorv) $);
	\coordinate (BB) at ($ (A) + (.5*\sectorh,-.5*\sectorv) $);
	\coordinate (C) at ($ (A) + (\sectorh,0) $);
	\coordinate (D) at ($ (A) + (\coordhhil,\coordvhil) $);
	\coordinate (E) at ($ (D) + (90-\halfidangle:.5*\halfidsep) $);
	\coordinate (F) at ($ (D) + (90-\halfidangle:-.5*\halfidsep) $);
	\coordinate (T) at ($ (A) + (.5*\sectorh,.5*\sectorv+.5*\tickv) $);
	\coordinate (S) at ($ (A) + (.5*\sectorh,.5*\sectorv-.5*\tickv) $);		
	\begin{scope}
		\clip (A) to [out=\sectorangle,in=180] (B) to [out=0,in=180-\sectorangle] (C) to [out=180+\sectorangle,in=0] (BB) to [out=180,in=-\sectorangle] (A);
		\draw (E) to +(180-\halfidangle:3);
		\draw (F) to +(180-\halfidangle:3);
	\end{scope}
	\sectortoparc;
	\sectorbotarc;
	\draw (T) to (S);
}

\def\hir{
	\begin{scope}[cm={-1,0,0,1,(\x+\x+\sectorh,0)}]
		\hil;
	\end{scope}
}

\def\hiltop{
	\coordinate (A) at (\x,\y);
	\coordinate (B) at ($ (A) + (.5*\sectorh,.5*\sectorv) $);
	\coordinate (C) at ($ (A) + (\sectorh,0) $);
	\coordinate (D) at ($ (A) + (\coordhhiltop,\coordvhiltop) $);
	\coordinate (E) at ($ (D) + (90-\halfidangle:.5*\halfidsep) $);
	\coordinate (F) at ($ (D) + (90-\halfidangle:-.5*\halfidsep) $);
	\coordinate (T) at ($ (A) + (.5*\sectorh,.5*\sectorv+.5*\tickv) $);
	\coordinate (S) at ($ (A) + (.5*\sectorh,.5*\sectorv-.5*\tickv) $);		
	\begin{scope}
		\clip (A) to [out=\sectorangle,in=180] (B) to [out=0,in=180-\sectorangle] (C) to (A);
		\draw (E) to +(180-\halfidangle:3);
		\draw (F) to +(180-\halfidangle:3);
	\end{scope}
	\sectortop;
	\draw (T) to (S);
}

\def\hirtop{
	\begin{scope}[cm={-1,0,0,1,(\x+\x+\sectorh,0)}]
		\hiltop;
	\end{scope}
}

\def\hilbotcut{
	\coordinate (A) at (\x,\y);
	\coordinate (B) at ($ (A) + (\defecth,0) $);
	\coordinate (C) at ($ (A) + (\defecth,-\defecth) $);
	\coordinate (D) at ($ (A) + (0,-\defecth) $);
	\coordinate (S) at ($ (A) + (\coordhhilbotcut,\coordvhilbotcut) $);
	\coordinate (E) at ($ (S) + (90-\halfidangle:.5*\halfidsep) $);
	\coordinate (F) at ($ (S) + (90-\halfidangle:-.5*\halfidsep) $);
	\begin{scope}
		\clip (A) to (B) to (C) to (D) to (A);
		\draw (E) to +(180-\halfidangle:3);
		\draw (F) to +(180-\halfidangle:3);
	\end{scope}
	\doublebothigh;
	\defect;
	\begin{scope}[shift={(\defecth,0)}]
		\arcd;
	\end{scope}
}

\def\hilbotcuttight{
	\coordinate (A) at (\x,\y);
	\coordinate (B) at ($ (A) + (\defecth,0) $);
	\coordinate (C) at ($ (A) + (\defecth,-\defecth) $);
	\coordinate (D) at ($ (A) + (0,-\defecth) $);
	\coordinate (S) at ($ (A) + (\coordhhilbotcut,\coordvhilbotcut) $);
	\coordinate (E) at ($ (S) + (90-\halfidangle:.5*\halfidsep) $);
	\coordinate (F) at ($ (S) + (90-\halfidangle:-.5*\halfidsep) $);
	\begin{scope}
		\clip (A) to (B) to (C) to (D) to (A);
		\draw (E) to +(180-\halfidangle:3);
		\draw (F) to +(180-\halfidangle:3);
	\end{scope}
	\doublebothigh;
	\defect;
	\begin{scope}[shift={(\defecth,0)}]
	\end{scope}
}

\def\hirbotcut{
	\begin{scope}[cm={-1,0,0,1,(\x+\x+\sectorh+\defecth,0)}]
		\hilbotcut;
	\end{scope}
}

\def\hirbotcuttight{
	\begin{scope}[cm={-1,0,0,1,(\x+\x+\sectorh+\defecth,0)}]
		\hilbotcuttight;
	\end{scope}
}

\def\hiltopthreecut{
	\coordinate (A) at (\x,\y);
	\coordinate (B) at ($ (A) + (.5*\dsectorh,.5*\dsectorvhigh) $);
	\coordinate (C) at ($ (A) + (\dsectorh,0) $);
	\coordinate (S) at ($ (A) + (\coordhhiltopthreecut,\coordvhiltopthreecut) $);
	\coordinate (E) at ($ (S) + (90-\halfidangleshallow:.5*\halfidsep) $);
	\coordinate (F) at ($ (S) + (90-\halfidangleshallow:-.5*\halfidsep) $);
	\begin{scope}
		\clip (A) to [out=\dsectorhighangle,in=180] (B) to [out=0,in=180-\dsectorhighangle] (C) to (A);
		\draw (E) to +(180-\halfidangleshallow:3);
		\draw (F) to +(180-\halfidangleshallow:3);
	\end{scope}
	\doubletophigh;
	\defect;
	\begin{scope}[shift={(\defecth,0)}]
		\arcu;
	\end{scope}
	\begin{scope}[shift={(.69*\dsectorh,\dtickh)}]
	\begin{scope}[rotate around={-\tickangle-10:(\x,\y)}]
		\tick;
	\end{scope}
	\end{scope}
	\begin{scope}[shift={(.225*\dsectorh,\dtickh-.265)}] %!%!%
	\begin{scope}[rotate around={\tickangle+10:(\x,\y)}]
		\caretlwhite;
	\end{scope}
	\end{scope}	
}

\def\hirtopthreecut{
	\begin{scope}[cm={-1,0,0,1,(\x+\x+\sectorh+\defecth,0)}]
		\hiltopthreecut;
	\end{scope}
}

\def\hilbotthreecut{
	\coordinate (A) at (\x,\y);
	\coordinate (B) at ($ (A) + (1.2*\defecth,0) $);
	\coordinate (C) at ($ (A) + (1.2*\defecth,-\defecth) $);
	\coordinate (D) at ($ (A) + (0,-\defecth) $);
	\coordinate (S) at ($ (A) + (\coordhhilbotthreecut,\coordvhilbotthreecut) $);
	\coordinate (E) at ($ (S) + (90-\halfidangle:.5*\halfidsep) $);
	\coordinate (F) at ($ (S) + (90-\halfidangle:-.5*\halfidsep) $);
	\begin{scope}
		\clip (A) to (B) to (C) to (D) to (A);
		\draw (E) to +(180-\halfidangle:3);
		\draw (F) to +(180-\halfidangle:3);
	\end{scope}
	\triplebothigh;
	\arrowl;
	\begin{scope}[shift={(\defecth,0)}]
		\defectshortl;
	\end{scope}
	\begin{scope}[shift={(\defecth+\defecth,0)}]
		\arcd;
	\end{scope}
	\begin{scope}[shift={(.5*\sectorh+.5*\defecth+.5*\defecth,-.5*\tsectorvhigh)}]
		\tick;
	\end{scope}
}

\def\hirbotthreecut{
	\begin{scope}[cm={-1,0,0,1,(\x+\x+\sectorh+\defecth+\defecth,0)}]
		\hilbotthreecut;
	\end{scope}
}

\def\hilbotthreecutmid{
	\coordinate (A) at (\x,\y);
	\coordinate (B) at ($ (A) + (1.2*\defecth,0) $);
	\coordinate (C) at ($ (A) + (1.2*\defecth,-\defecth) $);
	\coordinate (D) at ($ (A) + (0,-\defecth) $);
	\coordinate (S) at ($ (A) + (\coordhhilbotthreecut,\coordvhilbotthreecut) $);
	\coordinate (E) at ($ (S) + (90-\halfidangle:.5*\halfidsep) $);
	\coordinate (F) at ($ (S) + (90-\halfidangle:-.5*\halfidsep) $);
	\begin{scope}
		\clip (A) to (B) to (C) to (D) to (A);
		\draw (E) to +(180-\halfidangle:3);
		\draw (F) to +(180-\halfidangle:3);
	\end{scope}
	\triplebothigh;
	\arrowlgap;
	\begin{scope}[shift={(\defecth,0)}]
		\arcd;
	\end{scope}
	\begin{scope}[shift={(\defecth+\sectorh,0)}]
		\defect;
		\tickadddown;
	\end{scope}
	\begin{scope}[shift={(.5*\sectorh+.5*\defecth+.5*\defecth,-.5*\tsectorvhigh)}]
		\tick;
	\end{scope}
}

\def\hilbotthreecutmidshorttick{
	\coordinate (A) at (\x,\y);
	\coordinate (B) at ($ (A) + (1.2*\defecth,0) $);
	\coordinate (C) at ($ (A) + (1.2*\defecth,-\defecth) $);
	\coordinate (D) at ($ (A) + (0,-\defecth) $);
	\coordinate (S) at ($ (A) + (\coordhhilbotthreecut,\coordvhilbotthreecut) $);
	\coordinate (E) at ($ (S) + (90-\halfidangle:.5*\halfidsep) $);
	\coordinate (F) at ($ (S) + (90-\halfidangle:-.5*\halfidsep) $);
	\begin{scope}
		\clip (A) to (B) to (C) to (D) to (A);
		\draw (E) to +(180-\halfidangle:3);
		\draw (F) to +(180-\halfidangle:3);
	\end{scope}
	\triplebothigh;
	\arrowlgap;
	\begin{scope}[shift={(\defecth,0)}]
		\arcd;
	\end{scope}
	\begin{scope}[shift={(\defecth+\sectorh,0)}]
		\defect;
		\tickshortupadddown;
	\end{scope}
	\begin{scope}[shift={(.5*\sectorh+.5*\defecth+.5*\defecth,-.5*\tsectorvhigh)}]
		\tick;
	\end{scope}
}

\def\hilbotthreecutmidnotick{
	\coordinate (A) at (\x,\y);
	\coordinate (B) at ($ (A) + (1.2*\defecth,0) $);
	\coordinate (C) at ($ (A) + (1.2*\defecth,-\defecth) $);
	\coordinate (D) at ($ (A) + (0,-\defecth) $);
	\coordinate (S) at ($ (A) + (\coordhhilbotthreecut,\coordvhilbotthreecut) $);
	\coordinate (E) at ($ (S) + (90-\halfidangle:.5*\halfidsep) $);
	\coordinate (F) at ($ (S) + (90-\halfidangle:-.5*\halfidsep) $);
	\begin{scope}
		\clip (A) to (B) to (C) to (D) to (A);
		\draw (E) to +(180-\halfidangle:3);
		\draw (F) to +(180-\halfidangle:3);
	\end{scope}
	\triplebothigh;
	\arrowlgap;
	\begin{scope}[shift={(\defecth,0)}]
		\arcd;
	\end{scope}
	\begin{scope}[shift={(\defecth+\sectorh,0)}]
		\defect;
	\end{scope}
	\begin{scope}[shift={(.5*\sectorh+.5*\defecth+.5*\defecth,-.5*\tsectorvhigh)}]
		\tick;
	\end{scope}
}

\def\hirbotthreecutmidnotick{
	\begin{scope}[cm={-1,0,0,1,(\x+\x+\tsectorh,0)}]
		\hilbotthreecutmidnotick;
	\end{scope}
}

\def\hiltopthree{
	\coordinate (A) at (\x,\y);
	\coordinate (B) at ($ (A) + (.5*\tsectorh,.5*\sectorv) $);
	\coordinate (C) at ($ (A) + (\tsectorh,0) $);
	\coordinate (S) at ($ (A) + (\coordhhiltopthree,\coordvhiltopthree) $);
	\coordinate (E) at ($ (S) + (90-\halfidanglesupshallow:.5*\halfidsep) $);
	\coordinate (F) at ($ (S) + (90-\halfidanglesupshallow:-.5*\halfidsep) $);
	\begin{scope}
		\clip (A) to [out=\tsectorangle,in=180] (B) to [out=0,in=180-\tsectorangle] (C) to (A);
		\draw (E) to +(180-\halfidanglesupshallow:3);
		\draw (F) to +(180-\halfidanglesupshallow:3);
	\end{scope}
	\draw (C) to (A);
	\tripletopl;
}

\def\hirtopthree{
	\begin{scope}[cm={-1,0,0,1,(\x+\x+\tsectorh,0)}]
		\hiltopthree;
	\end{scope}
}

\def\hilbot{
	\coordinate (A) at (\x,\y);
	\coordinate (B) at ($ (A) + (\defecthsuplong,0) $);
	\coordinate (C) at ($ (A) + (\defecthsuplong,-\defecth) $);
	\coordinate (D) at ($ (A) + (0,-\defecth) $);
	\coordinate (H) at ($ (A) + (\tsectorh,0) $);
	\coordinate (S) at ($ (A) + (\coordhhilbot,\coordvhilbot) $);
	\coordinate (E) at ($ (S) + (90-\halfidangleshallow:.5*\halfidsep) $);
	\coordinate (F) at ($ (S) + (90-\halfidangleshallow:-.5*\halfidsep) $);
	\begin{scope}
		\clip (A) to (B) to (C) to (D) to (A);
		\draw (E) to +(180-\halfidangleshallow:3);
		\draw (F) to +(180-\halfidangleshallow:3);
	\end{scope}
	\triplebotarc;
	\draw (A) to (H);
	\begin{scope}[shift={(\defecthsuplong,0)}]
		\tick;
	\end{scope}
}

\def\hilbottight{
	\coordinate (A) at (\x,\y);
	\coordinate (B) at ($ (A) + (\defecthsuplong,0) $);
	\coordinate (C) at ($ (A) + (\defecthsuplong,-\defecth) $);
	\coordinate (D) at ($ (A) + (0,-\defecth) $);
	\coordinate (H) at ($ (A) + (\tsectorh,0) $);
	\coordinate (S) at ($ (A) + (\coordhhilbot,\coordvhilbot) $);
	\coordinate (E) at ($ (S) + (90-\halfidangleshallow:.5*\halfidsep) $);
	\coordinate (F) at ($ (S) + (90-\halfidangleshallow:-.5*\halfidsep) $);
	\begin{scope}
		\clip (A) to (B) to (C) to (D) to (A);
		\draw (E) to +(180-\halfidangleshallow:3);
		\draw (F) to +(180-\halfidangleshallow:3);
	\end{scope}
	\triplebotarc;
	\begin{scope}[shift={(\defecthsuplong,0)}]
	\end{scope}

}

\def\hirbot{
	\begin{scope}[cm={-1,0,0,1,(\x+\x+\tsectorh,0)}]
		\hilbot;
	\end{scope}
}

\def\hirbottight{
	\begin{scope}[cm={-1,0,0,1,(\x+\x+\tsectorh,0)}]
		\hilbottight;
	\end{scope}
}

\def\hilbotcutleft{
	\coordinate (A) at (\x,\y);
	\coordinate (BT) at ($ (A) + (.5*\sectorh,-.5*\sectorv) $); %(.5*\sectorh,-.5*\sectorv)
	\coordinate (CT) at ($ (A) + (\sectorh,0) $);
	\coordinate (DT) at ($ (A) + (\sectorh+\defecth,0) $);
	\coordinate (DDT) at ($ (A) + (.5*\sectorh+.5*\defecth,-.5*\dsectorvhigh) $);
	\coordinate (ST) at ($ (A) + (\coordhhilbotcutleft,\coordvhilbotcutleft) $);
	\coordinate (ET) at ($ (ST) + (90-\halfidangleshallow:.5*\halfidsep) $);
	\coordinate (FT) at ($ (ST) + (90-\halfidangleshallow:-.5*\halfidsep) $);
	\begin{scope}
		\clip (A) to [out=-\sectorangle,in=180] (BT) to [out=0,in=180+\sectorangle] (CT) to (DT) to [out=180+\sectorangle,in=0] (DDT) to [out=180,in=-\sectorangle] (A);
		\draw (ET) to +(180-\halfidangleshallow:3);
		\draw (FT) to +(180-\halfidangleshallow:3);
	\end{scope}
	\doublebothigh;
	\begin{scope}[shift={(\sectorh,0)}]
		\defect;
	\end{scope}
}

\def\ax#1-#2{
\ifthenelse{\equal{#1}{1}}{
			\axiombA{#2}}{}
\ifthenelse{\equal{#1}{2}}{
			\axiombB{#2}}{}
\ifthenelse{\equal{#1}{3}}{
			\axiombC{#2}}{}
\ifthenelse{\equal{#1}{4}}{
			\axiombD{#2}}{}
\ifthenelse{\equal{#1}{5}}{
			\axiombE{#2}}{}
\ifthenelse{\equal{#1}{6}}{
			\axiombF{#2}}{}
\ifthenelse{\equal{#1}{7}}{
			\axiombG{#2}}{}
\ifthenelse{\equal{#1}{8}}{
			\axiombH{#2}}{}
\ifthenelse{\equal{#1}{9}}{
			\axiombI{#2}}{}
\ifthenelse{\equal{#1}{10}}{
			\axiombJ{#2}}{}			
\ifthenelse{\equal{#1}{11}}{
			\axiombK{#2}}{}
\ifthenelse{\equal{#1}{12}}{
			\axiombL{#2}}{}			
\ifthenelse{\equal{#1}{13}}{
			\axiombM{#2}}{}
\ifthenelse{\equal{#1}{14}}{
			\axiombN{#2}}{}			
\ifthenelse{\equal{#1}{15}}{
			\axiombO{#2}}{}
\ifthenelse{\equal{#1}{16}}{
			\axiombP{#2}}{}
\ifthenelse{\equal{#1}{17}}{
			\axiombQ{#2}}{}
\ifthenelse{\equal{#1}{18}}{
			\axiombR{#2}}{}
\ifthenelse{\equal{#1}{19}}{
			\axiombS{#2}}{}
\ifthenelse{\equal{#1}{20}}{
			\axiombT{#2}}{}			
\ifthenelse{\equal{#1}{21}}{
			\axiombU{#2}}{}
\ifthenelse{\equal{#1}{22}}{
			\axiombV{#2}}{}			
\ifthenelse{\equal{#1}{23}}{
			\axiombW{#2}}{}
\ifthenelse{\equal{#1}{24}}{
			\axiombX{#2}}{}			
\ifthenelse{\equal{#1}{25}}{
			\axiombY{#2}}{}
\ifthenelse{\equal{#1}{26}}{
			\axiombZ{#2}}{}			
\ifthenelse{\equal{#1}{27}}{
			\axiombRefl{#2}}{}			
}

\def\axgrid#1{
\ifthenelse{\equal{#1}{1}}{
			\axgridA}{}
\ifthenelse{\equal{#1}{2}}{
			\axgridB}{}
\ifthenelse{\equal{#1}{3}}{
			\axgridC}{}
\ifthenelse{\equal{#1}{4}}{
			\axgridD}{}
\ifthenelse{\equal{#1}{5}}{
			\axgridE}{}
\ifthenelse{\equal{#1}{6}}{
			\axgridF}{}
\ifthenelse{\equal{#1}{7}}{
			\axgridG}{}
\ifthenelse{\equal{#1}{8}}{
			\axgridH}{}
\ifthenelse{\equal{#1}{9}}{
			\axgridI}{}
\ifthenelse{\equal{#1}{10}}{
			\axgridJ}{}			
\ifthenelse{\equal{#1}{11}}{
			\axgridK}{}
\ifthenelse{\equal{#1}{12}}{
			\axgridL}{}			
\ifthenelse{\equal{#1}{13}}{
			\axgridM}{}
\ifthenelse{\equal{#1}{14}}{
			\axgridN}{}			
\ifthenelse{\equal{#1}{15}}{
			\axgridO}{}
\ifthenelse{\equal{#1}{16}}{
			\axgridP}{}
\ifthenelse{\equal{#1}{17}}{
			\axgridQ}{}
\ifthenelse{\equal{#1}{18}}{
			\axgridR}{}
\ifthenelse{\equal{#1}{19}}{
			\axgridS}{}
\ifthenelse{\equal{#1}{20}}{
			\axgridT}{}			
\ifthenelse{\equal{#1}{21}}{
			\axgridU}{}
\ifthenelse{\equal{#1}{22}}{
			\axgridV}{}			
\ifthenelse{\equal{#1}{23}}{
			\axgridW}{}
\ifthenelse{\equal{#1}{24}}{
			\axgridX}{}			
\ifthenelse{\equal{#1}{25}}{
			\axgridY}{}
\ifthenelse{\equal{#1}{26}}{
			\axgridZ}{}			
\ifthenelse{\equal{#1}{27}}{
			\axgridRefl}{}			
}

\def\table#1-#2{
\ifthenelse{\equal{#1}{1}}{
			\tableAA{#2}}{}
\ifthenelse{\equal{#1}{2}}{
			\tableBB{#2}}{}
\ifthenelse{\equal{#1}{3}}{
			\tableCC{#2}}{}
\ifthenelse{\equal{#1}{4}}{
			\tableDD{#2}}{}
}

\def\tablegrid#1{
\ifthenelse{\equal{#1}{1}}{
			\tablegridAA}{}
\ifthenelse{\equal{#1}{2}}{
			\tablegridBB}{}
\ifthenelse{\equal{#1}{3}}{
			\tablegridCC}{}
\ifthenelse{\equal{#1}{4}}{
			\tablegridDD}{}
}

%% file: ib-code-forcniv.tex
\def\axiombI#1{
\ifnum#1=1
	\sh(0,0);
	\sector;
	\vid;
	\sh(0,-1);
	\strutv;
	\arrowlarcd;
	\sh(4,0);
	\defectshortr;
	\sh(7,0);
	\arrowr;
	\sh(4,-1);
	\defectshortl;
	\sh(7,-1);
	\tick;
	\defect;
	\sh(10,-1);
	\strutv;
\else \ifnum#1=2
	\sh(0,0);
	\sectorcaretr;
	\sh(-.5*\caretlength,0);
	\vid;
	\sh(4,0);
	\defect;
	\sh(7,0);
	\tick;
	\defect;
	\sh(10,0);
	\strutv;
	\sh(0,0);
	\strutv;
	\sh(0,1);
	\arcu;
	\sh(4,1);
	\tickaddup;
	\defectshortr;
	\sh(7,1);
	\arrowr;
\else \ifnum#1=3
	\sh(0,0);
	\tripletopl;
	\triplebotl;
	\vidt;
	\strutv;
	\sh(0,1);
	\tripletopr;
	\sh(10,0);
	\strutv;
\else \ifnum#1=4
	\sh(0,0);
	\arrowl;
	\strutv;
	\sh(3,0);
	\defectshortl;
	\sh(6,0);
	\defect;
	\tick;
	\sh(9,0);
	\strutv;
	\sh(0,1);
	\defect;
	\sh(3,1);
	\tick;
	\defectshortr;
	\sh(6,1);
	\arrowr;
		
\else \ifnum#1=5
	\sh(0,0);
	\tripletopr;
	\triplebotr;
	\vidt;
	\sh(0,-1);
	\triplebotl;
	\strutv;
	\sh(10,-1);
	\strutv;

\else \ifnum#1=6
	\sh(0,0);
	\vidd;
	\doubletoptick;
	\doublebottick;
	\sh(7,0);
	\defect;
	\sh(0,-1);
	\strutv;
	\doublebotl;
	\sh(7,-1);
	\tick;
	\defect;
	\sh(10,-1);
	\strutv;
	
\else \ifnum#1=0
	\begin{scope}[cm={-1,0,0,-1,(0,0)}]
		\axiombI{2}
	\end{scope}

\fi\fi\fi\fi\fi\fi\fi}

\def\axiombJ#1{
\ifnum#1=1
	\sh(0,0);
	\defect;
	\sh(3,0);
	\sector;
	\vid;
	\sh(7,0);
	\arrowrgap;
	\sh(10,-1);
	\strutv;
	\sh(7,-1);
	\tickadddown;
	\defect;
	\sh(3,-1);
	\arcd;
	\sh(0,-1);
	\strutv;
	\arrowlgap;

\else \ifnum#1=2
	\sh(0,0);
	\defect;
	\sh(3,0);
	\tickaddup;
	\arcu;
	\sh(7,0);
	\arrowrgap;
	\sh(10,-1);
	\strutv;
	\sh(7,-1);
	\defect;
	\sh(3,-1);
	\sector;
	\vid;
	\sh(0,-1);
	\strutv;
	\arrowlgap;

\else \ifnum#1=3
	\sh(0,0);
	\defect;
	\sh(3,0);
	\tick;
	\doubletopr;
	\sh(0,-1);
	\strutv;
	\defect;
	\sh(3,-1);
	\vidd;
	\doubletoptick;
	\doublebottick;
	\sh(10,-1);
	\strutv;
\else \ifnum#1=4
	\sh(0,0);
	\tripletopl;
	\triplebotl;
	\vidt;
	\strutv;
	\sh(0,1);
	\tripletopr;
	\sh(10,0);
	\strutv;
\else \ifnum#1=5
	\sh(0,0);
	\arrowl;
	\strutv;
	\sh(3,0);
	\defectshortl;
	\sh(6,0);
	\defect;
	\tick;
	\sh(9,0);
	\strutv;
	\sh(0,1);
	\defect;
	\sh(3,1);
	\tick;
	\defectshortr;
	\sh(6,1);
	\arrowr;
		
\else \ifnum#1=6
	\sh(0,0);
	\tripletopr;
	\triplebotr;
	\vidt;
	\sh(0,-1);
	\triplebotl;
	\strutv;
	\sh(10,-1);
	\strutv;

\else \ifnum#1=7
	\sh(0,0);
	\vidd;
	\doubletoptick;
	\doublebottick;
	\sh(7,0);
	\defect;
	\sh(0,-1);
	\strutv;
	\doublebotl;
	\sh(7,-1);
	\tick;
	\defect;
	\sh(10,-1);
	\strutv;

\fi\fi\fi\fi\fi\fi\fi}

\def\axiombK#1{
\ifnum#1=1
	\sh(0,0);
	\sectortop;
	\sectorbotarc;
	\sh(4,0);
	\arrowrgap;
	\sh(7,0);
	\arrowrgap;
	\sh(10,-1);
	\strutv;
	\sh(7,-1);
	\defect;
	\tick;
	\sh(4,-1);
	\defectshortl;
	\sh(0,-1);
	\arrowlarcd;
	\strutv;
	
\else \ifnum#1=2
	\sh(0,0);
	\sectortop;
	\sh(4,0);
	\defectshortr;
	\sh(0,-.5);
	\sectorbot;
	\sh(4,-.5);
	\defectshortr;
	\sh(7,-.25);
	\arrowr;
	\sh(0,-1.5);
	\strutv;
	\arrowlarcd;
	\sh(4,-1.5);
	\defectshortl;
	\sh(7,-1.5);
	\tick;
	\defect;
	\sh(10,-1.5);
	\strutlong;
	
\else \ifnum#1=3
	\sh(0,0);
	\strutshort;
	\arrowlarcd;
	\sh(4,0);
	\defectshortl;
	\sh(7,0);
	\tick;
	\defect;
	\sh(10,0);
	\strutv;
	\sh(0,1.5);
	\sectorbot;
	\sh(4,1.5);
	\defectshortr;
	\sh(0,2);
	\sectortop;
	\sh(4,2);
	\defectshortr;	
	\sh(7,1.5);
	\arrowr;
	\sh(7,2);
	\arrowr;
	\sh(0,1.5);
	\assoctickl;
	\sh(10,1.5);
	\assoctickr;

\else \ifnum#1=4
	\sh(0,0);
	\strutv;
	\arrowlarcd;
	\sh(4,0);
	\defectshortl;
	\sh(7,0);
	\tick;
	\defect;
	\sh(10,0);
	\strutv;
	\sh(0,1);
	\sectorbot;
	\sh(4,1);
	\defectshortr;
	\sh(0,1.5);
	\sectortop;
	\sh(4,1.5);
	\defectshortr;	
	\sh(7,1);
	\arrowr;
	\sh(7,1.5);
	\arrowr;

\else \ifnum#1=5
	\sh(0,0);
	\sectorbot;
	\strutv;
	\sh(4,0);
	\defectshortl;
	\sh(7,0);
	\tick;
	\defect;
	\sh(10,0);
	\strutv;
	\sh(0,1);
	\defectlong;
	\sh(4,1);
	\tickshort;
	\defectshortr;
	\sh(7,1);
	\arrowr;
	\sh(0,1.5);
	\sectortop;
	\sh(4,1.5);
	\defectshortr;
	\sh(7,1.5);
	\arrowr;
\else \ifnum#1=6
	\sh(0,0);
	\sectorbot;
	\sh(4,0);
	\defectshortl;
	\sh(7,0);
	\tick;
	\defect;
	\sh(0,.5);
	\strutv;
	\arrowllong;
	\sh(4,.5);
	\defectshortl;
	\sh(7,.5);
	\tick;
	\defect;
	\sh(10,.5);
	\strutv;
	\sh(0,1.5);
	\sectortop;
	\sh(4,1.5);
	\defectshortr;
	\sh(7,1.5);
	\arrowr;

\else \ifnum#1=7
	\sh(0,0);
	\sectorbot;
	\sh(4,0);
	\defectshortl;
	\sh(7,0);
	\tick;
	\defect;
	\sh(0,.5);
	\sectortop;
	\strutv;
	\sh(4,.5);
	\defectshortl;
	\sh(7,.5);
	\tick;
	\defect;
	\sh(10,.5);
	\strutv;
	\sh(0,1.5);
	\arcu;
	\sh(4,1.5);
	\tickaddup;
	\defectshortr;
	\sh(7,1.5);
	\arrowr;
	
\else \ifnum#1=8
	\sh(0,0);
	\sectorbot;
	\sh(4,0);
	\defectshortl;
	\sh(7,0);
	\tick;
	\defect;
	\sh(0,.5);
	\sectortop;
	\sh(4,.5);
	\defectshortl;
	\sh(7,.5);
	\tick;
	\defect;
	\sh(0,2);
	\strutshortdown;
	\sh(10,1);
	\strutv;
	\sh(0,2);
	\arcu;
	\sh(4,2);
	\tickaddup;
	\defectshortr;
	\sh(7,2);
	\arrowr;
	\sh(0,0);
	\assoctickl;
	\sh(10,0);
	\assoctickr;	

\else \ifnum#1=9
	\sh(0,0);
	\sectorbot;
	\sectortoparc;
	\strutv;
	\sh(4,0);
	\defectshortl;
	\sh(7,0);
	\tick;
	\defect;
	\sh(10,0);
	\strutv;
	\sh(0,1);
	\arcu;
	\sh(4,1);
	\tickaddup;
	\defectshortr;
	\sh(7,1);
	\arrowr;
	
\else \ifnum#1=0
	\begin{scope}[cm={-1,0,0,-1,(0,0)}]
		\axiombK{9}
	\end{scope}

\fi\fi\fi\fi\fi\fi\fi\fi\fi\fi}

\def\axiombL#1{
\ifnum#1=1
	\sh(0,0);
	\arrowlgap;
	\strutv;
	\sh(3,0);
	\arcd;
	\sh(7,0);
	\tickadddown;
	\defect;
	\sh(10,0);
	\strutv;
	\sh(0,1);
	\arrowlgap;
	\sh(3,1);
	\sectorbot;
	\sectortoparc;
	\sh(7,1);
	\arrowrgap;
	\sh(7.9,1);
	\caretr;
	
\else \ifnum#1=2
	\begin{scope}[cm={-1,0,0,-1,(0,0)}]
		\axiombL{1}
	\end{scope}
\else \ifnum#1=3
	\sh(0,0);
	\arrowlgapbig; 
	\strutlong;
	\sh(3,-.25);
	\sectorbot;
	\sh(3,.25);
	\sectortop;
	\sh(7,-.25);
	\defect;
	\sh(7,.25);
	\defect;
	\sh(10,.25);
	\strutv;
	\sh(0,1.25);
	\defect;
	\sh(3,1.25);
	\tickaddup;
	\arcu;
	\sh(7,1.25);
	\arrowrgap;
	
\else \ifnum#1=4
	\sh(0,0);
	\arrowlgap;
	\assoctickl;
	\sh(3,0);
	\sectorbot;
	\sh(7,0);
	\defect;
	\sh(10,0);
	\assoctickr;
	\sh(0,.5);
	\arrowlgap;
	\sh(3,.5);
	\sectortop;
	\sh(7,.5);
	\defect;
	\sh(0,1);
	\strutv;
	\sh(10,1);
	\strutv;
	\sh(0,2);
	\defect;
	\sh(3,2);
	\tickaddup;
	\arcu;
	\sh(7,2);
	\arrowrgap;
	
\else \ifnum#1=5
	\sh(0,0);
	\arrowlgap;
	\sh(3,0);
	\sectorbot;
	\sh(7,0);
	\defect;
	\sh(0,.5);
	\arrowlgap;
	\sh(3,.5);
	\sectortop;
	\sh(7,.5);
	\defect;
	\sh(0,.5);
	\strutv;
	\sh(10,.5);
	\strutv;
	\sh(0,1.5);
	\defect;
	\sh(3,1.5);
	\tickaddup;
	\arcu;
	\sh(7,1.5);
	\arrowrgap;

\else \ifnum#1=6
	\sh(0,0);
	\arrowlgap;
	\sh(3,0);
	\sectorbot;
	\sh(7,0);
	\defect;
	\sh(0,.5);
	\arrowlgap;
	\strutv;
	\sh(3,.5);
	\defectlong;
	\sh(7,.5);
	\defect;
	\tickshort;
	\sh(10,.5);
	\strutv;
	\sh(0,1.5);
	\defect;
	\sh(3,1.5);
	\sectortop;
	\sh(7,1.5);
	\arrowrgap;
	
\else \ifnum#1=7
	\sh(0,0);
	\arrowlgap;
	\strutv;
	\sh(3,0);
	\sectorbot;
	\sh(7,0);
	\defect;
	\sh(10,0);
	\strutv;
	\sh(0,1);
	\defect;
	\sh(3,1);
	\tickshort;
	\defectlong;
	\sh(7,1);
	\arrowrgap;
	\sh(0,1.5);
	\defect;
	\sh(3,1.5);
	\sectortop;
	\sh(7,1.5);
	\arrowrgap;
	
\else \ifnum#1=8
	\sh(0,0);
	\strutv;
	\arrowlgap;
	\sh(3,0);
	\arcd;
	\sh(7,0);
	\tickadddown;
	\defect;
	\sh(10,0);
	\strutv;
	\sh(0,1);
	\defect;
	\sh(3,1);
	\sectorbot;
	\sh(7,1);
	\arrowrgap;
	\sh(0,1.5);
	\defect;
	\sh(3,1.5);
	\sectortop;
	\sh(7,1.5);
	\arrowrgap;
	
\else \ifnum#1=9
	\sh(0,0);
	\strutv;
	\arrowlgap;
	\sh(3,0);
	\arcd;
	\sh(7,0);
	\tickadddown;
	\defect;
	\sh(10,0);
	\strutv;
	\sh(0,1.5);
	\defect;
	\assoctickl;
	\sh(3,1.5);
	\sectorbot;
	\sh(7,1.5);
	\arrowrgap;
	\sh(10,1.5);
	\assoctickr;
	\sh(0,2);
	\defect;
	\sh(3,2);
	\sectortop;
	\sh(7,2);
	\arrowrgap;
	
\else \ifnum#1=10
	\sh(0,0);
	\strutv;
	\arrowlgap;
	\sh(3,0);
	\arcd;
	\sh(7,0);
	\tickadddown;
	\defect;
	\sh(10,0);
	\strutlong;
	\sh(0,1);
	\defect;
	\sh(3,1);
	\sectorbot;
	\sh(7,1.25);
	\arrowrgapbig;
	\sh(0,1.5);
	\defect;
	\sh(3,1.5);
	\sectortop;

\fi\fi\fi\fi\fi\fi\fi\fi\fi\fi}

\def\axiombM#1{
\ifnum#1=1
	\sh(4,0);
	\arcdgapl;
	\sh(0,0);
	\strutv;
	\arrowlarcdgap;
	\sh(8,0);
	\tickadddown;
	\defect;
	\sh(11,0);
	\strutlong;
	\sh(0,1);
	\arcd;
	\sh(4,1);
	\sector;
	\sh(8,1.25);
	\arrowrgapbig;
	\sh(0,1.5);
	\sector;
	\sh(4,1.5);
	\arcu;
	
\else \ifnum#1=2
	\sh(4,0);
	\arcdgapl;
	\sh(0,0);
	\strutv;
	\arrowlarcdgap;
	\sh(8,0);
	\tickadddown;
	\defect;
	\sh(11,0);
	\strutlong;
	\sh(0,1);
	\sector;
	\sh(4,1);
	\arcd;
	\sh(8,1.25);
	\arrowrgapbig;
	\sh(0,1.5);
	\arcu;
	\sh(4,1.5);
	\sector;

\else \ifnum#1=3
	\sh(4,0);
	\arcdgapl;
	\sh(0,0);
	\strutshort;
	\arrowlarcdgap;
	\sh(8,0);
	\tickadddown;
	\defect;
	\sh(11,0);
	\strutv;
	\sh(0,1.5);
	\sector;
	\sh(4,1.5);
	\arcd;
	\sh(8,1.5);
	\arrowrgap;
	\sh(8,2);
	\arrowrgap;
	\sh(11,1.5);
	\assoctickr;
	\sh(0,1.5);
	\assoctickl;
	\sh(0,2);
	\arcu;
	\sh(4,2);
	\sector;

\else \ifnum#1=4
	\sh(4,0);
	\arcdgapl;
	\sh(0,0);
	\strutv;
	\arrowlarcdgap;
	\sh(8,0);
	\tickadddown;
	\defect;
	\sh(11,0);
	\strutv;
	\sh(0,1);
	\sector;
	\sh(4,1);
	\arcd;
	\sh(8,1);
	\arrowrgap;
	\sh(8,1.5);
	\arrowrgap;
	\sh(0,1.5);
	\arcu;
	\sh(4,1.5);
	\sector;
	
\else \ifnum#1=5
	\sh(0,0);
	\sectorcaretr;
	\strutv;
	\sh(4,0);
	\arcd;
	\sh(8,0);
	\tickadddown;
	\defect;
	\sh(11,0);
	\strutv;
	\sh(0,1);
	\arcu;
	\sh(4,1);
	\arcd;
	\tickshort;
	\sh(8,1);
	\arrowrgap;
	\sh(8,1.5);
	\arrowrgap;
	\sh(0,1.5);
	\arcu;
	\sh(4,1.5);
	\sector;

\else \ifnum#1=6
	\sh(0,0);
	\sectorcaretr;
	\sh(4,0);
	\arcd;
	\sh(8,0);
	\tickadddown;
	\defect;
	\sh(0,.5);
	\arrowlarcuturn;
	\strutv;
	\sh(4,.5);
	\arcd;
	\sh(8,.5);
	\tickadddown;
	\defect;
	\sh(11,.5);
	\strutv;
	\sh(8,1.5);
	\arrowrgap;
	\sh(0,1.5);
	\arcu;
	\sh(4,1.5);
	\sector;

\else \ifnum#1=7
	\sh(0,0);
	\sectorcaretr;
	\sh(4,0);
	\arcdlong;
	\sh(9,0);
	\tickshortupadddown;
	\defect;
	\sh(0,.5);
	\arrowlarcu;
	\strutv;
	\sh(4.5,.5);
	\sector;
	\sh(8.5,.5);
	\defectsemilong;
	\sh(12,.5);
	\strutv;
	\sh(0,1.5);
	\arcu;
	\sh(4,1.5);
	\tickaddup;
	\arculong;
	\sh(9,1.5);
	\arrowrgap;

\else \ifnum#1=8
	\sh(0,0);
	\sectorcaretr;
	\assoctickl;
	\sh(4,0);
	\arcdlong;
	\sh(9,0);
	\tickshortupadddown;
	\defect;
	\sh(12,0);
	\assoctickr;
	\sh(0,.5);
	\arrowlarcu;
	\sh(4.5,.5);
	\sector;
	\sh(8.5,.5);
	\defectsemilong;
	\sh(12,.5);
	\sh(0,2);
	\strutshortdown;
	\sh(12,1);
	\strutv;
	\sh(0,2);
	\arcu;
	\sh(4,2);
	\tickaddup;
	\arculong;
	\sh(9,2);
	\arrowrgap;

\else \ifnum#1=9
	\sh(0,.5);
	\sectorcaretr;
	\sh(4,.5);
	\arculong;
	\sh(9,.5);
	\tickshortdownaddup;
	\defect;
	\sh(12,0);
	\assoctickr;
	\sh(0,0);
	\arrowlarcd;
	\assoctickl;
	\sh(4.5,0);
	\sector;
	\sh(8.5,0);
	\defectsemilong;
	\sh(12,.5);
	\sh(0,2);
	\strutshortdown;
	\sh(12,1);
	\strutv;
	\sh(0,2);
	\arcu;
	\sh(4,2);
	\tickaddup;
	\arculong;
	\sh(9,2);
	\arrowrgap;
	
\else \ifnum#1=10
	\sh(0,.5);
	\sectorcaretr;
	\sh(4,.5);
	\arculong;
	\sh(9,.5);
	\tickshortdownaddup;
	\defect;
	\sh(0,0);
	\arrowlarcd;
	\sh(4.5,0);
	\sector;
	\sh(8.5,0);
	\defectsemilong;
	\sh(12,.5);
	\sh(0,.5);
	\strutv;
	\sh(12,.5);
	\strutv;
	\sh(0,1.5);
	\arcu;
	\sh(4,1.5);
	\tickaddup;
	\arculong;
	\sh(9,1.5);
	\arrowrgap;
	
\else \ifnum#1=11
	\sh(0,.5);
	\strutv;
	\arrowlarcdturn;
	\sh(4,.5);
	\arcu;
	\sh(8,.5);
	\tickaddup;
	\defect;
	\sh(11,.5);
	\strutv;
	\sh(0,1.5);
	\sector;
	\sh(4,1.5);
	\arcu;
	\sh(8,1.5);
	\arrowrgap;
	\sh(0,-.25);
	\arcdgapr;
	\sh(4,-.25);
	\begin{scope}[rotate=-135,shift={(-6.17,3.38)}] %!% insane hack
		\caretrmini;
	\end{scope}
	\sh(4,-.25);
	\sector;
	\sh(8,-.25);
	\defect;

\else \ifnum#1=12
	\sh(0,.65);
	\arcd;
	\sh(4,.65);
	\arcu;
	\tickshort;
	\sh(8,.65);
	\arrowrgap;
	\sh(0,1.15);
	\sector;
	\sh(4,1.15);
	\arcu;
	\sh(8,1.15);
	\arrowrgap;
	\sh(0,-.35);
	\arcdgapr;
	\strutv;
	\sh(4,-.35);
	\begin{scope}[rotate=-135,shift={(-6.10,3.55)}] %!% insane hack
		\caretrmini;
	\end{scope}
	\sh(4,-.35);
	\sector;
	\sh(8,-.35);
	\defect;
	\sh(11,-.35);
	\strutv;
	
\else \ifnum#1=13
	\sh(4,0);
	\arcdgapl;
	\sh(0,0);
	\strutv;
	\arrowlarcdgap;
	\sh(8,0);
	\tickadddown;
	\defect;
	\sh(11,0);
	\strutv;
	\sh(0,1);
	\arcd;
	\sh(4,1);
	\sector;
	\sh(8,1);
	\arrowrgap;
	\sh(8,1.5);
	\arrowrgap;
	\sh(0,1.5);
	\sector;
	\sh(4,1.5);
	\arcu;

\else \ifnum#1=14
	\sh(4,0);
	\arcdgapl;
	\sh(0,0);
	\strutshort;
	\arrowlarcdgap;
	\sh(8,0);
	\tickadddown;
	\defect;
	\sh(11,0);
	\strutv;
	\sh(0,1.5);
	\arcd;
	\sh(4,1.5);
	\sector;
	\sh(8,1.5);
	\arrowrgap;
	\sh(8,2);
	\arrowrgap;
	\sh(11,1.5);
	\assoctickr;
	\sh(0,1.5);
	\assoctickl;
	\sh(0,2);
	\sector;
	\sh(4,2);
	\arcu;
	
\else \ifnum#1=0
	\sh(0,0);
	\defect;
	\sh(3,0);
	\sector;
	\sh(7,0);
	\sectorcaretl;
	\sh(0,-1);
	\strutv;
	\arrowlgap;
	\sh(3,-1);
	\arcd;
	\sh(7,-1);
	\tickshortupadddown;
	\arcd;
	\sh(11,-1);
	\strutv;

\fi\fi\fi\fi\fi\fi\fi\fi\fi\fi\fi\fi\fi\fi\fi}

\def\axiombN#1{
\ifnum#1=1
	\sh(0,0);
	\strutshort;
	\arrowlarcd;
	\sh(4,0);
	\defectshortl;
	\sh(7,0);
	\tickadddown;
	\arcd;
	\sh(11,0);
	\strutshort;
	\sh(0,1.5);
	\assoctickl;
	\arcd;
	\sh(4,1.5);
	\tickshortupadddown;
	\defect;
	\sh(7,1.5);
	\sectorcaretl;
	\sh(11,1.5);
	\assoctickr;
	\sh(0,2);
	\sector;
	\sh(4,2);
	\defectshortr;
	\sh(7,2);
	\arcu;
	\caretrmini;

\else \ifnum#1=2
	\sh(0,0);
	\strutshort;
	\arrowlarcd;
	\sh(4,0);
	\defectshortl;
	\sh(7,0);
	\tickadddown;
	\arcd;
	\sh(11,0);
	\strutshort;
	\sh(0,1.5);
	\assoctickl;
	\sector;
	\sh(4,1.5);
	\defectshortr;
	\sh(7,1.5);
	\arcd;
	\caretrmini;
	\sh(11,1.5);
	\assoctickr;
	\sh(0,2);
	\arcu;
	\sh(4,2);
	\tickshortdownaddup;
	\defect;
	\sh(7,2);
	\sectorcaretl;

\else \ifnum#1=3
	\sh(0,0);
	\strutv;
	\arrowlarcd;
	\sh(4,0);
	\defectshortl;
	\sh(7,0);
	\tickadddown;
	\arcd;
	\sh(11,0);
	\strutv;
	\sh(0,1);
	\sector;
	\sh(4,1);
	\defectshortr;
	\sh(7,1);
	\arcd;
	\caretrmini;
	\sh(0,1.5);
	\arcu;
	\sh(4,1.5);
	\tickshortdownaddup;
	\defect;
	\sh(7,1.5);
	\sectorcaretl;

\else \ifnum#1=4
	\sh(0,0);
	\sectorcaretr;
	\strutv;
	\sh(4,0);
	\defect;
	\sh(7,0);
	\tickadddown;
	\arcd;
	\sh(11,0);
	\strutv;
	\sh(0,1);
	\arcu;
	\sh(4,1);
	\tickaddup;
	\defectshortr;
	\sh(7,1);
	\arcd;
	\caretrmini;
	\sh(0,1.5);
	\arcu;
	\sh(4,1.5);
	\tickaddup;
	\defect;
	\sh(7,1.5);
	\sectorcaretl;
	
\else \ifnum#1=5
	\begin{scope}[cm={-1,0,0,-1,(0,0)}]
		\axiombN{4}
	\end{scope}
\else \ifnum#1=6
	\begin{scope}[cm={-1,0,0,-1,(0,0)}]
		\axiombN{3}
	\end{scope}
\else \ifnum#1=7
	\begin{scope}[cm={-1,0,0,-1,(0,0)}]
		\axiombN{2}
	\end{scope}
\else \ifnum#1=8
	\begin{scope}[cm={-1,0,0,-1,(0,0)}]
		\axiombN{1}
	\end{scope}
\else \ifnum#1=9
	\begin{scope}[cm={-1,0,0,-1,(0,0)}]
		\axiombN{12}
	\end{scope}
\else \ifnum#1=10
	\begin{scope}[cm={-1,0,0,-1,(0,0)}]
		\axiombN{11}
	\end{scope}
\else \ifnum#1=11
	\sh(0,0);
	\strutv;
	\arrowlarcd;
	\sh(4,0);
	\defectshortl;
	\sh(7,0);
	\sector
	\sh(11,0);
	\strutv;
	\sh(0,1);
	\arcd;
	\sh(4,1);
	\tickshortupadddown;
	\defectshortr;
	\sh(7,1);
	\arrowrarcu;
	\sh(0,1.5);
	\sector;
	\sh(4,1.5);
	\defectshortr;
	\sh(7,1.5);
	\arrowrarcu;

\else \ifnum#1=12
	\sh(0,0);
	\strutv;
	\arrowlarcd;
	\sh(4,0);
	\defectshortl;
	\sh(7,0);
	\tickadddown;
	\arcd;
	\sh(11,0);
	\strutv;
	\sh(0,1);
	\arcd;
	\sh(4,1);
	\tickshortupadddown;
	\defect;
	\sh(7,1);
	\sectorcaretl;
	\sh(0,1.5);
	\sector;
	\sh(4,1.5);
	\defectshortr;
	\sh(7,1.5);
	\arcu;
	\caretrmini;

\fi\fi\fi\fi\fi\fi\fi\fi\fi\fi\fi\fi}

\def\axiombT#1{
\ifnum#1=1
	\sh(0,0);
	\hirbotcut;
	\sectortoparc;
	\sh(7,0);
	\arrowrgap;
\else \ifnum#1=2
	\sh(0,0);
	\hirbotcut;
	\sh(7,0);
	\arrowrgap;
	\sh(0,.5);
	\sector;
	\sh(4,.5);
	\defect;
	\iticks;
	\sh(7,.5);
	\arrowrgap;
\else \ifnum#1=3
	\sh(0,0);
	\arcd;
	\strutv;
	\sh(4,0);
	\hil;
	\sh(8,0);
	\strutv;
	\sh(0,1);
	\arcd;
	\sh(4,1);
	\tickminicc;
	\splittopr;
	\iticksrotlminifact;
	\sh(0,1.5);
	\sector;
	\sh(4,1.5);
	\splittopr;
	\iticksrotlminifact;
		
\else \ifnum#1=4
	\sh(0,.5);
	\strutv;
	\arrowlarcdturn;
	\sh(4,.5);
	\arcu;
	\iticksrotlminifact;	
	\sh(6,.5+\arcfactor*.5*\sectorv);
	\tick;
	\sh(8,.5);
	\strutv;
	\sh(0,1.5);
	\sector;
	\sh(4,1.5);
	\splittopr;
	\iticksrotlminifact;
	\sh(0,0);
	\arcd;
	\sh(4,0);
	\hil;
	
\else \ifnum#1=5
	\sh(0,0);
	\arcd;
	\sh(4,0);
	\hil;
	\sh(0,.5);
	\sectorcaretr;
	\strutv;
	\sh(4,.5);
	\arcu;
	\iticksrotlminifact;
	\sh(6,.5+\arcfactor*.5*\sectorv);
	\tick;
	\sh(8,.5);
	\strutv;
	\sh(0,1.5);
	\arcu;
	\sh(4,1.5);
	\tickshortdownaddup;
	\splittopr;
	\iticksrotlminifact;
	
\else \ifnum#1=6
	\sh(0,0);
	\arcd;
	\assoctickl;
	\sh(4,0);
	\hil;
	\sh(8,0);
	\assoctickr;
	\sh(0,.5);
	\sectorcaretr;
	\sh(4,.5);
	\arcu;
	\iticksrotlminifact;
	\sh(6,.5+\arcfactor*.5*\sectorv);
	\tick;
	\sh(0,2);
	\strutshortdown;
	\sh(8,2);
	\strutshortdown;
	\sh(0,2);
	\arcu;
	\sh(4,2);
	\tickshortdownaddup;
	\splittopr;
	\iticksrotlminifact;	

\else \ifnum#1=7
	\sh(0,0);
	\sector;
	\assoctickl;
	\sh(4,0);
	\arcd;
	\sh(8,0);
	\assoctickr;
	\sh(0,.5);
	\arcu;
	\sh(4,.5);
	\hil;
	\sh(0,2);
	\strutshortdown;
	\sh(8,2);
	\strutshortdown;
	\sh(0,2);
	\arcu;
	\sh(4,2);
	\tickshortdownaddup;
	\splittopr;
	\iticksrotlminifact;	

\else \ifnum#1=8
	\sh(0,0);
	\sector;
	\sh(4,0);
	\arcd;
	\sh(0,.5);
	\arcu;
	\sh(4,.5);
	\hil;
	\sh(0,.5);
	\strutv;
	\sh(8,.5);
	\strutv;
	\sh(0,1.5);
	\arcu;
	\sh(4,1.5);
	\tickshortdownaddup;
	\splittopr;
	\iticksrotlminifact;	

\else \ifnum#1=9
	\sh(0,0);
	\sectorbot;
	\sh(4,0);
	\defect;
	\sh(0,.5);
	\hirtop;
	\sh(4,.5);
	\defect;
	
\else \ifnum#1=10
	\sh(0,0);
	\hirtop;
	\sectorbotarc;
	\sh(4,0);
	\arrowrgap;
	
\else \ifnum#1=0
	\sh(0,0);
	\hir;
	\sh(4,0);
	\sector;
	
\fi\fi\fi\fi\fi\fi\fi\fi\fi\fi\fi}

\def\axiombU#1{
\ifnum#1=1
	\sh(0,0);
	\arcd;
	\assoctickl;
	\sh(4,0);
	\hir;
	\sh(8,0);
	\assoctickr;
	\sh(0,.5);
	\sectorcaretr;
	\sh(4,.5);
	\arcu;
	\iticksrotrminifact;
	\sh(6,.5+\arcfactor*.5*\sectorv);
	\tick;
	\sh(0,2);
	\strutshortdown;
	\sh(8,2);
	\strutshortdown;
	\sh(0,2);
	\arcu;
	\sh(4,2);
	\tickshortdownaddup;
	\splittopr;
	\iticksrotrminifact;

\else \ifnum#1=2
	\sh(0,0);
	\sector;
	\assoctickl;
	\sh(4,0);
	\arcd;
	\sh(8,0);
	\assoctickr;
	\sh(0,.5);
	\arcu;
	\sh(4,.5);
	\hir;
	\sh(0,2);
	\strutshortdown;
	\sh(8,2);
	\strutshortdown;
	\sh(0,2);
	\arcu;
	\sh(4,2);
	\tickshortdownaddup;
	\splittopr;
	\iticksrotrminifact;	

\else \ifnum#1=3
	\sh(0,0);
	\sector;
	\sh(4,0);
	\arcd;
	\sh(0,.5);
	\arcu;
	\sh(4,.5);
	\hir;
	\sh(0,.5);
	\strutv;
	\sh(8,.5);
	\strutv;
	\sh(0,1.5);
	\arcu;
	\sh(4,1.5);
	\tickshortdownaddup;
	\splittopr;
	\iticksrotrminifact;	

\else \ifnum#1=4
	\sh(0,0);
	\hirtopthreecut;
	\sectorbotarc;
	
\else \ifnum#1=5
	\sh(0,0);
	\hirbotthreecut;
	\sectortoparc;
	
\else \ifnum#1=6
	\sh(0,0);
	\arcd;
	\strutv;
	\sh(4,0);
	\hir;
	\sh(8,0);
	\strutv;
	\sh(0,1);
	\arcd;
	\sh(4,1);
	\tickminicc;
	\splittopr;
	\iticksrotrminifact;
	\sh(0,1.5);
	\sector;
	\sh(4,1.5);
	\splittopr;
	\iticksrotrminifact;

\else \ifnum#1=7
	\sh(0,.5);
	\arrowlarcdturn;
	\strutv;
	\sh(0,0);
	\arcd;
	\sh(4,0);
	\hir;
	\sh(4,.5);
	\arcu;
	\iticksrotrminifact;
	\sh(6,.5+\arcfactor*.5*\sectorv);
	\tick;
	\sh(8,.5);
	\strutv;
	\sh(0,1.5);
	\sector;
	\sh(4,1.5);
	\splittopr;
	\iticksrotrminifact;	

\else \ifnum#1=8
	\sh(0,0);
	\arcd;
	\sh(4,0);
	\hir;
	\sh(0,.5);
	\sectorcaretr;
	\sh(4,.5);
	\arcu;
	\iticksrotrminifact;
	\sh(6,.5+\arcfactor*.5*\sectorv);
	\tick;
	\sh(0,.5);
	\strutv;
	\sh(8,.5);
	\strutv;
	\sh(0,1.5);
	\arcu;
	\sh(4,1.5);
	\tickshortdownaddup;
	\splittopr;
	\iticksrotrminifact;
	
\else \ifnum#1=0
	\sh(0,0);
	\hil;
	\sh(4,0);
	\sector;		

\fi\fi\fi\fi\fi\fi\fi\fi\fi}

\def\axiombV#1{
\ifnum#1=1
	\sh(0,0);
	\hilbotcut;
	\sh(3,0);
	\sectortoparc;
	\sh(7,0);
	\arrowrgap;
	
\else \ifnum#1=2
	\sh(0,0);
	\hilbotcut;
	\sh(7,0);
	\arrowrgap;
	\sh(0,.5);
	\defect;
	\iticks;
	\sh(3,.5);
	\sector;
	\sh(7,.5);
	\arrowrgap;
	
\else \ifnum#1=3
	\sh(0,0);
	\hilbotthreecutmid;
	\strutv;
	\sh(10,0);
	\strutv;
	\sh(0,1);
	\defect;
	\iticks;
	\sh(3,1);
	\tickshortupadddown;
	\arcd;
	\sh(7,1);
	\arrowrgap;
	\sh(0,1.5);
	\defect;
	\iticks;
	\sh(3,1.5);
	\sector;
	\sh(7,1.5);
	\arrowrgap;
	
\else \ifnum#1=4
	\sh(0,0);
	\hilbotthreecutmid;
	\sh(0,.5);
	\strutv;
	\arrowlgap;
	\iticks;
	\sh(3,.5);
	\arcd;
	\sh(7,.5);
	\defect;
	\tickadddown;
	\sh(10,.5);
	\strutv;
	\sh(0,1.5);
	\defect;
	\iticks;
	\sh(3,1.5);
	\sector;
	\sh(7,1.5);
	\arrowrgap;

\else \ifnum#1=5
	\sh(0,0);
	\hilbotthreecutmidshorttick;
	\sh(0,.5);
	\strutv;
	\arrowlgap;
	\iticks;
	\sh(3,.5);
	\sector;
	\sh(7,.5);
	\defect;
	\sh(10,.5);
	\strutv;
	\sh(0,1.5);
	\defect;
	\iticks;
	\sh(3,1.5);
	\arcu;
	\tickaddup;
	\sh(7,1.5);
	\arrowrgap;

\else \ifnum#1=6
	\sh(0,0);
	\hilbotthreecutmidnotick;
	\sh(3,0);
	\sectortoparc;
	\sh(0,.5);
	\strutv;
	\sh(10,.5);
	\strutv;
	\sh(0,1.5);
	\defect;
	\iticks;
	\sh(3,1.5);
	\tickaddup;
	\arcu;
	\sh(7,1.5);
	\arrowrgap;
	
\else \ifnum#1=7
	\sh(0,0);
	\hiltopthree;
	\sectorbotarclong;
	\sh(0,.5);
	\strutv;
	\sh(10,.5);
	\strutv;
	\sh(0,1.5);
	\tripletopr;
	\itickstriprotlshallow;
	
\else \ifnum#1=8
	\sh(0,0);
	\sectorbotarclong;
	\defectsuplong;
	\sh(5,0);
	\defectsuplong;
	\sh(0,.5);
	\hiltopthree;
	\strutv;
	\sh(10,.5);
	\strutv;
	\sh(0,1.5);
	\tripletopr;
	\itickstriprotlshallow;

\else \ifnum#1=9
	\sh(0,0);
	\sectorbot;
	\sh(4,0);
	\defect;
	\sh(0,.5);
	\hiltop;
	\sh(4,.5);
	\defect;
	
\else \ifnum#1=10
	\sh(0,0);
	\hiltop;
	\sectorbotarc;
	\sh(4,0);
	\arrowrgap;
	
\else \ifnum#1=0
	\sh(0,0);
	\hirbotthreecutmidnotick;
	\sh(3,0);
	\sectortoparc;	

\fi\fi\fi\fi\fi\fi\fi\fi\fi\fi\fi}

\def\axiombW#1{
\ifnum#1=1
	\sh(0,0);
	\hilbotcuttight;
	\sh(3,0);
	\hil;
	\sh(0,.5);
	\strutv;
	\sh(7,1.5);
	\strutshortdown;
	\sh(0,1.5);
	\defect;
	\sh(3,1.5);
	\tickaddup;
	\splittopr;
	\sh(0,1.5);
	\iticks;
	\sh(3,1.5);
	\iticksrotlfact;
	
\else \ifnum#1=2
	\sh(0,0);
	\hilbottight;
	\hiltopthree;
	\sh(5,0);
	\tick;
	\sh(5,1.5);
	\itickshere;
	\sh(0,.5);
	\strutv;
	\sh(10,.5);
	\strutv;
	\sh(0,1.5);
	\tripletopr;
	\itickstriprotlshallow;
	\sh(5,3);
	\itickshere;
	
\else \ifnum#1=3
	\sh(0,0);
	\hilbot;
	\sh(0,.5);
	\hiltopthree;
	\sh(5,.5);
	\tick;
	\sh(5,2);
	\itickshere;
	\sh(0,.5);
	\strutv;
	\sh(10,.5);
	\strutv;
	\sh(0,1.5);
	\tripletopr;
	\itickstriprotlshallow;
	\sh(5,3);
	\itickshere;

\else \ifnum#1=4
	\sh(0,0);
	\hilbotcutleft;
	\hil;
	\iticksrotr;
	
\else \ifnum#1=5
	\sh(0,0);
	\hilbotcutleft;
	\hir;
	\iticksrotl;

\else \ifnum#1=6
	\sh(0,0);
	\hilbotcut;
	\sh(0,.5);
	\strutv;
	\defect;
	\iticks;
	\sh(3,.5);
	\hil;
	\sh(7,.5);
	\strutv;
	\sh(0,1.5);
	\defect;
	\iticks;
	\sh(3,1.5);
	\tickaddup;
	\splittopr;
	\iticksrotlfact;
	
\else \ifnum#1=0
	\sh(0,0);
	\hirbotcuttight;
	\hir;
	
\fi\fi\fi\fi\fi\fi\fi}

\def\axiombX#1{
\ifnum#1=1
	\sh(0,0);
	\hirbotcuttight;
	\hil;

\else \ifnum#1=2
	\sh(0,0);
	\hirbot;
	\sh(0,.5);
	\hiltopthree;
	\sh(5,.5);
	\tick;
	\sh(10,.5);
	\itickstriprotrshallow;
	\sh(0,.5);
	\strutv;
	\sh(10,.5);
	\strutv;
	\sh(0,1.5);
	\tripletopr;
	\itickstriprotlshallow;
	\sh(10,1.5);
	\itickstriprotrshallow;

\else \ifnum#1=3
	\sh(0,0);
	\hirbottight;
	\hiltopthree;
	\sh(5,0);
	\tick;
	\sh(0,.5);
	\strutv;
	\sh(10,.5);
	\strutv;
	\sh(0,1.5);
	\tripletopr;
	\itickstriprotlshallow;
	\sh(10,0);
	\itickstriprotrshallow;
	\sh(10,1.5);
	\itickstriprotrshallow;

\else \ifnum#1=4
	\sh(0,0);
	\hilbotcuttight;
	\sh(3,0);
	\hir;
	\sh(0,.5);
	\strutv;
	\sh(7,1.5);
	\strutshortdown;
	\sh(0,1.5);
	\defect;
	\sh(3,1.5);
	\tickaddup;
	\splittopr;
	\sh(0,1.5);
	\iticks;
	\sh(3,1.5);
	\iticksrotrfact;

\else \ifnum#1=5
	\sh(0,0);
	\hilbotcut;
	\sh(0,.5);
	\strutv;
	\defect;
	\iticks;
	\sh(3,.5);
	\hir;
	\sh(7,.5);
	\strutv;
	\sh(0,1.5);
	\defect;
	\iticks;
	\sh(3,1.5);
	\tickaddup;
	\splittopr;
	\iticksrotrfact;

\else \ifnum#1=6
	\sh(0,0);
	\hilbottight;
	\hirtopthree;
	\sh(5,0);
	\tick;
	\sh(0,0);
	\itickstriprotlshallow;

\fi\fi\fi\fi\fi\fi}

\def\axiombA#1{
\ifnum#1=1
\pgftext{\ingeps{3stuff-item1-1.eps}}
\else \ifnum#1=2
\pgftext{\ingeps{3stuff-item1-2.eps}} 
\else \ifnum#1=3 % This is the null node, though in fact it is not called
\fi\fi\fi}

\def\axiombB#1{
\ifnum#1=1
\pgftext{\ingeps{3stuff-item2-1.eps}}
\else \ifnum#1=2
\pgftext{\ingeps{3stuff-item2-2.eps}} 
\else \ifnum#1=3
\pgftext{\ingeps{3stuff-item2-3.eps}}
\else \ifnum#1=4
\else \ifnum#1=0
\pgftext{\ingeps{3stuff-item2a.eps}}
\fi\fi\fi\fi\fi}

\def\axiombC#1{
\ifnum#1=1
\pgftext{\ingeps{3stuff-item3-1.eps}}
\else \ifnum#1=2
\pgftext{\ingeps{3stuff-item3-2.eps}} 
\else \ifnum#1=3
\pgftext{\ingeps{3stuff-item3-3.eps}}
\else \ifnum#1=4
\fi\fi\fi\fi}

\def\axiombD#1{
\ifnum#1=1
\pgftext{\ingeps{3stuff-item4-1.eps}}
\else \ifnum#1=2
\pgftext{\ingeps{3stuff-item4-2.eps}} 
\else \ifnum#1=3
\pgftext{\ingeps{3stuff-item4-3.eps}}
\else \ifnum#1=4
\pgftext{\ingeps{3stuff-item4-4.eps}} 
\else \ifnum#1=5
\pgftext{\ingeps{3stuff-item4-5.eps}}
\else \ifnum#1=6
\fi\fi\fi\fi\fi\fi}

\def\axiombE#1{
\ifnum#1=1
\pgftext{\ingeps{3stuff-item5-1.eps}}
\else \ifnum#1=2
\pgftext{\ingeps{3stuff-item5-2.eps}} 
\else \ifnum#1=3
\pgftext{\ingeps{3stuff-item5-3.eps}}
\else \ifnum#1=4
\pgftext{\ingeps{3stuff-item5-4.eps}} 
\else \ifnum#1=5
\else \ifnum#1=0
\pgftext{\ingeps{3stuff-item5a.eps}}
\else \ifnum#1=-1
\pgftext{\ingeps{3stuff-item5b.eps}}
\else \ifnum#1=-2
\pgftext{\ingeps{3stuff-item5c.eps}}

\fi\fi\fi\fi\fi\fi\fi\fi}

\def\axiombF#1{
\ifnum#1=1
\pgftext{\ingeps{3stuff-item6-1.eps}}
\else \ifnum#1=2
\pgftext{\ingeps{3stuff-item6-2.eps}} 
\else \ifnum#1=3
\pgftext{\ingeps{3stuff-item6-3.eps}}
\else \ifnum#1=4
\pgftext{\ingeps{3stuff-item6-4.eps}} 
\else \ifnum#1=5
\pgftext{\ingeps{3stuff-item6-5.eps}}
\else \ifnum#1=6
\else \ifnum#1=0
\pgftext{\ingeps{3stuff-item6a.eps}}
\fi\fi\fi\fi\fi\fi\fi}

\def\axiombG#1{
\ifnum#1=1
\pgftext{\ingeps{3stuff-item7-1.eps}}
\else \ifnum#1=2
\pgftext{\ingeps{3stuff-item7-2.eps}} 
\else \ifnum#1=3
\pgftext{\ingeps{3stuff-item7-3.eps}}
\else \ifnum#1=4
\pgftext{\ingeps{3stuff-item7-4.eps}} 
\else \ifnum#1=5
\pgftext{\ingeps{3stuff-item7-5.eps}}
\else \ifnum#1=6
\pgftext{\ingeps{3stuff-item7-6.eps}}
\else \ifnum#1=7
\else \ifnum#1=0
\pgftext{\ingeps{3stuff-item7a.eps}}
\fi\fi\fi\fi\fi\fi\fi\fi}

\def\axiombH#1{
\ifnum#1=1
\pgftext{\ingeps{3stuff-item8-1.eps}}
\else \ifnum#1=2
\pgftext{\ingeps{3stuff-item8-2.eps}} 
\else \ifnum#1=3
\pgftext{\ingeps{3stuff-item8-3.eps}}
\else \ifnum#1=4
\pgftext{\ingeps{3stuff-item8-4.eps}} 
\else \ifnum#1=5
\pgftext{\ingeps{3stuff-item8-5.eps}}
\else \ifnum#1=6
\pgftext{\ingeps{3stuff-item8-6.eps}} 
\else \ifnum#1=7
\pgftext{\ingeps{3stuff-item8-7.eps}}
\else \ifnum#1=8
\pgftext{\ingeps{3stuff-item8-8.eps}}
\else \ifnum#1=9
\else \ifnum#1=0
\pgftext{\ingeps{3stuff-item8a.eps}}
\fi\fi\fi\fi\fi\fi\fi\fi\fi\fi}

\def\axiombO#1{
\ifnum#1=1
\pgftext{\ingeps{b3stuff-item15-1A.eps}}
\else \ifnum#1=2
\pgftext{\ingeps{b3stuff-item15-1B.eps}} 
\else \ifnum#1=3
\pgftext{\ingeps{b3stuff-item15-2A.eps}}
\else \ifnum#1=4
\pgftext{\ingeps{b3stuff-item15-2B.eps}} 
\else \ifnum#1=5
\pgftext{\ingeps{b3stuff-item15-3A.eps}}
\else \ifnum#1=6
\pgftext{\ingeps{b3stuff-item15-4A.eps}} 
\else \ifnum#1=7
\pgftext{\ingeps{b3stuff-item15-5A.eps}}
\else \ifnum#1=8
\pgftext{\ingeps{b3stuff-item15-5B.eps}}
\else \ifnum#1=9
\pgftext{\ingeps{b3stuff-item15-6A.eps}} 
\else \ifnum#1=10
\pgftext{\ingeps{b3stuff-item15-6B.eps}}
\else \ifnum#1=11
\pgftext{\ingeps{b3stuff-item15-7A.eps}} 
\else \ifnum#1=12
\pgftext{\ingeps{b3stuff-item15-7B.eps}}
\else \ifnum#1=13
\pgftext{\ingeps{b3stuff-item15-8A.eps}} 
\else \ifnum#1=14
\pgftext{\ingeps{b3stuff-item15-8B.eps}}
\else \ifnum#1=15
\pgftext{\ingeps{b3stuff-item15-9A.eps}}
\else \ifnum#1=16
\pgftext{\ingeps{b3stuff-item15-9B.eps}}
\else \ifnum#1=0
\pgftext{\ingeps{b3stuff-item15a.eps}}
\fi\fi\fi\fi\fi\fi\fi\fi\fi\fi\fi\fi\fi\fi\fi\fi\fi}

\def\axiombP#1{
\ifnum#1=1
\pgftext{\ingeps{b3stuff-item16-1A.eps}}
\else \ifnum#1=2
\pgftext{\ingeps{b3stuff-item16-1B.eps}} 
\else \ifnum#1=3
\pgftext{\ingeps{b3stuff-item16-2A.eps}}
\else \ifnum#1=4
\pgftext{\ingeps{b3stuff-item16-2B.eps}} 
\else \ifnum#1=5
\pgftext{\ingeps{b3stuff-item16-3A.eps}}
\else \ifnum#1=6
\pgftext{\ingeps{b3stuff-item16-4A.eps}} 
\else \ifnum#1=7
\pgftext{\ingeps{b3stuff-item16-5A.eps}}
\else \ifnum#1=8
\pgftext{\ingeps{b3stuff-item16-5B.eps}}
\else \ifnum#1=9
\pgftext{\ingeps{b3stuff-item16-6A.eps}} 
\else \ifnum#1=10
\pgftext{\ingeps{b3stuff-item16-6B.eps}}
\else \ifnum#1=11
\pgftext{\ingeps{b3stuff-item16-7A.eps}} 
\else \ifnum#1=12
\pgftext{\ingeps{b3stuff-item16-8A.eps}}
\else \ifnum#1=13
\pgftext{\ingeps{b3stuff-item16-9A.eps}} 
\else \ifnum#1=14
\pgftext{\ingeps{b3stuff-item16-9B.eps}}
\else \ifnum#1=15
\pgftext{\ingeps{b3stuff-item16-10A.eps}}
\else \ifnum#1=16
\pgftext{\ingeps{b3stuff-item16-10B.eps}}
\else \ifnum#1=17
\pgftext{\ingeps{b3stuff-item16-11A.eps}}
\else \ifnum#1=18
\pgftext{\ingeps{b3stuff-item16-11B.eps}}
\else \ifnum#1=0
\pgftext{\ingeps{b3stuff-item16a.eps}}
\fi\fi\fi\fi\fi\fi\fi\fi\fi\fi\fi\fi\fi\fi\fi\fi\fi\fi\fi}

\def\axiombQ#1{
\ifnum#1=1
\pgftext{\ingeps{b3stuff-item17-1A.eps}}
\else \ifnum#1=2
\pgftext{\ingeps{b3stuff-item17-13B.eps}} 
\else \ifnum#1=3
\pgftext{\ingeps{b3stuff-item17-13A.eps}}
\else \ifnum#1=4
\pgftext{\ingeps{b3stuff-item17-12A.eps}} 
\else \ifnum#1=5
\pgftext{\ingeps{b3stuff-item17-11A.eps}}
\else \ifnum#1=6
\pgftext{\ingeps{b3stuff-item17-10B.eps}} 
\else \ifnum#1=7
\pgftext{\ingeps{b3stuff-item17-10A.eps}}
\else \ifnum#1=8
\pgftext{\ingeps{b3stuff-item17-9B.eps}}
\else \ifnum#1=9
\pgftext{\ingeps{b3stuff-item17-9A.eps}} 
\else \ifnum#1=10
\pgftext{\ingeps{b3stuff-item17-8B.eps}}
\else \ifnum#1=11
\pgftext{\ingeps{b3stuff-item17-8A.eps}} 
\else \ifnum#1=12
\pgftext{\ingeps{b3stuff-item17-7B.eps}}
\else \ifnum#1=13
\pgftext{\ingeps{b3stuff-item17-7A.eps}} 
\else \ifnum#1=14
\pgftext{\ingeps{b3stuff-item17-6B.eps}}
\else \ifnum#1=15
\pgftext{\ingeps{b3stuff-item17-6A.eps}}
\else \ifnum#1=16
\pgftext{\ingeps{b3stuff-item17-5B.eps}}
\else \ifnum#1=17
\pgftext{\ingeps{b3stuff-item17-5A.eps}}
\else \ifnum#1=18
\pgftext{\ingeps{b3stuff-item17-4B.eps}}
\else \ifnum#1=19
\pgftext{\ingeps{b3stuff-item17-4A.eps}}
\else \ifnum#1=20
\pgftext{\ingeps{b3stuff-item17-3A.eps}}
\else \ifnum#1=21
\pgftext{\ingeps{b3stuff-item17-2A.eps}}
\else \ifnum#1=22
\pgftext{\ingeps{b3stuff-item17-1B.eps}}
\fi\fi\fi\fi\fi\fi\fi\fi\fi\fi\fi\fi\fi\fi\fi\fi\fi\fi\fi\fi\fi\fi}

\def\axiombR#1{
\ifnum#1=1
\pgftext{\ingeps{d3stuff-item18-1.eps}}
\else \ifnum#1=2
\pgftext{\ingeps{d3stuff-item18-2.eps}} 
\else \ifnum#1=3
\pgftext{\ingeps{d3stuff-item18-3.eps}}
\else \ifnum#1=4
\pgftext{\ingeps{d3stuff-item18-4.eps}} 
\else \ifnum#1=5
\else \ifnum#1=0
\pgftext{\ingeps{d3stuff-item18a.eps}} 
\fi\fi\fi\fi\fi\fi}

\def\axiombS#1{
\ifnum#1=1
\pgftext{\ingeps{d3stuff-item19-1.eps}}
\else \ifnum#1=2
\pgftext{\ingeps{d3stuff-item19-2.eps}} 
\else \ifnum#1=3
\pgftext{\ingeps{d3stuff-item19-3.eps}}
\else \ifnum#1=4
\pgftext{\ingeps{d3stuff-item19-4.eps}} 
\else \ifnum#1=5
\pgftext{\ingeps{d3stuff-item19-5.eps}} 
\else \ifnum#1=6
\pgftext{\ingeps{d3stuff-item19-6.eps}}
\else \ifnum#1=7
\pgftext{\ingeps{d3stuff-item19-7.eps}} 
\else \ifnum#1=8
\else \ifnum#1=0
\pgftext{\ingeps{d3stuff-item19a.eps}} 
\fi\fi\fi\fi\fi\fi\fi\fi\fi}

\def\axiombY#1{
\ifnum#1=1
\pgftext{\ingeps{b3stuff-item25-1.eps}}
\else \ifnum#1=2
\pgftext{\ingeps{b3stuff-item25-2.eps}} 
\else \ifnum#1=3
\pgftext{\ingeps{b3stuff-item25-3A.eps}}
\else \ifnum#1=4
\pgftext{\ingeps{b3stuff-item25-3B.eps}} 
\else \ifnum#1=5
\pgftext{\ingeps{b3stuff-item25-4A.eps}}
\else \ifnum#1=6
\pgftext{\ingeps{b3stuff-item25-4B.eps}} 
\else \ifnum#1=7
\pgftext{\ingeps{b3stuff-item25-5A.eps}}
\else \ifnum#1=8
\pgftext{\ingeps{b3stuff-item25-5B.eps}}
\else \ifnum#1=9
\pgftext{\ingeps{b3stuff-item25-6A.eps}}
\else \ifnum#1=10
\pgftext{\ingeps{b3stuff-item25-7.eps}}
\else \ifnum#1=11
\else \ifnum#1=0
\pgftext{\ingeps{b3stuff-item25a.eps}}
\fi\fi\fi\fi\fi\fi\fi\fi\fi\fi\fi\fi}

\def\axiombZ#1{
\ifnum#1=1
\pgftext{\ingeps{b3stuff-item26-1.eps}}
\else \ifnum#1=2
\pgftext{\ingeps{b3stuff-item26-8.eps}} 
\else \ifnum#1=3
\pgftext{\ingeps{b3stuff-item26-7B.eps}}
\else \ifnum#1=4
\pgftext{\ingeps{b3stuff-item26-7A.eps}} 
\else \ifnum#1=5
\pgftext{\ingeps{b3stuff-item26-6B.eps}}
\else \ifnum#1=6
\pgftext{\ingeps{b3stuff-item26-6A.eps}} 
\else \ifnum#1=7
\pgftext{\ingeps{b3stuff-item26-5A.eps}}
\else \ifnum#1=8
\pgftext{\ingeps{b3stuff-item26-4A.eps}}
\else \ifnum#1=9
\pgftext{\ingeps{b3stuff-item26-3B.eps}} 
\else \ifnum#1=10
\pgftext{\ingeps{b3stuff-item26-3A.eps}}
\else \ifnum#1=11
\pgftext{\ingeps{b3stuff-item26-2.eps}}
\else \ifnum#1=12
\else \ifnum#1=0
\pgftext{\ingeps{b3stuff-item26a.eps}}
\fi\fi\fi\fi\fi\fi\fi\fi\fi\fi\fi\fi\fi}

\def\axiombRefl#1{
\ifnum#1=2
\axiombB{0}
\else \ifnum#1=51
\axiombE{0}
\else \ifnum#1=52
\axiombE{-1}
\else \ifnum#1=53
\axiombE{-2}
\else \ifnum#1=6
\axiombF{0}
\else \ifnum#1=7
\axiombG{0}
\else \ifnum#1=8
\axiombH{0}
\else \ifnum#1=9
\axiombI{0}
\else \ifnum#1=11
\axiombK{0}
\else \ifnum#1=13
\axiombM{0}
\else \ifnum#1=15
\axiombO{0}
\else \ifnum#1=16
\axiombP{0}
\else \ifnum#1=18
\axiombR{0}
\else \ifnum#1=19
\axiombS{0}
\else \ifnum#1=20
\axiombT{0}
\else \ifnum#1=21
\axiombU{0}
\else \ifnum#1=22
\axiombV{0}
\else \ifnum#1=23
\axiombW{0}
\else \ifnum#1=25
\axiombY{0}
\else \ifnum#1=26
\axiombZ{0}
\fi\fi\fi\fi\fi\fi\fi\fi\fi\fi\fi\fi\fi\fi\fi\fi\fi\fi\fi\fi}

\def\tableAA#1{
\ifnum#1=1
	\sh(0,0);
	\draw[fill=black] (0,0) circle (.3ex);
\else \ifnum#1=2
	\sh(0,0);
	\defect;
\else \ifnum#1=3
	\sh(0,0);
	\sector;
\fi\fi\fi}

\def\tableBB#1{
\ifnum#1=1
	\sh(0,0);
	\defect;
	\iticks;
\else \ifnum#1=2
	\sh(0,0);
	\defect;
	\sh(3,0);
	\tick;
	\defect;
\else \ifnum#1=3
	\sh(0,0);
	\sector;
	\vid;
\else \ifnum#1=4
	\sh(0,0);
	\sectortop;
	\sectorbotarc;
\else \ifnum#1=5
	\sh(0,0);
	\sector;
	\sh(4,0);
	\defect;
\else \ifnum#1=6
	\sh(0,0);
	\defect;
	\sh(3,0);
	\sector;
\else \ifnum#1=7
	\sh(0,0);
	\defect;
	\iticks;
	\sh(3,0);
	\tick;
	\defect;
\else \ifnum#1=8
	\sh(0,0);
	\defect;
	\sh(3,0);
	\tick;
	\defect;
	\iticks;
\else \ifnum#1=9
	\sh(0,0);
	\defect;
	\sh(3,0);
	\tick;
	\defect;
	\sh(6,0);
	\tick;
	\defect;
\fi\fi\fi\fi\fi\fi\fi\fi\fi}

\def\tableCC#1{
\ifnum#1=1
\pgftext{\ingeps{2stuff-item1-1.eps}}
\else \ifnum#1=2
\pgftext{\ingeps{2stuff-item2-1.eps}}
\else \ifnum#1=3
	\sh(0,0);
	\sectortoparclow;
	\sectortoparchigh;
	\sectorbotarclow;
	\sectorbotarchigh;
\else \ifnum#1=4
	\sh(0,0);
	\sector;
	\vid;
	\sh(4,0);
	\defect;
\else \ifnum#1=5
	\sh(0,0);
	\defect;
	\sh(3,0);
	\sector;
	\vid;
\else \ifnum#1=6
	\sh(0,0);
	\sectortop;
	\sectorbotarc;
	\sh(4,0);
	\defect;
\else \ifnum#1=7
	\sh(0,0);
	\defect;
	\sh(3,0);
	\sectortop;
	\sectorbotarc;
\else \ifnum#1=8
	\sh(0,0);
	\sector;
	\sh(4,0);
	\sector;
\else \ifnum#1=9
	\sh(0,0);
	\sector;
	\sh(4,0);
	\defect;
	\sh(7,0);
	\tick;
	\defect;
\else \ifnum#1=10
	\sh(0,0);
	\defect;
	\sh(3,0);
	\tick;
	\defect;
	\sh(6,0);
	\sector;
\else \ifnum#1=11
	\sh(0,0);
	\defect;
	\sh(3,0);
	\sector;
	\sh(7,0);
	\defect;
\else \ifnum#1=12
	\sh(0,0);
	\defect;
	\sh(3,0);
	\tick;
	\defect;
	\sh(6,0);
	\tick;
	\defect;
	\sh(9,0);
	\tick;
	\defect;
\else \ifnum#1=13
	\sh(0,0);
	\defect;
	\iticks;
	\sh(3,0);
	\sector;
\else \ifnum#1=14
	\sh(0,0);
	\sector;
	\sh(4,0);
	\defect;
	\iticks;
\else \ifnum#1=15
	\sh(0,0);
	\defect;
	\iticks;
	\sh(3,0);
	\tick;
	\defect;
	\iticks;
\else \ifnum#1=16
	\sh(0,0);
	\defect;
	\iticks;
	\sh(3,0);
	\tick;
	\defect;
	\sh(6,0);
	\tick;
	\defect;
\else \ifnum#1=17
	\sh(0,0);
	\defect;
	\sh(3,0);
	\tick;
	\defect;
	\sh(6,0);
	\tick;
	\defect;
	\iticks;
\else \ifnum#1=18
	\sh(0,0);
	\defect;
	\sh(3,0);
	\tick;
	\defect;
	\iticks;
	\sh(6,0);
	\tick;
	\defect;
\fi\fi\fi\fi\fi\fi\fi\fi\fi\fi\fi\fi\fi\fi\fi\fi\fi\fi}

\def\tableDD#1{
\ifnum#1=1
	\pgftext{\ingeps{3stuff-item1-1.eps}}
\else \ifnum#1=2
	\begin{scope}[cm={1,0,0,-1,(0,0)}]
		\pgftext{\ingeps{3stuff-item2a.eps}}
	\end{scope}
\else \ifnum#1=3
	\pgftext{\ingeps{3stuff-item2a.eps}}
\else \ifnum#1=4
	\sh(0,0);
	\sectortoparclow;
	\sectortoparchigh;
	\sectorbotarclow;
	\sectorbotarchigh;
	\vids;
\else \ifnum#1=5
	\sh(0,0);
	\sectortoparcmid;
	\sectortoparchigh;
	\sectorbotarcmid;
	\sectorbotarchigh;
	\draw (A) to (C);
\else \ifnum#1=6
	\begin{scope}[cm={1,0,0,-1,(0,0)}]
		\pgftext{\ingeps{3stuff-item5a.eps}}
	\end{scope}
\else \ifnum#1=7
	\pgftext{\ingeps{3stuff-item5a.eps}}
\else \ifnum#1=8
	\pgftext{\ingeps{3stuff-item5b.eps}}
\else \ifnum#1=9
	\pgftext{\ingeps{3stuff-item5c.eps}}
\else \ifnum#1=10
	\sh(0,0);
	\sector;
	\sh(4,0);
	\sector;
	\vid;
\else \ifnum#1=11
	\sh(0,0);
	\sector;
	\vid;
	\sh(4,0);
	\sector;
\else \ifnum#1=12
	\begin{scope}[cm={-1,0,0,1,(0,0)}]
		\pgftext{\ingeps{3stuff-item7a.eps}}
	\end{scope}
\else \ifnum#1=13
	\pgftext{\ingeps{3stuff-item7a.eps}}
\else \ifnum#1=14
	\sh(0,0);
	\sectortop;
	\sectorbotarc;
	\sh(4,0);
	\sector;
\else \ifnum#1=15
	\sh(0,0);
	\sector;
	\sh(4,0);
	\sectortop;
	\sectorbotarc;
\else \ifnum#1=16
	\sh(0,0);
	\sector;
	\vid;
	\sh(4,0);
	\defect;
	\sh(7,0);
	\tick;
	\defect;
\else \ifnum#1=17
	\sh(0,0);
	\defect;
	\sh(3,0);
	\tick;
	\defect;
	\sh(6,0);
	\sector;
	\vid;
\else \ifnum#1=18
	\sh(0,0);
	\defect;
	\sh(3,0);
	\sector;
	\vid;
	\sh(7,0);
	\defect;
\else \ifnum#1=19
	\sh(0,0);
	\sectortop;
	\sectorbotarc;
	\sh(4,0);
	\defect;
	\sh(7,0);
	\tick;
	\defect;
\else \ifnum#1=20
	\sh(0,0);
	\defect;
	\sh(3,0);
	\tick;
	\defect;
	\sh(6,0);
	\sectortop;
	\sectorbotarc;
\else \ifnum#1=21
	\sh(0,0);
	\defect;
	\sh(3,0);
	\sectortop;
	\sectorbotarc;
	\sh(7,0);
	\defect;
\else \ifnum#1=22
	\sh(0,0);
	\sector;
	\sh(4,0);
	\sector;
	\sh(8,0);
	\defect;
\else \ifnum#1=23
	\sh(0,0);
	\defect;
	\sh(3,0);
	\sector;
	\sh(7,0);
	\sector;
\else \ifnum#1=24
	\sh(0,0);
	\sector;
	\sh(4,0);
	\defect;
	\sh(7,0);
	\sector;
\else \ifnum#1=25
	\sh(0,0);
	\sector;
	\sh(4,0);
	\defect;
	\sh(7,0);
	\tick;
	\defect;
	\sh(10,0);
	\tick;
	\defect;
\else \ifnum#1=26
	\sh(0,0);
	\defect;
	\sh(3,0);
	\tick;
	\defect;
	\sh(6,0);
	\tick;
	\defect;
	\sh(9,0);
	\sector;
\else \ifnum#1=27
	\sh(0,0);
	\defect;
	\sh(3,0);
	\sector;
	\sh(7,0);
	\defect;
	\sh(10,0);
	\tick;
	\defect;
\else \ifnum#1=28
	\sh(0,0);
	\defect;
	\sh(3,0);
	\tick;
	\defect;
	\sh(6,0);
	\sector;
	\sh(10,0);
	\defect;
\else \ifnum#1=29
	\sh(0,0);
	\defect;
	\sh(3,0);
	\tick;
	\defect;
	\sh(6,0);
	\tick;
	\defect;
	\sh(9,0);
	\tick;
	\defect;
	\sh(12,0);
	\tick;
	\defect;
\else \ifnum#1=30
	\sh(0,0);
	\sector;
	\vid;
	\sh(4,0);
	\defect;
	\iticks;
\else \ifnum#1=31
	\sh(0,0);
	\defect;
	\iticks;
	\sh(3,0);
	\sector;
	\vid;
\else \ifnum#1=32
	\sh(0,0);
	\sectortop;
	\sectorbotarc;
	\sh(4,0);
	\defect;
	\iticks;
\else \ifnum#1=33
	\sh(0,0);
	\defect;
	\iticks;
	\sh(3,0);
	\sectortop;
	\sectorbotarc;
\else \ifnum#1=34
	\sh(0,0);
	\sector;
	\sh(4,0);
	\defect;
	\iticks;
	\sh(7,0);
	\tick;
	\defect;
\else \ifnum#1=35
	\sh(0,0);
	\defect;
	\sh(3,0);
	\tick;
	\defect;
	\iticks;
	\sh(6,0);
	\sector;
\else \ifnum#1=36
	\sh(0,0);
	\sector;
	\sh(4,0);
	\defect;
	\sh(7,0);
	\tick;
	\defect;
	\iticks;
\else \ifnum#1=37
	\sh(0,0);
	\defect;
	\iticks;
	\sh(3,0);
	\tick;
	\defect;
	\sh(6,0);
	\sector;
\else \ifnum#1=38
	\sh(0,0);
	\defect;
	\iticks;
	\sh(3,0);
	\sector;
	\sh(7,0);
	\defect;
\else \ifnum#1=39
	\sh(0,0);
	\defect;
	\sh(3,0);
	\sector;
	\sh(7,0);
	\defect;
	\iticks;
\else \ifnum#1=40
	\sh(0,0);
	\defect;
	\iticks;
	\sh(3,0);
	\tick;
	\defect;
	\iticks;
	\sh(6,0);
	\tick;
	\defect;
\else \ifnum#1=41
	\sh(0,0);
	\defect;
	\sh(3,0);
	\tick;
	\defect;
	\iticks;
	\sh(6,0);
	\tick;
	\defect;
	\iticks;
\else \ifnum#1=42
	\sh(0,0);
	\defect;
	\iticks;
	\sh(3,0);
	\tick;
	\defect;
	\sh(6,0);
	\tick;
	\defect;
	\iticks;
\else \ifnum#1=43
	\sh(0,0);
	\defect;
	\iticks;
	\sh(3,0);
	\tick;
	\defect;
	\sh(6,0);
	\tick;
	\defect;
	\sh(9,0);
	\tick;
	\defect;
\else \ifnum#1=44
	\sh(0,0);
	\defect;
	\sh(3,0);
	\tick;
	\defect;
	\sh(6,0);
	\tick;
	\defect;
	\sh(9,0);
	\tick;
	\defect;
	\iticks;
\else \ifnum#1=45
	\sh(0,0);
	\defect;
	\sh(3,0);
	\tick;
	\defect;
	\iticks;
	\sh(6,0);
	\tick;
	\defect;
	\sh(9,0);
	\tick;
	\defect;
\else \ifnum#1=46
	\sh(0,0);
	\defect;
	\sh(3,0);
	\tick;
	\defect;
	\sh(6,0);
	\tick;
	\defect;
	\iticks;
	\sh(9,0);
	\tick;
	\defect;
\fi\fi\fi\fi\fi\fi\fi\fi\fi\fi\fi\fi\fi\fi\fi\fi\fi\fi\fi\fi\fi\fi\fi\fi\fi\fi\fi\fi\fi\fi\fi\fi\fi\fi\fi\fi\fi\fi\fi\fi\fi\fi\fi\fi\fi\fi}

\def\axgridA{
\path
	(0,0) coordinate (1)
	(0,-1) coordinate (2)
	(0,-1.75) coordinate (n); % This is a null node, may need to change to a number if need to put content in the node
\draw (1) -- (2) -- (n);
}

\def\axgridB{
\path
	(0,0) coordinate (1)
	(0,-1) coordinate (2)
	(0,-2) coordinate (3)
	(0,-2.75) coordinate (n)
	(0,-3.5) coordinate (0); % This is a reflected node
\draw (1) -- (2) -- (3) -- (n);
}

\def\axgridC{
\path
	(0,0) coordinate (1)
	(0,-1) coordinate (2)
	(0,-2) coordinate (3)
	(0,-2.75) coordinate (n);
\draw (1) -- (2) -- (3) -- (n);
}

\def\axgridD{
\path
	(0,0) coordinate (1)
	(0,-1) coordinate (2)
	(0,-2) coordinate (3)
	(0,-3) coordinate (4)
	(0,-4) coordinate (5)
	(0,-4.75) coordinate (n);
\draw (1) -- (2) -- (3) -- (4) -- (5) -- (n);
}

\def\axgridE{
\path
	(0,0) coordinate (1)
	(0,-1) coordinate (2)
	(0,-2) coordinate (3)
	(0,-3) coordinate (4)
	(0,-3.75) coordinate (n)
	(0,-4.5) coordinate (0)
	(0,-5.5) coordinate (-1)
	(0,-6.5) coordinate (-2);
\draw (1) -- (2) -- (3) -- (4) -- (n);
}

\def\axgridF{
\path
	(0,0) coordinate (1)
	(0,-1) coordinate (2)
	(0,-2) coordinate (3)
	(0,-3) coordinate (4)
	(0,-4) coordinate (5)
	(0,-4.75) coordinate (n)
	(0,-5.5) coordinate (0);
\draw (1) -- (2) -- (3) -- (4) -- (5) -- (n);
}

\def\axgridG{
\path
	(0,0) coordinate (1)
	(0,-1) coordinate (2)
	(0,-2) coordinate (3)
	(0,-3) coordinate (4)
	(0,-4) coordinate (5)
	(0,-5) coordinate (6)
	(0,-5.75) coordinate (n)
	(0,-6.5) coordinate (0);
\draw (1) -- (2) -- (3) -- (4) -- (5) -- (6) -- (n);
}

\def\axgridH{
\path
	(0,0) coordinate (1)
	(0,-1) coordinate (2)
	(0,-2) coordinate (3)
	(0,-3) coordinate (4)
	(0,-4) coordinate (5)
	(0,-5) coordinate (6)
	(0,-6) coordinate (7)
	(0,-7) coordinate (8)
	(0,-7.75) coordinate (n)
	(0,-8.5) coordinate (0);
\draw (1) -- (2) -- (3) -- (4) -- (5) -- (6) -- (7) -- (8) -- (n);
}

\def\axgridI{
\path
	(0,0) coordinate (1)
	(0,-1) coordinate (2)
	(0,-2) coordinate (3)
	(0,-3) coordinate (4)
	(0,-4) coordinate (5)
	(0,-5) coordinate (6)
	(0,-5.75) coordinate (n)
	(0,-6.5) coordinate (0);
\draw (1) -- (2) -- (3) -- (4) -- (5) -- (6) -- (n);
}

\def\axgridJ{
\path
	(0,0) coordinate (1)
	(0,-1) coordinate (2)
	(0,-2) coordinate (3)
	(0,-3) coordinate (4)
	(0,-4) coordinate (5)
	(0,-5) coordinate (6)
	(0,-6) coordinate (7)
	(0,-6.75) coordinate (n);
\draw (1) -- (2) -- (3) -- (4) -- (5) -- (6) -- (7) -- (n);
}

\def\axgridK{
\path
	(0,0) coordinate (1)
	(0,-1) coordinate (2)
	(0,-2) coordinate (3)
	(0,-3) coordinate (4)
	(0,-4) coordinate (5)
	(0,-5) coordinate (6)
	(0,-6) coordinate (7)
	(0,-7) coordinate (8)
	(0,-8) coordinate (9)
	(0,-8.75) coordinate (n)
	(0,-9.5) coordinate (0);
\draw (1) -- (2) -- (3) -- (4) -- (5) -- (6) -- (7) -- (8) -- (9) -- (n);
}

\def\axgridL{
\path
	(0,0) coordinate (1)
	(0,-1) coordinate (2)
	(0,-2) coordinate (3)
	(0,-3) coordinate (4)
	(0,-4) coordinate (5)
	(0,-5) coordinate (6)
	(0,-6) coordinate (7)
	(0,-7) coordinate (8)
	(0,-8) coordinate (9)
	(0,-9) coordinate (10)
	(0,-9.75) coordinate (n)
	(0,-10.5) coordinate (0);
\draw (1) -- (2) -- (3) -- (4) -- (5) -- (6) -- (7) -- (8) -- (9) -- (10) -- (n);
}

\def\axgridM{
\path
	(0,0) coordinate (1)
	(1,0) coordinate (2)
	(1,-1) coordinate (3)
	(1,-2) coordinate (4)
	(1,-3) coordinate (5)
	(1,-4) coordinate (6)
	(1,-5) coordinate (7)
	(1,-6) coordinate (8)	
	(0,-6) coordinate (9)
	(0,-5) coordinate (10)
	(0,-4) coordinate (11)
	(0,-3) coordinate (12)
	(0,-2) coordinate (13)
	(0,-1) coordinate (14);
\draw (1) -- (2) -- (3) -- (4) -- (5) -- (6) -- (7) -- (8) -- (9) -- (10) -- (11) -- (12) -- (13) -- (14) -- (1);
}

\def\axgridN{
\path
	(0,0) coordinate (1)
	(1,0) coordinate (2)
	(1,-1) coordinate (3)
	(1,-2) coordinate (4)
	(1,-3) coordinate (5)
	(1,-4) coordinate (6)
	(1,-5) coordinate (7)
	(0,-5) coordinate (8)	
	(0,-4) coordinate (9)
	(0,-3) coordinate (10)
	(0,-2) coordinate (11)
	(0,-1) coordinate (12);
\draw (1) -- (2) -- (3) -- (4) -- (5) -- (6) -- (7) -- (8) -- (9) -- (10) -- (11) -- (12) -- (1);
}

\def\axgridO{
\path
	(0,0) coordinate (1)
	(1,0) coordinate (2)
	(1,-1) coordinate (3)
	(1,-2) coordinate (4)
	(1,-3) coordinate (5)
	(1,-4) coordinate (6)
	(1,-5) coordinate (7)
	(1,-6) coordinate (8)	
	(1,-7) coordinate (9)	
	(0,-7) coordinate (10)
	(0,-6) coordinate (11)
	(0,-5) coordinate (12)
	(0,-4) coordinate (13)
	(0,-3) coordinate (14)
	(0,-2) coordinate (15)
	(0,-1) coordinate (16);
\draw (1) -- (2) -- (3) -- (4) -- (5) -- (6) -- (7) -- (8) -- (9) -- (10) -- (11) -- (12) -- (13) -- (14) -- (15) -- (16) -- (1);
}

\def\axgridP{
\path
	(0,0) coordinate (1)
	(1,0) coordinate (2)
	(1,-1) coordinate (3)
	(1,-2) coordinate (4)
	(1,-3) coordinate (5)
	(1,-4) coordinate (6)
	(1,-5) coordinate (7)
	(1,-6) coordinate (8)	
	(1,-7) coordinate (9)	
	(1,-8) coordinate (10)	
	(0,-8) coordinate (11)
	(0,-7) coordinate (12)
	(0,-6) coordinate (13)
	(0,-5) coordinate (14)
	(0,-4) coordinate (15)
	(0,-3) coordinate (16)
	(0,-2) coordinate (17)
	(0,-1) coordinate (18);
\draw (1) -- (2) -- (3) -- (4) -- (5) -- (6) -- (7) -- (8) -- (9) -- (10) -- (11) -- (12) -- (13) -- (14) -- (15) -- (16) -- (17) -- (18) -- (1);
}

\def\axgridQ{
\path
	(0,0) coordinate (1)
	(1,0) coordinate (2)
	(1,-1) coordinate (3)
	(1,-2) coordinate (4)
	(1,-3) coordinate (5)
	(1,-4) coordinate (6)
	(1,-5) coordinate (7)
	(1,-6) coordinate (8)	
	(1,-7) coordinate (9)	
	(1,-8) coordinate (10)	
	(1,-9) coordinate (11)	
	(1,-10) coordinate (12)	
	(0,-10) coordinate (13)
	(0,-9) coordinate (14)
	(0,-8) coordinate (15)
	(0,-7) coordinate (16)
	(0,-6) coordinate (17)
	(0,-5) coordinate (18)
	(0,-4) coordinate (19)
	(0,-3) coordinate (20)
	(0,-2) coordinate (21)
	(0,-1) coordinate (22);
\draw (1) -- (2) -- (3) -- (4) -- (5) -- (6) -- (7) -- (8) -- (9) -- (10) -- (11) -- (12) -- (13) -- (14) -- (15) -- (16) -- (17) -- (18) -- (19) -- (20) -- (21) -- (22) -- (1);
}

\def\axgridR{
\path
	(0,0) coordinate (1)
	(0,-1) coordinate (2)
	(0,-2) coordinate (3)
	(0,-3) coordinate (4)
	(0,-3.75) coordinate (n);
\draw (1) -- (2) -- (3) -- (4) -- (n);
}

\def\axgridS{
\path
	(0,0) coordinate (1)
	(0,-1) coordinate (2)
	(0,-2) coordinate (3)
	(0,-3) coordinate (4)
	(0,-4) coordinate (5)
	(0,-5) coordinate (6)
	(0,-6) coordinate (7)
	(0,-6.75) coordinate (n);
\draw (1) -- (2) -- (3) -- (4) -- (5) -- (6) -- (7) -- (n);
}

\def\axgridT{
\path
	(0,0) coordinate (1)
	(0,-1) coordinate (2)
	(0,-2) coordinate (3)
	(0,-3) coordinate (4)
	(0,-4) coordinate (5)
	(0,-5) coordinate (6)
	(0,-6) coordinate (7)
	(0,-7) coordinate (8)
	(0,-8) coordinate (9)
	(0,-9) coordinate (10)
	(0,-9.75) coordinate (n);
\draw (1) -- (2) -- (3) -- (4) -- (5) -- (6) -- (7) -- (8) -- (9) -- (10) -- (n);
}

\def\axgridU{
\path
	(0,0) coordinate (1)
	(0,-1) coordinate (2)
	(0,-2) coordinate (3)
	(0,-3) coordinate (4)
	(0,-4) coordinate (5)
	(0,-5) coordinate (6)
	(0,-6) coordinate (7)
	(0,-7) coordinate (8)
	(0,-7.75) coordinate (n);
\draw (1) -- (2) -- (3) -- (4) -- (5) -- (6) -- (7) -- (8) -- (n);
}

\def\axgridV{
\path
	(0,0) coordinate (1)
	(0,-1) coordinate (2)
	(0,-2) coordinate (3)
	(0,-3) coordinate (4)
	(0,-4) coordinate (5)
	(0,-5) coordinate (6)
	(0,-6) coordinate (7)
	(0,-7) coordinate (8)
	(0,-8) coordinate (9)
	(0,-9) coordinate (10)
	(0,-9.75) coordinate (n);
\draw (1) -- (2) -- (3) -- (4) -- (5) -- (6) -- (7) -- (8) -- (9) -- (10) -- (n);
}

\def\axgridW{
\path
	(0,0) coordinate (1)
	(0,-1) coordinate (2)
	(0,-2) coordinate (3)
	(0,-3) coordinate (4)
	(0,-4) coordinate (5)
	(0,-5) coordinate (6)
	(0,-5.75) coordinate (n);
\draw (1) -- (2) -- (3) -- (4) -- (5) -- (6) -- (n);
}

\def\axgridX{
\path
	(0,0) coordinate (1)
	(0,-1) coordinate (2)
	(0,-2) coordinate (3)
	(0,-3) coordinate (4)
	(0,-4) coordinate (5)
	(0,-5) coordinate (6)
	(0,-5.75) coordinate (n);
\draw (1) -- (2) -- (3) -- (4) -- (5) -- (6) -- (n);
}

\def\axgridY{
\path
	(0,0) coordinate (1)
	(0,-1) coordinate (2)
	(0,-2) coordinate (3)
	(0,-3) coordinate (4)
	(0,-4) coordinate (5)
	(0,-5) coordinate (6)
	(0,-6) coordinate (7)
	(0,-7) coordinate (8)
	(0,-8) coordinate (9)
	(0,-9) coordinate (10)
	(0,-9.75) coordinate (n);
\draw (1) -- (2) -- (3) -- (4) -- (5) -- (6) -- (7) -- (8) -- (9) -- (10) -- (n);
}

\def\axgridZ{
\path
	(0,0) coordinate (1)
	(0,-1) coordinate (2)
	(0,-2) coordinate (3)
	(0,-3) coordinate (4)
	(0,-4) coordinate (5)
	(0,-5) coordinate (6)
	(0,-6) coordinate (7)
	(0,-7) coordinate (8)
	(0,-8) coordinate (9)
	(0,-9) coordinate (10)
	(0,-10) coordinate (11)
	(0,-10.75) coordinate (n);
\draw (1) -- (2) -- (3) -- (4) -- (5) -- (6) -- (7) -- (8) -- (9) -- (10) -- (11) -- (n);
}

\def\axgridRefl{
\path
	(0,0) coordinate (2)
	(0,-1) coordinate (51)
	(0,-2) coordinate (52)
	(0,-3) coordinate (53)
	(0,-4) coordinate (6)
	(0,-5) coordinate (7)
	(0,-6) coordinate (8)
	(0,-7) coordinate (9)
	(0,-8) coordinate (11)
	(0,-9) coordinate (13)
	(0,-10) coordinate (15)
	(0,-11) coordinate (16)
	(0,-12) coordinate (18)
	(0,-13) coordinate (19)
	(0,-14) coordinate (20)
	(0,-15) coordinate (21)
	(0,-16) coordinate (22)
	(0,-17) coordinate (23)
	(0,-18) coordinate (25)
	(0,-19) coordinate (26);
	
}

\def\tablegridAA{
\path
	(0,0) coordinate (1)
	(0,-.75) coordinate (2)
	(0,-1.75) coordinate (3);
}

\def\tablegridBB{
\path
	(0,0) coordinate (1)
	(0,-.75) coordinate (2)
	(0,-1.75) coordinate (3)
	(0,-2.75) coordinate (4)
	(0,-3.75) coordinate (5)
	(0,-4.75) coordinate (6)
	(0,-5.75) coordinate (7)
	(0,-6.5) coordinate (8)
	(0,-7.25) coordinate (9);
}

\def\tablegridCC{
\path
	(0.02,0) coordinate (1)
	(0.02,-1) coordinate (2)
	(0,-2.25) coordinate (3)
	(0,-3.5) coordinate (4)
	(0,-4.5) coordinate (5)
	(0,-5.5) coordinate (6)
	(0,-6.5) coordinate (7)
	(0,-7.5) coordinate (8)
	(1,0) coordinate (9)
	(1,-1) coordinate (10)
	(1,-2) coordinate (11)
	(1,-3) coordinate (12)
	(1,-4) coordinate (13)
	(1,-5) coordinate (14)
	(1,-6) coordinate (15)
	(1,-6.75) coordinate (16)
	(1,-7.5) coordinate (17)
	(1,-8.25) coordinate (18);
}

\def\tablegridDD{
\path
	(0.01,0) coordinate (1)
	(0.01,-1.25) coordinate (2)
	(0.01,-2.625) coordinate (3)
	(0,-4) coordinate (4)
	(0,-5.25) coordinate (5)
	(0.01,-6.5) coordinate (6)
	(0.01,-7.5) coordinate (7)
	(0.01,-8.5) coordinate (8)
	(0.01,-9.5) coordinate (9)
	(1.05,0) coordinate (10)
	(1.05,-1) coordinate (11)
	(1.06,-2.375) coordinate (12)
	(1.06,-3.75) coordinate (13)
	(1.05,-5.125) coordinate (14)
	(1.05,-6.125) coordinate (15)
	(1.05,-7.125) coordinate (16)
	(1.05,-8.125) coordinate (17)
	(1.05,-9.125) coordinate (18)
	(2.25,0) coordinate (19)
	(2.25,-1) coordinate (20)
	(2.25,-2) coordinate (21)
	(2.25,-3) coordinate (22)
	(2.25,-4) coordinate (23)
	(2.25,-5) coordinate (24)
	(2.25,-6) coordinate (25)
	(2.25,-7) coordinate (26)
	(2.25,-8) coordinate (27)
	(2.25,-9) coordinate (28)
	(2.25,-10) coordinate (29)
	(3.75,0) coordinate (30)
	(3.75,-1) coordinate (31)
	(3.75,-2) coordinate (32)
	(3.75,-3) coordinate (33)
	(3.75,-4) coordinate (34)
	(3.75,-5) coordinate (35)
	(3.75,-6) coordinate (36)
	(3.75,-7) coordinate (37)
	(3.75,-8) coordinate (38)
	(3.75,-9) coordinate (39)
	(5,0) coordinate (40)
	(5,-1) coordinate (41)
	(5,-2) coordinate (42)
	(5,-3) coordinate (43)
	(5,-4) coordinate (44)
	(5,-5) coordinate (45)
	(5,-6) coordinate (46);
}

%% file: cnii-defs.tex
\def\common{1} % Comm on (1) or off (0).
\def\labon{1} % Labels on (1) or off (0).

\pgfdeclarelayer{background} 
\pgfdeclarelayer{foreground}  % put bullets in this layer
\pgfsetlayers{background,main,foreground}

\def\commsize{\footnotesize} % Text size in comm.
\tikzset{c/.style={font=\commsize}}

\def\labelsize{\footnotesize} % Text size in edge labels.
\tikzset{la/.style={midway,above,sloped,font=\labelsize}}
\tikzset{lb/.style={midway,below,sloped,font=\labelsize}}

\def\gap{.2} % This is the gap
\def\fudge{0.9} % This is the factor by which hook boxes are scaled when they look too tall
\def\myinnersep{.15cm}

\def\hexangle{62}
\def\hexangleA{50}
\def\gapfactorA{1.45} %60 --> 2, 45 --> 1.41, 50 --> 1.56, 55 --> 1.74
\def\hexangleB{55}
\def\gapfactorB{1.74} %properly 1.74
\def\hexangleC{50}
\def\gapfactorC{1.45} %properly 1.56

\def\hidist{.15} % Distance between half identity ticks
\def\hilength{.4} % Length of half identity ticks
\def\hiangle{40} % Angle of half identity ticks
\def\ticklength{.2} % Length of vertical ticks
\def\tickfactor{.6} % Factor by which the ticks are shortened in pictures where they double up.

\def\cusplength{.3}%.3 is good %Total distance over which the bezier curves occur
\def\cuspwidth{.18}%.18 is good %Total width of diffeomorphism ballons.
\def\cuspfactor{.5} %.5 is good %Controls the bezier curve steepness.
\def\cw{.5*\cuspwidth}
\def\cuspstretch{1.5} % Amount by which the cuspwidth stretches when it is too close to a dot.
\def\bulletsize{.06} %.06 reasonable % Note that for a few pictures, .05 would be better.
\def\bulletcolor{black!80!gray}

\def\dfctangle{60} % Angle from horizontal of crosses in defect x
\def\dfctlength{.16} % Length of defect; Calculate as 2*bulletsize/sin(dfctangle); .14 is for dfct angle 60, bulletsize .06. 
\def\idfctdist{.07} % Distance between the identity defect ticks; This must be visually the same as the distance between the end points of the defect x, ie dfctlength*cos(dfctangle); for dfct angle 60, is is .5*dfctlength, but that's a touch big visually. 
\def\idfctlength{.14} % Length of identity defect ticks; Use 2*bulletsize as base.
\tikzset{littlepic/.style={scale=.1ex, draw=black}}  %originally .08ex 
\tikzset{littlepicbull/.style={scale=.125ex, draw=black}}  %
\tikzset{littlenode/.style={fill=white,inner sep=\myinnersep}}
\tikzset{displaypic/.style={scale=.15ex, draw=black}} 
\def\innertikz{\tikz[littlepic]}
\def\innertikzb{\tikz[littlepicbull]}
\def\displaytikz{\tikz[displaypic]}

\def\mycolor#1{
\ifnum #1=1
	black!0 %white
\else \ifnum #1=2
	gray!40 %40?
\else \ifnum #1=3
	black!50 %60?
\fi\fi\fi
}

\def\bx(#1,#2)#3{
\def\myclr{\mycolor{#3}}
	\draw[fill=\myclr] (\x,\y) rectangle +(#1,#2);
}

\def\hex(#1,#2)#3{
\def\myclr{\mycolor{#3}}
	\coordinate (A) at (\x,\y);
	\coordinate (B) at ($ (A) + (0,.5*#2) $) ;
	\coordinate (C) at ($ (A) + (0,#2) $) ;
	\coordinate (F) at ($ (A) + (#1,0) $) ;
	\coordinate (E) at ($ (F) + (0,.5*#2) $) ;
	\coordinate (D) at ($ (F) + (0,#2) $) ;
	\coordinate (G) at ($ (B) + (360-\hexangleA:1) $) ;
	\coordinate (J) at ($ (B) + (\hexangleA:1) $) ;
	\coordinate (H) at ($ (E) + (180+\hexangleA:1) $) ;
	\coordinate (I) at ($ (E) + (180-\hexangleA:1) $) ;
	\draw[fill=\myclr] (B) -- 
		(intersection of B--G and A--F) -- 
		(intersection of A--F and H--E) --
		(E) --
		(intersection of E--I and D--C) --
		(intersection of D--C and J--B) --
		cycle;
}

\def\hexb(#1,#2)#3{ % "Hexagon" with bullets, uses \hexangle.
\def\myclr{\mycolor{#3}}
	\coordinate (A) at (\x,\y);
	\coordinate (B) at ($ (A) + (0,.5*#2) $) ;
	\coordinate (C) at ($ (A) + (0,#2) $) ;
	\coordinate (F) at ($ (A) + (#1,0) $) ;
	\coordinate (E) at ($ (F) + (0,.5*#2) $) ;
	\coordinate (D) at ($ (F) + (0,#2) $) ;
	\coordinate (G) at ($ (B) + (360-\hexangle:1) $) ;
	\coordinate (J) at ($ (B) + (\hexangle:1) $) ;
	\coordinate (H) at ($ (E) + (180+\hexangle:1) $) ;
	\coordinate (I) at ($ (E) + (180-\hexangle:1) $) ;
	\draw[fill=\myclr] (B) -- 
		(intersection of B--G and A--F) -- 
		(intersection of A--F and H--E) --
		(E) --
		(intersection of E--I and D--C) --
		(intersection of D--C and J--B) --
		cycle;
	\begin{pgfonlayer}{foreground}
			\fill [color=\bulletcolor] (B) circle (\bulletsize);
			\fill [color=\bulletcolor] (intersection of B--G and A--F) circle (\bulletsize);
			\fill [color=\bulletcolor] (intersection of A--F and H--E) circle (\bulletsize);
			\fill [color=\bulletcolor] (E) circle (\bulletsize);
			\fill [color=\bulletcolor] (intersection of E--I and D--C) circle (\bulletsize);
			\fill [color=\bulletcolor] (intersection of D--C and J--B) circle (\bulletsize);
	\end{pgfonlayer}
}

\def\hkur(#1,#2)#3,#4,#5{ %length, height, half spike (1) or full spike (2), no gap (0) or gap (1), color
\def\myclr{\mycolor{#5}}
\ifnum#3=2
	\def\hexangleH{\hexangleB}
	\def\gapfactorH{\gapfactorB}
\else \ifnum#3=1
	\def\hexangleH{\hexangleC}
	\def\gapfactorH{\gapfactorC}
\fi\fi

	\coordinate (A) at (\x,\y);
	\coordinate (B) at ($ (A) + (0,#2) $);
	\coordinate (C) at ($ (A) + (#1,0) $);
	\coordinate (D) at ($ (A) + (#1,#2) $);
	\coordinate (E) at ($ (D) + (0,.5*#3-\gapfactorH*\gap*#4+\gap*#4) $);
	\coordinate (F) at ($ (E) + (180+\hexangleH:1) $);
	\draw[fill=\myclr] 
		(A) --
		(B) --
		(intersection of B--D and F--E) --
		(E) --
		(C) --
		cycle;
}

\def\hkul(#1,#2)#3,#4,#5{ 
	\begin{scope}[cm={-1,0,0,1,(\x+\x+#1,0)}]
		\hkur(#1,#2)#3,#4,#5;
	\end{scope}
}

\def\hkdr(#1,#2)#3,#4,#5{ 
	\begin{scope}[cm={1,0,0,-1,(0,\y+\y+#2)}]
		\hkur(#1,#2)#3,#4,#5;
	\end{scope}
}

\def\hkdl(#1,#2)#3,#4,#5{ 
	\begin{scope}[cm={-1,0,0,-1,(\x+\x+#1,\y+\y+#2)}]
		\hkur(#1,#2)#3,#4,#5;
	\end{scope}
}

\def\dhku(#1,#2)#3,#4,#5{ %length, height, half spike (1) or full spike (2), no gap (0) or gap (1), color
\def\myclr{\mycolor{#5}}
\ifnum#3=2
	\def\hexangleH{\hexangleB}
	\def\gapfactorH{\gapfactorB}
\else \ifnum#3=1
	\def\hexangleH{\hexangleC}
	\def\gapfactorH{\gapfactorC}
\fi\fi

	\coordinate (A) at (\x,\y);
	\coordinate (B) at ($ (A) + (0,#2) $);
	\coordinate (C) at ($ (A) + (#1,0) $);
	\coordinate (D) at ($ (A) + (#1,#2) $);
	\coordinate (E) at ($ (D) + (0,.5*#3-\gapfactorH*\gap*#4+\gap*#4) $);
	\coordinate (F) at ($ (B) + (0,.5*#3-\gapfactorH*\gap*#4+\gap*#4) $);
	\coordinate (G) at ($ (E) + (180+\hexangleH:1) $);
	\coordinate (H) at ($ (F) + (360-\hexangleH:1) $);
	\draw[fill=\myclr] 
		(A) --
		(F) --
		(intersection of B--D and F--H) --
		(intersection of B--D and G--E) --
		(E) --
		(C) --
		cycle;
}

\def\dhkd(#1,#2)#3,#4,#5{ 
	\begin{scope}[cm={1,0,0,-1,(0,\y+\y+#2)}]
		\dhku(#1,#2)#3,#4,#5;
	\end{scope}
}

\def\halfhexr(#1,#2)#3{ 
\def\myclr{\mycolor{#3}}
	\def\hexangleH{\hexangleC}
	\def\gapfactorH{\gapfactorC}

	\coordinate (A) at (\x,\y);
	\coordinate (B) at ($ (A) + (0,#2) $);
	\coordinate (C) at ($ (A) + (#1,0) $);
	\coordinate (D) at ($ (A) + (#1,#2) $);
	\coordinate (E) at ($ (C) + (0,.5*#2) $);
	\coordinate (F) at ($ (E) + (180-\hexangleH:1) $);
	\coordinate (G) at ($ (E) + (180+\hexangleH:1) $);
	\draw[fill=\myclr] 
		(A) --
		(B) --
		(intersection of B--D and E--F) --
		(E) --
		(intersection of A--C and E--G) --
		cycle;
}

\def\halfhexl(#1,#2)#3{ 
	\begin{scope}[cm={-1,0,0,1,(\x+\x+#1,0)}]
		\halfhexr(#1,#2)#3;
	\end{scope}
}

\def\qhexur(#1,#2)#3{ 
\def\myclr{\mycolor{#3}}
	\def\hexangleH{\hexangleB}
	\def\gapfactorH{\gapfactorB}

	\coordinate (A) at (\x,\y);
	\coordinate (B) at ($ (A) + (0,#2) $);
	\coordinate (C) at ($ (A) + (#1,0) $);
	\coordinate (D) at ($ (A) + (#1,#2) $);
	\coordinate (E) at ($ (D) + (180+\hexangleH:1) $);
	\draw[fill=\myclr] 
		(A) --
		(B) --
		(D) --
		(intersection of D--E and C--A) --
		cycle;
}

\def\qhexul(#1,#2)#3{ 
	\begin{scope}[cm={-1,0,0,1,(\x+\x+#1,0)}]
		\qhexur(#1,#2)#3;
	\end{scope}
}

\def\qhexdr(#1,#2)#3{ 
	\begin{scope}[cm={1,0,0,-1,(0,\y+\y+#2)}]
		\qhexur(#1,#2)#3;
	\end{scope}
}

\def\qhexdl(#1,#2)#3{ 
	\begin{scope}[cm={-1,0,0,-1,(\x+\x+#1,\y+\y+#2)}]
		\qhexur(#1,#2)#3;
	\end{scope}
}

\def\halfidr(#1,#2)#3{ %length and height of box, percentage along
	\coordinate (A) at ($ (\x,\y) + (0,#2) $);
	\coordinate (B) at ($ (\x,\y) + (#1,#2) $);
	\coordinate (C) at ($ (A)!#3!(B) - (.5*\hidist,0) $);
	\coordinate (D) at ($ (C) +(270-\hiangle:\hilength) $);
	\coordinate (E) at ($ (A)!#3!(B) + (.5*\hidist,0) $);
	\coordinate (F) at ($ (E) +(270-\hiangle:\hilength) $);
	\draw (C) -- (D);	
	\draw (E) -- (F);
}

\def\halfidl(#1,#2)#3{ %length and height of box, percentage along
	\coordinate (A) at ($ (\x,\y) + (0,#2) $);
	\coordinate (B) at ($ (\x,\y) + (#1,#2) $);
	\coordinate (C) at ($ (A)!#3!(B) - (.5*\hidist,0) $);
	\coordinate (D) at ($ (C) +(270+\hiangle:\hilength) $);
	\coordinate (E) at ($ (A)!#3!(B) + (.5*\hidist,0) $);
	\coordinate (F) at ($ (E) +(270+\hiangle:\hilength) $);
	\draw (C) -- (D);	
	\draw (E) -- (F);
}

\def\tick(#1,#2){
	\draw ($ (#1,#2) + (0,.5*\ticklength) $) -- ($ (#1,#2) - (0,.5*\ticklength) $);
}

\def\cusp(#1,#2,#3,#4){ %angle (0/90/180/270), left/right (-1/1), first, second length
	(\x,\y) -- ++(#1:#3-.5*\cusplength) .. controls +(#1:\cuspfactor*\cusplength) and +(180+#1:\cuspfactor*\cusplength) .. ++($ (#1:\cusplength) + (#1-#2*90:.5*\cuspwidth) $) -- ++ (#1:#4-.5*\cusplength)
}

\def\bend(#1,#2,#3,#4,#5){ %angle, direction, lengths side middle side
	\cusp(#1,#2,#3,.5*#4)
	{[shift={($ (#1:#3+.5*#4) + (#1-#2*90:.5*\cuspwidth) $)}] -- \cusp(#1,-#2,.5*#4,#5)}
}

\def\bubble(#1,#2,#3)#4{
\pgfmathsetmacro{\one}{#1}
\pgfmathsetmacro{\two}{#2}
\pgfmathsetmacro{\three}{#3}
\pgfmathsetmacro{\four}{#4}
	\draw[fill=\mycolor{#4}]
		\bend(0,1,\one,\two,\three) 
		{[shift={(#1+#2+#3,0)}] -- \bend(180,1,\three,\two,\one)} -- cycle;
}

\def\bubblev(#1,#2,#3)#4{
\pgfmathsetmacro{\one}{#1}
\pgfmathsetmacro{\two}{#2}
\pgfmathsetmacro{\three}{#3}
\pgfmathsetmacro{\four}{#4}
	\draw[fill=\mycolor{#4},even odd rule]
		\bend(90,1,\one,\two,\three) 
		{[shift={(0,#1+#2+#3)}] -- \bend(270,1,\three,\two,\one)} -- cycle;
}

\def\cornerur(#1,#2,#3,#4)#5{
	\draw[fill=\mycolor{#5},even odd rule]
		{[shift={(0,#3+#4)}]
			\cusp(0,1,#1,#2-.5*\cuspwidth)}
		{[shift={(#1+#2-.5*\cuspwidth,#3+#4-.5*\cuspwidth)}] 
			-- \cusp(270,-1,#4-.5*\cuspwidth,#3)}
		{[shift={(#1+#2,0)}] 
			-- \cusp(90,1,#3,#4+.5*\cuspwidth)}
		{[shift={(#1+#2+.5*\cuspwidth,#3+#4+.5*\cuspwidth)}] 
			-- \cusp(180,-1,#2+.5*\cuspwidth,#1)} -- cycle;
	\begin{pgfonlayer}{foreground}	
		\bull(\x+#1+#2-.5*\cuspwidth,\y-.5*\cuspwidth+#3+#4);
		\bull(\x+#1+#2+.5*\cuspwidth,\y+.5*\cuspwidth+#3+#4);
	\end{pgfonlayer}
}

\def\cornerul(#1,#2,#3,#4)#5{
	\begin{scope}[cm={-1,0,0,1,(\x+\x+#1+#2,0)}]
		\cornerur(#1,#2,#3,#4)#5;
	\end{scope}
}

\def\cornerdr(#1,#2,#3,#4)#5{
	\begin{scope}[cm={1,0,0,-1,(0,\y+\y+#3+#4)}]
		\cornerur(#1,#2,#3,#4)#5;
	\end{scope}
}

\def\cornerdl(#1,#2,#3,#4)#5{
	\begin{scope}[cm={-1,0,0,-1,(\x+\x+#1+#2,\y+\y+#3+#4)}]
		\cornerur(#1,#2,#3,#4)#5;
	\end{scope}
}

\def\boxt(#1,#2,#3,#4)#5{ %before, during, after cusp lengths, height, color
\pgfmathsetmacro{\one}{#1}
\pgfmathsetmacro{\two}{#2}
\pgfmathsetmacro{\three}{#3}
\pgfmathsetmacro{\four}{#4}
	\draw[fill=\mycolor{#5},even odd rule]
		(\x,\y) -- ++(#1+#2+#3,0) -- ++ (0,#4)
		{[shift={(#1+#2+#3,#4)}] -- \bend(180,-1,\three,\two,\one)} -- cycle;
	\begin{pgfonlayer}{foreground}	
		\bull(\x,\y);
		\bull(\x+#1+#2+#3,\y);
		\bull(\x+#1+#2+#3,\y+#4);
		\bull(\x,\y+#4);
	\end{pgfonlayer}
} 

\def\boxb(#1,#2,#3,#4)#5{
	\begin{scope}[cm={1,0,0,-1,(0,\y+\y+#4)}]
		\boxt(#1,#2,#3,#4)#5;
	\end{scope}
}

\def\boxslopeur(#1,#2,#3)#4{ %before, during cusp lengths, height, color
	\draw[fill=\mycolor{#4},even odd rule]
		(\x,\y) -- ++(#1+#2,0) -- ++ (0,#3-.5*\cuspwidth)
		{[shift={(#1+#2,#3-.5*\cuspwidth)}] -- \cusp(180,1,#2,#1)} -- cycle;
	\begin{pgfonlayer}{foreground}
		\bull(\x,\y);
		\bull(\x+#1+#2,\y);
		\bull(\x+#1+#2,\y+#3-.5*\cuspwidth);
		\bull(\x,\y+#3);
	\end{pgfonlayer}
}

\def\boxslopeul(#1,#2,#3)#4{
	\begin{scope}[cm={-1,0,0,1,(\x+\x+#1+#2,0)}]
		\boxslopeur(#1,#2,#3)#4;
	\end{scope}
}

\def\boxslopedr(#1,#2,#3)#4{
	\begin{scope}[cm={1,0,0,-1,(0,\y+\y+#3)}]
		\boxslopeur(#1,#2,#3)#4;
	\end{scope}
}

\def\boxslopedl(#1,#2,#3)#4{
	\begin{scope}[cm={-1,0,0,-1,(\x+\x+#1+#2,\y+\y+#3)}]
		\boxslopeur(#1,#2,#3)#4;
	\end{scope}
}

\def\boxsloperu(#1,#2,#3)#4{ %length, before, after cusp height, color
	\begin{scope}[cm={1,0,0,1,(0,#2+#3)},rotate around={270:(\x,\y)}]
		\boxslopeul(#2,#3,#1)#4;
	\end{scope}
}

\def\boxsloperd(#1,#2,#3)#4{ 
	\begin{scope}[cm={1,0,0,1,(0,#2+#3)},rotate around={270:(\x,\y)}]
		\boxslopeur(#2,#3,#1)#4;
	\end{scope}
}

\def\boxslopelu(#1,#2,#3)#4{
	\begin{scope}[cm={1,0,0,1,(#1,0)},rotate around={90:(\x,\y)}]
		\boxslopeur(#2,#3,#1)#4;
	\end{scope}
}

\def\boxslopeld(#1,#2,#3)#4{
	\begin{scope}[cm={1,0,0,1,(#1,0)},rotate around={90:(\x,\y)}]
		\boxslopeul(#2,#3,#1)#4;
	\end{scope}
}

\def\boxcornur(#1,#2,#3,#4)#5{ %before and after cusp horiz, before and after vert, color
	\draw[fill=\mycolor{#5}, even odd rule]
		(\x,\y) -- ++(#1+#2,0)
		{[shift={(#1+#2,0)}] 
			-- \cusp(90,-1,#3,#4-.5*\cuspwidth)}
		{[shift={(#1+#2-.5*\cuspwidth,#3+#4-.5*\cuspwidth)}] 
			-- \cusp(180,1,#2-.5*\cuspwidth,#1)} -- cycle;
	\begin{pgfonlayer}{foreground}
		\bull(\x,\y);
		\bull(\x+#1+#2,\y);
		\bull(\x+#1+#2-.5*\cuspwidth,\y+#3+#4-.5*\cuspwidth);
		\bull(\x,\y+#3+#4);
	\end{pgfonlayer}
}

\def\boxcornul(#1,#2,#3,#4)#5{
	\begin{scope}[cm={-1,0,0,1,(\x+\x+#1+#2,0)}]
		\boxcornur(#1,#2,#3,#4)#5;
	\end{scope}
}

\def\boxcorndr(#1,#2,#3,#4)#5{
	\begin{scope}[cm={1,0,0,-1,(0,\y+\y+#3+#4)}]
		\boxcornur(#1,#2,#3,#4)#5;
	\end{scope}
}

\def\boxcorndl(#1,#2,#3,#4)#5{
	\begin{scope}[cm={-1,0,0,-1,(\x+\x+#1+#2,\y+\y+#3+#4)}]
		\boxcornur(#1,#2,#3,#4)#5;
	\end{scope}
}

\def\sbx(#1,#2)#3{
\pgfmathsetmacro{\one}{#1}
\pgfmathsetmacro{\two}{#2}
\def\myclr{\mycolor{#3}}
	\draw[fill=\myclr] (\x,\y) rectangle +(#1,#2);
}

\def\bull(#1,#2){
	\fill [color=\bulletcolor] (#1,#2) circle (\bulletsize);
}

\def\bxb(#1,#2)#3{
\pgfmathsetmacro{\one}{#1}
\pgfmathsetmacro{\two}{#2}
\def\myclr{\mycolor{#3}}
	\draw[fill=\myclr] (\x,\y) rectangle +(#1,#2);
\begin{pgfonlayer}{foreground}
	\bull(\x,\y);\bull(\x,\y+.5*\two);\bull(\x,\y+#2);
	\bull(\x+#1,\y);\bull(\x+#1,\y+.5*\two);\bull(\x+#1,\y+#2);
\end{pgfonlayer}
}

\def\bxbt(#1,#2)#3{
\pgfmathsetmacro{\one}{#1}
\pgfmathsetmacro{\two}{#2}
\def\myclr{\mycolor{#3}}
	\draw[fill=\myclr] (\x,\y) rectangle +(#1,#2);
\begin{pgfonlayer}{foreground}
	\bull(\x,\y);;\bull(\x,\y+#2);
	\bull(\x+#1,\y);\bull(\x+#1,\y+#2);
	\bull(\x+.33*\one,\y+#2);\bull(\x+.67*\one,\y+#2);
\end{pgfonlayer}
}

\def\bxbb(#1,#2)#3{
\pgfmathsetmacro{\one}{#1}
\pgfmathsetmacro{\two}{#2}
\def\myclr{\mycolor{#3}}
	\draw[fill=\myclr] (\x,\y) rectangle +(#1,#2);
\begin{pgfonlayer}{foreground}
	\bull(\x,\y);;\bull(\x,\y+#2);
	\bull(\x+#1,\y);\bull(\x+#1,\y+#2);
	\bull(\x+.33*\one,\y);\bull(\x+.67*\one,\y);
\end{pgfonlayer}
}

\def\bxbur(#1,#2)#3{
\pgfmathsetmacro{\one}{#1}
\pgfmathsetmacro{\two}{#2}
\def\myclr{\mycolor{#3}}
	\draw[fill=\myclr] (\x,\y) rectangle +(#1,#2);
\begin{pgfonlayer}{foreground}
	\bull(\x,\y);;\bull(\x,\y+#2);
	\bull(\x+#1,\y);\bull(\x+#1,\y+#2);
	\bull(\x+.5*\one,\y+#2);\bull(\x+#1,\y+.5*\two);
\end{pgfonlayer}
}

\def\bxbdr(#1,#2)#3{
\pgfmathsetmacro{\one}{#1}
\pgfmathsetmacro{\two}{#2}
\def\myclr{\mycolor{#3}}
	\draw[fill=\myclr] (\x,\y) rectangle +(#1,#2);
\begin{pgfonlayer}{foreground}
	\bull(\x,\y);;\bull(\x,\y+#2);
	\bull(\x+#1,\y);\bull(\x+#1,\y+#2);
	\bull(\x+.5*\one,\y);\bull(\x+#1,\y+.5*\two); 
\end{pgfonlayer}
}

\def\bxbul(#1,#2)#3{
\pgfmathsetmacro{\one}{#1}
\pgfmathsetmacro{\two}{#2}
\def\myclr{\mycolor{#3}}
	\draw[fill=\myclr] (\x,\y) rectangle +(#1,#2);
\begin{pgfonlayer}{foreground}
	\bull(\x,\y);;\bull(\x,\y+#2);
	\bull(\x+#1,\y);\bull(\x+#1,\y+#2);
	\bull(\x+.5*\one,\y+#2);\bull(\x,\y+.5*\two);
\end{pgfonlayer}
}

\def\bxbdl(#1,#2)#3{
\pgfmathsetmacro{\one}{#1}
\pgfmathsetmacro{\two}{#2}
\def\myclr{\mycolor{#3}}
	\draw[fill=\myclr] (\x,\y) rectangle +(#1,#2);
\begin{pgfonlayer}{foreground}
	\bull(\x,\y);;\bull(\x,\y+#2);
	\bull(\x+#1,\y);\bull(\x+#1,\y+#2);
	\bull(\x+.5*\one,\y);\bull(\x,\y+.5*\two);
\end{pgfonlayer}
}

\def\sh(#1,#2){\def\x{#1}\def\y{#2}}

\def\dfct(#1,#2){
	\coordinate (A) at ($ (#1,#2) +(180+\dfctangle:.5*\dfctlength) $);
	\coordinate (B) at ($ (#1,#2) +(360-\dfctangle:.5*\dfctlength) $);
	\coordinate (C) at ($ (#1,#2) +(180-\dfctangle:.5*\dfctlength) $);
	\coordinate (D) at ($ (#1,#2) +(\dfctangle:.5*\dfctlength) $);
	\draw (A) -- (D);
	\draw (B) -- (C);
}

\def\idfct(#1,#2){
	\coordinate (A) at ($ (#1,#2) + (-.5*\idfctdist,-.5*\idfctlength) $);
	\coordinate (B) at ($ (#1,#2) + (-.5*\idfctdist,.5*\idfctlength) $);
	\coordinate (C) at ($ (#1,#2) + (.5*\idfctdist,-.5*\idfctlength) $);
	\coordinate (D) at ($ (#1,#2) + (.5*\idfctdist,.5*\idfctlength) $);
	\draw (A) -- (B);
	\draw (C) -- (D);
}

\def\expand{expand}
\def\whisk{whisker}
\def\ident{identity}
\def\switch{switch}
\def\pass{pass}
\def\swap{swap}
\def\flip{flip}
\def\hswitch{half switch}

\def\quasiidentity{quasi-id}
\def\identity{id}

\def\onetimesone{$1 \boxtimes 1$}

\def\ax#1-#2{
\ifthenelse{\equal{#1}{1}}{
			\AxiomA{#2}}{}
\ifthenelse{\equal{#1}{2}}{
			\AxiomB{#2}}{}
\ifthenelse{\equal{#1}{3}}{
			\AxiomC{#2}}{}		
\ifthenelse{\equal{#1}{4}}{
			\AxiomD{#2}}{}
\ifthenelse{\equal{#1}{5}}{
			\AxiomE{#2}}{}
\ifthenelse{\equal{#1}{6}}{
			\AxiomF{#2}}{}
\ifthenelse{\equal{#1}{6a}}{
			\AxiomFa{#2}}{}
\ifthenelse{\equal{#1}{7}}{
			\AxiomG{#2}}{}
\ifthenelse{\equal{#1}{8}}{
			\AxiomH{#2}}{}
\ifthenelse{\equal{#1}{8s}}{
			\AxiomHs{#2}}{}		
\ifthenelse{\equal{#1}{8a}}{
			\AxiomHa{#2}}{}
\ifthenelse{\equal{#1}{8b}}{
			\AxiomHb{#2}}{}
\ifthenelse{\equal{#1}{11}}{
			\AxiomK{#2}}{}
\ifthenelse{\equal{#1}{12}}{
			\AxiomL{#2}}{}
\ifthenelse{\equal{#1}{13}}{
			\AxiomM{#2}}{}
\ifthenelse{\equal{#1}{13a}}{
			\AxiomMa{#2}}{}		
\ifthenelse{\equal{#1}{14}}{
			\AxiomN{#2}}{}
\ifthenelse{\equal{#1}{14a}}{
			\AxiomNa{#2}}{}
\ifthenelse{\equal{#1}{15}}{
			\AxiomO{#2}}{}
\ifthenelse{\equal{#1}{16}}{
			\AxiomP{#2}}{}
\ifthenelse{\equal{#1}{18}}{
			\AxiomR{#2}}{}
\ifthenelse{\equal{#1}{18new}}{
			\AxiomRnew{#2}}{}
\ifthenelse{\equal{#1}{19}}{
			\AxiomS{#2}}{}	
\ifthenelse{\equal{#1}{19new}}{
			\AxiomSnew{#2}}{}	
\ifthenelse{\equal{#1}{20}}{
			\AxiomT{#2}}{}
\ifthenelse{\equal{#1}{20new}}{
			\AxiomTnew{#2}}{}
\ifthenelse{\equal{#1}{20a}}{
			\AxiomTa{#2}}{}
\ifthenelse{\equal{#1}{21}}{
			\AxiomU{#2}}{}
\ifthenelse{\equal{#1}{21a}}{
			\AxiomUa{#2}}{}
\ifthenelse{\equal{#1}{21new}}{
			\AxiomUnew{#2}}{}
\ifthenelse{\equal{#1}{21anew}}{
			\AxiomUanew{#2}}{}
\ifthenelse{\equal{#1}{22}}{
			\AxiomV{#2}}{}
\ifthenelse{\equal{#1}{22new}}{
			\AxiomVnew{#2}}{}
\ifthenelse{\equal{#1}{23}}{
			\AxiomW{#2}}{}		
\ifthenelse{\equal{#1}{23new}}{
			\AxiomWnew{#2}}{}		
\ifthenelse{\equal{#1}{23a}}{
			\AxiomWa{#2}}{}
\ifthenelse{\equal{#1}{24}}{
			\AxiomX{#2}}{}
\ifthenelse{\equal{#1}{24new}}{
			\AxiomXnew{#2}}{}
\ifthenelse{\equal{#1}{25}}{
			\AxiomY{#2}}{}
\ifthenelse{\equal{#1}{26}}{
			\AxiomZ{#2}}{}			
}

\def\axgrid#1{
\ifthenelse{\equal{#1}{1}}{
			\axgridA}{}
\ifthenelse{\equal{#1}{2}}{
			\axgridB}{}
\ifthenelse{\equal{#1}{3}}{
			\axgridC}{}		
\ifthenelse{\equal{#1}{4}}{
			\axgridD}{}
\ifthenelse{\equal{#1}{5}}{
			\axgridE}{}
\ifthenelse{\equal{#1}{6}}{
			\axgridF}{}
\ifthenelse{\equal{#1}{6a}}{
			\axgridFa}{}
\ifthenelse{\equal{#1}{7}}{
			\axgridG}{}
\ifthenelse{\equal{#1}{8}}{
			\axgridH}{}
\ifthenelse{\equal{#1}{8s}}{
			\axgridHs}{}		
\ifthenelse{\equal{#1}{8a}}{
			\axgridHa}{}
\ifthenelse{\equal{#1}{8b}}{
			\axgridHb}{}
\ifthenelse{\equal{#1}{11}}{
			\axgridK}{}
\ifthenelse{\equal{#1}{12}}{
			\axgridL}{}
\ifthenelse{\equal{#1}{13}}{
			\axgridM}{}
\ifthenelse{\equal{#1}{13a}}{
			\axgridMa}{}		
\ifthenelse{\equal{#1}{14}}{
			\axgridN}{}
\ifthenelse{\equal{#1}{14a}}{
			\axgridNa}{}
\ifthenelse{\equal{#1}{15}}{
			\axgridO}{}
\ifthenelse{\equal{#1}{16}}{
			\axgridP}{}
\ifthenelse{\equal{#1}{18}}{
			\axgridR}{}
\ifthenelse{\equal{#1}{18new}}{
			\axgridRnew}{}
\ifthenelse{\equal{#1}{19}}{
			\axgridS}{}
\ifthenelse{\equal{#1}{19new}}{
			\axgridSnew}{}		
\ifthenelse{\equal{#1}{20}}{
			\axgridT}{}
\ifthenelse{\equal{#1}{20new}}{
			\axgridTnew}{}
\ifthenelse{\equal{#1}{20a}}{
			\axgridTa}{}
\ifthenelse{\equal{#1}{21}}{
			\axgridU}{}
\ifthenelse{\equal{#1}{21a}}{
			\axgridUa}{}
\ifthenelse{\equal{#1}{21new}}{
			\axgridUnew}{}
\ifthenelse{\equal{#1}{21anew}}{
			\axgridUanew}{}
\ifthenelse{\equal{#1}{22}}{
			\axgridV}{}
\ifthenelse{\equal{#1}{22new}}{
			\axgridVnew}{}
\ifthenelse{\equal{#1}{23}}{
			\axgridW}{}
\ifthenelse{\equal{#1}{23new}}{
			\axgridWnew}{}		
\ifthenelse{\equal{#1}{23a}}{
			\axgridWa}{}
\ifthenelse{\equal{#1}{24}}{
			\axgridX}{}
\ifthenelse{\equal{#1}{24new}}{
			\axgridXnew}{}
\ifthenelse{\equal{#1}{25}}{
			\axgridY}{}
\ifthenelse{\equal{#1}{26}}{
			\axgridZ}{}			
}

\def\comm#1{
\ifthenelse{\equal{#1}{1}}{
			\commA}{}
\ifthenelse{\equal{#1}{2}}{
			\commB}{}
\ifthenelse{\equal{#1}{3}}{
			\commC}{}		
\ifthenelse{\equal{#1}{4}}{
			\commD}{}
\ifthenelse{\equal{#1}{5}}{
			\commE}{}
\ifthenelse{\equal{#1}{6}}{
			\commF}{}
\ifthenelse{\equal{#1}{6a}}{
			\commFa}{}
\ifthenelse{\equal{#1}{7}}{
			\commG}{}
\ifthenelse{\equal{#1}{8}}{
			\commH}{}
\ifthenelse{\equal{#1}{8s}}{
			\commHs}{}		
\ifthenelse{\equal{#1}{8a}}{
			\commHa}{}
\ifthenelse{\equal{#1}{8b}}{
			\commHb}{}
\ifthenelse{\equal{#1}{11}}{
			\commK}{}
\ifthenelse{\equal{#1}{12}}{
			\commL}{}
\ifthenelse{\equal{#1}{13}}{
			\commM}{}
\ifthenelse{\equal{#1}{13a}}{
			\commMa}{}		
\ifthenelse{\equal{#1}{14}}{
			\commN}{}
\ifthenelse{\equal{#1}{14a}}{
			\commNa}{}
\ifthenelse{\equal{#1}{15}}{
			\commO}{}
\ifthenelse{\equal{#1}{16}}{
			\commP}{}
\ifthenelse{\equal{#1}{18}}{
			\commR}{}
\ifthenelse{\equal{#1}{18new}}{
			\commRnew}{}
\ifthenelse{\equal{#1}{19}}{
			\commS}{}		
\ifthenelse{\equal{#1}{19new}}{
			\commSnew}{}		
\ifthenelse{\equal{#1}{20}}{
			\commT}{}
\ifthenelse{\equal{#1}{20new}}{
			\commTnew}{}
\ifthenelse{\equal{#1}{20a}}{
			\commTa}{}
\ifthenelse{\equal{#1}{21}}{
			\commU}{}
\ifthenelse{\equal{#1}{21a}}{
			\commUa}{}
\ifthenelse{\equal{#1}{21new}}{
			\commUnew}{}
\ifthenelse{\equal{#1}{21anew}}{
			\commUanew}{}
\ifthenelse{\equal{#1}{22}}{
			\commV}{}
\ifthenelse{\equal{#1}{22new}}{
			\commVnew}{}
\ifthenelse{\equal{#1}{23}}{
			\commW}{}	
\ifthenelse{\equal{#1}{23new}}{
			\commWnew}{}		
\ifthenelse{\equal{#1}{23a}}{
			\commWa}{}
\ifthenelse{\equal{#1}{24}}{
			\commX}{}
\ifthenelse{\equal{#1}{24new}}{
			\commXnew}{}
\ifthenelse{\equal{#1}{25}}{
			\commY}{}
\ifthenelse{\equal{#1}{26}}{
			\commZ}{}			
}			

\def\labelout#1{
\ifthenelse{\equal{#1}{1}}{
			\labelA}{}
\ifthenelse{\equal{#1}{2}}{
			\labelB}{}
\ifthenelse{\equal{#1}{3}}{
			\labelC}{}		
\ifthenelse{\equal{#1}{4}}{
			\labelD}{}
\ifthenelse{\equal{#1}{5}}{
			\labelE}{}
\ifthenelse{\equal{#1}{6}}{
			\labelF}{}
\ifthenelse{\equal{#1}{6a}}{
			\labelFa}{}
\ifthenelse{\equal{#1}{7}}{
			\labelG}{}
\ifthenelse{\equal{#1}{8}}{
			\labelH}{}
\ifthenelse{\equal{#1}{8s}}{
			\labelHs}{}		
\ifthenelse{\equal{#1}{8a}}{
			\labelHa}{}
\ifthenelse{\equal{#1}{8b}}{
			\labelHb}{}
\ifthenelse{\equal{#1}{11}}{
			\labelK}{}
\ifthenelse{\equal{#1}{12}}{
			\labelL}{}
\ifthenelse{\equal{#1}{13}}{
			\labelM}{}
\ifthenelse{\equal{#1}{13a}}{
			\labelMa}{}		
\ifthenelse{\equal{#1}{14}}{
			\labelN}{}
\ifthenelse{\equal{#1}{14a}}{
			\labelNa}{}
\ifthenelse{\equal{#1}{15}}{
			\labelO}{}
\ifthenelse{\equal{#1}{16}}{
			\labelP}{}
\ifthenelse{\equal{#1}{18}}{
			\labelR}{}
\ifthenelse{\equal{#1}{18new}}{
			\labelRnew}{}
\ifthenelse{\equal{#1}{19}}{
			\labelS}{}		
\ifthenelse{\equal{#1}{19new}}{
			\labelSnew}{}
\ifthenelse{\equal{#1}{20}}{
			\labelT}{}
\ifthenelse{\equal{#1}{20new}}{
			\labelTnew}{}
\ifthenelse{\equal{#1}{20a}}{
			\labelTa}{}
\ifthenelse{\equal{#1}{21}}{
			\labelU}{}
\ifthenelse{\equal{#1}{21new}}{
			\labelUnew}{}
\ifthenelse{\equal{#1}{21a}}{
			\labelUa}{}
\ifthenelse{\equal{#1}{21anew}}{
			\labelUanew}{}
\ifthenelse{\equal{#1}{22}}{
			\labelV}{}
\ifthenelse{\equal{#1}{22new}}{
			\labelVnew}{}
\ifthenelse{\equal{#1}{23}}{
			\labelW}{}
\ifthenelse{\equal{#1}{23new}}{
			\labelWnew}{}		
\ifthenelse{\equal{#1}{23a}}{
			\labelWa}{}
\ifthenelse{\equal{#1}{24}}{
			\labelX}{}
\ifthenelse{\equal{#1}{24new}}{
			\labelXnew}{}
\ifthenelse{\equal{#1}{25}}{
			\labelY}{}
\ifthenelse{\equal{#1}{26}}{
			\labelZ}{}			
}

%% file: cnii-code.tex
\def\AxiomE#1{ %5
	\ifnum#1=1
		\sh(0,0);
		\bx(1,1)3;
		\sh(0,1);
		\bx(1,1)2;
		\sh(1+\gap,0);
		\bx(1,2)2;
\else \ifnum#1=2
		\sh(0,0);
		\bx(1,1)3;
		\sh(1,0);
		\bx(1,1)2;
\else \ifnum#1=3
		\sh(0,0);
		\bx(2,1)2;
\else \ifnum#1=4
		\sh(0,0);
		\bx(2,1)2;
		\sh(0,1);
		\bx(2,1)2;
\else \ifnum#1=5
		\sh(0,0);
		\bx(1,1)3;
		\sh(1,0);
		\bx(1,1)2;		
		\sh(0,1+\gap);
		\bx(1,1)2;
		\sh(1,1+\gap);
		\bx(1,1)2;		
\else \ifnum#1=6
		\sh(0,0);
		\bx(1,1)3;
		\sh(1,0);
		\bx(1,1)2;
		\sh(0,1+\gap);
		\bx(2,1)2;
\fi\fi\fi\fi\fi\fi}

\def\AxiomF#1{ %6
	\ifnum#1=1
		\sh(0,0);
		\bx(1,1)1;
		\sh(1,0);
		\bx(1,1)2;		
		\sh(0,1+\gap);
		\bx(1,1)2;
		\sh(1,1+\gap);
		\bx(1,1)2;
\else \ifnum#1=2
		\sh(0,0);
		\bx(1,1)1;
		\sh(1,0);
		\bx(1,1)2;		
		\sh(0,1+\gap);
		\bx(2,1)2;
\else \ifnum#1=3
		\sh(0,0);
		\bx(1,1)1;
		\sh(1,0);
		\bx(1,1)2;		
\else \ifnum#1=4
		\sh(0,0);
		\bx(1,1)2;
		\sh(1,0);
		\bx(1,1)2;		
		\sh(0,1+\gap);
		\bx(1,1)1;
		\sh(1,1+\gap);
		\bx(1,1)2;
\else \ifnum#1=5
		\sh(0,0);
		\bx(2,1)2;
		\sh(0,1+\gap);
		\bx(1,1)1;
		\sh(1,1+\gap);
		\bx(1,1)2;
\fi\fi\fi\fi\fi}

\def\AxiomFa#1{ %6a
	\ifnum#1=1
		\sh(0,0);
		\bx(1,1)3;
		\sh(1,0);
		\bx(1,1)2;		
		\sh(0,1+\gap);
		\bx(1,1)2;
		\sh(1,1+\gap);
		\bx(1,1)2;		
\else \ifnum#1=2
		\sh(0,0);
		\bx(1,1)3;
		\sh(1,0);
		\bx(1,1)2;
		\sh(0,1+\gap);
		\bx(2,1)2;
\else \ifnum#1=3
		\sh(0,0);
		\bx(2,1)2;
		\sh(0,1);
		\bx(2,1)2;
\else \ifnum#1=4
		\sh(0,0);
		\bx(2,1)2;
\else \ifnum#1=5
		\sh(0,0);
		\bx(1,1)3;
		\sh(1,0);
		\bx(1,1)2;
\fi\fi\fi\fi\fi}

\def\AxiomG#1{ %7
	\ifnum#1=1
		\sh(0,0);
		\bx(1,1)3;
		\sh(0,1);
		\qhexdl(1,1)1;
		\sh(0,2+\gap);
		\hkdl(1+.5*\gap,1*\fudge)2,1,3;
		\sh(1+\gap,0);
		\bx(1,2)2;
		\sh(1+.5*\gap,2+\gap);
		\bx(1+.5*\gap,1*\fudge)2;
\else \ifnum#1=2
	\def\hexangleCtmp{\hexangleC}
	\def\gapfactorCtmp{\gapfactorC}
	\def\hexangleAtmp{\hexangleA}
	\def\gapfactorAtmp{\gapfactorA}
	\def\hexangleC{55}
	\def\gapfactorC{1.64}
	\def\hexangleA{55}
	\def\gapfactorA{1.74}
		\sh(0,0);
		\dhku(1,1*\fudge)1,1,3;
		\sh(0,1*\fudge+\gap);
		\hex(1,1)1;
		\sh(0,1+1*\fudge+2*\gap);
		\dhkd(1,1*\fudge)1,1,3;
		\sh(1+\gap,0);
		\bx(1,2*\fudge+1+2*\gap)2;
	\def\hexangleC{\hexangleCtmp}
	\def\gapfactorC{\gapfactorCtmp}
	\def\hexangleA{\hexangleAtmp}
	\def\gapfactorA{\gapfactorAtmp}
\else \ifnum#1=3
		\sh(0,0);
		\bx(1,1)2;
		\sh(0,1);
		\qhexdl(1,1)1;
		\sh(0,2+\gap);
		\hkdl(2+\gap,1*\fudge)2,1,2;
		\sh(1+\gap,0);
		\bx(1,2)2;
\else \ifnum#1=4
		\sh(0,0);
		\bx(1,1)2;
		\sh(0,1);
		\bx(1,1)1;		
		\sh(1+\gap,0);
		\bx(1,2)2;
\else \ifnum#1=5
		\sh(0,0);
		\bx(1,1)2;
		\sh(1,0);
		\bx(1,1)2;		
		\sh(0,1+\gap);
		\bx(1,1)1;
		\sh(1,1+\gap);
		\bx(1,1)2;
\else \ifnum#1=6
		\sh(0,0);
		\bx(1,1)1;
		\sh(1,0);
		\bx(1,1)2;		
\else \ifnum#1=7
		\sh(0,0);
		\bx(2,1)2;
		\sh(0,1+\gap);
		\bx(1,1)1;
		\sh(1,1+\gap);
		\bx(1,1)2;
\else \ifnum#1=8
		\sh(0,0);
		\hkul(1,1*\fudge)1,1,2;
		\sh(1,0);
		\hkur(1,1*\fudge)1,1,2;
		\sh(0,1*\fudge+\gap);
		\halfhexl(1,1)1;
		\sh(1,1*\fudge+\gap);
		\halfhexr(1,1)2;
		\sh(0,1*\fudge+1+2*\gap);
		\dhkd(2,1*\fudge)1,1,2;
\else \ifnum#1=9
		\sh(0,0);
		\bx(1,1)1;
		\sh(0,1);
		\bx(1,1)2;		
		\sh(1+\gap,0);
		\bx(1,2)2;
\else \ifnum#1=10
		\sh(0,0);
		\bx(1,1)1;
		\sh(1,0);
		\bx(1,1)2;		
		\sh(0,1+\gap);
		\bx(2,1)2;
\else \ifnum#1=11
		\sh(0,0);
		\dhku(2,1*\fudge)1,1,2;
		\sh(0,1*\fudge+\gap);
		\halfhexl(1,1)1;
		\sh(1,1*\fudge+\gap);
		\halfhexr(1,1)2;
		\sh(0,1*\fudge+1+2*\gap);
		\dhkd(2,1*\fudge)1,1,2;
\else \ifnum#1=12
		\sh(0,0);
		\bx(1,1)1;
		\sh(1,0);
		\bx(1,1)2;		
		\sh(0,1+\gap);
		\bx(1,1)2;
		\sh(1,1+\gap);
		\bx(1,1)2;
\else \ifnum#1=13
		\sh(0,0);
		\dhku(2,1*\fudge)1,1,2;
		\sh(0,1*\fudge+\gap);
		\halfhexl(1,1)1;
		\sh(1,1*\fudge+\gap);
		\halfhexr(1,1)2;
		\sh(0,1*\fudge+1+2*\gap);
		\hkdl(1,1*\fudge)1,1,2;
		\sh(1,1*\fudge+1+2*\gap);
		\hkdr(1,1*\fudge)1,1,2;
\else \ifnum#1=14
		\sh(0,0);
		\hkul(2+\gap,1*\fudge)2,1,2;
		\sh(0,1*\fudge+\gap);
		\qhexul(1,1)1;
		\sh(0,1*\fudge+1+\gap);
		\bx(1,1)2;
		\sh(1+\gap,1*\fudge+\gap);
		\bx(1,2)2;
\else \ifnum#1=15
		\sh(0,0);
		\hkul(1,1*\fudge)1,1,3;
		\sh(0,1*\fudge+\gap);
		\halfhexl(1,1)1;
		\sh(0,1*\fudge+1+2*\gap);
		\hkdl(1,1*\fudge)1,1,3;
		\sh(1,0);
		\hkur(1,1*\fudge)1,1,2;
		\sh(1,1*\fudge+\gap);
		\halfhexr(1,1)2;
		\sh(1,1*\fudge+1+2*\gap);
		\hkdr(1,1*\fudge)1,1,2;
\else \ifnum#1=16
		\sh(0,0);
		\hkul(1+.5*\gap,1*\fudge)2,1,3;
		\sh(0,1*\fudge+\gap);
		\qhexul(1,1)1;
		\sh(0,1+1*\fudge+\gap);
		\bx(1,1)3;
		\sh(1+.5*\gap,0);
		\bx(1+.5*\gap,1*\fudge)2;
		\sh(1+\gap,1*\fudge+\gap);
		\bx(1,2)2;
\fi\fi\fi\fi\fi\fi\fi\fi\fi\fi\fi\fi\fi\fi\fi\fi}

\def\AxiomH#1{
	\ifnum#1=1
		\sh(0,0);
		\hkul(1,1*\fudge)1,1,1;
		\sh(0,1*\fudge+\gap);
		\halfhexl(1,1)2;
		\sh(0,1*\fudge+1+2*\gap);
		\hkdl(1,1*\fudge)1,1,1;
		\sh(1,0);
		\hkur(1,1*\fudge)1,1,2;
		\sh(1,1*\fudge+\gap);
		\halfhexr(1,1)1;
		\sh(1,1*\fudge+1+2*\gap);
		\hkdr(1,1*\fudge)1,1,2;
\else \ifnum#1=2
		\sh(0,0);
		\hkul(1,1*\fudge)1,1,2;
		\sh(0,1*\fudge+\gap);
		\halfhexl(1,1)1;
		\sh(0,1*\fudge+1+2*\gap);
		\hkdl(1,1*\fudge)1,1,1;
		\sh(1,0);
		\hkur(1,1*\fudge)1,1,1;
		\sh(1,1*\fudge+\gap);
		\halfhexr(1,1)2;
		\sh(1,1*\fudge+1+2*\gap);
		\hkdr(1,1*\fudge)1,1,2;
\else \ifnum#1=3
		\sh(0,0);
		\bx(1,1)1;
		\sh(1,0);
		\bx(1,1)1;		
		\sh(0,1+\gap);
		\bx(1,1)1;
		\sh(1,1+\gap);
		\bx(1,1)2;
\else \ifnum#1=4
		\sh(0,0);
		\bx(1,1)1;
		\sh(1,0);
		\bx(1,1)2;		
		\sh(0,1+\gap);
		\bx(1,1)1;
		\sh(1,1+\gap);
		\bx(1,1)1;
\else \ifnum#1=5
		\sh(0,0);
		\hkul(1,1*\fudge)1,1,1;
		\sh(0,1*\fudge+\gap);
		\halfhexl(1,1)1;
		\sh(0,1*\fudge+1+2*\gap);
		\hkdl(1,1*\fudge)1,1,2;
		\sh(1,0);
		\hkur(1,1*\fudge)1,1,2;
		\sh(1,1*\fudge+\gap);
		\halfhexr(1,1)2;
		\sh(1,1*\fudge+1+2*\gap);
		\hkdr(1,1*\fudge)1,1,1;
\else \ifnum#1=6
		\sh(0,0);
		\bx(1,1)1;
		\sh(1+\gap,0);
		\bx(1,2)1;		
		\sh(0,1);
		\bx(1,1)1;
\else \ifnum#1=7
		\sh(0,0);
		\hkul(1+.5*\gap,1*\fudge)2,1,2;
		\sh(0,1*\fudge+\gap);
		\qhexul(1,1)1;
		\sh(0,1+1*\fudge+\gap);
		\bx(1,1)1;
		\sh(1+.5*\gap,0);
		\bx(1+.5*\gap,1*\fudge)1;
		\sh(1+\gap,1*\fudge+\gap);
		\bx(1,2)2;
\else \ifnum#1=8
		\sh(0,0);
		\bx(1,1)1;
		\sh(0,1);
		\qhexdl(1,1)1;
		\sh(0,2+\gap);
		\hkdl(1+.5*\gap,1*\fudge)2,1,2;
		\sh(1+\gap,0);
		\bx(1,2)2;
		\sh(1+.5*\gap,2+\gap);
		\bx(1+.5*\gap,1*\fudge)1;
\fi\fi\fi\fi\fi\fi\fi\fi}

\def\AxiomHs#1{
	\ifnum#1=1
		\sh(0,0);
		\bx(1,1)1;
		\sh(1,0);
		\bx(1,1)1;		
		\sh(0,1+\gap);
		\bx(1,1)3;
		\sh(1,1+\gap);
		\bx(1,1)2;
\else \ifnum#1=2
		\sh(0,0);
		\bx(1,1)1;
		\sh(1,0);
		\bx(1,1)1;		
		\sh(0,1+\gap);
		\bx(2,1)2;
\else \ifnum#1=3
		\sh(0,0);
		\bx(1,1)1;
		\sh(1,0);
		\bx(1,1)1;		
\else \ifnum#1=4
		\sh(0,0);
		\bx(1,1)1;
		\sh(0,1);
		\bx(1,1)3;
		\sh(1+\gap,0);
		\bx(1,2)1;
\fi\fi\fi\fi}

\def\AxiomHa#1{
	\ifnum#1=1
		\sh(0,0);
		\hkul(1,1*\fudge)1,1,2;
		\sh(0,1*\fudge+\gap);
		\halfhexl(1,1)1;
		\sh(0,1*\fudge+1+2*\gap);
		\hkdl(1,1*\fudge)1,1,3;
		\sh(1,0);
		\hkur(1,1*\fudge)1,1,3;
		\sh(1,1*\fudge+\gap);
		\halfhexr(1,1)2;
		\sh(1,1*\fudge+1+2*\gap);
		\hkdr(1,1*\fudge)1,1,2;
\else \ifnum#1=2
		\sh(0,0);
		\bx(1,1)1;
		\sh(1,0);
		\bx(1,1)3;		
		\sh(0,1+\gap);
		\bx(1,1)3;
		\sh(1,1+\gap);
		\bx(1,1)2;
\else \ifnum#1=3
		\sh(0,0);
		\hkul(1,1*\fudge)1,1,2;
		\sh(1,0);
		\hkur(1,1*\fudge)1,1,2;
		\sh(0,1*\fudge+\gap);
		\halfhexl(1,1)1;
		\sh(1,1*\fudge+\gap);
		\halfhexr(1,1)2;
		\sh(0,1*\fudge+1+2*\gap);
		\dhkd(2,1*\fudge)1,1,2;
\else \ifnum#1=4
		\sh(0,0);
		\bx(1,1)1;
		\sh(1,0);
		\bx(1,1)2;		
		\sh(0,1+\gap);
		\bx(2,1)2;
\else \ifnum#1=5
		\sh(0,0);
		\bx(1,1)2;
		\sh(1,0);
		\bx(1,1)2;		
		\sh(0,1+\gap);
		\bx(1,1)1;
		\sh(1,1+\gap);
		\bx(1,1)2;
\else \ifnum#1=6
		\sh(0,0);
		\dhku(2,1*\fudge)1,1,2;
		\sh(0,1*\fudge+\gap);
		\halfhexl(1,1)1;
		\sh(1,1*\fudge+\gap);
		\halfhexr(1,1)2;
		\sh(0,1*\fudge+1+2*\gap);
		\dhkd(2,1*\fudge)1,1,2;
\else \ifnum#1=7
		\sh(0,0);
		\bx(2,1)2;
		\sh(0,1+\gap);
		\bx(1,1)1;
		\sh(1,1+\gap);
		\bx(1,1)2;
\else \ifnum#1=8
		\sh(0,0);
		\dhku(2,1*\fudge)1,1,2;
		\sh(0,1*\fudge+\gap);
		\halfhexl(1,1)1;
		\sh(1,1*\fudge+\gap);
		\halfhexr(1,1)2;
		\sh(0,1*\fudge+1+2*\gap);
		\hkdl(1,1*\fudge)1,1,2;
		\sh(1,1*\fudge+1+2*\gap);
		\hkdr(1,1*\fudge)1,1,2;
\else \ifnum#1=9
		\sh(0,0);
		\bx(1,1)1;
		\sh(1,0);
		\bx(1,1)2;		
		\sh(0,1+\gap);
		\bx(2,1)2;
\else \ifnum#1=10
		\sh(0,0);
		\bx(1,1)1;
		\sh(1,0);
		\bx(1,1)2;		
\else \ifnum#1=11
		\sh(0,0);
		\bx(1,1)1;
		\sh(1,0);
		\bx(1,1)2;		
		\sh(0,1+\gap);
		\bx(1,1)2;
		\sh(1,1+\gap);
		\bx(1,1)2;
\else \ifnum#1=12
		\sh(0,0);
		\hkul(2+\gap,1*\fudge)2,1,2;
		\sh(0,1*\fudge+\gap);
		\qhexul(1,1)1;
		\sh(0,1*\fudge+1+\gap);
		\bx(1,1)2;
		\sh(1+\gap,1*\fudge+\gap);
		\bx(1,2)2;
\else \ifnum#1=13
		\sh(0,0);
		\hkul(1+.5*\gap,1*\fudge)2,1,2;
		\sh(0,1*\fudge+\gap);
		\qhexul(1,1)1;
		\sh(0,1+1*\fudge+\gap);
		\bx(1,1)3;
		\sh(1+.5*\gap,0);
		\bx(1+.5*\gap,1*\fudge)3;
		\sh(1+\gap,1*\fudge+\gap);
		\bx(1,2)2;
\else \ifnum#1=14
		\sh(0,0);
		\bx(1,1)1;
		\sh(0,1);
		\bx(1,1)3;		
		\sh(1+\gap,0);
		\bx(1,2)3;
\fi\fi\fi\fi\fi\fi\fi\fi\fi\fi\fi\fi\fi\fi}

\def\AxiomHb#1{
	\ifnum#1=1
		\sh(0,0);
		\hkul(1,1*\fudge)1,1,3;
		\sh(0,1*\fudge+\gap);
		\halfhexl(1,1)2;
		\sh(0,1*\fudge+1+2*\gap);
		\hkdl(1,1*\fudge)1,1,3;
		\sh(1,0);
		\hkur(1,1*\fudge)1,1,2;
		\sh(1,1*\fudge+\gap);
		\halfhexr(1,1)1;
		\sh(1,1*\fudge+1+2*\gap);
		\hkdr(1,1*\fudge)1,1,2;
\else \ifnum#1=2
		\sh(0,0);
		\dhku(2,1*\fudge)1,1,2;
		\sh(0,1*\fudge+\gap);
		\halfhexl(1,1)2;
		\sh(1,1*\fudge+\gap);
		\halfhexr(1,1)1;
		\sh(0,1*\fudge+1+2*\gap);
		\hkdl(1,1*\fudge)1,1,2;
		\sh(1,1*\fudge+1+2*\gap);
		\hkdr(1,1*\fudge)1,1,2;
\else \ifnum#1=3
		\sh(0,0);
		\hkul(1,1*\fudge)1,1,2;
		\sh(1,0);
		\hkur(1,1*\fudge)1,1,2;
		\sh(0,1*\fudge+\gap);
		\halfhexl(1,1)2;
		\sh(1,1*\fudge+\gap);
		\halfhexr(1,1)1;
		\sh(0,1*\fudge+1+2*\gap);
		\dhkd(2,1*\fudge)1,1,2;
\else \ifnum#1=4
		\sh(0,0);
		\bx(1,1)3;
		\sh(1,0);
		\bx(1,1)1;		
		\sh(0,1+\gap);
		\bx(1,1)3;
		\sh(1,1+\gap);
		\bx(1,1)2;
\else \ifnum#1=5
		\sh(0,0);
		\dhku(2,1*\fudge)1,1,2;
		\sh(0,1*\fudge+\gap);
		\halfhexl(1,1)2;
		\sh(1,1*\fudge+\gap);
		\halfhexr(1,1)1;
		\sh(0,1*\fudge+1+2*\gap);
		\dhkd(2,1*\fudge)1,1,2;
\else \ifnum#1=6
		\sh(0,0);
		\bx(1,1)3;
		\sh(1,0);
		\bx(1,1)2;		
		\sh(0,1+\gap);
		\bx(1,1)3;
		\sh(1,1+\gap);
		\bx(1,1)1;
\else \ifnum#1=7
		\sh(0,0);
		\bx(1,1)2;
		\sh(1,0);
		\bx(1,1)1;		
		\sh(0,1+\gap);
		\bx(2,1)2;
\else \ifnum#1=8
		\sh(0,0);
		\bx(2,1)2;
		\sh(0,1+\gap);
		\bx(1,1)2;
		\sh(1,1+\gap);
		\bx(1,1)1;
\else \ifnum#1=9
		\sh(0,0);
		\bx(1,1)2;
		\sh(1,0);
		\bx(1,1)1;		
\else \ifnum#1=10
		\sh(0,0);
		\bx(1,1)3;
		\sh(0,1);
		\bx(1,1)3;		
		\sh(1+\gap,0);
		\bx(1,2)1;
\fi\fi\fi\fi\fi\fi\fi\fi\fi\fi}

\def\AxiomK#1{
	\ifnum#1=1
		\sh(0,0);
		\bx(1,1)3;
		\sh(1,0);
		\bx(1,1)2;		
		\sh(2,0);
		\bx(1,1)2;
		\sh(0,1+\gap);
		\bx(1,1)3;		
		\sh(1,1+\gap);
		\bx(1,1)2;
		\sh(2,1+\gap);
		\bx(1,1)2;		
\else \ifnum#1=2
		\sh(0,0);
		\bx(2,1)2;
		\sh(2,0);
		\bx(1,1)2;		
		\sh(0,1+\gap);
		\bx(2,1)2;
		\sh(2,1+\gap);
		\bx(1,1)2;		
\else \ifnum#1=3
		\sh(0,0);
		\bx(1,1)3;
		\sh(1,0);
		\bx(2,1)2;		
		\sh(0,1+\gap);
		\bx(1,1)3;
		\sh(1,1+\gap);
		\bx(2,1)2;		
\else \ifnum#1=4
		\sh(0,0);
		\bx(3,1)2;
		\sh(0,1);
		\bx(3,1)2;		
\else \ifnum#1=5
		\sh(0,0);
		\bx(2,1)2;
		\sh(0,1);
		\bx(2,1)2;		
		\sh(2+\gap,0);
		\bx(1,2)2;
\else \ifnum#1=6
		\sh(0,0);
		\bx(1,1)3;
		\sh(1,0);
		\bx(1,1)2;		
		\sh(0,1+\gap);
		\bx(1,1)3;
		\sh(1,1+\gap);
		\bx(1,1)2;		
		\sh(2+\gap,0);
		\bx(1,2+\gap)2;
\else \ifnum#1=7
		\sh(0,0);
		\bx(3,1)2;
\else \ifnum#1=8
		\sh(0,0);
		\bx(2,1)2;
		\sh(2,0);
		\bx(1,1)2;		
\else \ifnum#1=9
		\sh(0,0);
		\bx(1,1)2;
		\sh(1,0);
		\bx(2,1)2;		
\else \ifnum#1=10
		\sh(0,0);
		\bx(1,1)2;
		\sh(1,0);
		\bx(1,1)2;		
		\sh(2,0);
		\bx(1,1)2;
\else \ifnum#1=11
		\sh(0,0);
		\bx(1,1)3;
		\sh(0,1);
		\bx(1,1)3;		
		\sh(1+\gap,0);
		\bx(2,2)2;
\else \ifnum#1=12
		\sh(0,0);
		\bx(1,1)3;
		\sh(0,1);
		\bx(1,1)3;		
		\sh(1+\gap,0);
		\bx(1,2)2;
		\sh(2+2*\gap,0);
		\bx(1,2)2;
\fi\fi\fi\fi\fi\fi\fi\fi\fi\fi\fi\fi}

\def\AxiomL#1{
	\ifnum#1=1
		\sh(0,0);
		\bx(1,2)2;
		\sh(1+\gap,0);
		\bx(1,1)3;		
		\sh(1+\gap,1);
		\bx(1,1)3;
		\sh(2+2*\gap,0);
		\bx(1,2)2;		
\else \ifnum#1=2
		\sh(0,0);
		\bx(1,1)2;
		\sh(1,0);
		\bx(1,1)3;		
		\sh(0,1+\gap);
		\bx(1,1)2;
		\sh(1,1+\gap);
		\bx(1,1)3;		
		\sh(2+\gap,0);
		\bx(1,2+\gap)2;
\else \ifnum#1=3
		\sh(0,0);
		\bx(1,1)2;
		\sh(1,0);
		\bx(1,1)2;		
		\sh(2,0);
		\bx(1,1)2;
\else \ifnum#1=4
		\sh(0,0);
		\bx(2,1)2;
		\sh(2,0);
		\bx(1,1)2;		
\else \ifnum#1=5
		\sh(0,0);
		\bx(2,1)2;
		\sh(0,1);
		\bx(2,1)2;		
		\sh(2+\gap,0);
		\bx(1,2)2;
\else \ifnum#1=6
		\sh(0,0);
		\bx(1,1)2;
		\sh(1,0);
		\bx(2,1)2;		
\else \ifnum#1=7
		\sh(0,0);
		\bx(3,1)2;
\else \ifnum#1=8
		\sh(0,0);
		\bx(3,1)2;
		\sh(0,1);
		\bx(3,1)2;		
\else \ifnum#1=9
		\sh(0,0);
		\bx(2,1)2;
		\sh(2,0);
		\bx(1,1)2;		
		\sh(0,1+\gap);
		\bx(2,1)2;
		\sh(2,1+\gap);
		\bx(1,1)2;		
\else \ifnum#1=10
		\sh(0,0);
		\bx(1,2)2;
		\sh(1+\gap,0);
		\bx(2,1)2;		
		\sh(1+\gap,1);
		\bx(2,1)2;
\else \ifnum#1=11
		\sh(0,0);
		\bx(1,1)2;
		\sh(1,0);
		\bx(2,1)2;		
		\sh(0,1+\gap);
		\bx(1,1)2;
		\sh(1,1+\gap);
		\bx(2,1)2;		
\else \ifnum#1=12
		\sh(0,0);
		\bx(1,2+\gap)2;
		\sh(1+\gap,0);
		\bx(1,1)3;		
		\sh(2+\gap,0);
		\bx(1,1)2;
		\sh(1+\gap,1+\gap);
		\bx(1,1)3;		
		\sh(2+\gap,1+\gap);
		\bx(1,1)2;
\else \ifnum#1=13
		\sh(0,0);
		\bx(1,1)2;
		\sh(1,0);
		\bx(1,1)3;		
		\sh(2,0);
		\bx(1,1)2;
		\sh(0,1+\gap);
		\bx(1,1)2;		
		\sh(1,1+\gap);
		\bx(1,1)3;
		\sh(2,1+\gap);
		\bx(1,1)2;		
\fi\fi\fi\fi\fi\fi\fi\fi\fi\fi\fi\fi\fi}

\def\AxiomM#1{
	\ifnum#1=1
		\sh(0,0);
		\bx(1,1)2;
		\sh(1,0);
		\bx(1,1)1;		
		\sh(2,0);
		\bx(1,1)2;
		\sh(0,1+\gap);
		\bx(1,1)1;		
		\sh(1,1+\gap);
		\bx(1,1)2;
		\sh(2,1+\gap);
		\bx(1,1)2;		
\else \ifnum#1=2
		\sh(0,0);
		\bx(1,1)2;
		\sh(1,0);
		\bx(1,1)1;		
		\sh(0,1+\gap);
		\bx(1,1)1;
		\sh(1,1+\gap);
		\bx(1,1)2;		
		\sh(2+\gap,0);
		\bx(1,2+\gap)2;
\else \ifnum#1=3
		\sh(0,0);
		\bx(1,1)2;
		\sh(1,0);
		\bx(1,1)1;		
		\sh(2,0);
		\bx(1,1)2;
		\sh(0,1+\gap);
		\bx(1,1)1;		
		\sh(1,1+\gap);
		\bx(2,1)2;
\else \ifnum#1=4
		\sh(0,0);
		\bx(1,1)1;
		\sh(1,0);
		\bx(1,1)1;		
		\sh(2,0);
		\bx(1,1)2;
\else \ifnum#1=5
		\sh(0,0);
		\bx(1,1)1;
		\sh(1,0);
		\bx(1,1)2;		
		\sh(0,1+\gap);
		\bx(1,1)2;
		\sh(1,1+\gap);
		\bx(1,1)1;		
		\sh(2+\gap,0);
		\bx(1,2+\gap)2;
\else \ifnum#1=6
		\sh(0,0);
		\bx(1,1)1;
		\sh(1,0);
		\bx(2,1)2;		
		\sh(0,1+\gap);
		\bx(1,1)2;		
		\sh(1,1+\gap);
		\bx(1,1)1;
		\sh(2,1+\gap);
		\bx(1,1)2;		
\else \ifnum#1=7
		\sh(0,0);
		\bx(1,1)1;
		\sh(1,0);
		\bx(1,1)2;		
		\sh(2,0);
		\bx(1,1)2;
		\sh(0,1+\gap);
		\bx(1,1)2;		
		\sh(1,1+\gap);
		\bx(1,1)1;
		\sh(2,1+\gap);
		\bx(1,1)2;		
\fi\fi\fi\fi\fi\fi\fi}

\def\AxiomMa#1{
	\ifnum#1=1
		\sh(0,0);
		\bx(1,1)2;
		\sh(1,0);
		\bx(1,1)1;		
		\sh(2,0);
		\bx(1,1)2;
		\sh(0,1+\gap);
		\bx(1,1)3;		
		\sh(1,1+\gap);
		\bx(1,1)2;
		\sh(2,1+\gap);
		\bx(1,1)2;		
\else \ifnum#1=2
		\sh(0,0);
		\bx(1,1)2;
		\sh(1,0);
		\bx(1,1)1;		
		\sh(2,0);
		\bx(1,1)2;
		\sh(0,1+\gap);
		\bx(2,1)2;		
		\sh(2,1+\gap);
		\bx(1,1)2;		
\else \ifnum#1=3
		\sh(0,0);
		\bx(1,1)2;
		\sh(1,0);
		\bx(1,1)1;		
		\sh(0,1+\gap);
		\bx(1,1)3;
		\sh(1,1+\gap);
		\bx(1,1)2;		
		\sh(2+\gap,0);
		\bx(1,2+\gap)2;
\else \ifnum#1=4
		\sh(0,0);
		\bx(1,1)2;
		\sh(1,0);
		\bx(1,1)1;		
		\sh(0,1+\gap);
		\bx(2,1)2;
		\sh(2+\gap,0);
		\bx(1,2+\gap)2;
\else \ifnum#1=5
		\sh(0,0);
		\bx(1,1)2;
		\sh(1,0);
		\bx(1,1)1;		
		\sh(2,0);
		\bx(1,1)2;
		\sh(0,1+\gap);
		\bx(3,1)2;		
\else \ifnum#1=6
		\sh(0,0);
		\bx(1,1)2;
		\sh(1,0);
		\bx(1,1)1;		
		\sh(2,0);
		\bx(1,1)2;
		\sh(0,1+\gap);
		\bx(1,1)3;		
		\sh(1,1+\gap);
		\bx(2,1)2;
\else \ifnum#1=7
		\sh(0,0);
		\bx(1,1)3;
		\sh(1,0);
		\bx(1,1)1;		
		\sh(2,0);
		\bx(1,1)2;
\fi\fi\fi\fi\fi\fi\fi}

\def\AxiomN#1{
	\ifnum#1=1
		\sh(0,0);
		\bx(2,1)2;
		\sh(2,0);
		\bx(1,1)1;
		\sh(0,1+\gap);
		\bx(1,1)1;		
		\sh(1,1+\gap);
		\bx(1,1)2;
		\sh(2,1+\gap);
		\bx(1,1)2;		
\else \ifnum#1=2
		\sh(0,0);
		\bx(1,1)2;
		\sh(1,0);
		\bx(1,1)2;		
		\sh(2,0);
		\bx(1,1)1;
		\sh(0,1+\gap);
		\bx(1,1)1;		
		\sh(1,1+\gap);
		\bx(2,1)2;
\else \ifnum#1=3
		\sh(0,0);
		\bx(1,1)1;
		\sh(1,0);
		\bx(1,1)2;		
		\sh(2,0);
		\bx(1,1)1;
\else \ifnum#1=4
		\sh(0,0);
		\bx(1,1)1;
		\sh(1,0);
		\bx(1,1)2;		
		\sh(2,0);
		\bx(1,1)2;
		\sh(0,1+\gap);
		\bx(2,1)2;		
		\sh(2,1+\gap);
		\bx(1,1)1;		
\else \ifnum#1=5
		\sh(0,0);
		\bx(1,1)1;
		\sh(1,0);
		\bx(2,1)2;		
		\sh(0,1+\gap);
		\bx(1,1)2;		
		\sh(1,1+\gap);
		\bx(1,1)2;
		\sh(2,1+\gap);
		\bx(1,1)1;		
\fi\fi\fi\fi\fi}

\def\AxiomNa#1{
	\ifnum#1=1
		\sh(0,0);
		\bx(2,1)2;
		\sh(2,0);
		\bx(1,1)3;
		\sh(0,1+\gap);
		\bx(1,1)3;		
		\sh(1,1+\gap);
		\bx(1,1)2;
		\sh(2,1+\gap);
		\bx(1,1)2;		
\else \ifnum#1=2
		\sh(0,0);
		\bx(1,1)2;
		\sh(1,0);
		\bx(1,1)2;		
		\sh(2,0);
		\bx(1,1)3;
		\sh(0,1+\gap);
		\bx(1,1)3;		
		\sh(1,1+\gap);
		\bx(2,1)2;
\else \ifnum#1=3
		\sh(0,0);
		\bx(2,1)2;
		\sh(2,0);
		\bx(1,1)2;
		\sh(0,1+\gap);
		\bx(1,1)2;		
		\sh(1,1+\gap);
		\bx(2,1)2;
\else \ifnum#1=4
		\sh(0,0);
		\bx(3,1)2;
		\sh(0,1+\gap);
		\bx(1,1)2;		
		\sh(1,1+\gap);
		\bx(1,1)2;
		\sh(2,1+\gap);
		\bx(1,1)2;		
\else \ifnum#1=5
		\sh(0,0);
		\bx(3,1)2;
		\sh(0,1+\gap);
		\bx(1,1)2;		
		\sh(1,1+\gap);
		\bx(2,1)2;
\else \ifnum#1=6
		\sh(0,0);
		\bx(3,1)2;
		\sh(0,1+\gap);
		\bx(3,1)2;		
\else \ifnum#1=7
		\sh(0,0);
		\bx(2,1)2;
		\sh(2,0);
		\bx(1,1)2;
		\sh(0,1+\gap);
		\bx(3,1)2;		
\else \ifnum#1=8
		\sh(0,0);
		\bx(1,1)2;
		\sh(1,0);
		\bx(1,1)2;		
		\sh(2,0);
		\bx(1,1)2;
		\sh(0,1+\gap);
		\bx(3,1)2;		
\else \ifnum#1=9
		\sh(0,0);
		\bx(3,1)2;
\else \ifnum#1=10
		\sh(0,0);
		\bx(1,1)2;
		\sh(1,0);
		\bx(2,1)2;		
\else \ifnum#1=11
		\sh(0,0);
		\bx(2,1)2;
		\sh(2,0);
		\bx(1,1)2;
\else \ifnum#1=12
		\sh(0,0);
		\bx(1,1)3;
		\sh(1,0);
		\bx(1,1)2;		
		\sh(2,0);
		\bx(1,1)3;
\fi\fi\fi\fi\fi\fi\fi\fi\fi\fi\fi\fi}

\def\AxiomO#1{
	\ifnum#1=1
		\sh(0,0);
		\bx(1,1)1;
		\sh(1,0);
		\bx(1,1)2;		
		\sh(2,0);
		\bx(1,1)2;
		\sh(3,0);
		\bx(1,1)2;
\else \ifnum#1=2
		\sh(0,0);
		\bx(1,1)1;
		\sh(1,0);
		\bx(1,1)2;		
		\sh(2,0);
		\bx(2,1)2;
\else \ifnum#1=3
		\sh(0,0);
		\bx(1,1)1;
		\sh(1,0);
		\bx(2,1)2;		
		\sh(3,0);
		\bx(1,1)2;
\else \ifnum#1=4
		\sh(0,0);
		\bx(1,1)1;
		\sh(1,0);
		\bx(3,1)2;		
\fi\fi\fi\fi}

\def\AxiomP#1{
	\ifnum#1=1
		\sh(0,0);
		\bx(1,1)2;
		\sh(1,0);
		\bx(1,1)2;		
		\sh(2,0);
		\bx(1,1)2;
		\sh(3,0);
		\bx(1,1)2;
\else \ifnum#1=2
		\sh(0,0);
		\bx(1,1)2;
		\sh(1,0);
		\bx(1,1)2;		
		\sh(2,0);
		\bx(2,1)2;
\fi\fi}

\def\AxiomR#1{ %18
	\ifnum#1=1
		\def\x{0}\def\y{0}
		\boxt(.75,1,.25,1)2;
		\def\x{0}\def\y{1}		
		\bubble(.75,1,.25)2;
		\def\x{0}\def\y{1+\gap}
		\boxslopedr(.75,.25,1)2;
		\def\x{1}\def\y{1+\gap}
		\boxslopedl(.25,.75,1)2;
	\begin{pgfonlayer}{foreground}
		\bull(.5,1);\bull(.5,1+\gap);
		\bull(1.25,1-\cw);  
		\bull(1.5,1+\cw);\bull(1.5,1+\cw+\gap);
		\bull(1,1.5+\gap);
		\bull(1,1+\cw);
	\end{pgfonlayer}
		
\else \ifnum#1=2
		\def\x{0}\def\y{0}
		\bxbt(1.5,1)2;
		\def\x{0}\def\y{1}
		\boxcornur(.5,1,.5,.5)2;
		\cornerur(.5,1,.5,.5)2;
	\begin{pgfonlayer}{foreground}
		\bull(.75,2+\cw);
	\end{pgfonlayer}
\else \ifnum#1=3
	\def\x{0}\def\y{0}
	\boxslopedr(.5,.5,1)2;
	\def\x{1}\def\y{\cw}
	\boxsloperd(1,.75,.25-\cw)2;
	\def\x{0}\def\y{0}
	\cornerdr(.5,1.5,.75,.25)2;
	\begin{pgfonlayer}{foreground}
		\bull(0,.5);
		\bull(1,.5);
		\bull(2,.5);
	\end{pgfonlayer}
\else \ifnum#1=4%
	\def\x{0}\def\y{0}
	\boxcorndr(.5,1,.75,.25)2;
	\cornerdr(.5,1,.75,.25)2;
	\begin{pgfonlayer}{foreground}
		\bull(0,.5);\bull(1.5,.5);\bull(.75,\cw);
	\end{pgfonlayer}

\else \ifnum#1=5%!% Cuspstretched.  (Would look better with .05 bullets.)
		
	\sh(0,0)
	\sbx(1.5,1+\cw)2;
	\draw[fill=\mycolor{2},even odd rule]
		{[shift={(0,0)}]
			\cusp(0,1,1,.5+.5*\cuspwidth+.5*\cuspwidth*\cuspstretch)}
		{[shift={(1+.5+.5*\cuspwidth+.5*\cuspwidth*\cuspstretch,-.5*\cuspwidth)}] 
			-- (\x,\y) -- ++(90:.52-.5*\cusplength) .. controls +(90:\cuspfactor*\cusplength) and +(180+90:\cuspfactor*\cusplength) .. ++($ (90:\cusplength) + (90+1*90:.5*\cuspwidth*\cuspstretch) $) -- ++ (90:.15+.5*\cuspwidth-.5*\cusplength)	}		
		{[shift={(1+.5+.5*\cuspwidth,.15+.52)}] 
			-- \cusp(270,1,.15,.52-.5*\cuspwidth)}
		{[shift={(1+.5,.5*\cuspwidth)}] 
			-- \cusp(180,-1,.5,.5)} -- cycle;
	\sh(1.5,\cw);
	\bubblev(.5*\cusplength,1-\cusplength,.5*\cusplength)2;
	\begin{pgfonlayer}{foreground}
		\bull(1.5+\cw+\cw*\cuspstretch,-\cw);\bull(1.5,\cw);
		\bull(1.5,1+\cw);
		\bull(1.5-\cw,.5+\cw);
		\bull(1.5+\cw,.25+\cw);\bull(1.5+\cw,.5+\cw);\bull(1.5+\cw,.75+\cw);
		\bull(0,0);\bull(0,.5+.5*\cw);\bull(0,1+\cw);
	\end{pgfonlayer}

\else \ifnum#1=6%
	\def\x{0}\def\y{0}
	\boxt(.67,1,.33,1)2;
	\def\x{0}\def\y{1}
	\bubble(.67,1,.33)2;
	\boxb(.67,1,.33,1)2;	
	\begin{pgfonlayer}{foreground}
		\bull(.5,1);\bull(1.25,1-\cw);\bull(1,1+\cw);\bull(1.5,1+\cw);
		\bull(1,2);
	\end{pgfonlayer}
	
\else \ifnum#1=7%
	\def\x{0}\def\y{0}
	\boxcornur(.5,1,.75,.25)2;
	\cornerur(.5,1,.75,.25)2;
	\begin{pgfonlayer}{foreground}
		\bull(0,.5);\bull(1.5,.5);\bull(.75,1+\cw);
	\end{pgfonlayer}

\fi\fi\fi\fi\fi\fi\fi}

\def\AxiomRnew#1{ %18new
	\ifnum#1=1
		\sh(0,0);
		\boxt(.875,.5,.625,1)2;
		\sh(0,1);
		\bubble(.875,.5,.625)2;
		\sh(0,1+\gap);
		\boxslopedr(.875,.125,1)2;
		\sh(1,1+\gap);
		\boxslopedl(.625,.375,1)2;
		\begin{pgfonlayer}{foreground}
			\bull(.5,1);\bull(1.5,1);\bull(.5,1+\gap);
                        \bull(1,1+\cw);
			\bull(1.5,1+\gap);
                        \bull(1,1.5+\gap);
                        \dfct(.5,2+\gap);\dfct(.75,1+\gap);\dfct(.75,1);
                        \dfct(1,0);\idfct(1.25,1+\cw);\idfct(1.25,1+\gap+\cw);\idfct(1.5,2+\gap);
		\end{pgfonlayer}
\else \ifnum#1=2
  	\sh(0,0);
	\bx(1.5,1)2;
	\sh(0,1);
	\boxt(.75*.75,.75,.25*.75,1)2;
	\sh(0,2);
	\bubble(.75*.75,.75,.25*.75)2;
		\begin{pgfonlayer}{foreground}
			\bull(0,0);\bull(1.5,0);\bull(.5,1);\bull(1,1);\bull(.75,2+\cw);
			\dfct(.375,2);\dfct(.75,0);\dfct(.75,1);\idfct(1.125,2+\cw);
		\end{pgfonlayer}
\else \ifnum#1=3
		\sh(0,0);
		\bubble(.875,.5,.625)2;
		\boxslopedr(.875,.125,1)2;
		\sh(1,0);
		\boxslopedl(.625,.375,1)2;
		\begin{pgfonlayer}{foreground}
			\bull(.5,0);
			\bull(1.5,0);
                        \bull(1,.5);
                        \dfct(.5,1);\dfct(.75,0);
                        \idfct(1.25,\cw);\idfct(1.5,1);
		\end{pgfonlayer}
\else \ifnum#1=4%
		\sh(0,0);
		\bubble(.875,.5,.625)2;
		\boxb(.875,.5,.625,1)2;
		\begin{pgfonlayer}{foreground}
			\bull(.5,0); \bull(1,1);\bull(1,\cw);
			\bull(1.5,0);
                        \dfct(.5,1);\dfct(.75,0);
                        \idfct(1.25,\cw);\idfct(1.5,1);
		\end{pgfonlayer}
\else \ifnum#1=5%
		\sh(0,0);
		\boxt(.875,.5,.625,1)2;
		\sh(0,1);
		\bubble(.875,.5,.625)2;
		\sh(0,1+\gap);
		\boxb(.875,.5,.625,1)2;
		\begin{pgfonlayer}{foreground}
			\bull(.5,1);\bull(1.5,1);\bull(.5,1+\gap);
                        \bull(1,1+\cw);
			\bull(1.5,1+\gap);
                        \bull(1,1+\gap+\cw);
                        \bull(1,2+\gap);
                        \dfct(.5,2+\gap);\dfct(.75,1+\gap);\dfct(.75,1);
                        \dfct(1,0);\idfct(1.25,1+\cw);\idfct(1.25,1+\gap+\cw);\idfct(1.5,2+\gap);
		\end{pgfonlayer}
\else \ifnum#1=6%
	\sh(0,0);
	\boxt(.75*.75,.75,.25*.75,1)2;
	\sh(0,1);
	\bubble(.75*.75,.75,.25*.75)2;
		\begin{pgfonlayer}{foreground}
			\bull(.75,1+\cw);\bull(0,.5);\bull(1.5,.5);
			\dfct(.75,0);\dfct(.375,1);\idfct(1.125,1+\cw);
		\end{pgfonlayer}
	
\fi\fi\fi\fi\fi\fi}

\def\AxiomS#1{ %19

	\def\ticklengthtmp{\ticklength}
	\pgfmathsetmacro{\ticklengthnew}{\tickfactor*\ticklengthtmp} 
	\def\ticklength{\ticklengthnew}

	\ifnum#1=1
		\sh(0,0);
		\dhku(2,1*\fudge)1,1,1;
		\sh(0,1*\fudge+\gap);
		\hex(2,1)2;
		\halfidr(2,1){.67};
		\sh(0,1+1*\fudge+2*\gap);
		\hkdl(1,1*\fudge)1,1,3;
		\sh(1,1+1*\fudge+2*\gap);
		\hkdr(1,1*\fudge)1,1,2;
		\tick(1,1+1*\fudge+\gap);
		
\else \ifnum#1=2
		\sh(0,0);
		\dhku(2,1*\fudge)1,1,2;
		\halfidr(2,1*\fudge){.67};
		\sh(0,1*\fudge+\gap);
		\halfhexl(1,1)1;
		\sh(1,1*\fudge+\gap);
		\halfhexr(1,1)2;
		\sh(0,1*\fudge+1+2*\gap);
		\hkdl(1,1*\fudge)1,1,3;
		\sh(1,1*\fudge+1+2*\gap);
		\hkdr(1,1*\fudge)1,1,2;
		\tick(1,1*\fudge);
\else \ifnum#1=3
		\sh(0,0);
		\dhku(2,1*\fudge)1,1,1;
		\sh(0,1*\fudge+\gap);
		\hex(2,1)2;
		\halfidr(2,1){.67};
		\sh(0,1+1*\fudge+2*\gap);
		\dhkd(2,1*\fudge)1,1,2;
		\tick(1,1+1*\fudge+\gap);\tick(1,1+1*\fudge+2*\gap);\tick(1,1+2*\fudge+2*\gap);
		
\else \ifnum#1=4
		\sh(0,0);
		\dhku(2,1*\fudge)1,1,2;
		\halfidr(2,1*\fudge){.67};
		\sh(0,1*\fudge+\gap);
		\halfhexl(1,1)1;
		\sh(1,1*\fudge+\gap);
		\halfhexr(1,1)2;
		\sh(0,1*\fudge+1+2*\gap);
		\dhkd(2,1*\fudge)1,1,2;
		\tick(1,1+1*\fudge+2*\gap);\tick(1,1+2*\fudge+2*\gap);\tick(1,1*\fudge);		
		
\else \ifnum#1=5
		\sh(0,0);
		\bx(2,1)1;
		\sh(0,1);
		\bx(2,1)2;
		\halfidr(2,1){.75};
		\tick(1,2);
\else \ifnum#1=6
		\sh(0,0);
		\bx(2,1)2;
		\halfidr(2,1){.75};
		\sh(0,1+\gap);
		\bx(1,1)1;
		\sh(1,1+\gap);
		\bx(1,1)2;
		\tick(1,1);
\else \ifnum#1=7
		\sh(0,0);
		\dhku(2,1*\fudge)1,1,1;
		\sh(0,1*\fudge+\gap);
		\hex(2,1)3;
		\sh(0,1+1*\fudge+2*\gap);
		\dhkd(2,1*\fudge)1,1,2;
		\sh(0,1+1*\fudge+2*\gap);
		\halfidr(2,1*\fudge){.75};
		\tick(1,1+2*\fudge+2*\gap)
\else \ifnum#1=8
		\sh(0,0);
		\hkul(2+\gap,1*\fudge)2,1,2;
		\sh(0,0);
		\halfidr(2,1*\fudge){.85};
		\sh(0,1*\fudge+\gap);
		\qhexul(1,1)1;
		\sh(0,1*\fudge+1+\gap);
		\bx(1,1)3;
		\sh(1+\gap,1*\fudge+\gap);
		\bx(1,2)2;
		\tick(1+.5*\gap,1*\fudge);
		
\def\ticklength{\ticklengthtmp}
\fi\fi\fi\fi\fi\fi\fi\fi}

\def\AxiomSnew#1{ %19new

	\ifnum#1=1
		\sh(0,0);
		\dhku(2,1*\fudge)1,1,1;
		\sh(0,1*\fudge+\gap);
		\hex(2,1)2;
		\halfidr(2,1){.67};
		\sh(0,1+1*\fudge+2*\gap);
		\hkdl(1,1*\fudge)1,1,1;
		\sh(1,1+1*\fudge+2*\gap);
		\hkdr(1,1*\fudge)1,1,2;
		\begin{pgfonlayer}{foreground}
		\end{pgfonlayer}

\else \ifnum#1=2
		\sh(0,0);
		\dhku(2,1*\fudge)1,1,2;
		\halfidr(2,1*\fudge){.67};
		\sh(0,1*\fudge+\gap);
		\halfhexl(1,1)1;
		\sh(1,1*\fudge+\gap);
		\halfhexr(1,1)2;
		\sh(0,1*\fudge+1+2*\gap);
		\hkdl(1,1*\fudge)1,1,1;
		\sh(1,1*\fudge+1+2*\gap);
		\hkdr(1,1*\fudge)1,1,2;
		\tick(1,1*\fudge);
\else \ifnum#1=3
		\sh(0,0);
		\boxt(.875,.5,.625,1)1;
		\sh(0,1+\gap);
		\bubble(.875,.5,.625)2;
		\boxslopedr(.875,.125,1)1;
		\sh(1,1+\gap);
		\boxslopedl(.625,.375,1)2;
		\begin{pgfonlayer}{foreground}
			\bull(.5,1);\bull(1.5,1);\bull(.5,1+\gap);
			\bull(1.5,1+\gap);
                        \bull(1,1.5+\gap);
                        \dfct(.5,2+\gap);\dfct(.75,1+\gap);\dfct(.75,1);
                        \dfct(1,0);\idfct(1.25,1+\gap+\cw);\idfct(1.5,2+\gap);
		\end{pgfonlayer}
\else \ifnum#1=4
  	\sh(0,0);
	\bx(1.5,1)1;
	\sh(0,1);
	\boxt(.75*.75,.75,.25*.75,1)1;
	\sh(0,2);
	\bubble(.75*.75,.75,.25*.75)1;
		\begin{pgfonlayer}{foreground}
			\bull(0,0);\bull(1.5,0);\bull(.5,1);\bull(1,1);\bull(.75,2+\cw);
			\dfct(.375,2);\dfct(.75,0);\dfct(.75,1);\idfct(1.125,2+\cw);
		\end{pgfonlayer}
\else \ifnum#1=5
	\sh(0,0);
	\boxslopedr(.75,.25,1)1;
	\bubble(.75,1,.25)2;
	\sh(0,1);
	\bx(1,1)1;
	\sh(1,0);
	\boxslopedl(.25,.75,2)2;
	\begin{pgfonlayer}{foreground}
		\bull(.25,1);\bull(.75,1);\bull(0,2);\bull(2,1);
		\dfct(.5,0);\dfct(.5,1);\dfct(.5,2);\idfct(1.5,\cw);\idfct(1.5,2);
	\end{pgfonlayer}
\else \ifnum#1=6
		\sh(0,0);
		\dhku(2,1*\fudge)1,1,1;
		\sh(0,1*\fudge+\gap);
		\hex(2,1)1;
		\sh(0,1+1*\fudge+2*\gap);
		\dhkd(2,1*\fudge)1,1,2;
		\sh(0,1+1*\fudge+2*\gap);
		\halfidr(2,1*\fudge){.75};
		\tick(1,1+2*\fudge+2*\gap);
\else \ifnum#1=7
		\sh(0,0);
		\hkul(2+\gap,1*\fudge)2,1,2;
		\sh(0,0);
		\halfidr(2,1*\fudge){.85};
		\sh(0,1*\fudge+\gap);
		\qhexul(1,1)1;
		\sh(0,1*\fudge+1+\gap);
		\bx(1,1)1;
		\sh(1+\gap,1*\fudge+\gap);
		\bx(1,2)2;
		\tick(1+.5*\gap,1*\fudge);
\fi\fi\fi\fi\fi\fi\fi}

\def\AxiomT#1{ %20
\ifnum#1=1
	\def\x{0}\def\y{0}
	\bxbur(1,1)1;
	\def\x{1}\def\y{0}
	\bxbul(1,1)2;
	\def\x{0}\def\y{1+\gap} 
	\boxsloperu(1,.75,.25+\cw)2;
	\def\x{1}\def\y{1+\gap}
	\cornerul(.25,.75,.75,.25)2;
	\def\x{1}\def\y{1+\gap}
	\boxcornul(.25,.75,.75,.25)2;
	\begin{pgfonlayer}{foreground}	
		\bull(.5,1+\gap);\bull(1.5,1+\gap);
		\bull(1,1+\gap+.5);\bull(1.5,2+\gap+\cw);
	\end{pgfonlayer}

\else \ifnum#1=2
	\def\x{0}\def\y{0}
	\boxsloperu(1,.75,.25+\cw)2;
	\def\x{1}\def\y{0}
	\cornerul(.75,.75,.75,.25)2;
	\def\x{1}\def\y{0}
	\boxcornul(.75,.75,.75,.25)2;
	\def\x{0}\def\y{1+\gap+\cw}
	\bxbdr(1-\cw,1)1;
	\def\x{1-\cw}\def\y{1+\gap}
	\boxslopedl(.75,.75+\cw,1+\cw)2;
	\begin{pgfonlayer}{foreground}	
		\bull(1,.5);\bull(.5,1+\cw);\bull(1.5,1+\cw);\bull(2,1);
		\bull(1.5,1+\cw+\gap);\bull(2,1+\gap);
		\bull(1.75,2+\cw+\gap);
	\end{pgfonlayer}

\else \ifnum#1=3
	\def\x{0}\def\y{0}
	\bxbur(1,1)1;
	\def\x{1}\def\y{0}
	\bxbul(1,1)2;
	\def\x{0}\def\y{1+\gap}
	\bxbdr(1,1)2;
	\def\x{1}\def\y{1+\gap}
	\boxt(.25,.5,.25,1)2;
	\def\x{1}\def\y{2+\gap}
	\bubble(.25,.5,.25)2;
	\begin{pgfonlayer}{foreground}	
		\bull(.5,1+\gap);\bull(1.5,1+\gap);
		\bull(1,1+\gap+.5);\bull(1.5,2+\gap+\cw);
	\end{pgfonlayer}
	
\else \ifnum#1=4
	\def\x{0}\def\y{0}
	\bxbur(1,1)2;
	\def\x{1}\def\y{0}
	\boxt(.25,.5,.75,1)2;
	\def\x{1}\def\y{1}
	\bubble(.25,.5,.75)2;
	\def\x{0}\def\y{1+\gap}
	\bxbdr(1,1)1;
	\def\x{1}\def\y{1+\gap}
	\boxb(.25,.5,.75,1)2;
	\begin{pgfonlayer}{foreground}	
		\bull(1.5,1+\cw);\bull(2,1);
		\bull(1.5,1+\cw+\gap);\bull(2,1+\gap);
		\bull(1.75,2+\gap);
	\end{pgfonlayer}

\else \ifnum#1=5
	\def\x{0}\def\y{0}
	\bxbur(1,1)1;
	\def\x{1}\def\y{0}
	\bxbul(1,1)2;
	\def\x{0}\def\y{1+\gap}
	\boxt(1.25,.5,.25,1)2;
	\def\x{0}\def\y{2+\gap}
	\bubble(1.25,.5,.25)2;	
	\begin{pgfonlayer}{foreground}	
		\bull(.5,1+\gap);\bull(1,1+\gap);\bull(1.5,1+\gap);
		\bull(1,2+\gap);\bull(1.5,2+\gap+\cw);
	\end{pgfonlayer}

\else \ifnum#1=6
	\def\x{0}\def\y{0}
	\bxb(1,1)1;
	\def\x{1}\def\y{0}
	\boxt(.25,.5,.25,1)2;
	\def\x{1}\def\y{1}
	\bubble(.25,.5,.25)2;
	\begin{pgfonlayer}{foreground}	
		\bull(1.5,1+\cw);\bull(2,.5);
	\end{pgfonlayer}
	
\else \ifnum#1=7
	\def\x{0}\def\y{0}
	\bxbur(1,1)1;
	\def\x{1}\def\y{0}
	\bxbul(1,1)2;
	\def\x{0}\def\y{1+\gap} 
	\boxcornur(.25,.75,.75,.25)2;
	\def\x{0}\def\y{1+\gap}
	\cornerur(.25,.75,.75,.25)2;
	\def\x{1}\def\y{1+\gap}
	\boxslopelu(1,.75,.25+\cw)2;
	\begin{pgfonlayer}{foreground}	
		\bull(.5,1+\gap);\bull(1.5,1+\gap);
		\bull(1,1+\gap+.5);\bull(.5,2+\gap+\cw);
	\end{pgfonlayer}

\else \ifnum#1=8
	\def\x{0}\def\y{0}
	\bxbur(1,1)1;
	\def\x{1}\def\y{0}
	\bxbul(1,1)2;
	\def\x{0}\def\y{1+\gap}
	\boxt(.25,.5,.25,1)2;
	\def\x{0}\def\y{2+\gap}
	\bubble(.25,.5,.25)2;
	\def\x{1}\def\y{1+\gap}
	\bxbdl(1,1)2;
	\begin{pgfonlayer}{foreground}	
		\bull(.5,1+\gap);\bull(.5,2+\gap+\cw);
	\end{pgfonlayer}
	
\else \ifnum#1=9
	\def\x{0}\def\y{0}
	\bxbur(1,1)1;
	\def\x{1}\def\y{0}
	\bxbul(1,1)2;
	\def\x{0}\def\y{1+\gap}
	\boxt(.25,.5,1.25,1)2;
	\def\x{0}\def\y{2+\gap}
	\bubble(.25,.5,.25)2;
	\begin{pgfonlayer}{foreground}	
		\bull(.5,1+\gap);\bull(1,1+\gap);\bull(1.5,1+\gap);\bull(.5,2+\gap+\cw);\bull(1,2+\gap);
	\end{pgfonlayer}

\else \ifnum#1=10
	\def\x{0}\def\y{0}
	\boxcornur(.75,1.25,.75,.25)2;
	\cornerur(.75,1.25,.75,.25)2;
	\def\x{2}\def\y{0}
	\boxslopelu(1,.75,.25+\cw)2;
	\def\x{0}\def\y{1+\gap}
	\boxslopedr(.75,.25,1)1;
	\def\x{1}\def\y{1+\gap+\cw}
	\bxb(1,1-\cw)2;
	\def\x{2}\def\y{1+\gap+\cw}
	\bxbdl(1,1-\cw)2;
	\begin{pgfonlayer}{foreground}	
		\bull(.5,1);\bull(.5,1+\gap);\bull(1,1+\cw);\bull(2,.5);\bull(2.5,1+\cw);
	\end{pgfonlayer}
	
\else \ifnum#1=11
	\def\x{0}\def\y{0}
	\boxt(.25,.5,.25,1)1;
	\def\x{0}\def\y{1}
	\bubble(.25,.5,.25)1;
	\def\x{1}\def\y{0}
	\bxb(1,1)2;
	\begin{pgfonlayer}{foreground}	
		\bull(.5,1+\cw);\bull(0,.5);
	\end{pgfonlayer}
	
\else \ifnum#1=12
	\def\x{0}\def\y{0}
	\bxbt(1,1)1;
	\def\x{0}\def\y{1}
	\boxt(.25,.5,.25,1)2;
	\def\x{0}\def\y{2}
	\bubble(.25,.5,.25)2;
	\def\x{1+\gap}\def\y{0}
	\bxb(1,2)2;
	\begin{pgfonlayer}{foreground}	
		\bull(.5,2+\cw);
	\end{pgfonlayer}
	
\else \ifnum#1=13
	\def\x{0}\def\y{0}
	\boxcornur(.25,.75,.75,.25)1;
	\cornerur(.25,.75,.75,.25)1;
	\def\x{1}\def\y{0}
	\boxslopelu(1,.75,.25+\cw)2;
	\begin{pgfonlayer}{foreground}	
		\bull(.5,1+\cw);\bull(0,.5);\bull(1,.5)\bull(2,.5);
	\end{pgfonlayer}
	
\else \ifnum#1=14 %!% This has two predefined cuspwidths, both of which are different from the usual cuspwidth and from the cuspstretched cuspwidth.  This picture is likely to break if there are parameter changes.

	\def\cuspwidthtmp{\cuspwidth}
	\def\cuspwidth{.26}

	\def\x{0}\def\y{0}
	\sbx(1,1)1;

\begin{scope}[cm={1,0,0,-1,(0,\y+\y+.6)}]

	\draw[fill=\mycolor{1},even odd rule]
		{[shift={(0,.6)}]
			\cusp(0,1,.5,.5+\cw-.5*\cuspwidth)}
		{[shift={(.5+.5+\cw-.5*\cuspwidth,.6-.5*\cuspwidth)}] 
			-- \cusp(270,-1,.55-.5*\cuspwidth,.05)}
		{[shift={(.5+.5+\cw,0)}] 
			-- \cusp(90,1,.05,.55+.5*\cuspwidth)}
		{[shift={(.5+.5+\cw+.5*\cuspwidth,.6+.5*\cuspwidth)}] 
			-- \cusp(180,-1,.5+\cw+.5*\cuspwidth,.5)} -- cycle;

\end{scope}
	
	\def\x{0}\def\y{\cw}
	
	\def\cuspwidth{.2}
	
	\cornerur(.48+.5*\cusplength,.52-.5*\cusplength,.5*\cusplength,1-.5*\cusplength-\cw)1;
	
	\def\cuspwidth{.26}
	
	\def\x{1+.1}\def\y{-.13}
	\boxslopeld(1,.45+.1,.55+.13)2;
	
	\def\cuspwidth{.2}
	
	\begin{pgfonlayer}{foreground}	
		\bull(0,0);\bull(0,.5);\bull(0,1);
		\bull(2+\cw,.5);
		\bull(.5,1); \bull(1+\cw,.4);\bull(1+\cw,.8);
		\bull(1,.13);
	\end{pgfonlayer}

	\def\cuspwidth{\cuspwidthtmp}

\else \ifnum#1=15
	\def\x{0}\def\y{0}
	\boxslopedr(.5,.5,1)1;
	\cornerdr(.5,1.5,.75,.25)2;
	\def\x{1}\def\y{\cw}
	\boxsloperd(1,.75,.25-\cw)2;
	\def\x{2}\def\y{-\cw}
	\boxslopeld(1,.75,.25+\cw)2;
	\begin{pgfonlayer}{foreground}	
		\bull(0,.5);\bull(1,.5);\bull(2,.5);\bull(3,.5-.5*\cw);
	\end{pgfonlayer}
	
\else \ifnum#1=16
	\def\x{0}\def\y{0}
	\bxbt(1,1)1;
	\def\x{0}\def\y{1}
	\boxcornur(.25,.75,.5,.5)2;
	\cornerur(.25,.75,.5,.5)2;
	\def\x{1+\gap}\def\y{0}
	\boxslopelu(1,1.5,.5+\cw)2;
	\begin{pgfonlayer}{foreground}	
		\bull(.5,2+\cw);\bull(1+\gap,1);\bull(2+\gap,1);
	\end{pgfonlayer}
	
\else \ifnum#1=17
	\def\x{0}\def\y{0}
	\boxt(.75,1,.25,1)2;
	\def\x{0}\def\y{1}
	\bubble(.75,1,.25)2;
	\def\x{0}\def\y{1+\gap}
	\boxslopedr(.75,.25,1)1;
	\def\x{1}\def\y{1+\gap}
	\boxslopedl(.25,.75,1)2;
	\def\x{2+\gap}\def\y{0}
	\bxb(1,2+\gap)2;
	\begin{pgfonlayer}{foreground}	
		\bull(.5,1);\bull(1.25,1-\cw);\bull(1,1+\cw);\bull(1.5,1+\cw);
		\bull(.5,1+\gap);\bull(1.5,1+\cw+\gap);
		\bull(1,1.5+\gap);
	\end{pgfonlayer}
	
\fi\fi\fi\fi\fi\fi\fi\fi\fi\fi\fi\fi\fi\fi\fi\fi\fi}

\def\AxiomTnew#1{ %20new
\ifnum#1=1
	\sh(0,0);
	\bxbur(1,1)1;
	\sh(1,0);
	\bxbul(1,1)2;
	\sh(0,1+\gap);
	\bxbdr(1,1)2;
	\sh(1,1+\gap);
	\boxt(.125,.5,.375,1)2;
	\sh(1,2+\gap);
	\bubble(.125,.5,.375)2;
	\begin{pgfonlayer}{foreground}
		\bull(1.5,1+\gap); \bull(1.5,2+\gap+\cw);
		\dfct(.5,0);\dfct(1.5,0);\dfct(.75,1);\dfct(1.25,1);\dfct(.75,1+\gap);\dfct(1.25,1+\gap);\dfct(.5,2+\gap);\dfct(1.75,2+\gap);\idfct(1.25,2+\gap+\cw);
	\end{pgfonlayer}

\else \ifnum#1=2
	\sh(0,0);
	\bxbur(1,1)2;
	\sh(1,0);
	\boxt(.125,.5,.875,1)2;
	\sh(1,1);
	\bubble(.125,.5,.875)2;
	\sh(0,1+\gap);
	\bxbdr(1,1)1;
	\sh(1,1+\gap);
	\boxb(.125,.5,.875,1)2;
	\begin{pgfonlayer}{foreground}
		\bull(1.5,1+\cw);\bull(2,1);\bull(1.5,1+\cw+\gap);\bull(2,1+\gap);\bull(1.5,2+\gap);
		\dfct(.5,0);\dfct(.75,1);\dfct(.75,1+\gap);\dfct(.5,2+\gap);\dfct(1.75,0);\dfct(1.75,1);\dfct(1.75,1+\gap);\dfct(1.75,2+\gap);\idfct(1.25,1+\cw);\idfct(1.25,1+\cw+\gap);\idfct(1.25,2+\gap);
	\end{pgfonlayer}

\else \ifnum#1=3
	\sh(0,0);
	\bxbur(1,1)1;
	\sh(1,0);
	\bxbul(1,1)2;
	\sh(0,1+\gap);
	\boxt(1.125,.5,.375,1)2;
	\sh(0,2+\gap);
	\bubble(1.125,.5,.375)2;
	\begin{pgfonlayer}{foreground}
		\bull(1,2+\gap);\bull(.5,1+\gap);\bull(1,1+\gap);\bull(1.5,1+\gap); \bull(1.5,2+\gap+\cw);
		\dfct(.5,0);\dfct(1.5,0);\dfct(.75,1);\dfct(1.25,1);\dfct(.75,1+\gap);\dfct(1.25,1+\gap);\dfct(.5,2+\gap);\dfct(1.75,2+\gap);\idfct(1.25,2+\gap+\cw);
	\end{pgfonlayer}	
\else \ifnum#1=4
	\sh(0,0);
	\bxb(1,1)1;
	\sh(1,0);
	\boxt(.125,.5,.375,1)2;
	\sh(1,1);
	\bubble(.125,.5,.375)2;
	\begin{pgfonlayer}{foreground}
		\bull(2,.5);\bull(1.5,1+\cw);\dfct(.5,0);\dfct(.5,1);\dfct(1.5,0);\dfct(1.75,1);\idfct(1.25,1+\cw);
	\end{pgfonlayer}
\else \ifnum#1=5
	\sh(0,0);
	\bxbur(1,1)1;
	\sh(1,0);
	\bxbul(1,1)2;
	\sh(1,1+\gap);
	\bxbdl(1,1)2;
	\sh(0,1+\gap);
	\boxt(.375,.5,.125,1)2;
	\sh(0,2+\gap);
	\bubble(.375,.5,.125)2;
	\begin{pgfonlayer}{foreground}
		\bull(.5,1+\gap); \bull(.5,2+\gap+\cw);
		\dfct(.5,0);\dfct(1.5,0);\dfct(.75,1);\dfct(1.25,1);\dfct(.75,1+\gap);\dfct(1.25,1+\gap);\dfct(1.5,2+\gap);\dfct(.25,2+\gap);\idfct(.75,2+\gap+\cw);
	\end{pgfonlayer}

\else \ifnum#1=6
	\sh(0,0);
	\bxbur(1,1)1;
	\sh(1,0);
	\bxbul(1,1)2;
	\sh(0,1+\gap);
	\boxt(.375,.5,1.125,1)2;
	\sh(0,2+\gap);
	\bubble(.375,.5,1.125)2;
	\begin{pgfonlayer}{foreground}
		\bull(1,2+\gap);\bull(.5,1+\gap);\bull(1,1+\gap);\bull(1.5,1+\gap); \bull(.5,2+\gap+\cw);
		\dfct(.5,0);\dfct(1.5,0);\dfct(.75,1);\dfct(1.25,1);\dfct(.75,1+\gap);\dfct(1.25,1+\gap);\dfct(1.5,2+\gap);\dfct(.25,2+\gap);\idfct(.75,2+\gap+\cw);
	\end{pgfonlayer}
			
\else \ifnum#1=7
	\def\x{0}\def\y{0}
	\boxt(.875,.875,.25,1)2;
	\def\x{0}\def\y{1}
	\bubble(.875,.875,.25)2;
	\def\x{0}\def\y{1+\gap}
	\boxslopedr(.875,.125,1)1;
	\def\x{1}\def\y{1+\gap}
	\boxslopedl(.25,.75,1)2;
	\sh(2,0);
	\bxbul(1,1)2;
	\sh(2,1+\gap);
	\bxbdl(1,1)2;
	\begin{pgfonlayer}{foreground}	
		\bull(.5,1);\bull(1,1+\cw);\idfct(1.5,1+\cw);
		\bull(.5,1+\gap);\idfct(1.5,1+\cw+\gap);
		\bull(1,1.5+\gap);
		\dfct(1,0);\dfct(2.5,0);\dfct(.75,1);\dfct(2.25,1);\dfct(2.25,1+\gap);\dfct(.75,1+\gap);\dfct(.5,2+\gap);\dfct(2.5,2+\gap);\idfct(1.5,2+\gap);
	\end{pgfonlayer}

\else \ifnum#1=8
	\sh(1,0);
	\bxb(1,1)2;
	\sh(0,0);
	\boxt(.375,.5,.125,1)1;
	\sh(0,1);
	\bubble(.375,.5,.125)1;
	\begin{pgfonlayer}{foreground}
		\bull(2,.5);\bull(.5,1+\cw);\dfct(.5,0);\dfct(1.5,1);\dfct(1.5,0);\dfct(.25,1);\idfct(.75,1+\cw);\bull(0,.5);
	\end{pgfonlayer}	
	
\else \ifnum#1=9
	\def\x{0}\def\y{0}
	\bubble(.875,.875,.25)2;
	\def\x{0}\def\y{0}
	\boxslopedr(.875,.125,1)1;
	\def\x{1}\def\y{0}
	\boxslopedl(.25,.75,1)2;
	\sh(2,0);
	\bxb(1,1)2;
	\begin{pgfonlayer}{foreground}	
		\bull(0,.5);\bull(1,.5);\dfct(.5,0);\dfct(2.5,0);\dfct(.5,1);\dfct(2.5,1);\idfct(1.5,\cw);\idfct(1.5,1);
	\end{pgfonlayer}
	
\else \ifnum#1=10
	\sh(0,0);
	\bxbt(1,1)1;
	\sh(0,1);
	\boxt(.375,.5,.125,1)2;
	\sh(0,2);
	\bubble(.375,.5,.125)2;
	\sh(1+\gap,0);
	\bxb(1,2)2;
	\begin{pgfonlayer}{foreground}
		\bull(.5,2+\cw);\dfct(.5,0);\dfct(.5,1);\dfct(.25,2);\dfct(1.5+\gap,0);\dfct(1.5+\gap,2);\idfct(.75,2+\cw);
	\end{pgfonlayer}
	
\else \ifnum#1=11
	\def\x{0}\def\y{0}
	\boxt(.875,.5,.625,1)2;
	\def\x{0}\def\y{1}
	\bubble(.875,.5,.625)2;
	\def\x{0}\def\y{1+\gap}
	\boxslopedr(.875,.125,1)1;
	\def\x{1}\def\y{1+\gap}
	\boxslopedl(.625,.375,1)2;
	\def\x{2+\gap}\def\y{0}
	\bxb(1,2+\gap)2;
	\begin{pgfonlayer}{foreground}	
		\bull(.5,1);\bull(1,1+\cw);
		\bull(.5,1+\gap);
		\bull(1,1.5+\gap); \bull(1.5,1);\bull(1.5,1+\gap);
		\dfct(1,0);\dfct(.75,1);\dfct(.75,1+\gap);\dfct(.5,2+\gap);\dfct(2.5+\gap,0);\dfct(2.5+\gap,2+\gap);\idfct(1.25,1+\cw);\idfct(1.25,1+\cw+\gap);\idfct(1.5,2+\gap);
	\end{pgfonlayer}
	
\fi\fi\fi\fi\fi\fi\fi\fi\fi\fi\fi}

\def\AxiomTa#1{
	\def\ticklengthtmp{\ticklength}
	\pgfmathsetmacro{\ticklengthnew}{\tickfactor*\ticklengthtmp} 
	\def\ticklength{\ticklengthnew}

	\ifnum#1=1
		\sh(0,0);
		\bx(1,1)3;
		\sh(1,0);
		\bx(2,1)2;
		\sh(0,1+\gap);
		\bx(1,1)2;
		\sh(1,1+\gap);
		\bx(2,1)2;
		\halfidl(2,1){.25};
		\tick(2,2+\gap);
\else \ifnum#1=2
		\sh(0,0);
		\bx(3,1)2;
		\sh(0,1+\gap);
		\bx(1,1)2;
		\sh(1,1+\gap);
		\bx(2,1)2;
		\halfidl(2,1){.25};
		\tick(2,2+\gap);
\else \ifnum#1=3
		\sh(0,0);
		\bx(2,1)3;
		\sh(2,0);
		\bx(1,1)2;
		\sh(0,1+\gap);
		\bx(2,1)2;
		\halfidr(2,1){.75};
		\sh(2,1+\gap);
		\bx(1,1)2;
		\tick(1,2+\gap);\tick(1,0);\tick(1,1);
\else \ifnum#1=4
		\sh(0,0);
		\bx(3,1)2;
		\sh(0,1+\gap);
		\bx(2,1)2;
		\halfidr(2,1){.75};
		\sh(2,1+\gap);
		\bx(1,1)2;
		\tick(1,2+\gap);\tick(2,0);\tick(2,1);
\else \ifnum#1=5
		\sh(0,0);
		\bx(1,1)2;
		\sh(1,0);
		\bx(2,1)2;
		\halfidl(2,1){.25};
		\tick(2,1);
\else \ifnum#1=6
		\sh(0,0);
		\bx(1,1)2;
		\sh(1,0);
		\bx(2,1)2;
		\halfidl(2,1){.25};
		\sh(0,1+\gap);
		\bx(3,1)2;
		\tick(1,1+\gap);\tick(1,2+\gap);\tick(2,1);\tick(2,1+\gap);\tick(2,2+\gap);
\else \ifnum#1=7
		\sh(0,0);
		\bx(1,1)2;
		\sh(1,0);
		\bx(2,1)2;
		\halfidl(2,1){.25};
		\sh(0,1+\gap);
		\bx(1,1)3;
		\sh(1,1+\gap);
		\bx(2,1)2;
		\tick(2,1);\tick(2,1+\gap);\tick(2,2+\gap);
\else \ifnum#1=8
		\sh(0,0);
		\bx(2,1)2;
		\halfidr(2,1){.75};
		\sh(2,0);
		\bx(1,1)2;
		\sh(0,1+\gap);
		\bx(3,1)2;
		\tick(1,1);\tick(1,1+\gap);\tick(1,2+\gap);\tick(2,1+\gap);\tick(2,2+\gap);
\else \ifnum#1=9
		\sh(0,0);
		\bx(1,1)2;
		\sh(1,0);
		\bx(2,1)2;
		\halfidl(2,1){.25};
		\sh(0,1+\gap);
		\bx(2,1)2;
		\sh(2,1+\gap);
		\bx(1,1)2;
		\tick(2,1);\tick(1,1+\gap);\tick(1,2+\gap);
\else \ifnum#1=10
		\sh(0,0);
		\bx(1,1)2;
		\sh(1,0);
		\bx(2,1)2;
		\halfidl(2,1){.25};
		\sh(0,1+\gap);
		\bx(1,1)2;
		\sh(1,1+\gap);
		\bx(1,1)2;
		\sh(2,1+\gap);
		\bx(1,1)2;
		\tick(2,1);
\else \ifnum#1=11
		\sh(0,0);
		\bx(2,1)3;
		\sh(0,1);
		\bx(2,1)2;
		\halfidr(2,1){.75};
		\sh(2+\gap,0);
		\bx(1,2)2;
		\tick(1,2);
\else \ifnum#1=12
		\sh(0,0);
		\bx(2,1)2;
		\halfidr(2,1){.75};
		\sh(2,0);
		\bx(1,1)2;
		\tick(1,1);
\else \ifnum#1=13
		\sh(0,0);
		\bx(2,1)2;
		\halfidr(2,1){.75};
		\sh(0,1);
		\bx(2,1)2;
		\sh(2+\gap,0);
		\bx(1,2)2;
		\tick(1,1);\tick(1,2);
\else \ifnum#1=14
		\sh(0,0);
		\bx(2,1)2;
		\halfidr(2,1){.75};
		\sh(2,0);
		\bx(1,1)2;
		\sh(0,1+\gap);
		\bx(2,1)2;
		\sh(2,1+\gap);
		\bx(1,1)2;
		\tick(1,1);\tick(1,1+\gap);\tick(1,2+\gap);
\else \ifnum#1=15
		\sh(0,0);
		\bx(2,1)2;
		\halfidr(2,1){.75};
		\sh(2,0);
		\bx(1,1)2;
		\sh(0,1+\gap);
		\bx(1,1)3;
		\sh(1,1+\gap);
		\bx(1,1)2;
		\sh(2,1+\gap);
		\bx(1,1)2;
		\tick(1,1);
\else \ifnum#1=16
		\sh(0,0);
		\bx(2,1)2;
		\halfidr(2,1){.75};
		\sh(0,1+\gap);
		\bx(1,1)3;
		\sh(1,1+\gap);
		\bx(1,1)2;
		\sh(2+\gap,0);
		\bx(1,2+\gap)2;
		\tick(1,1);
\def\ticklength{\ticklengthtmp}
		
\fi\fi\fi\fi\fi\fi\fi\fi\fi\fi\fi\fi\fi\fi\fi\fi}

\def\AxiomU#1{ %21
\ifnum#1=1
	\def\x{0}\def\y{0}
	\bxbur(1,1)1;
	\def\x{1}\def\y{0}
	\bxbul(1,1)2;
	\def\x{0}\def\y{1+\gap}
	\bxbdr(1,1)2;
	\def\x{1}\def\f{1+\gap}
	\boxcornur(.25,.75,.5,.5)2
	\cornerur(.25,.75,.5,.5)2;
	\begin{pgfonlayer}{foreground}
		\bull(1.5,1+\gap);\bull(1.5,2+\gap+\cw);
	\end{pgfonlayer}
\else \ifnum#1=2
	\sh(0,0)
	\bxbur(1,1)1;
	\sh(1,0)
	\bxbul(1,1)2;
	\sh(0,1+\gap)
	\boxcornur(1.25,.75,.5,.5)2;
	\cornerur(1.25,.75,.5,.5)2;
	\begin{pgfonlayer}{foreground}
		\bull(.5,1+\gap);\bull(1,1+\gap);\bull(1.5,1+\gap);
		\bull(1,2+\gap);\bull(1.5,2+\gap+\cw);
	\end{pgfonlayer}

\else \ifnum#1=3
	\def\x{0}\def\y{0}
	\bxb(1,1)1;
	\def\x{1}\def\y{0}
	\boxcornur(.5,1,.75,.25)2;
	\cornerur(.5,1,.75,.25)2;
	\begin{pgfonlayer}{foreground}
		\bull(2.5,.5);\bull(1.75,1+\cw);
	\end{pgfonlayer}

\else \ifnum#1=4  %!% Cuspstretched a la 18-5.
	\sh(0,0)
	\bxb(1,1+\cw)1;
	\sh(1,0)
	\sbx(1,1+\cw)2;
	\draw[fill=\mycolor{2},even odd rule]
		{[shift={(0,0)}]
			\cusp(0,1,.5,.5+.5*\cuspwidth+.5*\cuspwidth*\cuspstretch)}
		{[shift={(.5+.5+.5*\cuspwidth+.5*\cuspwidth*\cuspstretch,-.5*\cuspwidth)}] 
			-- (\x,\y) -- ++(90:.52-.5*\cusplength) .. controls +(90:\cuspfactor*\cusplength) and +(180+90:\cuspfactor*\cusplength) .. ++($ (90:\cusplength) + (90+1*90:.5*\cuspwidth*\cuspstretch) $) -- ++ (90:.15+.5*\cuspwidth-.5*\cusplength)	}		
		{[shift={(.5+.5+.5*\cuspwidth,.15+.52)}] 
			-- \cusp(270,1,.15,.52-.5*\cuspwidth)}
		{[shift={(.5+.5,.5*\cuspwidth)}] 
			-- \cusp(180,-1,.5,.5)} -- cycle;
	\sh(2,\cw);
	\bubblev(.5*\cusplength,1-\cusplength,.5*\cusplength)2;
	\begin{pgfonlayer}{foreground}
		\bull(2+\cw+\cw*\cuspstretch,-\cw);\bull(2,\cw);
		\bull(2,1+\cw);
		\bull(2-\cw,.5+\cw);
		\bull(2+\cw,.25+\cw);\bull(2+\cw,.5+\cw);\bull(2+\cw,.75+\cw);
	\end{pgfonlayer}

\else \ifnum#1=5
	\def\x{0}\def\y{0}
	\bxb(1,2)1;
	\def\x{1+\gap}\def\y{0}
	\boxt(.75,.75,.5,1)2;
	\def\x{1+\gap}\def\y{1}
	\bubble(.75,1,.25)2;
	\boxb(.75,1,.25,1)2;	
	\begin{pgfonlayer}{foreground}
		\bull(2+\gap,2);
		\bull(1.5+\gap,1);\bull(2.25+\gap,1-\cw);\bull(2+\gap,1+\cw);\bull(2.5+\gap,1+\cw);
	\end{pgfonlayer}

\else \ifnum#1=6
	\sh(0,0)
	\bxb(1,1)1;
	\sh(1,0);
	\boxcorndr(.5,1,.75,.25)2;
	\cornerdr(.5,1,.75,.25)2;
	\begin{pgfonlayer}{foreground}
		\bull(1.75,\cw);\bull(1.75,1);\bull(2.5,.5);
	\end{pgfonlayer}
	
\else \ifnum#1=7
	\sh(0,0)
	\bxb(1,1)1;
	\sh(1,0)
	\boxslopedr(.75,.25,1)2;
	\bubble(.75,1,.25)2;
	\sh(2,0)
	\boxslopedl(.25,.75,1)2;
	\begin{pgfonlayer}{foreground}
		\bull(2,.5);\bull(2.25,-\cw);\bull(2.5,\cw);
	\end{pgfonlayer}
	
\else \ifnum#1=8
	\sh(0,0)
	\bxbur(1,1)2;
	\sh(1,0)
	\boxt(.75,1,.25,1)2;
	\sh(1,1)
	\bubble(.75,1,.25)2;
	\sh(0,1+\gap)
	\bxbdr(1,1)1;
	\sh(1,1+\gap)
	\boxb(.75,1,.25,1)2;	
	\begin{pgfonlayer}{foreground}
		\bull(2.25,1-\cw);\bull(2,1+\cw);\bull(2.5,1+\cw);
		\bull(2,1+\cw+\gap);\bull(2.5,1+\cw+\gap);
		\bull(2,2+\gap);
	\end{pgfonlayer}

\else \ifnum#1=9 %!% The bullets here will break if hexangleB is changed.
	\def\x{0}\def\y{0}
	\bxbur(1,2)2;
	\def\x{1+\gap}\def\y{0}
	\boxt(.75,1,.25,1)2;
	\def\x{1+\gap}\def\y{1}
	\sh(1+\gap,1);
	\qhexdr(2,1)2;
	\bubble(.75,1,.25)2;
	\def\x{0}\def\y{2+\gap}
	\bxbdr(1+.5*\gap,1)1;
	\def\x{1+.5*\gap}\def\y{2+\gap}
	\sh(1+.5*\gap,2+\gap);
	\hkdr(2+.5*\gap,1)2,1,2;
	\begin{pgfonlayer}{foreground}
		\bull(2.25+\gap,1-\cw);\bull(1.5+\gap,1);\bull(2+\gap,1+\cw);\bull(2.5+\gap,1+\cw);
		\bull(1+\gap,2);
		\bull(1+.5*\gap+1+.25*\gap,3+\gap);\bull(3+\gap,3+\gap);
		\bull(3+\gap,1+\gap*\gapfactorB);
		\bull(1+\gap+1.3,2); %!% By hand.
		\bull(1+\gap+1.4,2+\gap); %!% By hand.
		\bull(1+.5*\gap+.75,2);
		\bull(1+.5*\gap+.75,2+\gap);
	\end{pgfonlayer}

\else \ifnum#1=10
	\def\x{0}\def\y{0}
	\bxbur(1,1)2;
	\def\x{1}\def\y{0}
	\boxcorndr(.75,1.25,.5,.5)2;
	\cornerdr(.75,1.25,.5,.5)2;
	\def\x{0}\def\y{1+\gap}
	\bxbdr(1,1)1;
	\def\x{1}\def\y{1+\gap}
	\bxbur(2,1)2;
	\begin{pgfonlayer}{foreground}
		\bull(2,\cw);\bull(1+.67,1);\bull(1+1.33,1);\bull(1+.67,1+\gap);\bull(1+1.33,1+\gap);
	\end{pgfonlayer}

\else \ifnum#1=11
	\sh(0,0)
	\boxcorndr(1.75,1.25,.5,.5)2;
	\cornerdr(1.75,1.25,.5,.5)2;
	\sh(0,1+\gap)
	\bxbdr(1,1)1;
	\sh(1,1+\gap)
	\bxbur(2,1)2;	
	\begin{pgfonlayer}{foreground}
		\bull(2,\cw);\bull(1+.67,1);\bull(1+1.33,1);\bull(1+.67,1+\gap);\bull(1+1.33,1+\gap);
		\bull(1,0);\bull(.5,1);\bull(1,1);
	\end{pgfonlayer}

\else \ifnum#1=12 %!% Bullets will break if hexangleB is changed.
	\sh(0,0);
	\def\hexangleAtmp{\hexangleA}
	\pgfmathsetmacro{\hexangleA}{\hexangleB} 

	\hex(3,2)2;

	\def\hexangleA{\hexangleAtmp}

	\sh(0,0)
	\boxt(1.75,1,.25,1)2;
	\sh(0,1)
	\bubble(1.75,1,.25)2;
	\sh(0,2+\gap);
	\hkdl(1,1)2,1,1;
	\sh(1,2+\gap);
	\hkdr(2,1)2,1,2;
	\begin{pgfonlayer}{foreground}
		\bull(1,0);\bull(.5,1);\bull(1,1);\bull(2.25,1-\cw);\bull(2,1+\cw);\bull(2.5,1+\cw);
		\bull(0,1+\gap*\gapfactorB);\bull(3,1+\gap*\gapfactorB);
		\bull(0,3+\gap);\bull(3,3+\gap);
		\bull(1,2);\bull(1,2+\gap);\bull(1,2.5+\gap);\bull(1,3+\gap);
		\bull(1+.67,2);\bull(1+.67,2+\gap);
		\bull(2,3+\gap);
		\bull(1-.3,2);\bull(1-.4,2+\gap);\bull(2+.3,2);\bull(2+.4,2+\gap);%!% By hand.
	\end{pgfonlayer}
	
\else \ifnum#1=13
	\sh(0,0)
	\boxt(1.75,1,.25,1)2;
	\sh(0,1)
	\bubble(1.75,1,.25)2;
	\sh(0,1+\gap)
	\bxbdr(1,1)1;
	\sh(1,1+\gap)
	\boxb(.75,1,.25,1)2;
	\begin{pgfonlayer}{foreground}
		\bull(1,0); \bull(.5,1);\bull(1,1);\bull(2.25,1-\cw);\bull(2,1+\cw);\bull(2.5,1+\cw);
		\bull(2,1+\cw+\gap);\bull(2.5,1+\cw+\gap);\bull(2,2+\gap);
	\end{pgfonlayer}
\else \ifnum#1=14
	\sh(0,0)
	\boxt(1.75,1,.25,1)2;
	\sh(0,1)
	\bubble(1.75,1,.25)2;
	\sh(0,1+\gap)
	\bxbdr(1,1)1;
	\sh(1,1+\gap)
	\boxslopedr(.75,.25,1)2;
	\sh(2,1+\gap)
	\boxslopedl(.25,.75,1)2;
	\begin{pgfonlayer}{foreground}
		\bull(1,0);\bull(.5,1);\bull(1,1);\bull(2,1+\cw);\bull(2.25,1-\cw);\bull(2.5,1+\cw);
		\bull(2.5,1+\cw+\gap);\bull(2,1.5+\gap);
	\end{pgfonlayer}
	
\fi\fi\fi\fi\fi\fi\fi\fi\fi\fi\fi\fi\fi\fi}

\def\AxiomUnew#1{ %21new
\ifnum#1=1
	\sh(0,0);
	\bxbur(1,1)1;
	\sh(1,0);
	\bxbul(1,1)2;
	\sh(0,1+\gap);
	\bxbdr(1,1)2;
	\sh(1,1+\gap);
	\boxt(.375,.5,.125,1)2;
	\sh(1,2+\gap);
	\bubble(.375,.5,.125)2;
	\begin{pgfonlayer}{foreground}
		\bull(1.5,1+\gap); \bull(1.5,2+\gap+\cw);
		\dfct(.5,0);\dfct(1.5,0);\dfct(.75,1);\dfct(1.25,1);\dfct(.75,1+\gap);\dfct(1.25,1+\gap);\dfct(.5,2+\gap);\dfct(1.25,2+\gap);\idfct(1.75,2+\gap+\cw);
	\end{pgfonlayer}

\else \ifnum#1=2
	\sh(0,0);
	\bxbur(1,1)2;
	\sh(1,0);
	\boxt(.375,.5,.625,1)2;
	\sh(1,1);
	\bubble(.375,.5,.625)2;
	\sh(0,1+\gap);
	\bxbdr(1,1)1;
	\sh(1,1+\gap);
	\boxb(.375,.5,.625,1)2;
	\begin{pgfonlayer}{foreground}
		\bull(1.5,1+\cw);\bull(2,1);\bull(1.5,1+\cw+\gap);\bull(2,1+\gap);\bull(1.5,2+\gap);
		\dfct(.5,0);\dfct(.75,1);\dfct(.75,1+\gap);\dfct(.5,2+\gap);\dfct(1.75,0);\dfct(1.25,1);\dfct(1.25,1+\gap);\dfct(1.25,2+\gap);\idfct(1.75,1+\cw);\idfct(1.75,1+\cw+\gap);\idfct(1.75,2+\gap);
	\end{pgfonlayer}

\else \ifnum#1=3
	\sh(0,0);
	\bxb(1,1)1;
	\sh(1,0);
	\boxt(.375,.5,.125,1)2;
	\sh(1,1);
	\bubble(.375,.5,.125)2;
	\begin{pgfonlayer}{foreground}
		\bull(2,.5);\bull(1.5,1+\cw);\dfct(.5,0);\dfct(.5,1);\dfct(1.5,0);\dfct(1.25,1);\idfct(1.75,1+\cw);
	\end{pgfonlayer}
	
\else \ifnum#1=4  
	\def\x{1}\def\y{0}
	\bubble(.875,.875,.25)2;
	\def\x{1}\def\y{0}
	\boxslopedr(.875,.125,1)2;
	\def\x{2}\def\y{0}
	\boxslopedl(.25,.75,1)2;
	\sh(0,0);
	\bxb(1,1)1;
	\begin{pgfonlayer}{foreground}	
		\bull(3,.5);\bull(2,.5);\dfct(.5,0);\dfct(1.5,0);\dfct(.5,1);\dfct(1.5,1);\idfct(2.5,\cw);\idfct(2.5,1);
	\end{pgfonlayer}
		
\else \ifnum#1=5
	\sh(0,0);
	\bxbur(1,1)1;
	\sh(1,0);
	\bxbul(1,1)2;
	\sh(0,1+\gap);
	\boxt(1.375,.5,.125,1)2;
	\sh(0,2+\gap);
	\bubble(1.375,.5,.125)2;
	\begin{pgfonlayer}{foreground}
		\bull(1,2+\gap);\bull(.5,1+\gap);\bull(1,1+\gap);\bull(1.5,1+\gap); \bull(1.5,2+\gap+\cw);
		\dfct(.5,0);\dfct(1.5,0);\dfct(.75,1);\dfct(1.25,1);\dfct(.75,1+\gap);\dfct(1.25,1+\gap);\dfct(1.25,2+\gap);\dfct(.5,2+\gap);\idfct(1.75,2+\gap+\cw);
	\end{pgfonlayer}
	
\else \ifnum#1=6
	\sh(0,0);
	\boxt(1.75,.625,.625,1)2;
	\def\x{0}\def\y{1}
	\bubble(1.75,.625,.625)2;
	\def\x{1}\def\y{1+\gap}
	\boxslopedr(.75,.25,1)2;
	\def\x{2}\def\y{1+\gap}
	\boxslopedl(.625,.375,1)2;
	\sh(0,1+\gap);
	\bxbdr(1,1)1;
	\begin{pgfonlayer}{foreground}	
		\bull(1,0);\bull(1,1);\dfct(.5,0);\dfct(2,0);\bull(.5,1);\dfct(.75,1);\dfct(.75,1+\gap);\dfct(1.5,1);\dfct(1.5,1+\gap);\dfct(.5,2+\gap);\dfct(1.5,2+\gap);\bull(2,1.5+\gap);\bull(2,1+\cw);\bull(2.5,1);\bull(2.5,1+\gap);\idfct(2.25,1+\cw);\idfct(2.25,1+\cw+\gap);\idfct(2.5,2+\gap);
	\end{pgfonlayer}
	
\fi\fi\fi\fi\fi\fi}

\def\AxiomUa#1{ %21a
	\ifnum#1=1 %!% Bullets will break if hexangleB is changed.
	\sh(0,1)
	\sbx(2,1)2;
	\sh(0,0);
	\dhku(2,1)2,0,1;
	\sh(0,2)
	\sbx(2,1)2;
	\bubble(.75,1,.25)2;
	\begin{pgfonlayer}{foreground}
		\bull(0,0);\bull(2,0);\bull(0,2);\bull(0,3);\bull(2,2);\bull(2,3);
		\bull(.5,2);\bull(1,2+\cw);\bull(1.5,2+\cw);\bull(1.25,2-\cw);
		\bull(1-.3,1);\bull(1+.3,1);
	\end{pgfonlayer}

\else \ifnum#1=2
	\def\x{0}\def\y{0}
	\boxt(.75,1,.25,1)1;
	\def\x{0}\def\y{1}
	\bubble(.75,1,.25)2;
	\boxb(.75,1,.25,1)2;
	\begin{pgfonlayer}{foreground}
		\bull(.5,1);\bull(1.25,1-\cw);\bull(1,1+\cw);\bull(1.5,1+\cw);
	\end{pgfonlayer}
	
\else \ifnum#1=3
	\def\x{0}\def\y{0}
	\bxbt(1,1)1;
	\def\x{0}\def\y{1}
	\boxcornur(.25,.75,.5,.5)2;
	\cornerur(.25,.75,.5,.5)2;
	\begin{pgfonlayer}{foreground}
		\bull(.5,2+\cw);
	\end{pgfonlayer}
	
\else \ifnum#1=4
	\def\x{0}\def\y{0}
	\boxcornur(.5,1,.75,.25)1;
	\cornerur(.5,1,.75,.25)1;
	\begin{pgfonlayer}{foreground}
		\bull(0,.5);\bull(1.5,.5);\bull(.75,1+\cw);
	\end{pgfonlayer}
	
\fi\fi\fi\fi}

\def\AxiomUanew#1{ %21anew
	\ifnum#1=1
	\sh(0,0);
	\bxb(1,1)3;
	\sh(1,0);
	\boxt(.375,.5,.125,1)2;
	\sh(1,1);
	\bubble(.375,.5,.125)2;
	\begin{pgfonlayer}{foreground}
		\bull(2,.5);\bull(1.5,1+\cw);\dfct(.5,0);\dfct(.5,1);\dfct(1.5,0);\dfct(1.25,1);\idfct(1.75,1+\cw);
	\end{pgfonlayer}
		
\else \ifnum#1=2
	\sh(0,0);
	\bxbur(1,1)2;
	\sh(1,0);
	\boxt(.375,.5,.625,1)2;
	\sh(1,1);
	\bubble(.375,.5,.625)2;
	\sh(0,1+\gap);
	\bxbdr(1,1)3;
	\sh(1,1+\gap);
	\boxb(.375,.5,.625,1)2;
	\begin{pgfonlayer}{foreground}
		\bull(1.5,1+\cw);\bull(2,1);\bull(1.5,1+\cw+\gap);\bull(2,1+\gap);\bull(1.5,2+\gap);
		\dfct(.5,0);\dfct(.75,1);\dfct(.75,1+\gap);\dfct(.5,2+\gap);\dfct(1.75,0);\dfct(1.25,1);\dfct(1.25,1+\gap);\dfct(1.25,2+\gap);\idfct(1.75,1+\cw);\idfct(1.75,1+\cw+\gap);\idfct(1.75,2+\gap);
	\end{pgfonlayer}
		
\else \ifnum#1=3
	\sh(0,0);
	\boxt(1.375,.5,.125,1)2;
	\sh(0,1);
	\bubble(1.375,.5,.125)2;
	\begin{pgfonlayer}{foreground}
		\bull(2,.5);\bull(1.5,1+\cw);\dfct(.5,0);\dfct(.5,1);\dfct(1.5,0);\dfct(1.25,1);\idfct(1.75,1+\cw);\bull(1,0);\bull(1,1);\bull(0,.5);
	\end{pgfonlayer}
		
\else \ifnum#1=4
	\sh(0,0);
	\bxbur(1,1)2;
	\sh(1,0);
	\boxt(.375,.5,.625,1)2;
	\sh(1,1);
	\bubble(.375,.5,.625)2;
	\sh(0,1+\gap);
	\boxb(1.375,.5,.625,1)2;
	\begin{pgfonlayer}{foreground}
		\bull(1.5,1+\cw);\bull(2,1);\bull(1.5,1+\cw+\gap);\bull(2,1+\gap);\bull(1.5,2+\gap);
		\dfct(.5,0);\dfct(.75,1);\dfct(.75,1+\gap);\dfct(.5,2+\gap);\dfct(1.75,0);\dfct(1.25,1);\dfct(1.25,1+\gap);\dfct(1.25,2+\gap);\idfct(1.75,1+\cw);\idfct(1.75,1+\cw+\gap);\idfct(2,2+\gap);\bull(.5,1+\gap);\bull(1,1+\gap);\bull(1,2+\gap);
	\end{pgfonlayer}
	
\else \ifnum#1=5
	\def\x{0}\def\y{0}
	\bubble(1.375,.875,.25)2;
	\boxslopedr(1.375,.125,1)2;
	\def\x{1.5}\def\y{0}
	\boxslopedl(.25,.75,1)2;
	\begin{pgfonlayer}{foreground}	
		\bull(.75,0);\bull(.75,1);\bull(0,.5);\bull(1.5,.5);\bull(2.5,.5);\dfct(.375,0);\dfct(.375,1);\dfct(1.125,0);\dfct(1.125,1);\idfct(2,\cw);\idfct(2,1);
	\end{pgfonlayer}
	
\else \ifnum#1=6
	\def\x{0}\def\y{0}
	\bubble(1.175,.65,.175)2;
	\boxb(1.175,.65,.175,1)2;
	\begin{pgfonlayer}{foreground}	
		\bull(.67,0);\bull(.67,1);\bull(1.33,\cw);\bull(1.33,1);\bull(0,.5);\bull(2,.5);\dfct(.33,0);\dfct(.33,1);\dfct(1,0);\dfct(1,1);\idfct(1.67,\cw);\idfct(1.67,1);
	\end{pgfonlayer}
	
\else \ifnum#1=7
	\sh(0,0);
	\boxt(1.75,.625,.625,1)2;
	\sh(0,1);
	\bubble(1.75,.625,.625)2;
	\sh(0,1+\gap);
	\boxb(1.75,.625,.625,1)2;
	\begin{pgfonlayer}{foreground}	
		\bull(1,0);\bull(1,1);\bull(1,1+\gap);\bull(1,2+\gap);
		\bull(.5,1);\bull(.5,1+\gap);\bull(2,1+\cw);\bull(2,1+\cw+\gap);\bull(2,2+\gap);\bull(2.5,1);\bull(2.5,1+\gap);\dfct(.5,0);\dfct(.75,1);\dfct(.75,1+\gap);\dfct(.5,2+\gap);\dfct(2,0);\dfct(1.5,1);\dfct(1.5,1+\gap);\dfct(1.5,2+\gap);\idfct(2.25,1+\cw);\idfct(2.25,1+\cw+\gap);\idfct(2.5,2+\gap);
	\end{pgfonlayer}
	
\else \ifnum#1=8
	\def\x{1}\def\y{0}
	\bubble(.875,.875,.25)2;
	\def\x{1}\def\y{0}
	\boxslopedr(.875,.125,1)2;
	\def\x{2}\def\y{0}
	\boxslopedl(.25,.75,1)2;
	\sh(0,0);
	\bxb(1,1)3;
	\begin{pgfonlayer}{foreground}	
		\bull(3,.5);\bull(2,.5);\dfct(.5,0);\dfct(1.5,0);\dfct(.5,1);\dfct(1.5,1);\idfct(2.5,\cw);\idfct(2.5,1);
	\end{pgfonlayer}
		
\else \ifnum#1=9
	\sh(0,0);
	\boxt(1.75,.625,.625,1)2;
	\def\x{0}\def\y{1}
	\bubble(1.75,.625,.625)2;
	\def\x{1}\def\y{1+\gap}
	\boxslopedr(.75,.25,1)2;
	\def\x{2}\def\y{1+\gap}
	\boxslopedl(.625,.375,1)2;
	\sh(0,1+\gap);
	\bxbdr(1,1)3;
	\begin{pgfonlayer}{foreground}	
		\bull(1,0);\bull(1,1);\dfct(.5,0);\dfct(2,0);\bull(.5,1);\dfct(.75,1);\dfct(.75,1+\gap);\dfct(1.5,1);\dfct(1.5,1+\gap);\dfct(.5,2+\gap);\dfct(1.5,2+\gap);\bull(2,1.5+\gap);\bull(2,1+\cw);\bull(2.5,1);\bull(2.5,1+\gap);\idfct(2.25,1+\cw);\idfct(2.25,1+\cw+\gap);\idfct(2.5,2+\gap);
	\end{pgfonlayer}
	
\fi\fi\fi\fi\fi\fi\fi\fi\fi}

\def\AxiomV#1{ %22
	\ifnum#1=1
	\sh(0,0)
	\bxbur(1,1)1;
	\sh(1,0)	
	\bxbul(1,1)2;
	\sh(0,1+\gap)
	\boxcornul(1.25,.75,.5,.5)2;
	\cornerul(1.25,.75,.5,.5)2;
	\begin{pgfonlayer}{foreground}
		\bull(.5,1+\gap);\bull(1,1+\gap);\bull(1.5,1+\gap);
		\bull(.5,2+\gap+\cw);\bull(1,2+\gap);
	\end{pgfonlayer}

\else \ifnum#1=2
	\sh(0,0)
	\boxt(.25,1,1.75,1)2;
	\sh(0,1+\gap)
	\boxslopedr(.25,.75,1)2;
	\bubble(.25,1,.75)2;
	\sh(1,1+\gap)
	\boxslopedl(.75,.25,1)1;
	\sh(2,1+\gap)
	\bxbdl(1,1)2;
	\begin{pgfonlayer}{foreground}
		\bull(2,0);
		\bull(.75,1-\cw);\bull(.75,1-\cw+\gap);\bull(2,1);\bull(2.5,1);
		\bull(.5,1+\gap+\cw);\bull(1,1.5+\gap);\bull(2,1.5+\gap);
	\end{pgfonlayer}

\else \ifnum#1=3
	\sh(0,0)
	\bxbur(1,1)1;
	\sh(1,0)	
	\bxbul(1,1)2;
	\sh(0,1+\gap)
	\cornerul(.25,.75,.75,.25)2;
	\boxcornul(.25,.75,.75,.25)2;
	\sh(1,1+\gap)
	\bxb(1,1)2;
	\begin{pgfonlayer}{foreground}
		\bull(.5,1+\gap);\bull(.5,2+\gap+\cw);\bull(1.5,1+\gap);
	\end{pgfonlayer}

\else \ifnum#1=4

	\sh(0,0)
	\cornerul(.25,.75,.75,.25)1;
	\boxcornul(.25,.75,.75,.25)1;
	\sh(1,0)
	\bxb(1,1)2;
	\begin{pgfonlayer}{foreground}
		\bull(0,.5);\bull(.5,1+\cw);
	\end{pgfonlayer}

\else \ifnum#1=5 %!% Cuspstretched a la 18-5.
	\begin{scope}[cm={-1,0,0,1,(2,0)}] %%% The shift is necessary for the correct bounding box.

	\sh(0,0)
	\bxb(1,1+\cw)2;
	\sh(1,0)
	\sbx(1,1+\cw)1;
	\draw[fill=\mycolor{1},even odd rule]
		{[shift={(0,0)}]
			\cusp(0,1,.5,.5+.5*\cuspwidth+.5*\cuspwidth*\cuspstretch)}
		{[shift={(.5+.5+.5*\cuspwidth+.5*\cuspwidth*\cuspstretch,-.5*\cuspwidth)}] 
			-- (\x,\y) -- ++(90:.52-.5*\cusplength) .. controls +(90:\cuspfactor*\cusplength) and +(180+90:\cuspfactor*\cusplength) .. ++($ (90:\cusplength) + (90+1*90:.5*\cuspwidth*\cuspstretch) $) -- ++ (90:.15+.5*\cuspwidth-.5*\cusplength)	}		
		{[shift={(.5+.5+.5*\cuspwidth,.15+.52)}] 
			-- \cusp(270,1,.15,.52-.5*\cuspwidth)}
		{[shift={(.5+.5,.5*\cuspwidth)}] 
			-- \cusp(180,-1,.5,.5)} -- cycle;
	\sh(2,\cw);
	\bubblev(.5*\cusplength,1-\cusplength,.5*\cusplength)1;
	\begin{pgfonlayer}{foreground}
		\bull(2+\cw+\cw*\cuspstretch,-\cw);\bull(2,\cw);
		\bull(2,1+\cw);
		\bull(2-\cw,.5+\cw);
		\bull(2+\cw,.25+\cw);\bull(2+\cw,.5+\cw);\bull(2+\cw,.75+\cw);
	\end{pgfonlayer}
	
	\end{scope}
	
\else \ifnum#1=6
	\sh(0,0)
	\boxslopedr(.25,.75,1)2;
	\bubble(.25,1,.75)2;
	\sh(1,0)
	\boxslopedl(.75,.25,1)1;
	\sh(2,0)
	\bxb(1,1)2;
	\begin{pgfonlayer}{foreground}
		\bull(.75,-\cw);\bull(.5,\cw);\bull(1,.5);
	\end{pgfonlayer}

\else \ifnum#1=7
	\sh(0,0)
	\boxt(.25,1,.75,1)2;
	\sh(2,0)
	\bxbul(1,1)2;
	\sh(0,1+\gap)
	\boxslopedr(.25,.75,1)2;
	\bubble(.25,1,.75)2;
	\sh(1,1+\gap)
	\boxslopedl(.75,.25,1)1;
	\sh(2,1+\gap)
	\bxb(1,1)2;
	\begin{pgfonlayer}{foreground}
		\bull(.75,1-\cw);\bull(.75,1-\cw+\gap);\bull(.5,1+\gap+\cw);
		\bull(1,1.5+\gap);\bull(2.5,1+\gap);
	\end{pgfonlayer}

\else \ifnum#1=8
	\sh(0,0)
	\bxbt(1,1)1;
	\sh(0,1)
	\cornerul(.25,.75,.75,.25)2;
	\boxcornul(.25,.75,.75,.25)2;
	\sh(1+\gap,0)
	\bxb(1,2)2;	
	\begin{pgfonlayer}{foreground}
		\bull(.5,2+\cw);
	\end{pgfonlayer}

\else \ifnum#1=9
	\sh(0,0)
	\boxt(.25,1,.75,1)2;
	\sh(2+\gap,0)
	\bxb(1,2+\gap)2;
	\sh(0,1+\gap)
	\boxslopedr(.25,.75,1)2;
	\bubble(.25,1,.75)2;
	\sh(1,1+\gap)
	\boxslopedl(.75,.25,1)1;
	\begin{pgfonlayer}{foreground}
		\bull(.75,1-\cw);\bull(.75,1-\cw+\gap);
		\bull(.5,1+\gap+\cw);\bull(1.5,1);\bull(1.5,1+\gap);
		\bull(1,1.5+\gap);
	\end{pgfonlayer}

\fi\fi\fi\fi\fi\fi\fi\fi\fi}

\def\AxiomVnew#1{ %22new
	\ifnum#1=1
	\sh(0,0);
	\bxbur(1,1)1;
	\sh(1,0);
	\bxbul(1,1)2;
	\sh(0,1+\gap);
	\boxt(.125,.5,1.375,1)2;
	\sh(0,2+\gap);
	\bubble(.125,.5,1.375)2;
	\begin{pgfonlayer}{foreground}
		\bull(1,2+\gap);\bull(.5,1+\gap);\bull(1,1+\gap);\bull(1.5,1+\gap); \bull(.5,2+\gap+\cw);
		\dfct(.5,0);\dfct(1.5,0);\dfct(.75,1);\dfct(1.25,1);\dfct(.75,1+\gap);\dfct(1.25,1+\gap);\dfct(1.5,2+\gap);\dfct(.75,2+\gap);\idfct(.25,2+\gap+\cw);
	\end{pgfonlayer}
	
\else \ifnum#1=2
	\sh(0,0);
	\boxt(.625,.625,1.75,1)2;
	\def\x{0}\def\y{1}
	\bubble(.625,.625,1.75)2;
	\def\x{1}\def\y{1+\gap}
	\boxslopedl(.75,.25,1)1;
	\def\x{0}\def\y{1+\gap}
	\boxslopedr(.625,.375,1)2;
	\sh(2,1+\gap);
	\bxbdl(1,1)2;
	\begin{pgfonlayer}{foreground}	
		\bull(2,0);\bull(2,1);\dfct(2.5,0);\dfct(1,0);\bull(2.5,1);\dfct(2.25,1);\dfct(2.25,1+\gap);\dfct(1.5,1);\dfct(1.5,1+\gap);\dfct(2.5,2+\gap);\dfct(1.5,2+\gap);\bull(1,1.5+\gap);\bull(1,1+\cw);\bull(.5,1);\bull(.5,1+\gap);\idfct(.75,1+\cw);\idfct(.75,1+\cw+\gap);\idfct(.5,2+\gap);
	\end{pgfonlayer}
	
\else \ifnum#1=3
	\sh(0,0);
	\bxbur(1,1)1;
	\sh(1,0);
	\bxbul(1,1)2;
	\sh(1,1+\gap);
	\bxbdl(1,1)2;
	\sh(0,1+\gap);
	\boxt(.125,.5,.375,1)2;
	\sh(0,2+\gap);
	\bubble(.125,.5,.375)2;
	\begin{pgfonlayer}{foreground}
		\bull(.5,1+\gap); \bull(.5,2+\gap+\cw);
		\dfct(.5,0);\dfct(1.5,0);\dfct(.75,1);\dfct(1.25,1);\dfct(.75,1+\gap);\dfct(1.25,1+\gap);\dfct(1.5,2+\gap);\dfct(.75,2+\gap);\idfct(.25,2+\gap+\cw);
	\end{pgfonlayer}

\else \ifnum#1=4
	\sh(1,0);
	\bxb(1,1)2;
	\sh(0,0);
	\boxt(.125,.5,.375,1)1;
	\sh(0,1);
	\bubble(.125,.5,.375)1;
	\begin{pgfonlayer}{foreground}
		\bull(2,.5);\bull(.5,1+\cw);\dfct(.5,0);\dfct(1.5,1);\dfct(1.5,0);\dfct(.75,1);\idfct(.25,1+\cw);\bull(0,.5);
	\end{pgfonlayer}	

\else \ifnum#1=5 
	\def\x{0}\def\y{0}
	\bubble(.25,.875,.875)2;
	\def\x{0}\def\y{0}
	\boxslopedr(.25,.75,1)2;
	\def\x{1}\def\y{0}
	\boxslopedl(.875,.125,1)1;
	\sh(2,0);
	\bxb(1,1)2;
	\begin{pgfonlayer}{foreground}	
		\bull(0,.5);\bull(1,.5);\dfct(1.5,0);\dfct(2.5,0);\dfct(1.5,1);\dfct(2.5,1);\idfct(.5,\cw);\idfct(.5,1);
	\end{pgfonlayer}
	
\else \ifnum#1=6
	\sh(0,0);
	\boxt(.625,.625,.75,1)2;
	\def\x{0}\def\y{1}
	\bubble(.625,.625,.75)2;
	\def\x{1}\def\y{1+\gap}
	\boxslopedl(.75,.25,1)1;
	\sh(2,0);
	\bxbul(1,1)2;
	\def\x{0}\def\y{1+\gap}
	\boxslopedr(.625,.375,1)2;
	\sh(2,1+\gap);
	\bxbdl(1,1)2;
	\begin{pgfonlayer}{foreground}	
		\bull(2,0);\bull(2,1);\dfct(2.5,0);\dfct(1,0);\bull(2.5,1);\dfct(2.25,1);\dfct(2.25,1+\gap);\dfct(1.5,1);\dfct(1.5,1+\gap);\dfct(2.5,2+\gap);\dfct(1.5,2+\gap);\bull(1,1.5+\gap);\bull(1,1+\cw);\bull(.5,1);\bull(.5,1+\gap);\idfct(.75,1+\cw);\idfct(.75,1+\cw+\gap);\idfct(.5,2+\gap);
	\end{pgfonlayer}
	
\else \ifnum#1=7
	\sh(0,0);
	\bxbt(1,1)1;
	\sh(0,1);
	\boxt(.125,.5,.375,1)2;
	\sh(0,2);
	\bubble(.125,.5,.375)2;
	\sh(1+\gap,0);
	\bxb(1,2)2;
	\begin{pgfonlayer}{foreground}
		\bull(.5,2+\cw);\dfct(.5,0);\dfct(.5,1);\dfct(.75,2);\dfct(1.5+\gap,0);\dfct(1.5+\gap,2);\idfct(.25,2+\cw);
	\end{pgfonlayer}

\else \ifnum#1=8
	\def\x{0}\def\y{0}
	\boxt(.625,.5,.875,1)2;
	\def\x{0}\def\y{1}
	\bubble(.625,.5,.875)2;
	\def\x{0}\def\y{1+\gap}
	\boxslopedr(.625,.375,1)2;
	\def\x{1}\def\y{1+\gap}
	\boxslopedl(.875,.125,1)1;
	\def\x{2+\gap}\def\y{0}
	\bxb(1,2+\gap)2;
	\begin{pgfonlayer}{foreground}	
		\bull(1.5,1);\bull(1,1+\cw);
		\bull(1.5,1+\gap);
		\bull(1,1.5+\gap); \bull(.5,1);\bull(.5,1+\gap);
		\dfct(1,0);\dfct(1.25,1);\dfct(1.25,1+\gap);\dfct(1.5,2+\gap);\dfct(2.5+\gap,0);\dfct(2.5+\gap,2+\gap);\idfct(.75,1+\cw);\idfct(.75,1+\cw+\gap);\idfct(.5,2+\gap);
	\end{pgfonlayer}
	
\fi\fi\fi\fi\fi\fi\fi\fi}

\def\AxiomW#1{ %23
	\ifnum#1=1
	\sh(0,0)
	\boxt(.25,1,.75,1)2;
	\sh(0,1)
	\bubble(.25,1,.75)2;
	\sh(0,1+\gap)
	\draw[fill=\mycolor{2},even odd rule]
		\cusp(0,-1,.25,.75)
		{[shift={(1,\cw)}] -- \cusp(90,-1,.75,.25) 
		-- ++(-1+\cw,0)}
		-- cycle;
	\sh(1,1+\gap+\cw)
	\cornerul(.25,.75,.75-\cw,.25)2;
	\draw[fill=\mycolor{2},even odd rule]
		\cusp(0,1,.25,.75) -- ++(0,1)
		{[shift={(1,1-\cw)}] -- \cusp(180,-1,.25,.75-\cw)}
		{[shift={(\cw,1-\cw-\cw)}] -- \cusp(270,1,.25-\cw,.75-\cw)}
		-- cycle;
	\begin{pgfonlayer}{foreground}
		\bull(.75,1-\cw);\bull(.5,1+\cw);\bull(1,1+\cw);\bull(1.5,1);
		\bull(0,1+\gap);\bull(.5,1+\cw+\gap);\bull(1,1+\cw+\gap);\bull(1.5,1+\gap);\bull(2,1+\gap);
		\bull(1,1.5+\gap);\bull(0,2+\cw+\gap);\bull(1.5,2+\cw+\gap);\bull(2,2+\gap);
	\end{pgfonlayer}
		
\else \ifnum#1=2
	\sh(0,0)
	\boxt(.25,1,.75,1)2;
	\sh(0,1)
	\bubble(.25,1,.75)2;
	\sh(0,1+\gap)
	\draw[fill=\mycolor{2},even odd rule]
		\cusp(0,-1,.25,.75)
		{[shift={(1,\cw)}] -- \cusp(90,-1,.75-\cw,.25-\cw)}
		{[shift={(1-\cw,1-\cw)}] -- \cusp(180,1,.75-\cw,.25)} 
		-- cycle;
	\sh(0,1+\gap+\cw)
	\cornerur(.25,.75,.75-\cw,.25)2;
	\sh(1,1+\gap+\cw)
	\draw[fill=\mycolor{2},even odd rule]
		\cusp(0,1,.25,.75) -- ++(0,1+\cw) -- ++(-1+\cw,0)
		{[shift={(\cw,1)}] -- \cusp(270,1,.25+\cw,.75-\cw)}
		-- cycle;
	\begin{pgfonlayer}{foreground}
		\bull(.75,1-\cw);\bull(.5,1+\cw);\bull(1,1+\cw);\bull(1.5,1);
		\bull(0,1+\gap);\bull(.5,1+\cw+\gap);\bull(1,1+\cw+\gap);\bull(1.5,1+\gap);\bull(2,1+\gap);
		\bull(1,1.5+\gap);\bull(0,2+\gap);\bull(.5,2+\cw+\gap);\bull(2,2+\cw+\gap);
	\end{pgfonlayer}

\else \ifnum#1=3
	\sh(0,0)
	\bubble(.25,1,.75)2;
	\sh(0,0)
	\draw[fill=\mycolor{2},even odd rule]
		\cusp(0,-1,.25,.75)
		{[shift={(1,\cw)}] -- \cusp(90,-1,.75,.25) 
		-- ++(-1+\cw,0)}
		-- cycle;
	\sh(1,\cw)
	\cornerul(.25,.75,.75-\cw,.25)2;
	\draw[fill=\mycolor{2},even odd rule]
		\cusp(0,1,.25,.75) -- ++(0,1)
		{[shift={(1,1-\cw)}] -- \cusp(180,-1,.25,.75-\cw)}
		{[shift={(\cw,1-\cw-\cw)}] -- \cusp(270,1,.25-\cw,.75-\cw)}
		-- cycle;
	\begin{pgfonlayer}{foreground}
		\bull(0,0);\bull(.75,-\cw);\bull(.5,\cw);\bull(1,\cw);\bull(2,0);
		\bull(0,.5+.5*\cw);\bull(0,1+\cw);\bull(1,.5);\bull(2,.5);\bull(1.5,1+\cw);\bull(2,1);
	\end{pgfonlayer}

\else \ifnum#1=4
	\sh(0,0)
	\cornerdl(.75,1.25,.75,.25)2;
	\sh(\cw,\cw)
	\draw[fill=\mycolor{2},even odd rule]
		(\x,\y) -- ++(1-\cw,0)
		{[shift={(1,0)}] -- \cusp(90,-1,.75-\cw,.25-\cw)}
		{[shift={(1-\cw-\cw,1-\cw-\cw)}] -- \cusp(180,1,.75-\cw,.25)} 
		{[shift={(-\cw,1-\cw)}] -- \cusp(270,-1,.75,.25-\cw)}
		-- cycle;
	\sh(0,\cw)
	\cornerur(.25,.75,.75-\cw,.25)2;
	\sh(1,\cw)
	\draw[fill=\mycolor{2},even odd rule]
		\cusp(0,1,.25,.75) -- ++(0,1+\cw) -- ++(-1+\cw,0)
		{[shift={(\cw,1)}] -- \cusp(270,1,.25+\cw,.75-\cw)}
		-- cycle;
	\begin{pgfonlayer}{foreground}
		\bull(1,\cw);\bull(2,0);
		\bull(0,.5);\bull(0,1);\bull(1,.5);\bull(2,.5+.5*\cw);\bull(.5,1+\cw);\bull(2,1+\cw);
	\end{pgfonlayer}

\else \ifnum#1=5 %!% Cuspstretched a la 18-5, though with slight distance modifications.  Handbubbled as well.

	\sh(0,0);
	\bx(1.5,1)2;
	\sh(-1.5,0);
	\begin{scope}[cm={-1,0,0,1,(0,0)}]
	\draw[fill=\mycolor{2},even odd rule]
		{[shift={(0,0)}]
			\cusp(0,1,1+.13,-.13+.5+.5*\cuspwidth+.5*\cuspwidth*\cuspstretch)}
		{[shift={(1+.5+.5*\cuspwidth+.5*\cuspwidth*\cuspstretch,-.5*\cuspwidth)}] 
			-- (\x,\y) -- ++(90:.52-.5*\cusplength) .. controls +(90:\cuspfactor*\cusplength) and +(180+90:\cuspfactor*\cusplength) .. ++($ (90:\cusplength) + (90+1*90:.5*\cuspwidth*\cuspstretch) $) -- ++ (90:.15+.5*\cuspwidth-.5*\cusplength)	}		
		{[shift={(1+.5+.5*\cuspwidth,.15+.52)}] 
			-- \cusp(270,1,.15,.52-.5*\cuspwidth)}
		{[shift={(1+.5,.5*\cuspwidth)}] 
			-- \cusp(180,-1,.5-.13,.13+.5)} -- cycle;
	\end{scope}

	\sh(\cw,.5)
	\cornerul(.25,1,.25,.25)2;
	\sh(0,\cw)
\pgfmathsetmacro{\three}{3*\cw}
\pgfmathsetmacro{\four}{1-6*\cw}
\pgfmathsetmacro{\five}{3*\cw}
\pgfmathsetmacro{\threep}{.5*\cusplength}
\pgfmathsetmacro{\fourp}{1-\cusplength}
\pgfmathsetmacro{\fivep}{.5*\cusplength}

	\draw[fill=\mycolor{2}]
		\bend(90,1,\three,\four,\five) 
		{[shift={(0,1)}] -- \bend(270,1,\threep,\fourp,\fivep)} -- cycle;
	
	\begin{pgfonlayer}{foreground}
		\bull(.75,0);\bull(1.5,0);\bull(1.5,.5);\bull(1.5,1);\bull(.75,1+\cw);
		\bull(\cw,.5+\cw);\bull(-\cw,.5+\cw);\bull(-\cw,.25+\cw);\bull(-\cw,.75+\cw);
		\bull(0,\cw);\bull(-\cw-\cw*\cuspstretch,-\cw);
	\end{pgfonlayer}

\else \ifnum#1=6
	\sh(0,0)
	\sbx(1,1)2;
	\sh(1,0)
	\sbx(1,1)2;
	\sh(0,0)
	\cornerdl(.75,1.25,.25,.25)2;
	\sh(0,.5)
	\cornerul(.25,.75,.25,.25)2;
	\begin{pgfonlayer}{foreground}
		\bull(0,.5);\bull(1,\cw);\bull(2,0);
		\bull(1,.5);\bull(2,.5);\bull(.5,1+\cw);
		\bull(1,1);\bull(2,1);	
	\end{pgfonlayer}
	
\else \ifnum#1=7 %!% Cuspwidth defined by hand.  Could break.
	\def\cuspwidthtmp{\cuspwidth}
	\def\cuspwidth{.26}
	\sh(0,0)
	\bxb(1.5,1)2;
	\draw[fill=\mycolor{2},even odd rule]
		{[shift={(0,0)}]
			\cusp(90,1,.75,.25-\cw)}
		{[shift={(\cw,1-\cw)}] 
			-- \cusp(0,-1,.75+.20*.75-\cw,.75-.25*.75)}
		{[shift={(1.5,1)}] 
			-- \cusp(180,1,.75+.25*.75,.75-.25*.75+\cw)}
		{[shift={(-\cw,1+\cw)}] 
			-- \cusp(270,-1,.25+\cw,.75)} -- cycle;
	\draw (0,0) -- (0,1) -- (1.5,1);
	\begin{pgfonlayer}{foreground}
		\bull(-\cw,1+\cw);\bull(\cw,1-\cw);\bull(.375,1+\cw);\bull(.75,1);\bull(1.125,1);
	\end{pgfonlayer}

	\def\cuspwidth{\cuspwidthtmp}

\else \ifnum#1=8
	\sh(0,0)
	\sbx(1.5,1)2;
	\cornerdl(.25,1.25,.25,.25)2;
	\sh(0,.5)
	\cornerul(.25,.75,.25,.25)2;
	\begin{pgfonlayer}{foreground}
		\bull(0,.5); \bull(1,\cw);\bull(.5,1+\cw);\bull(1,1);\bull(1.5,0);\bull(1.5,.5);\bull(1.5,1);
	\end{pgfonlayer}

\else \ifnum#1=9 
	\sh(0,0)
	\boxt(.25,1,.75,1)2;
	\sh(0,1)
	\boxcornul(1.25,.75,.5,.5)2;
	\bubble(.25,1,.75)2;
	\cornerul(1.25,.75,.5,.5)2;
	\begin{pgfonlayer}{foreground}
		\bull(.75,1-\cw);\bull(.5,1+\cw);\bull(1,1+\cw);\bull(1.5,1);\bull(.5,2+\cw);\bull(1,2);
	\end{pgfonlayer}

\else \ifnum#1=10
	\sh(0,0)
	\boxt(.25,1,.75,1)2;
	\sh(0,1+\gap)
	\sbx(1,1)2;
	\bubble(.25,1,.75)2;
	\cornerul(.25,.75,.5,.5)2;
	\sh(1,1+\gap)
	\boxslopedl(.75,.25,1)2;	
	\begin{pgfonlayer}{foreground}
		\bull(.75,1-\cw);\bull(.5,1+\gap+\cw);
		\bull(.75,1+\gap-\cw);
		\bull(0,1+\gap);
		\bull(1,1.5+\gap);\bull(.5,2+\gap+\cw);
		\bull(1.5,1); \bull(1.5,1+\gap);
	\end{pgfonlayer}

\fi\fi\fi\fi\fi\fi\fi\fi\fi\fi}

\def\AxiomWnew#1{ %23new
\ifnum#1=1
	\sh(0,0);
	\boxt(.625,.5,.875,1)2;
	\sh(0,1);
	\bubble(.625,.5,.875)2;
	\sh(0,1+\gap);
	\boxslopedr(.625,.375,1)2;
	\sh(1,1+\gap);
	\boxslopedl(.875,.125,1)2;
	\sh(1,2+\gap);
	\bubble(.125,.5,.375)2;
	\begin{pgfonlayer}{foreground}
		\bull(.5,1);\bull(.5,1+\gap),\bull(1,1.5+\gap);\bull(1,1+\cw);\bull(1.5,1);\bull(1.5,1+\gap);\bull(1.5,2+\gap+\cw);\dfct(1,0);\dfct(1.25,1);\dfct(1.25,1+\gap);\dfct(1.75,2+\gap);\idfct(.75,1+\cw);\idfct(.75,1+\gap+\cw);\idfct(.5,2+\gap);\idfct(1.25,2+\cw+\gap);
	\end{pgfonlayer}
		
\else \ifnum#1=2
	\sh(0,0);
	\boxt(.625,.5,.875,1)2;
	\sh(0,1);
	\bubble(.625,.5,.875)2;
	\sh(0,1+\gap);
	\boxslopedr(.625,.375,1)2;
	\sh(1,1+\gap);
	\boxslopedl(.875,.125,1)2;
	\sh(0,2+\gap);
	\bubble(.375,.5,.125)2;
	\begin{pgfonlayer}{foreground}
		\bull(.5,1);\bull(.5,1+\gap),\bull(1,1.5+\gap);\bull(1,1+\cw);\bull(1.5,1);\bull(1.5,1+\gap);\bull(.5,2+\gap+\cw);\dfct(1,0);\dfct(1.25,1);\dfct(1.25,1+\gap);\dfct(1.5,2+\gap);\idfct(.75,1+\cw);\idfct(.75,1+\gap+\cw);\idfct(.2,2+\gap);\idfct(.75,2+\cw+\gap);
	\end{pgfonlayer}

\else \ifnum#1=3
	\sh(0,0);
	\bubble(.625,.5,.875)2;
	\sh(0,0);
	\boxslopedr(.625,.375,1)2;
	\sh(1,0);
	\boxslopedl(.875,.125,1)2;
	\sh(1,1);
	\bubble(.125,.5,.375)2;
	\begin{pgfonlayer}{foreground}
		\bull(0,.5),\bull(1,.5);\bull(2,.5);\bull(1.5,1+\cw);\dfct(1.5,0);\dfct(1.75,1);\idfct(.75,\cw);\idfct(.5,1);\idfct(1.25,1+\cw);
	\end{pgfonlayer}

\else \ifnum#1=4
	\sh(0,0);
	\bubble(.625,.5,.875)2;
	\sh(0,0);
	\boxslopedr(.625,.375,1)2;
	\sh(1,0);
	\boxslopedl(.875,.125,1)2;
	\sh(0,1);
	\bubble(.375,.5,.125)2;
	\begin{pgfonlayer}{foreground}
		\bull(0,.5),\bull(1,.5);\bull(2,.5);\bull(.5,1+\cw);\dfct(1.5,0);\dfct(1.5,1);\idfct(.75,\cw);\idfct(.2,1);\idfct(.75,1+\cw);
	\end{pgfonlayer}

\else \ifnum#1=5 
	\def\cuspwidth{.28};	%c%
\sh(0,1);
\bubble(.25,.5,1.75)2;
\sh(0,0);
\bxb(2.5,1)2;
\sh(0,1);
\draw \bend(0,1,.25,1.5,.75);
\def\cuspwidth{.18}; %c% coded by hand!!!
	\begin{pgfonlayer}{foreground}
		\dfct(1.25,0);
		\dfct(2,1); \bull(.615,1+.14);\idfct(.385,1+.14); %c%
		\idfct(1,1);\bull(1.5,1);
	\end{pgfonlayer}

\else \ifnum#1=6
	\sh(0,0);
	\bubble(.625,.5,.875)2;
	\sh(0,0);
	\boxslopedr(.625,.375,1)2;
	\sh(1,0);
	\boxslopedl(.875,.125,1)2;
	\sh(0,1);
	\bubble(.125,.5,.375)2;
	\begin{pgfonlayer}{foreground}
		\bull(0,.5),\bull(1,.5);\bull(2,.5);\bull(.5,1+\cw);\dfct(1.5,0);\dfct(1.5,1);\idfct(.75,\cw);\idfct(.8,1);\idfct(.25,1+\cw);
	\end{pgfonlayer}
		
\else \ifnum#1=7 
	\sh(0,0);
	\bxb(2,1)2;
	\bubble(.25,1,.75)2;
	\sh(0,1);
	\bubble(.25,.5,1.25)2;
	\begin{pgfonlayer}{foreground}
		\bull(1,\cw);\bull(.615,1+\cw);\bull(1.33,1);\dfct(1.5,0);\dfct(1.67,1);\idfct(.5,\cw);\idfct(1,1);\idfct(.385,1+\cw);
	\end{pgfonlayer}
	
\else \ifnum#1=8
	\sh(0,0);
	\boxt(.625,.5,.875,1)2;
	\sh(0,1+\gap);
	\boxb(.625,.5,.875,1)2;
	\sh(0,1);
	\bubble(.625,.5,.875)2;
	\sh(0,2+\gap);
	\bubble(.25,.5,1.25)2;
	\begin{pgfonlayer}{foreground}
		\bull(1,1+\gap+\cw);\bull(.615,2+\gap+\cw);\bull(1.33,2+\gap);\dfct(1.25,1+\gap);\dfct(1.67,2+\gap);\idfct(.75,1+\gap+\cw);\idfct(1,2+\gap);\idfct(.385,2+\gap+\cw);
		\bull(.5,1);\bull(.5,1+\gap);\bull(1.5,1);\bull(1.5,1+\gap);\bull(1,1+\cw);\idfct(.75,1+\cw);\dfct(1,0);\dfct(1.25,1);
	\end{pgfonlayer}

\else \ifnum#1=9 
	\sh(0,0);
	\boxt(.625,.5,.875,1)2;
	\sh(0,1);
	\bubble(.625,.5,.875)2;
	\sh(0,1+\gap);
	\boxslopedr(.625,.375,1)2;
	\sh(1,1+\gap);
	\boxslopedl(.875,.125,1)2;
	\sh(0,2+\gap);
	\bubble(.125,.5,.375)2;
	\begin{pgfonlayer}{foreground}
		\bull(1,1.5+\gap);\bull(.5,2+\gap+\cw);\dfct(1.5,2+\gap);\idfct(.75,1+\gap+\cw);\idfct(.8,2+\gap);\idfct(.25,2+\gap+\cw);
		\dfct(1,0);\dfct(1.25,1);\idfct(.75,1+\cw);\bull(.5,1);\bull(1.5,1);\bull(1.5,1+\gap);\bull(1,1+\cw);\dfct(1.25,1+\gap);\bull(.5,1+\gap);
	\end{pgfonlayer}

\fi\fi\fi\fi\fi\fi\fi\fi\fi}

\def\AxiomWa#1{ %23a
	\ifnum#1=1
		\sh(0,0);
		\bx(3,1)2;
		\halfidr(3,1){.75};
		\sh(0,1);
		\bx(3,1)2;
		\halfidr(3,1){.84};
		\tick(1,2);\tick(2,2);\tick(1.5,1);
\else \ifnum#1=2
		\sh(0,0);
		\bx(3,1)2;
		\halfidr(3,1){.67};
		\sh(0,1);
		\bx(1,1)2;
		\sh(1,1);
		\bx(2,1)2;
		\halfidr(2,1){.75};
		\tick(2,2);
\else \ifnum#1=3
		\sh(0,0);
		\bx(3,1)2;
		\halfidl(3,1){.25};
		\sh(0,1);
		\bx(3,1)2;
		\halfidr(3,1){.84};
		\tick(1,2);\tick(2,2);\tick(1.5,1);
\else \ifnum#1=4
		\sh(0,0);
		\bx(3,1)2;
		\halfidl(3,1){.16};
		\sh(0,1);
		\bx(1,1)2;
		\sh(1,1);
		\bx(2,1)2;
		\halfidr(2,1){.75};
		\tick(2,2);
\else \ifnum#1=5
		\sh(0,0);
		\bx(3,1)2;
		\halfidr(3,1){.84};
		\sh(0,1);
		\bx(2,1)2;
		\halfidl(2,1){.25};
		\sh(2,1);
		\bx(1,1)2;
		\tick(1,2);
\else \ifnum#1=6
		\sh(0,0);
		\bx(3,1)2;
		\halfidr(3,1){.75};
		\sh(0,1);
		\bx(3,1)2;
		\halfidl(3,1){.16};
		\tick(1,2);\tick(2,2);\tick(1.5,1);
\else \ifnum#1=7
		\sh(0,0);
		\bx(3,1)2;
		\halfidl(3,1){.33};
		\sh(0,1);
		\bx(2,1)2;
		\halfidl(2,1){.25};
		\sh(2,1);
		\bx(1,1)2;
		\tick(1,2);
\else \ifnum#1=8
		\sh(0,0);
		\bx(3,1)2;
		\halfidl(3,1){.25};
		\sh(0,1);
		\bx(3,1)2;
		\halfidl(3,1){.16};
		\tick(1,2);\tick(2,2);\tick(1.5,1);
\else \ifnum#1=9
		\sh(0,0);
		\bx(3,1)2;
		\halfidl(3,1){.33};
		\sh(0,1);
		\bx(2,1)2;
		\halfidr(2,1){.75};
		\sh(2,1);
		\bx(1,1)2;
		\tick(1,2);
\else \ifnum#1=10
		\sh(0,0);
		\bx(3,1)2;
		\halfidl(3,1){.16};
		\sh(0,1);
		\bx(1,1)2;
		\sh(1,1);
		\bx(2,1)2;
		\halfidl(2,1){.25};
		\tick(2,2);
\else \ifnum#1=11
		\sh(0,0);
		\bx(3,1)2;
		\halfidr(3,1){.84};
		\sh(0,1);
		\bx(2,1)2;
		\halfidr(2,1){.75};
		\sh(2,1);
		\bx(1,1)2;
		\tick(1,2);
\else \ifnum#1=12
		\sh(0,0);
		\bx(3,1)2;
		\halfidr(3,1){.67};
		\sh(0,1);
		\bx(1,1)2;
		\sh(1,1);
		\bx(2,1)2;
		\halfidl(2,1){.25};
		\tick(2,2)
\fi\fi\fi\fi\fi\fi\fi\fi\fi\fi\fi\fi}

\def\AxiomX#1{ %24
	\ifnum#1=1
	\sh(0,0)
	\boxt(.25,1,.75,1)2;
	\sh(1,1+\gap)
	\sbx(1,1)2;
	\sh(0,1+\gap)
	\boxslopedr(.25,.75,1)2;
	\bubble(.25,1,.75)2;
	\sh(1,1+\gap)
	\cornerur(.25,.75,.5,.5)2;	
	\begin{pgfonlayer}{foreground}	
		\bull(.75,1-\cw); \bull(1.5,1);
		\bull(.75,1-\cw+\gap); \bull(1.5,1+\gap);
		\bull(2,1+\gap);
		\bull(.5,1+\gap+\cw); \bull(1,1+\gap+\cw);
		\bull(1,1.5+\gap); \bull(1.5,2+\gap+\cw);
	\end{pgfonlayer}

\else \ifnum#1=2
	\sh(0,0)
	\boxt(.25,1,.75,1)2;
	\sh(0,1)
	\boxcornur(1.25,.75,.5,.5)2;
	\bubble(.25,1,.75)2;
	\cornerur(1.25,.75,.5,.5)2;
	\begin{pgfonlayer}{foreground}	
		\bull(.75,1-\cw); \bull(1.5,1);
		\bull(.5,1+\cw); \bull(1,1+\cw);
		\bull(1,2); \bull(1.5,2+\cw);
	\end{pgfonlayer}

\else \ifnum#1=3
	\sh(0,0)
	\sbx(1,1)2;
	\sh(1,0)
	\sbx(1,1)2;
	\sh(0,0)
	\cornerdl(.75,1.25,.75,.25)2;
	\sh(0,\cw)
	\boxslopeld(1,.75,.25-\cw)2;
	\sh(1,0)
	\cornerur(.25,.75,.75,.25)2;
	\begin{pgfonlayer}{foreground}
		\bull(2,0);\bull(0,.5);\bull(1,.5);\bull(2,.5);\bull(1.5,1+\cw);
	\end{pgfonlayer}
	
\else \ifnum#1=4
	\sh(0,0)
	\sbx(1.5,1)2;
	\cornerdl(.5,1,.75,.25)2;
	\cornerur(1,.5,.75,.25)2;
	\begin{pgfonlayer}{foreground}
		\bull(0,.5);\bull(1.5,.5);\bull(.75,\cw);\bull(0,1);\bull(.75,1);\bull(1.25,1+\cw);
	\end{pgfonlayer}

\else \ifnum#1=5 %!% Cuspstretched a la 18-5.
	\begin{scope}[cm={-1,0,0,1,(1.5,0)}] % Shift necessary for bounding box
	\sh(0,0)
	\sbx(1.5,1+\cw)2;
	\draw[fill=\mycolor{2},even odd rule]
		{[shift={(0,0)}]
			\cusp(0,1,1,.5+.5*\cuspwidth+.5*\cuspwidth*\cuspstretch)}
		{[shift={(1+.5+.5*\cuspwidth+.5*\cuspwidth*\cuspstretch,-.5*\cuspwidth)}] 
			-- (\x,\y) -- ++(90:.52-.5*\cusplength) .. controls +(90:\cuspfactor*\cusplength) and +(180+90:\cuspfactor*\cusplength) .. ++($ (90:\cusplength) + (90+1*90:.5*\cuspwidth*\cuspstretch) $) -- ++ (90:.15+.5*\cuspwidth-.5*\cusplength)	}		
		{[shift={(1+.5+.5*\cuspwidth,.15+.52)}] 
			-- \cusp(270,1,.15,.52-.5*\cuspwidth)}
		{[shift={(1+.5,.5*\cuspwidth)}] 
			-- \cusp(180,-1,.5,.5)} -- cycle;
	\sh(1.5,\cw);
	\bubblev(.5*\cusplength,1-\cusplength,.5*\cusplength)2;
	\sh(0,0);
	\cornerul(.75,.75,.75+\cw,.25)2;
	\begin{pgfonlayer}{foreground}
		\bull(1.5+\cw+\cw*\cuspstretch,-\cw);\bull(1.5,\cw);
		\bull(1.5,1+\cw);
		\bull(1.5-\cw,.5+\cw);
		\bull(1.5+\cw,.25+\cw);\bull(1.5+\cw,.5+\cw);\bull(1.5+\cw,.75+\cw);
		\bull(0,0);\bull(0,.5+.5*\cw);\bull(.5,1+\cw+\cw);
	\end{pgfonlayer}
	\end{scope}

\else \ifnum#1=6
	\sh(0,0)
	\sbx(1.5,1)2;
	\cornerur(1,.5,.75,.25)2;
	\cornerul(.25,.5,.75,.25)2;
	\begin{pgfonlayer}{foreground}
		\bull(0,0);\bull(0,.5);\bull(1.5,0);\bull(1.5,.5);\bull(.35,1+\cw);\bull(1.5-.35,1+\cw);
	\end{pgfonlayer}

\else \ifnum#1=7 %!% Cuspstretched a la 18-5.
	\sh(0,0)
	\sbx(1.5,1+\cw)2;
	\draw[fill=\mycolor{2},even odd rule]
		{[shift={(0,0)}]
			\cusp(0,1,1,.5+.5*\cuspwidth+.5*\cuspwidth*\cuspstretch)}
		{[shift={(1+.5+.5*\cuspwidth+.5*\cuspwidth*\cuspstretch,-.5*\cuspwidth)}] 
			-- (\x,\y) -- ++(90:.52-.5*\cusplength) .. controls +(90:\cuspfactor*\cusplength) and +(180+90:\cuspfactor*\cusplength) .. ++($ (90:\cusplength) + (90+1*90:.5*\cuspwidth*\cuspstretch) $) -- ++ (90:.15+.5*\cuspwidth-.5*\cusplength)	}		
		{[shift={(1+.5+.5*\cuspwidth,.15+.52)}] 
			-- \cusp(270,1,.15,.52-.5*\cuspwidth)}
		{[shift={(1+.5,.5*\cuspwidth)}] 
			-- \cusp(180,-1,.5,.5)} -- cycle;
	\sh(1.5,\cw);
	\bubblev(.5*\cusplength,1-\cusplength,.5*\cusplength)2;
	\sh(0,0);
	\cornerul(.75,.75,.75+\cw,.25)2;
	\begin{pgfonlayer}{foreground}
		\bull(1.5+\cw+\cw*\cuspstretch,-\cw);\bull(1.5,\cw);
		\bull(1.5,1+\cw);
		\bull(1.5-\cw,.5+\cw);
		\bull(1.5+\cw,.25+\cw);\bull(1.5+\cw,.5+\cw);\bull(1.5+\cw,.75+\cw);
		\bull(0,0);\bull(0,.5+.5*\cw);\bull(.5,1+\cw+\cw);
	\end{pgfonlayer}
\else \ifnum#1=8
	\sh(0,0)
	\sbx(1.5,1)2;
	\cornerdr(.5,1,.75,.25)2;
	\cornerul(1,.5,.75,.25)2;
	\begin{pgfonlayer}{foreground}
		\bull(0,0);\bull(0,.5);\bull(1.5,.5);\bull(.75,\cw);\bull(1.5,1);\bull(.75,1);\bull(.25,1+\cw);
	\end{pgfonlayer}

\else \ifnum#1=9
	\sh(0,0)
	\sbx(1,1)2;
	\sh(1,0)
	\sbx(1,1)2;
	\sh(0,0)
	\cornerdr(.75,1.25,.75,.25)2;
	\sh(1,\cw)
	\boxsloperd(1,.75,.25-\cw)2;
	\sh(0,0)
	\cornerul(.25,.75,.75,.25)2;
	\begin{pgfonlayer}{foreground}
		\bull(0,0);\bull(0,.5);\bull(1,.5);\bull(2,.5);\bull(.5,1+\cw);
	\end{pgfonlayer}

\else \ifnum#1=10
	\sh(0,0)
	\boxt(.75,1,.25,1)2;
	\sh(0,1)
	\boxcornul(1.25,.75,.5,.5)2;
	\bubble(.75,1,.25)2;
	\cornerul(1.25,.75,.5,.5)2;
	\begin{pgfonlayer}{foreground}	
		\bull(1.25,1-\cw); \bull(.5,1);
		\bull(1.5,1+\cw); \bull(1,1+\cw);
		\bull(1,2); \bull(.5,2+\cw);
	\end{pgfonlayer}

\else \ifnum#1=11
	\sh(0,0)
	\boxt(.75,1,.25,1)2;
	\sh(0,1+\gap)
	\sbx(1,1)2;
	\sh(1,1+\gap)
	\boxslopedl(.25,.75,1)2;
	\sh(0,1+\gap)
	\bubble(.75,1,.25)2;
	\sh(0,1+\gap)
	\cornerul(.25,.75,.5,.5)2;	
	\begin{pgfonlayer}{foreground}	
		\bull(1.25,1-\cw); \bull(.5,1);
		\bull(1.25,1-\cw+\gap); \bull(.5,1+\gap);
		\bull(2,1+\gap); \bull(0,1+\gap);
		\bull(1.5,1+\gap+\cw); \bull(1,1+\gap+\cw);
		\bull(1,1.5+\gap); \bull(.5,2+\gap+\cw);
	\end{pgfonlayer}

\fi\fi\fi\fi\fi\fi\fi\fi\fi\fi\fi}

\def\AxiomXnew#1{ %24new
	\ifnum#1=1
	\sh(0,0);
	\boxt(.625,.5,.875,1)2;
	\sh(0,1);
	\bubble(.625,.5,.875)2;
	\sh(0,1+\gap);
	\boxslopedr(.625,.375,1)2;
	\sh(1,1+\gap);
	\boxslopedl(.875,.125,1)2;
	\sh(1,2+\gap);
	\bubble(.375,.5,.125)2;
	\begin{pgfonlayer}{foreground}
		\bull(.5,1);\bull(.5,1+\gap),\bull(1,1.5+\gap);\bull(1,1+\cw);\bull(1.5,1);\bull(1.5,1+\gap);\bull(1.5,2+\gap+\cw);\dfct(1,0);\dfct(1.25,1);\dfct(1.25,1+\gap);\dfct(1.25,2+\gap);\idfct(.75,1+\cw);\idfct(.75,1+\gap+\cw);\idfct(.5,2+\gap);\idfct(1.75,2+\cw+\gap);
	\end{pgfonlayer}

\else \ifnum#1=2
	\sh(0,0);
	\boxt(.625,.5,.875,1)2;
	\sh(0,1+\gap);
	\boxb(.625,.5,.875,1)2;
	\sh(0,1);
	\bubble(.625,.5,.875)2;
	\sh(0,2+\gap);
	\bubble(1.25,.5,.25)2;
	\begin{pgfonlayer}{foreground}
		\bull(1,1+\gap+\cw);\bull(2-.615,2+\gap+\cw);\bull(2-1.33,2+\gap);\dfct(1.25,1+\gap);\idfct(.33,2+\gap);\idfct(.75,1+\gap+\cw);\dfct(1,2+\gap);\idfct(2-.385,2+\gap+\cw);
		\bull(.5,1);\bull(.5,1+\gap);\bull(1.5,1);\bull(1.5,1+\gap);\bull(1,1+\cw);\idfct(.75,1+\cw);\dfct(1,0);\dfct(1.25,1);
	\end{pgfonlayer}

\else \ifnum#1=3
	\sh(0,0);
	\bubble(.625,.5,.875)2;
	\sh(0,0);
	\boxslopedr(.625,.375,1)2;
	\sh(1,0);
	\boxslopedl(.875,.125,1)2;
	\sh(1,1);
	\bubble(.375,.5,.125)2;
	\begin{pgfonlayer}{foreground}
		\bull(0,.5),\bull(1,.5);\bull(2,.5);\bull(1.5,1+\cw);\dfct(1.5,0);\dfct(1.25,1);\idfct(.75,\cw);\idfct(.5,1);\idfct(1.75,1+\cw);
	\end{pgfonlayer}
	
\else \ifnum#1=4
	\sh(0,0);
	\bxb(2,1)2;
	\bubble(.25,1,.75)2;
	\sh(0,1);
	\bubble(1.25,.5,.25)2;
	\begin{pgfonlayer}{foreground}
		\bull(1,\cw);\bull(2-.615,1+\cw);\bull(2-1.33,1);\dfct(1.5,0);\idfct(2-1.67,1);\idfct(.5,\cw);\dfct(1,1);\idfct(2-.385,1+\cw);
	\end{pgfonlayer}

\else \ifnum#1=5 
	\sh(0,0);
	\bxb(2,1)2;
	\sh(0,1);
	\bubble(.25,.5,.25)2;
	\sh(1,1);
	\bubble(.25,.5,.25)2;
	\begin{pgfonlayer}{foreground}
		\idfct(.385,1+\cw);\bull(.615,1+\cw);\dfct(1,0);\dfct(1,1);
		\idfct(2-.385,1+\cw);\bull(2-.615,1+\cw);
	\end{pgfonlayer}

\else \ifnum#1=6
	\sh(0,0);
	\bxb(2,1)2;
	\bubble(.75,1,.25)2;
	\sh(0,1);
	\bubble(.25,.5,1.25)2;
	\begin{pgfonlayer}{foreground}
		\bull(1,\cw);\bull(.615,1+\cw);\bull(1.33,1);\dfct(.5,0);\idfct(1.67,1);\idfct(1.5,\cw);\dfct(1,1);\idfct(.385,1+\cw);
	\end{pgfonlayer}

\else \ifnum#1=7 
	\sh(0,0);
	\bubble(.875,.5,.625)2;
	\sh(1,0);
	\boxslopedl(.625,.375,1)2;
	\sh(0,0);
	\boxslopedr(.875,.125,1)2;
	\sh(0,1);
	\bubble(.125,.5,.375)2;
	\begin{pgfonlayer}{foreground}
		\bull(2,.5),\bull(1,.5);\bull(0,.5);\bull(.5,1+\cw);\dfct(.5,0);\dfct(.75,1);\idfct(1.25,\cw);\idfct(1.5,1);\idfct(.25,1+\cw);
	\end{pgfonlayer}

\else \ifnum#1=8
	\sh(0,0);
	\boxt(.875,.5,.625,1)2;
	\sh(0,1+\gap);
	\boxb(.875,.5,.625,1)2;
	\sh(0,1);
	\bubble(.875,.5,.625)2;
	\sh(0,2+\gap);
	\bubble(.25,.5,1.25)2;
	\begin{pgfonlayer}{foreground}
		\bull(1,1+\gap+\cw);\bull(.615,2+\gap+\cw);\bull(1.33,2+\gap);\dfct(2-1.25,1+\gap);\idfct(2-.33,2+\gap);\idfct(2-.75,1+\gap+\cw);\dfct(1,2+\gap);\idfct(.385,2+\gap+\cw);
		\bull(1.5,1);\bull(1.5,1+\gap);\bull(.5,1);\bull(.5,1+\gap);\bull(1,1+\cw);\idfct(2-.75,1+\cw);\dfct(1,0);\dfct(2-1.25,1);
	\end{pgfonlayer}

\else \ifnum#1=9
	\sh(0,0);
	\boxt(.875,.5,.625,1)2;
	\sh(0,1);
	\bubble(.875,.5,.625)2;
	\sh(0,1+\gap);
	\boxslopedr(.875,.125,1)2;
	\sh(1,1+\gap);
	\boxslopedl(.625,.375,1)2;
	\sh(0,2+\gap);
	\bubble(.125,.5,.375)2;
	\begin{pgfonlayer}{foreground}
		\bull(1.5,1);\bull(1.5,1+\gap),\bull(1,1.5+\gap);\bull(1,1+\cw);\bull(.5,1);\bull(.5,1+\gap);\bull(.5,2+\gap+\cw);\dfct(1,0);\dfct(2-1.25,1);\dfct(2-1.25,1+\gap);\dfct(2-1.25,2+\gap);\idfct(2-.75,1+\cw);\idfct(2-.75,1+\gap+\cw);\idfct(1.5,2+\gap);\idfct(2-1.75,2+\cw+\gap);
	\end{pgfonlayer}

\fi\fi\fi\fi\fi\fi\fi\fi\fi}

\def\AxiomY#1{
	\ifnum#1=1
	\sh(0,0);
	\boxt(.25*.75,.75,2+.75*.75,1)2;
	\sh(0,1);
	\bubble(.25*.75,.75,2+.75*.75)2;
	\begin{pgfonlayer}{foreground}	
		\bull(0,.5);\bull(3.5,.5);\bull(.75,1+\cw);
		\bull(1.5,0);\bull(2.5,0);\bull(1.5,1);\bull(2.5,1);
		\dfct(.75,0);\dfct(2,0);\dfct(3,0);
		\idfct(.375,1+\cw);\dfct(1.125,1);\dfct(2,1);\dfct(3,1);
	\end{pgfonlayer}

\else \ifnum#1=2
	\sh(0,0);
	\boxt(.25*.75,.75,.75*.75,1)2;
	\sh(0,1);
	\bubble(.25*.75,.75,.75*.75)2;
	\sh(1.5,0);
	\bxb(2,1)2;
	\begin{pgfonlayer}{foreground}	
		\bull(0,.5);\bull(.75,1+\cw);
		\bull(2.5,0);\bull(2.5,1);
		\dfct(.75,0);\dfct(2,0);\dfct(3,0);
		\idfct(.375,1+\cw);\dfct(1.125,1);\dfct(2,1);\dfct(3,1);
	\end{pgfonlayer}

\else \ifnum#1=3
	\sh(0,0);
	\boxt(.25*.75,.75,1+.75*.75,1)2;
	\sh(0,1);
	\bubble(.25*.75,.75,1+.75*.75)2;
	\sh(2.5,0);
	\bxb(1,1)2;
	\begin{pgfonlayer}{foreground}	
		\bull(0,.5);\bull(.75,1+\cw);
		\bull(1.5,0);\bull(1.5,1);
		\dfct(.75,0);\dfct(2,0);\dfct(3,0);
		\idfct(.375,1+\cw);\dfct(1.125,1);\dfct(2,1);\dfct(3,1);
	\end{pgfonlayer}

\else \ifnum#1=4
	\sh(0,0);
	\boxt(.25*.75,.75,.75*.75,1)2;
	\sh(0,1);
	\bubble(.25*.75,.75,.75*.75)2;
	\sh(1.5,0);
	\bxb(1,1)2;
	\sh(2.5,0);
	\bxb(1,1)2;
	\begin{pgfonlayer}{foreground}	
		\bull(0,.5);\bull(.75,1+\cw);
		\dfct(.75,0);\dfct(2,0);\dfct(3,0);
		\idfct(.375,1+\cw);\dfct(1.125,1);\dfct(2,1);\dfct(3,1);
	\end{pgfonlayer}

\fi\fi\fi\fi}

\def\AxiomZ#1{
	\ifnum#1=1
	\sh(0,0);
	\boxt(.75*.75,.75,.25*.75,1)2;
	\sh(0,1);
	\bubble(.75*.75,.75,.25*.75)2;
	\sh(1.5,0);
	\bxb(1,1)2;
	\sh(2.5,0);
	\bxb(1,1)2;
	\begin{pgfonlayer}{foreground}	
		\bull(0,.5);\bull(1.5,.5);\bull(.75,1+\cw);
		\dfct(.75,0);\dfct(2,0);\dfct(3,0);
		\dfct(.375,1);\idfct(1.125,1+\cw);\dfct(2,1);\dfct(3,1);
	\end{pgfonlayer}

\else \ifnum#1=2
	\sh(1,0);
	\boxt(.25*.75,.75,.75*.75,1)2;
	\sh(1,1);
	\bubble(.25*.75,.75,.75*.75)2;
	\sh(0,0);
	\bxb(1,1)2;
	\sh(2.5,0);
	\bxb(1,1)2;
	\begin{pgfonlayer}{foreground}	
		\bull(1.75,1+\cw);
		\dfct(1.75,0);\dfct(.5,0);\dfct(3,0);
		\dfct(2.125,1);\idfct(1.375,1+\cw);\dfct(.5,1);\dfct(3,1);
	\end{pgfonlayer}

\else \ifnum#1=3
	\sh(0,0);
	\boxt(.75*.75,.75,.25*.75,1)2;
	\sh(0,1);
	\bubble(.75*.75,.75,.25*.75)2;
	\sh(1.5,0);
	\bxb(2,1)2;
	\begin{pgfonlayer}{foreground}	
		\bull(0,.5);\bull(1.5,.5);\bull(.75,1+\cw);\bull(2.5,0);\bull(2.5,1);
		\dfct(.75,0);\dfct(2,0);\dfct(3,0);
		\dfct(.375,1);\idfct(1.125,1+\cw);\dfct(2,1);\dfct(3,1);
	\end{pgfonlayer}

\else \ifnum#1=4
	\sh(1,0);
	\boxt(.25*.75,.75,1+.75*.75,1)2;
	\sh(1,1);
	\bubble(.25*.75,.75,1+.75*.75)2;
	\sh(0,0);
	\bxb(1,1)2;
	\begin{pgfonlayer}{foreground}	
		\bull(1.75,1+\cw);\bull(2.5,0);\bull(2.5,1);\bull(3.5,.5);
		\dfct(1.75,0);\dfct(.5,0);\dfct(3,0);
		\dfct(2.125,1);\idfct(1.375,1+\cw);\dfct(.5,1);\dfct(3,1);
	\end{pgfonlayer}

\fi\fi\fi\fi}

\def\axgridE{
\path
	(0,3) 	coordinate (1)
	(3,3) 	coordinate (2)
	(2,2) 	coordinate (3)
	(1,1) 	coordinate (4)
	(0,0) 	coordinate (5)
	(3,0) 	coordinate (6);
\draw (5) -- (1) -- (2) -- (6) -- (5) -- (4) -- (3) -- (2)   (6) -- (4);
}

\def\axgridF{
\path
	(0,2) 	coordinate (1)
	(2,2) 	coordinate (2)
	(2,1) 	coordinate (3)
	(0,0) 	coordinate (4)
	(2,0) 	coordinate (5);
\draw  (1) -- (3) -- (2) -- (1) -- (4) -- (5) -- (3) -- (4);
}

\def\axgridFa{
\path
	(-1,3) 	coordinate (1)
	(4,3) 	coordinate (2)
	(2,2) 	coordinate (3)
	(3,.8) 	coordinate (4)
	(4,-1) 	coordinate (5);
\draw  (1) -- (2) -- (3) -- (4) -- (5) -- (1) -- (3) (2) -- (5);
}

\def\axgridG{
\path
	(0,10) 	coordinate (1)
	(10,10) 	coordinate (2)
	(1.5,8.5) 	coordinate (3)
	(5,9) 	coordinate (4)
	(3,7) 	coordinate (5)
	(7,7) 	coordinate (6)
	(4.3,5.7) 	coordinate (7)
	(1,5) 	coordinate (8)
	(9,5) 	coordinate (9)
	(5.7,4.3) 	coordinate (10)
	(3,3) 	coordinate (11)
	(7,3) 	coordinate (12)
	(5,1) 	coordinate (13)
	(8.5,1.5) 	coordinate (14)
	(0,0) 	coordinate (15)
	(10,0) 	coordinate (16);
\draw (1) -- (2) -- (16) -- (15) -- (1) -- (3) -- (4) -- (2) -- (9) -- (6) -- (4) -- (5) -- (8) -- (11) -- (7) -- (6) -- (10) -- (11) -- (13) -- (12) -- (9) -- (14) -- (13) -- (15) -- (8) -- (3) (5) -- (7) (10) -- (12) (14) -- (16);
}

\def\axgridH{
\path
	(3,4) 	coordinate (1)
	(0.2,2.5) 	coordinate (2)
	(2.1,2.3) 	coordinate (3)
	(3.9,2.3) 	coordinate (4)
	(5.8,2.5) 	coordinate (5)
	(3,.8) 	coordinate (6)
	(1.4,0) 	coordinate (7)
	(4.6,0) 	coordinate (8);
\draw (2) -- (1) -- (5) -- (4) -- (1) -- (3) -- (2) -- (7) -- (6) -- (8) -- (5) (3) -- node[above right,font=\footnotesize] {*} (6) -- node[above left,font=\footnotesize] {*} (4) (7) -- (8);
}

\def\axgridHs{
\path
	(0,0) 	coordinate (1)
	(1,0) 	coordinate (2)
	(2,0) 	coordinate (3)
	(3,0) 	coordinate (4);
\draw  (1) -- (2) -- (3) -- (4);
}

\def\axgridHa{
\path
	(0,7) 	coordinate (1)
	(7,7) 	coordinate (2)
	(3,6) 	coordinate (3)
	(6,6) 	coordinate (4)
	(4.8,5) 	coordinate (5)
	(2,4) 	coordinate (6)
	(3.8,4) 	coordinate (7)
	(1,3) 	coordinate (8)
	(4,2.5) 	coordinate (9)
	(6,2.2) 	coordinate (10)
	(3,1.5) 	coordinate (11)
	(1,1) 	coordinate (12)
	(0,0) 	coordinate (13)
	(7,0) 	coordinate (14);
\draw  (1) -- (2) -- (4) -- (3) -- (1) -- (8) -- (6) -- (3) -- (5) -- (10) -- (7) -- (6) -- (9) -- (10) -- (14) -- (11) -- (8) -- (1) -- (13) -- (12) -- (14) -- (13) (2) -- (14) (4) -- (10) (5) -- (7) (8) -- (12) (9) -- (11);
}

\def\axgridHb{
\path
	(3,6) 	coordinate (1)
	(2,3.8) 	coordinate (2)
	(4,3.8) 	coordinate (3)
	(0,3) 	coordinate (4)
	(3,2.9) 	coordinate (5)
	(6,3) 	coordinate (6)
	(2,2) 	coordinate (7)
	(4,2) 	coordinate (8)
	(3,1.1) 	coordinate (9)
	(3,0) 	coordinate (10);
\draw (10) -- (4) -- (1) -- (6) -- (10) -- (9) -- (7) -- (4) -- (2) -- (1) -- (3) -- (6) -- (8) -- (9)  (2) -- (5) -- (8)  (3) -- (5) -- (7);
}

\def\axgridK{
\path
	(3,5.2) 	coordinate (1)
	(3,4) 	coordinate (2)
	(-.2,3) 	coordinate (3)
	(1.5,3) 	coordinate (4)
	(4.5,3) 	coordinate (5)
	(6.2,3) 	coordinate (6)
	(2.2,2) 	coordinate (7)
	(3.8,2) 	coordinate (8)
	(2.2,1) 	coordinate (9)
	(3.8,1) 	coordinate (10)
	(1,0) 	coordinate (11)
	(5,0) 	coordinate (12);
\draw (1) -- (6) -- (12) -- (10) -- (8) -- (5) -- (2) -- (4) -- (7) -- (9) -- (11) -- (3) -- (1) -- (2) (3) -- (4) (5) -- (6) (7) -- (8) (9) -- (10) (11) -- (12);
}

\def\axgridL{
\path
	(0,3.9) 	coordinate (1)
	(4,3.9) 	coordinate (2)
	(.8,3.1) 	coordinate (3)
	(1.9,2.9) 	coordinate (4)
	(3.2,3) 	coordinate (5)
	(.9,2.1) 	coordinate (6)
	(1.7,2.2) 	coordinate (7)
	(2.4,1.7) 	coordinate (8)
	(3.2,1.9) 	coordinate (9)
	(0.7,1) 	coordinate (10)
	(2.1,1) 	coordinate (11)
	(0,0.1) 	coordinate (12)
	(4,0.1) 	coordinate (13);
\draw (1) -- (3) -- (6) -- (10) -- (12) -- (1) -- (2) -- (5) -- (4) -- (3) (4) -- (7) -- (6) (8) -- (7) (5) -- (9) -- (8) -- (11) -- (10) (2) -- (13) -- (9) (12) -- (13) -- (11);
}

\def\axgridM{
\path
	(2,4) 	coordinate (1)
	(0,3) 	coordinate (2)
	(4,3) 	coordinate (3)
	(2,2) 	coordinate (4)
	(0,1) 	coordinate (5)
	(4,1) 	coordinate (6)
	(2,0) 	coordinate (7);
\draw (2) -- (4) -- (3) -- (1) -- (2) -- (5) -- (4) -- (6) -- (3) (5) -- (7) -- (6);
}

\def\axgridMa{
\path
	(2,4) 	coordinate (1)
	(2,3) 	coordinate (2)
	(0,2) 	coordinate (3)
	(1.2,2) 	coordinate (4)
	(2.8,2) 	coordinate (5)
	(4,2) 	coordinate (6)
	(2,.5) 	coordinate (7);
\draw (1) -- (3) -- (4) -- (2) -- (1) -- (6) -- (5) -- (2) (3) -- (7) -- (4) (5) -- (7) -- (6);
}

\def\axgridN{
\path
	(0,1) 	coordinate (1)
	(1,1) 	coordinate (2)
	(.5,.5) 	coordinate (3)
	(0,0) 	coordinate (4)
	(1,0) 	coordinate (5);
\draw (1) -- (2) -- (3) -- (4) -- (5) -- (3) -- (1) -- (4) (2) -- (5);
}

\def\axgridNa{
\path
	(0,7) 	coordinate (1)
	(6,7) 	coordinate (2)
	(3,6.5) 	coordinate (3)
	(1.2,4.5) 	coordinate (4)
	(2.2,5.5) 	coordinate (5)
	(3,4.5) 	coordinate (6)
	(3.8,5.5) 	coordinate (7)
	(4.8,4.5) 	coordinate (8)
	(3,3.5) 	coordinate (9)
	(2.5,2.2) 	coordinate (10)
	(3.5,2.2) 	coordinate (11)
	(3,0) 	coordinate (12);
\draw (4) -- (1) to [bend left=15] (2) -- (3) -- (1) to [bend right=20] (12) -- (10) to [bend left=9] (5) -- (3) -- (7) -- (8) -- (2) to [bend left=20] (12) -- (11) -- (9) -- (6) (9) -- (10) (7) -- (6) -- (5) -- (4) to [bend right=10] (12) to [bend right=10] (8) (7) to [bend left=9] (11);
}

\def\axgridO{
\path
	(-.15,1) 	coordinate (1)
	(1.65,1) 	coordinate (2)
	(-.15,0) 	coordinate (3)
	(1.65,0) 	coordinate (4);
\draw (1) -- (2) -- (4) -- (3) -- (1);
}

\def\axgridP{
\path
	(-.1,0) 	coordinate (1)
	(1.6,0) 	coordinate (2);
\draw (1) -- (2);
}

\def\axgridR{
\path
	(0,3) 	coordinate (1)
	(3,3) 	coordinate (2)
	(1,2) 	coordinate (3)
	(1,1) 	coordinate (4)
	(2,1) 	coordinate (5)
	(0,0) 	coordinate (6)
	(3,0) 	coordinate (7);
\draw (1) -- (2) -- (7) -- (5) -- (3) -- (1) -- (6) -- (7) (3) -- (4) -- (5) (6) -- (4);
}

\def\axgridRnew{
\path
	(0,3) 	coordinate (1)
	(3,3) 	coordinate (2)
	(2,2) 	coordinate (3)
	(1,1) 	coordinate (4)
	(0,0) 	coordinate (5)
	(3,0) 	coordinate (6);
\draw (1) -- (2) -- (6) -- (5) -- (1) -- (3) -- (6) (3) -- (4) -- (5) ;
}

\def\axgridS{
\path
	(0,3) 	coordinate (1)
	(3,3) 	coordinate (2)
	(1,2.1) 	coordinate (3)
	(2,2.1) 	coordinate (4)
	(1,.9) 	coordinate (5)
	(2,.9) 	coordinate (6)
	(0,0) 	coordinate (7)
	(3,0) 	coordinate (8);
\draw (1) -- (2) -- (4) -- (3) -- (1) -- (7) -- (8) -- (2) (3) -- (5) -- (6) -- (4) (5) -- (7) (6) -- (8);
}

\def\axgridSnew{
\path
	(0,3) 	coordinate (1)
	(3,3) 	coordinate (2)
	(1,2.1) 	coordinate (3)
	(1,.9) 	coordinate (4)
	(2,.9) 	coordinate (5)
	(0,0) 	coordinate (6)
	(3,0) 	coordinate (7);
\draw (1) -- (2) -- (7) -- (6) -- (1) -- (3) -- (4) -- (5) -- (7) 
   (6) -- (4);
}

\def\axgridT{
\path
	(.8,4) 	coordinate (1)
	(3,4) 	coordinate (2)
	(1.2,3.2) 	coordinate (3)
	(2.7,3.2) 	coordinate (4)
	(1.58,2.55) 	coordinate (5)
	(2.5,2.7) 	coordinate (6)
	(0,2) 	coordinate (7)
	(.7,2) 	coordinate (8)
	(1.4,2) 	coordinate (9)
	(3.8,2) 	coordinate (10)
	(2,1.5) 	coordinate (11)
	(1.1,1.2) 	coordinate (12)
	(1.9,.9) 	coordinate (13)
	(2.6,.6) 	coordinate (14)
	(3.3,.3) 	coordinate (15)
	(1,0) 	coordinate (16)
	(4,0) 	coordinate (17);
\draw (1) -- (2) (7) -- (1) -- (3) -- (4) -- (2) -- (10) -- (17) (3) -- (6) -- (4) (12) -- (8) -- (3) -- (5) -- (6) -- (11) -- (9) (5) -- (9) -- (8) -- (7) -- (16) -- (12) -- (11) -- (13) -- (16) -- (17) -- (15) -- (14) -- (13);
}

\def\axgridTnew{
\path
	(0,4) 	coordinate (1)
	(3,4) 	coordinate (2)
	(.9,3) 	coordinate (3)
	(2.1,3.4) 	coordinate (4)
	(0,2) 	coordinate (5)
	(.9,2) 	coordinate (6)
	(3,2) 	coordinate (7)
	(1,1) 	coordinate (8)
	(2,1) 	coordinate (9)
	(0,0) 	coordinate (10)
	(3,0) 	coordinate (11);
\draw (1) -- (2) -- (7) -- (11) -- (10) -- (5) -- (1) -- (3) -- (4) -- (1) 
 (4) -- (8) -- (6) -- (3)  (5) -- (6) (10) -- (8) -- (9) -- (11) (2) -- (4);
}

\def\axgridTa{
\path
	(1.9,4.5) 	coordinate (1)
	(1.9,3.8) 	coordinate (2)
	(.4,3.1) 	coordinate (3)
	(1,2.8) 	coordinate (4)
	(2.1,3.1) 	coordinate (5)
	(2.9,2.9) 	coordinate (6)
	(4,3.2) 	coordinate (7)
	(2,2) 	coordinate (8)
	(2.8,2) 	coordinate (9)
	(3.5,2) 	coordinate (10)
	(.4,1) 	coordinate (11)
	(1,1.4) 	coordinate (12)
	(2,1) 	coordinate (13)
	(2.8,1.2) 	coordinate (14)
	(4,1) 	coordinate (15)
	(2,0) 	coordinate (16);
\draw (1) -- (7) -- (15) -- (16) -- (11) -- (3) -- (1) -- (2) -- (5) -- (6) -- (7) -- (10) -- (15) -- (14) -- (13) -- (16) (13) -- (12) -- (11) (12) -- (4) -- (3) (4) -- (2) (5) -- (12) -- (8) -- (6) -- (9) -- (10) (9) -- (14) -- (8);
}

\def\axgridU{
\path
	(0,6) 	coordinate (1)
	(0,0) 	coordinate (2)
	(.9,4.6) 	coordinate (3)
	(1.1,1) 	coordinate (4)
	(1.4,3.9) 	coordinate (5)
	(1.8,2.7) 	coordinate (6)
	(2.3,.6) 	coordinate (7)
	(4,6) 	coordinate (8)
	(3.5,5) 	coordinate (9)
	(3.1,4) 	coordinate (10)
	(3,3) 	coordinate (11)
	(3.1,2) 	coordinate (12)
	(3.4,1) 	coordinate (13)
	(4,0) 	coordinate (14);
\draw (1) -- (2) -- (3) -- (1) -- (8) -- (3) -- (4) -- (7) -- (14) -- (2) (3) -- (5) -- (6) -- (7) (5) -- (9) -- (8) -- (14) (4) -- (6) -- (10) -- (9) (6) -- (11) -- (10) (11) -- (12) -- (13) -- (8) (7) -- (6) -- (13) -- (14);
}

\def\axgridUnew{
\path
	(0,2) 	coordinate (1)
	(2,2) 	coordinate (2)
	(.9,1.3) 	coordinate (3)
	(1.2,.7) 	coordinate (4)
	(0,0) 	coordinate (5)
	(2,0) 	coordinate (6);
\draw (1) -- (2) -- (6) -- (5) -- (1) -- (3) -- (4) -- (6) (2) -- (3) -- (5); 
}

\def\axgridUa{
\path
	(0,1) 	coordinate (1)
	(1.5,1) 	coordinate (2)
	(0,0) 	coordinate (3)
	(1.5,0) 	coordinate (4);
\draw (1) -- (2) -- (4) -- (3) -- (1);
}

\def\axgridUanew{
\path
	(0,2.6) 	coordinate (1)
	(3.6,2.6) 	coordinate (2)
	(.6,1.8) 	coordinate (3)
	(2.6,1.8) 	coordinate (4)
        (.6,.8) 	coordinate (5)
        (1.6,.8) 	coordinate (6)
        (2.6,.8) 	coordinate (7)
        (0,0) 	coordinate (8)
        (3.6,0) 	coordinate (9)

;
\draw (1) -- (2) -- (9) -- (8) -- (1) -- (4) -- (7) -- (9)
     (1) -- (3) -- (7) -- (6) -- (5) -- (8) (3) -- (5) (2) -- (4);
}

\def\axgridV{
\path
	(0,2) 	coordinate (1)
	(1.6,2) 	coordinate (2)
	(-.9,1) 	coordinate (3)
	(.1,1) 	coordinate (4)
	(.8,1) 	coordinate (5)
	(1.5,1) 	coordinate (6)
	(2.5,1) 	coordinate (7)
	(0,0) 	coordinate (8)
	(1.6,0) 	coordinate (9);
\draw (1) -- (2) -- (6) -- (5) -- (4) -- (1) -- (3) -- (8) -- (4) (8) -- (9) -- (7) -- (2) (6) -- (9);
}

\def\axgridVnew{
\path
	(0,2) 	coordinate (1)
	(1.3,2) 	coordinate (2)
	(-.9,1) 	coordinate (3)
	(.1,1) 	coordinate (4)
	(1.2,1) 	coordinate (5)
	(2.2,1) 	coordinate (6)
	(0,0) 	coordinate (7)
	(1.3,0) 	coordinate (8);
\draw (1) -- (2) -- (6) -- (8) -- (7) -- (3) -- (1) -- (4) -- (5) -- (8) 
   (2) -- (5)  (4) -- (7);
}

\def\axgridW{
\path
	(0,4) 	coordinate (1)
	(3,4) 	coordinate (2)
	(1,3) 	coordinate (3)
	(2,3) 	coordinate (4)
	(1,2) 	coordinate (5)
	(2.2,2) 	coordinate (6)
	(1,1) 	coordinate (7)
	(1.8,1) 	coordinate (8)
	(0,0) 	coordinate (9)
	(3,0) 	coordinate (10);
\draw (1) -- (3) (2) -- (4) (8) -- (9) -- (7) (6) -- (8) -- (7) -- (5) -- (3) -- (4) -- (6) -- (10) -- (2) -- (1) -- (9) -- (10);
}

\def\axgridWnew{
\path
	(0,4) 	coordinate (1)
	(3,4) 	coordinate (2)
	(1,3) 	coordinate (3)
	(2,3) 	coordinate (4)
	(1,2) 	coordinate (5)
	(2,2) 	coordinate (6)
	(1.6,1.2) 	coordinate (7)
	(0,.6) 	coordinate (8)
	(3,.6) 	coordinate (9);
\draw (1) -- (2) -- (9) -- (8) -- (1) -- (3) -- (4) -- (6) -- (7) -- (8)
   -- (5) -- (3) (2) -- (4) (5) -- (6) -- (9) ;
}

\def\axgridWa{
\path
	(0,5) 	coordinate (1)
	(5,5) 	coordinate (2)
	(1,4) 	coordinate (3)
	(4,4) 	coordinate (4)
	(1,3) 	coordinate (5)
	(3,3) 	coordinate (6)
	(2,2) 	coordinate (7)
	(3,2) 	coordinate (8)
	(2,1) 	coordinate (9)
	(4,1) 	coordinate (10)
	(0,0) 	coordinate (11)
	(5,0) 	coordinate (12);
\draw (1) -- (2) -- (12) -- (11) -- (1) -- (3) -- (4) -- (2) (3) -- (5) -- (6) -- (4) -- (10) -- (12) (6) -- (8) -- (10) -- (9) -- (11) -- (5) -- (7) -- (9) (7) -- (8);
}

\def\axgridX{
\path
	(0,4) 	coordinate (1)
	(4,4) 	coordinate (2)
	(0.6,3.1) 	coordinate (3)
	(2.1,3.3) 	coordinate (4)
	(1.2,2.4) 	coordinate (5)
	(2,2) 	coordinate (6)
	(2.8,1.6) 	coordinate (7)
	(1.9,.7) 	coordinate (8)
	(3.4,.9) 	coordinate (9)
	(0,0) 	coordinate (10)
	(4,0) 	coordinate (11);
\draw (1) -- (2) -- (11) -- (10) -- (1) -- (3) -- (4) -- (2) -- (6) -- (5) -- (3) (4) -- (5) (7) -- (8) -- (10) -- (6) -- (7) -- (9) -- (11) (8) -- (9);
}

\def\axgridXnew{
\path
	(0,3.6) 	coordinate (1)
	(4,3.6) 	coordinate (2)
	(0.6,2.7) 	coordinate (3)
	(2.1,2.9) 	coordinate (4)
	(2,1.8) 	coordinate (5)
	(1.9,.7) 	coordinate (6)
	(3.4,.9) 	coordinate (7)
	(0,0) 	coordinate (8)
	(4,0) 	coordinate (9);
\draw (1) -- (2) -- (9) -- (8) -- (1) -- (3) -- (4) -- (2) -- 
     (5) -- (7) -- (6) -- (8)  -- (5) -- (3) (7) -- (9);
}

\def\axgridY{
\path
	(-.1,1) 	coordinate (1)
	(1.6,1) 	coordinate (2)
	(-.1,0) 	coordinate (3)
	(1.6,0) 	coordinate (4);
\draw (1) -- (2) -- (4) -- (3) -- (1);
}

\def\axgridZ{
\path
	(-.1,1) 	coordinate (1)
	(1.6,1) 	coordinate (2)
	(-.1,0) 	coordinate (3)
	(1.6,0) 	coordinate (4);
\draw (1) -- (2) -- (4) -- (3) -- (1);
}

\def\commE{ %3x3 %size of diagram
\ifnum\common=1
	\path (1,2) node[c] {\byDef};
	\path (2.2,1.1) node[c] {\byNat};
	\path (1.2,.35) node[c] {\byComOp};
\fi}

\def\labelE{
\ifnum\labon=1
	\path (1) -- node[la] {\ident} (2);
	\path (1) -- node[lb] {\whisk} (5);
	\path (5) -- node[lb] {\expand} (6);
	\path (2) -- node[la] {\ident} (6);
\fi}

\def\commF{ %2x2
\ifnum\common=1
	\path (.8,1) node[c] {\byDef};
	\path (1.3,1.7) node[c] {\fbox{3-6a}};
	\path (1.3,.3) node[c] {\fbox{3-6b}};
\fi}

\def\labelF{
\ifnum\labon=1
	\path (1) -- (2) node[la] {\expand};
	\path (1) -- (4) node[lb] {\switch};
	\path (2) -- (3) node[la,font=\tiny] {\ident}; %!%
	\path (3) -- (5) node[la,font=\tiny] {\ident}; %!%
	\path (4) -- (5) node[lb] {\expand};
\fi}

\def\commFa{ 
\ifnum\common=1
	\path (1.6,2.7) node[c] {\byComOp};
	\path (3.3,1.8) node[c] {\byNat};
	\path (2.0,1.2) node[c] {\byDef};
\fi}

\def\labelFa{
\ifnum\labon=1
	\path (1) -- (2) node[la] {\expand};
	\path (1) -- (5) node[lb] {\hswitch};
	\path (2) -- (5) node[la] {\ident};
\fi}

\def\commG{ %10x10
\ifnum\common=1
	\path (7,.7) node[c] {\byComOp};
	\path (3,9.3) node[c] {3-5};
	\path (1.7,1.7) node[c] {\byComOp};
	\path (8.3,8.3) node[c] {\byNat};
	\path (.7,7) node[c] {\byComOp};
	\path(9.3,3) node[c] {3-5};
	\path (7.75,2.25) node[c] {\byNat};
	\path (2.25,7.75) node[c] {\byNat};
	\path (5.2,2.9) node[c] {\byNat};
	\path (4.8,7.1) node[c] {3-5};
	\path (2.7,5) node[c] {\byNat};
	\path (7.3,5) node[c] {3-5};
	\path (5,5) node[c] {\byNat};
\fi}

\def\labelG{
\ifnum\labon=1
	\path (1) -- (2) node[la] {\whisk};
	\path (1) -- (15) node[lb] {\whisk};
	\path (2) -- (16) node[la] {\whisk};
	\path (15) -- (16) node[lb] {\whisk};
\fi}

\def\commH{
\ifnum\common=1
	\path (3,2.3) node[c] {\fbox{3-8b}};
	\path (1.9,3) node[c] {\byDef};
	\path (4.1,3) node[c] {\byDef};
	\path (1.6,1.5) node[c] {\fbox{3-8a}};
	\path(4.4,1.5) node[c] {\fbox{3-8c}};
	\path(3,.2) node[c] {\byDef};
\fi}

\def\labelH{
\ifnum\labon=1
	\path (1) -- (2) node[la] {\switch};
	\path (1) -- (5) node[la] {\switch};
	\path (2) -- (7) node[lb] {\whisk};
	\path (5) -- (8) node[lb] {\whisk};
	\path (7) -- (8) node[lb] {\switch};
\fi}

\def\commHs{
\ifnum\common=1

\fi}

\def\labelHs{
\ifnum\labon=1
	\path (1) -- (2) node[la] {\expand};
	\path (2) -- (3) node[la] {\ident};
	\path (3) -- (4) node[la] {\ident};
\fi}

\def\commHa{ %7x7
\ifnum\common=1
	\path(5,6.5) node[c] {\byComOp};
	\path(1.5,5.5) node[c] {\byComOp};
	\path(3.5,4.8) node[c] {\byComOp};
	\path(5.4,5.4) node[c] {\byNat};
	\path(4.7,4) node[c] {3-6a};
	\path (6.5,3.8) node[c] {\byDef};
	\path(3.7,3.3) node[c] {\byComOp};
	\path(2.5,2.7) node[c] {\byComOp};
	\path(5,1.5) node[c] {3-5};
	\path(.5,2) node[c] {\byComOp};
	\path (2,1.5) node[c] {\byNat};
	\path(2,.4) node[c] {3-6a};
\fi}

\def\labelHa{
\ifnum\labon=1
	\path (1) -- (2) node[la] {\hswitch};
	\path (1) -- (13) node[lb] {\whisk};
	\path (2) -- (14) node[la] {*};
	\path (13) -- (14) node[lb] {\hswitch};

\fi}

\def\commHb{
\ifnum\common=1 
	\path (3,2) node[c] {\byComOp};
	\path (3,4.2) node[c] {\byComOp};
	\path (1.8,2.9) node[c] {\byComOp};
	\path (4.2,2.9) node[c] {\byComOp};
	\path (2.25,1.25) node[c] {\byDef};
	\path (3.75,1.25) node[c] {\byDef};
	\path (2,4.55) node[c] {3-6a};
	\path (4,4.55) node[c] {3-6a};
\fi}

\def\labelHb{
\ifnum\labon=1
	\path (1) -- (4) node[la] {\hswitch};
	\path (1) -- (6) node[la] {\hswitch};
	\path (4) -- (10) node[lb] {*};
	\path (6) -- (10) node[lb] {*};
\fi}

\def\commK{
\ifnum\common=1
	\path (3,.5) node[c] {\byComOp};
	\path (3,1.5) node[c] {\byLemmaL};
	\path (3,2.8) node[c] {\byDef};
	\path (1.2,1.9) node[c] {\byDef};
	\path (4.8,1.9) node[c] {\byDef};
	\path (1.8,3.8) node[c] {\byLemmaL};
	\path (4.2,3.8) node[c] {\byNat};
\fi}

\def\labelK{
\ifnum\labon=1
	\path (1) -- node[la] {\expand} (3);
	\path (1) -- node[la] {\whisk} (6);
	\path (3) -- node[lb] {\whisk} (11);
	\path (11) -- node[lb] {\expand} (12);
	\path (12) -- node[lb] {\whisk} (6);
\fi}

\def\commL{ %4x4
\ifnum\common=1
	\path (1.5,.5) node[c] {\byNat};
	\path (2,3.5) node[c] {\byDef};
	\path (.3,2) node[c] {\byDef};
	\path (3.6,2.5) node[c] {\byNat};
	\path (1.3,2.6) node[c] {\byLemmaL};
	\path (2.5,2.4) node[c] {\byDef};
	\path (1.55,1.55) node[c] {\byDef};
	\path (3,1.1) node[c] {\byLemmaL};
\fi}

\def\labelL{
\ifnum\labon=1
	\path (1) -- (2) node[la] {\whisk};
	\path (1) -- (12) node[lb] {\whisk};
	\path (2) -- (13) node[la] {\whisk};
	\path (12) -- (13) node[lb] {\whisk};
\fi}

\def\commM{
\ifnum\common=1
	\path (2,3) node[c] {\fbox{3-13a}};
	\path (2,1) node[c] {\fbox{3-13b}};
	\path (.7,2) node[c] {\byDef};
	\path (3.3,2) node[c] {\byDef};
\fi}

\def\labelM{
\ifnum\labon=1
	\path (1) -- (2) node[la] {\whisk};
	\path (1) -- (3) node[la] {\expand}; %!% Notice that some labels are sketchy
	\path (2) -- (5) node[lb] {\switch};
	\path (3) -- (6) node[la] {\switch};
	\path (5) -- (7) node[lb] {\whisk};
	\path (6) -- (7) node[lb] {\expand};
\fi}

\def\commMa{
\ifnum\common=1
	\path (1.2,2.7) node[c] {\byNat};
	\path (2.8,2.7) node[c] {\byLemmaL};
	\path (2,2) node[c] {3-5};
	\path (1.2,1.4) node[c] {3-6a};
	\path (2.8,1.4) node[c] {3-6a};
\fi}

\def\labelMa{
\ifnum\labon=1
	\path (1) -- (3) node[la] {\whisk};
	\path (1) -- (6) node[la] {\expand};
	\path (3) -- (7) node[lb] {\hswitch};
	\path (6) -- (7) node[lb] {\hswitch};
\fi}

\def\commN{
\ifnum\common=1
	\path (.15,.5) node[c] {\byDef};
	\path (.85,.5) node[c] {\byDef};
	\path (.5,.8) node[c] {\fbox{3-14a}};
	\path (.5,.2) node[c] {\fbox{3-14b}};
\fi}

\def\labelN{
\ifnum\labon=1
	\path (1) -- (2) node[la] {\expand};
	\path (1) -- (4) node[lb] {\switch};
	\path (2) -- (5) node[la] {\switch};
	\path (4) -- (5) node[lb] {\expand};
\fi}

\def\commNa{
\ifnum\common=1
	\path (3,7.1) node[c] {\byComOp};
	\path (1.2,6.2) node[c] {\byComOp};
	\path (4.8,6.2) node[c] {\byComOp};
	\path (3,5.5) node[c] {\byComOp};
	\path (2.5,4) node[c] {\byNat};
	\path (3.5,4) node[c] {\byNat};
	\path (1.75,3.75) node[c] {\byNat};
	\path (4.25,3.75) node[c] {\byNat};
	\path (3,1.7) node[c] {\byLemmaL};
	\path (1,3.7) node[c] {3-6a};
	\path (5,3.7) node[c] {3-6a};
\fi}

\def\labelNa{
\ifnum\labon=1
	\path (1) -- (2) node[la,yshift=22pt] {\expand};
	\path (1) -- (12) node[lb,yshift=-37pt] {\hswitch};
	\path (2) -- (12) node[lb,yshift=-37pt] {\hswitch};
\fi}

\def\commO{ %15
\ifnum\common=1
	\path (.75,.5) node[c] {\byLemmaL};
\fi}

\def\labelO{
\ifnum\labon=1
	\path (1) -- (2) node[la] {\expand};
	\path (1) -- (3) node[lb] {\expand};
	\path (2) -- (4) node[la] {\expand};
	\path (3) -- (4) node[lb] {\expand};
\fi}

\def\commP{
\ifnum\common=1

\fi}

\def\labelP{
\ifnum\labon=1
	\path (1) -- (2) node[la] {\expand};
\fi}

\def\commR{ 
\ifnum\common=1
	\path (2.25,2.25) node[c] {\byDef};
	\path (1.3,1.3) node[c] {\byDef};
	\path (.5,1.5) node[c] {\byNat};
	\path (1.5,.5) node[c] {K};
\fi}

\def\commRnew{ 
\ifnum\common=1
	\path (.9,2) node[c] {\byNat};
	\path (2,.9) node[c] {\byLemmaM};
	\path (2.5,2.5) node[c] {\byDef};
\fi}

\def\labelR{
\ifnum\labon=1
	\path (1) -- (2) node[la] {\pass};
	\path (1) -- (6) node[lb] {\expand};
	\path (2) -- (7) node[la] {\ident};
	\path (6) -- (7) node[lb] {\ident};
\fi}

\def\labelRnew{
\ifnum\labon=1
	\path (1) -- (2) node[la] {\pass};
	\path (1) -- (5) node[lb] {\expand};
	\path (2) -- (6) node[la] {\ident};
	\path (5) -- (6) node[lb] {\ident};
\fi}

\def\commS{
\ifnum\common=1
	\path (1.5,.4) node[c] {\byNat};
	\path (1.5,2.6) node[c] {\byComOp};
	\path (1.5,1.5) node[c] {\byNat};
	\path (.5,1.5) node[c] {3-18};
	\path (2.5,1.5) node[c] {3-5};
\fi}

\def\commSnew{
\ifnum\common=1
	\path (1.5,.4) node[c] {\byDef};
	\path (.5,1.5) node[c] {\byDef};
\fi}

\def\labelS{
\ifnum\labon=1
	\path (1) -- (2) node[la] {\pass};
	\path (1) -- (7) node[lb] {\pass};
	\path (2) -- (8) node[la] {\whisk};
	\path (7) -- (8) node[lb] {\pass};
\fi}

\def\labelSnew{
\ifnum\labon=1
	\path (1) -- (2) node[la] {\pass};
	\path (1) -- (6) node[lb] {\pass};
	\path (2) -- (7) node[la] {\whisk};
	\path (6) -- (7) node[lb] {\pass};
\fi}

\def\commT{ %20
\ifnum\common=1
	\path (2,.4) node[c] {\byDef};
	\path (1.5,1) node[c] {\byNat};
	\path (.6,1.4) node[c] {\byNat};
	\path (1.25,1.6) node[c] {\byDef};
	\path (1.9,2.2) node[c] {\byNat};
	\path (1.2,2.4) node[c] {\byNat};
	\path (.6,2.8) node[c] {\byDef};
	\path (1.65,2.87) node[c] {3-6a};
	\path (2.25,3) node[c] {\byDef};
	\path (2,3.7) node[c] {\byNat};
\fi}

\def\labelT{
\ifnum\labon=1
	\path (1) -- (2) node[la] {\switch};
	\path (1) -- (7) node[la] {\flip};
	\path (7) -- (16) node[lb] {\whisk};
	\path (2) -- (10) node[la] {\expand,\flip};
	\path (10) -- (17) node[la] {\whisk};
	\path (16) -- (17) node[lb] {\pass};
\fi}

\def\commTnew{ %20new
\ifnum\common=1
	\path (1.8,3.75) node[c] {\byDef};
	\path (1,3.5) node[c] {3-6a};
	\path (.4,2.6) node[c] {\byNat};
	\path (1.3,2.5) node[c] {\byNat};
	\path (.4,1.3) node[c] {3-5};
	\path (1.5,.4) node[c] {\byDef};
\fi}

\def\labelTnew{
\ifnum\labon=1
	\path (1) -- (2) node[la] {\switch};
	\path (1) -- (5) node[lb] {\flip};
	\path (5) -- (10) node[lb] {\whisk};
	\path (2) -- (7) node[la] {\expand,\flip};
	\path (7) -- (11) node[la] {\whisk};
	\path (10) -- (11) node[lb] {\pass};
\fi}

\def\commTa{
\ifnum\common=1
	\path (1.25,.75) node[c] {3-18};
	\path (2.6,.75) node[c] {\byNat};
	\path (2,1.5) node[c] {3-5};
	\path (2.4,2) node[c] {\byComOp};
	\path (3.3,1.5) node[c] {\byComOp};
	\path (3.3,2.5) node[c] {\byLemmaL};
	\path (3.85,2.3) node[c] {\byComOp};
	\path (2,2.5) node[c] {\byComOp};
	\path (1.5,2.75) node[c] {\byComOp};
	\path (.7,2) node[c] {3-5};
	\path (1.2,3.5) node[c] {\byComOp};
	\path (2.6,3.5) node[c] {3-6};
\fi}

\def\labelTa{
\ifnum\labon=1
	\path (1) -- (3) node[la] {\flip};
	\path (1) -- (7) node[la] {\switch};
	\path (3) -- (11) node[lb] {\whisk};
	\path (11) -- (16) node[lb] {\pass};
	\path (7) -- (15) node[la] {\expand,\flip};
	\path (15) -- (16) node[lb] {\whisk};
\fi}

\def\commU{
\ifnum\common=1
	\path (.75,.5) node[c] {\byDef};
	\path (.3,4.25) node[c] {3-6a};
	\path (2.65,1.2) node[c] {\byNat};
	\path (1.75,1.25) node[c] {\byDef};
	\path (1.3,3) node[c] {K};
	\path (2.5,2.5) node[c] {3-21a};
	\path (2.5,3.15) node[c] {3-6a};
	\path (3.75,1.5) node[c] {\byComOp};
	\path (3.35,3.5) node[c] {J};
	\path (2.25,3.75) node[c] {\byNat};
	\path (2.25,4.75) node[c] {\byNat};
	\path (1.25,5.5) node[c] {\byDef};
\fi}

\def\commUnew{
\ifnum\common=1
	\path (1,1.7) node[c] {\byDef};
	\path (.4,1.1) node[c] {3-6a};
	\path (1.6,1.1) node[c] {\fbox{3-21a}};
	\path (.9,.3) node[c] {\byDef};
\fi}

\def\labelU{
\ifnum\labon=1
	\path (1) -- (8) node[la] {\switch};
	\path (1) -- (2) node[lb] {\expand};
	\path (8) -- (14) node[la] {\expand,\expand};
	\path (2) -- (14) node[lb] {\pass};
\fi}

\def\labelUnew{
\ifnum\labon=1
	\path (1) -- (2) node[la] {\switch};
	\path (1) -- (5) node[lb] {\expand};
	\path (2) -- (6) node[la] {\expand,\expand};
	\path (5) -- (6) node[lb] {\pass};
\fi}

\def\commUa{
\ifnum\common=1

\fi}

\def\labelUa{
\ifnum\labon=1
	\path (1) -- (2) node[la] {\ident};
	\path (1) -- (3) node[lb] {\ident};
	\path (2) -- (4) node[la] {\ident};
	\path (3) -- (4) node[lb] {\ident};
\fi}

\def\commUanew{
\ifnum\common=1
	\path (2.1,2.3) node[c] {3-6a};
	\path (1.6,1.7) node[c] {\byNat};
	\path (.3,1.3) node[c] {\byLemmaX};
	\path (1.2,1.1) node[c] {\byLemmaM};
        \path (3.1,1.3) node[c] {\byComOp};
	\path (1.6,.3) node[c] {\byNat};
\fi}

\def\labelUanew{
\ifnum\labon=1
	\path (1) -- (2) node[la] {\hswitch};
	\path (1) -- (8) node[lb] {\quasiidentity};
	\path (2) -- (9) node[la] {\expand,\expand};
	\path (8) -- (9) node[lb] {\ident};
\fi}

\def\commV{
\ifnum\common=1
	\path (-.4,1) node[c] {3-5};
	\path (2.0,1) node[c] {3-5};
	\path (.75,.5) node[c] {\byDef};
	\path (.75,1.5) node[c] {\byDef};
\fi}

\def\labelV{
\ifnum\labon=1
	\path (1) -- (2) node[la] {\pass};
	\path (1) -- (3) node[la] {\expand};
	\path (3) -- (8) node[lb] {\whisk};
	\path (8) -- (9) node[lb] {\pass};
	\path (2) -- (7) node[la] {\expand};
	\path (7) -- (9) node[lb] {\whisk};
\fi}

\def\commVnew{
\ifnum\common=1
	\path (-.4,1) node[c] {3-5};
	\path (1.7,1) node[c] {3-5};
	\path (.6,.5) node[c] {\byDef};
	\path (.6,1.5) node[c] {\byDef};
\fi}

\def\labelVnew{
\ifnum\labon=1
	\path (1) -- (2) node[la] {\pass};
	\path (1) -- (3) node[la] {\expand};
	\path (3) -- (7) node[lb] {\whisk};
	\path (7) -- (8) node[lb] {\pass};
	\path (2) -- (6) node[la] {\expand};
	\path (6) -- (8) node[lb] {\whisk};
\fi}

\def\commW{
\ifnum\common=1
	\path (2.1,.5) node[c] {\byNat};
	\path (1,.75) node[c] {\byComOp};
	\path (.5,2) node[c] {\byDef};
	\path (2.6,2.5) node[c] {\byNat};
	\path (1.5,3.5) node[c] {\byNat};
\fi}

\def\labelW{
\ifnum\labon=1
	\path (1) -- (2) node[la] {\flip};
	\path (1) -- (9) node[lb] {\pass};
	\path (2) -- (10) node[la] {\swap};
	\path (9) -- (10) node[lb] {\expand};
\fi}

\def\commWnew{
\ifnum\common=1
	\path (2.1,1) node[c] {\byNat};
	\path (1.2,1.5) node[c] {\byLemmaM};
	\path (.5,2.5) node[c] {\byDef};
	\path (1.5,2.5) node[c] {\byLemmaZ};
	\path (2.6,2.5) node[c] {\byNat};
	\path (1.5,3.5) node[c] {\byNat};
\fi}

\def\labelWnew{
\ifnum\labon=1
	\path (1) -- (2) node[la] {\flip};
	\path (1) -- (8) node[lb] {\pass};
	\path (2) -- (9) node[la] {\swap};
	\path (8) -- (9) node[lb] {\expand};
\fi}

\def\commWa{
\ifnum\common=1
	\path (3,.5) node[c] {\byComOp};
	\path (1.25,1.5) node[c] {\byComOp};
	\path (2.7,1.45) node[c] {3-23L};
	\path (2.25,2.5) node[c] {\byComOp};
	\path (3.5,2.5) node[c] {\byComOp};
	\path (4.5,2.5) node[c] {\byComOp};
	\path (2.2,3.5) node[c] {3-24};
	\path (.45,3.5) node[c] {\byNat};
	\path (2.5,4.5) node[c] {\byComOp};
\fi}

\def\labelWa{
\ifnum\labon=1
	\path (1) -- (2) node[la] {\expand};
	\path (1) -- (11) node[lb] {\pass};
	\path (2) -- (12) node[la] {\swap};
	\path (11) -- (12) node[lb] {\flip};
\fi}

\def\commX{
\ifnum\common=1
	\path (2.75,.4) node[c] {\byNat};
	\path (1.25,3.6) node[c] {\byNat};
	\path (1.9,1.3) node[c] {K};
	\path (2.1,2.7) node[c] {K};
	\path (2.75,1.15) node[c] {\byDef};
	\path (1.25,2.85) node[c] {\byDef};
	\path (.7,1.75) node[c] {\byDef};
	\path (3.3,2.25) node[c] {\byDef};
\fi}

\def\labelX{
\ifnum\labon=1
	\path (1) -- (2) node[la] {\expand};
	\path (1) -- (10) node[lb] {\pass};
	\path (2) -- (11) node[la] {\pass};
	\path (10) -- (11) node[lb] {\expand};
\fi}

\def\commXnew{
\ifnum\common=1
	\path (2.75,.4) node[c] {\byNat};
	\path (1.25,3.2) node[c] {\byNat};
	\path (1.9,1.1) node[c] {\byLemmaM};
	\path (2.1,2.5) node[c] {\byLemmaM};
	\path (.7,1.55) node[c] {\byDef};
	\path (3.3,1.85) node[c] {\byDef};
\fi}

\def\labelXnew{
\ifnum\labon=1
	\path (1) -- (2) node[la] {\expand};
	\path (1) -- (8) node[lb] {\pass};
	\path (2) -- (9) node[la] {\pass};
	\path (8) -- (9) node[lb] {\expand};
\fi}

\def\commY{
\ifnum\common=1
	\path (.75,.5) node[c] {\byLemmaL};
\fi}

\def\labelY{
\ifnum\labon=1
	\path (1) -- (2) node[la] {\expand};
	\path (1) -- (3) node[lb] {\expand};
	\path (2) -- (4) node[la] {\expand};
	\path (3) -- (4) node[lb] {\expand};
\fi}

\def\commZ{
\ifnum\common=1
	\path (.75,.5) node[c] {\byNat};
\fi}

\def\labelZ{
\ifnum\labon=1
	\path (1) -- (2) node[la] {\flip};
	\path (1) -- (3) node[lb] {\expand};
	\path (2) -- (4) node[la] {\expand};
	\path (3) -- (4) node[lb] {\flip};
\fi}

%% file: cn-proofs.bbl
\begin{thebibliography}{}
\providecommand\bibmarginpar{\leavevmode\marginpar}
\def\urlstyle#1{{\tt #1}}

\bibitem{BDH(nets)}
\textbf{A Bartels}, \textbf{C\,L Douglas}, \textbf{A Henriques},
  \emph{Conformal nets {I}: {C}oordinate-free nets}, Int. Math. Res. Not. 13
  (2015) 4975--5052 \xox{arXiv}{1302.2604v2}

\bibitem{BDH(modularity)}
\textbf{A Bartels}, \textbf{C\,L Douglas}, \textbf{A Henriques},
  \emph{Conformal nets {II}: {C}onformal blocks}, Comm. Math. Phys. 354 (2017)
  393--458 \xox{arXiv}{1409.8672}

\bibitem{BDH(1*1)}
\textbf{A Bartels}, \textbf{C\,L Douglas}, \textbf{A Henriques},
  \emph{Conformal nets {III}: {F}usion of defects}, Mem. Amer. Math. Soc.
  (2017) \xox{arXiv}{1310.8263}

\bibitem{Bischoff-Kawahigashi-Longo-Rehren(Phase-boundaries-in-algebraic-conformal-QFT)}
\textbf{M Bischoff}, \textbf{Y Kawahigashi}, \textbf{R Longo}, \textbf{K-H
  Rehren}, \emph{Phase
  boundaries in algebraic conformal {QFT}}, Comm. Math. Phys. 342 (2016) 1--45

\bibitem{Brunetti-Guido-Longo(1993modular+duality-in-CQFT)}
\textbf{R Brunetti}, \textbf{D Guido}, \textbf{R Longo},
  \emph{Modular structure and duality in conformal quantum field theory},
  Comm. Math. Phys. 156 (1993) 201--219

\bibitem{Buchholz-Mack-Todorov(1988current-alg)}
\textbf{D Buchholz}, \textbf{G Mack}, \textbf{I Todorov},
  \emph{The current
  algebra on the circle as a germ of local field theories}, Nuclear Phys. B
  Proc. Suppl. 5B (1988) 20--56

\bibitem{Buchholz-Schulz-Mirbach(Haag-duality-in-conformal-quantum-field-theory)}
\textbf{D Buchholz}, \textbf{H Schulz-Mirbach},
  \emph{Haag duality in
  conformal quantum field theory}, Rev. Math. Phys. 2 (1990) 105--125

\bibitem{Douglas-Henriques(Internal-bicategories)}
\textbf{C\,L Douglas}, \textbf{A Henriques}, \emph{Internal bicategories}
  (2012) \xox{arXiv}{1206.4284}

\bibitem{Gabbiani-Froehlich(OperatorAlg-CFT)}
\textbf{F Gabbiani}, \textbf{J Fr{\"o}hlich},
  \emph{Operator algebras and conformal field theory}, Comm. Math. Phys. 155
  (1993) 569--640

\bibitem{Kawahigashi-Longo(2004classification)}
\textbf{Y Kawahigashi}, \textbf{R Longo},
  \emph{Classification of
  local conformal nets. {C}ase {$c<1$}}, Ann. of Math. (2) 160 (2004) 493--522

\bibitem{Kawahigashi-Longo-Mueger(2001multi-interval)}
\textbf{Y Kawahigashi}, \textbf{R Longo}, \textbf{M M{\"u}ger},
  \emph{Multi-interval subfactors and modularity of representations in
  conformal field theory}, Comm. Math. Phys. 219 (2001) 631--669

\bibitem{Mather(1974commutators)}
\textbf{J\,N Mather}, \emph{Commutators of diffeomorphisms}, Comment. Math.
  Helv. 49 (1974) 512--528

\bibitem{Wassermann(Operator-algebras-and-conformal-field-theory)}
\textbf{A Wassermann}, \emph{Operator algebras and conformal field theory. {III}. {F}usion of
  positive energy representations of {${\rm LSU}(N)$} using bounded
  operators}, Invent. Math. 133 (1998) 467--538

\end{thebibliography}
